\documentclass{amsart}
\usepackage[all]{xy}
\usepackage{latexsym}
\usepackage{amssymb}
\usepackage{amsmath}
\usepackage{amsthm}
\usepackage{amscd}
\usepackage{fancyhdr}
\usepackage[dvips]{graphicx}
\usepackage{a4wide}
\usepackage{mathrsfs}

\theoremstyle{plain}
\newtheorem{pro}{Proposition}
\newtheorem{thm}[pro]{Theorem}
\newtheorem{lem}[pro]{Lemma}
\newtheorem{cor}[pro]{Corollary}
\newtheorem{pro-def}[pro]{Definition-Proposition}
\newtheorem{fact}[pro]{Fact}
\theoremstyle{definition}
\newtheorem*{dei}{Definition}
\newtheorem*{Rq}{\sc Remark}
\newtheorem*{Ex}{Example}

\newcommand{\X}{\mathcal{X}}
\newcommand{\I}{\mathcal{I}}
\newcommand{\oD}{\overline{\Delta}}
\newcommand{\KK}{\mathbb{K}}
\newcommand{\bc}{\boxtimes}

\newcommand{\bbc}{\boxtimes_{(1,1)}}
\newcommand{\R}{\mathcal{R}}
\newcommand{\End}{\mathrm{End}}
\newcommand{\oc}{\circ}
\newcommand{\cp}{\star}
\newcommand{\As}{\mathrm{As}}
\newcommand{\ad}{\mathrm{ad}}

\newcommand{\B}{\mathcal{B}}
\newcommand{\oC}{\overline{\mathcal{C}}}

\newcommand{\oPo}{\bar{\Po}}

\newcommand{\BLi}{\mathcal{B}i\mathcal{L}ie}
\newcommand{\IBi}{\varepsilon\mathcal{B}i}
\newcommand{\G}{\mathcal{G}}

\newcommand{\NN}{\mathbb{N}}
\newcommand{\ZZ}{\mathbb{Z}}

\newcommand{\kj}{{\bar{k},\, \bar{\jmath}}}
\newcommand{\oi}{{\bar{\imath}}}
\newcommand{\oj}{{\bar{\jmath}}}
\newcommand{\ok}{{\bar{k}}}
\newcommand{\ol}{{\bar{l}}}

\newcommand{\Ker}{\mathop{\rm Ker }}

\newcommand{\Sy}{\mathbb{S}}

\newcommand{\A}{\mathcal{A}s}

\newcommand{\Po}{\mathcal{P}}
\newcommand{\F}{\mathcal{F}}

\newcommand{\oF}{\overline{\mathcal{F}}}
\newcommand{\ac}{\scriptstyle \textrm{!`}}

\newcommand{\II}{I}
\newcommand{\Id}{\mathrm{Id}}
\newcommand{\Qo}{\mathcal{Q}}
\newcommand{\Ro}{\mathcal{R}}
\newcommand{\Li}{\mathcal{L}ie}
\newcommand{\Sc}{\mathbb{S}^c}

\newcommand{\Frob}{\mathcal{F}rob}
\newcommand{\NCFrob}{\mathcal{NC}$-$\mathcal{F}rob}
\newcommand{\Co}{\mathcal{C}}
\newcommand{\Do}{\mathcal{D}}

\newcommand{\UU}{\mathcal{U}}

\newcommand{\ot}{\otimes}

\newcommand{\Hom}{\mathrm{Hom}}
\newcommand{\Der}{\mathrm{Der}}
\newcommand{\CoDer}{\mathrm{CoDer}}
\newcommand{\sgn}{\mathrm{sgn}}

\newcommand{\epi}{\twoheadrightarrow}
\newcommand{\mono}{\rightarrowtail}

\newcommand{\Y}{\vcenter{\xymatrix@M=0pt@R=6pt@C=6pt{
\ar@{-}[dr] &  &\ar@{-}[dl]  \\
 &\ar@{-}[d] &  \\  & &}}}
\newcommand{\YY}{\vcenter{\xymatrix@M=0pt@R=6pt@C=6pt{
\ar@{-}[dr] &  &\ar@2{-}[dl]  \\
 &\ar@2{-}[d] &  \\  & &}}}
\newcommand{\YYY}{\vcenter{\xymatrix@M=0pt@R=6pt@C=6pt{
\ar@{-}[dr] &  &\ar@3{-}[dl]  \\
 &\ar@3{-}[d] &  \\  & &}}}
\newcommand{\cop}{\vcenter{\xymatrix@M=0pt@R=6pt@C=6pt{
 & \ar@{-}[d] & \\
 &\ar@{-}[dr] \ar@{-}[dl] &  \\  & &}}}

\newcommand{\copL}{\xymatrix@M=0pt@R=6pt@C=6pt{
 & \ar@{-}[d] & \\
 &\ar@{-}[dr] \ar@{-}[dl] &  \\  & &\\  & &\\  & &}}
\newcommand{\YL}{\vcenter{\xymatrix@M=0pt@R=6pt@C=6pt{
\ar@{-}[dr] &  &\ar@{-}[dl]  \\
 &\ar@{-}[d] &   \\  & &\\  & &\\  & &}}}
\newcommand{\YYL}{\vcenter{\xymatrix@M=0pt@R=6pt@C=6pt{
\ar@{-}[dr] &  &\ar@2{-}[dl]  \\
 &\ar@2{-}[d] &  \\  & &\\  & &\\  & &}}}
\newcommand{\LYY}{\vcenter{\xymatrix@M=0pt@R=6pt@C=6pt{
\ar@{-}[dr] &  &\ar@2{-}[dl]  \\
 &\ar@2{-}[d] &  \\  & &\\  & &\\  & &}}}

\newcommand{\XX}{\vcenter{\xymatrix@M=0pt@R=6pt@C=6pt{\ar@{-}[ddrr]&&\ar@{-}[ddll] \\ && \\ &&   }}}

\newcommand{\Ta}{\vcenter{\xymatrix@M=0pt@R=6pt@C=6pt{ \ar@{-}[dddrrr] && \ar@{-}[dl] &&  \\
&&& \ar@{-}[dl]  &  \\ &&&&  \ar@{-}[dl]  \\&&&  \ar@{-}[d] &
\\&&&& }}}
\newcommand{\Tb}{\vcenter{\xymatrix@M=0pt@R=6pt@C=6pt{  & \ar@{-}[dr]&&\ar@{-}[dl] \\
\ar@{-}[dr]&&\ar@{-}[dl]& \\&\ar@{-}[dr]&&\ar@{-}[dl]
\\&&\ar@{-}[d]& \\&&& }}}
\newcommand{\Tc}{\vcenter{\xymatrix@M=0pt@R=6pt@C=6pt{   \ar@{-}[dr]&&\ar@{-}[dl]& \\
&\ar@{-}[dr]&& \ar@{-}[dl] \\\ar@{-}[dr]&&\ar@{-}[dl]&
\\&\ar@{-}[d]&& \\&&& }}}
\newcommand{\Td}{\vcenter{\xymatrix@M=0pt@R=6pt@C=6pt{ && \ar@{-}[dr]&&\ar@{-}[dddlll] \\
 &\ar@{-}[dr]&&& \\ \ar@{-}[dr]&&&& \\& \ar@{-}[d]&&& \\&&&&  }}}
\newcommand{\Te}{\vcenter{\xymatrix@R=3pt@C=3pt{\ar@{-}[drdr] &&\ar@{-}[dl]  *=0{}
\ar@{-}[dr]&& \ar@{-}[ddll] \\ &&& *=0{}& \\&& *=0{} \ar@{-}[d]&&
\\&&&& }}}
\newcommand{\TaC}{\vcenter{\xymatrix@M=0pt@R=6pt@C=6pt{ \ar@{-}[ddddddrrrrrr] && \ar@{-}[dl] && && \\
&&& \ar@{-}[dl]  &&&  \\ &&&&  \ar@{-}[dl]&&  \\&&& &&&
\\&&&&\ar@{-}[dl]&& \\&&&&&\ar@{-}[dl]&\\&&&&&& }}}
\newcommand{\TreeL}{\vcenter{\xymatrix@M=0pt@R=5pt@C=5pt{ \ar@{-}[dr] &
&\ar@{-}[dl] & &  \\
& \ar@{-}[dr] & &\ar@{-}[dl]  & \\
& &\ar@{-}[d] & & \\
& & \\ & & }}}
\newcommand{\TreeR}{\vcenter{\xymatrix@M=0pt@R=5pt@C=5pt{
 & &\ar@{-}[dr] & & \ar@{-}[dl]  \\
& \ar@{-}[dr] & &\ar@{-}[dl]  & \\
& &\ar@{-}[d] & & \\
& & \\ & & }}}

 \newcommand{\draftnote}[1]{}

\def\one{{\mbox{1 \hskip -8pt 1}}}
 \newcommand{\Ass}{{\mathcal A}ss}
\newcommand{\Hm}{\Hom^\Sy_\bullet}
\newcommand{\Com}{\mathcal Com}
\newcommand{\sip}{\smallskip}
\newcommand{\bip}{\bigskip}
 \newcommand{\lon}{\longrightarrow}
 \newcommand{\rar}{\rightarrow}
 \newcommand{\Mor}{{\mathrm M\mathrm o\mathrm r}}

 \newcommand{\bZ}{{\mathbb Z}}
 \newcommand{\bS}{{\mathbb S}}
 \newcommand{\bK}{{\mathbb K}}
 \newcommand{\p}{{\partial}}

  \newcommand{\bN}{{\mathbb N}}
 \newcommand{\fg}{{\mathfrak g}}

 \newcommand{\f}{{\mathcal O}}

 \newcommand{\cC}{{\mathcal C}}
 \newcommand{\caD}{{\mathcal D}}
 \newcommand{\cE}{{\mathcal E}}
 \newcommand{\cF}{{\mathcal F}}

 \newcommand{\cP}{{\mathcal P}}
\newcommand{\cQ}{{\mathcal Q}}

\newcommand{\sB}{{\mathsf B}}
\newcommand{\ass}{{\mathsf A  \mathsf s  \mathsf s}}
\newcommand{\ab}{{\mathcal Ass\mathcal Bi}}
\newcommand{\asb}{{\mathcal A  s  s \mathcal B}}

\newcommand{\hsB}{{\mathsf B}}
\newcommand{\hp}{{\p_{ind}}}
\newcommand{\sE}{{\mathsf E}}

\newcommand{\sd}{{\mathsf d}}
\newcommand{\Lieb}{{\mathcal LieBi}}

\newcommand{\Beq}{\begin{equation}}
 \newcommand{\Eeq}{\end{equation}}
 \newcommand{\Beqr}{\begin{eqnarray}}
 \newcommand{\Eeqr}{\end{eqnarray}}
 \newcommand{\Beqrn}{\begin{eqnarray*}}
 \newcommand{\Eeqrn}{\end{eqnarray*}}
 \newcommand{\Ba}{\begin{array}}
 \newcommand{\Ea}{\end{array}}
 \newcommand{\Bit}{\begin{itemize}}
 \newcommand{\Eit}{\end{itemize}}

 \newcommand{\al}{\alpha}
 \newcommand{\be}{\beta}
 \newcommand{\ga}{\gamma}
 \newcommand{\Ga}{\Gamma}
 
\newenvironment{proo}{\begin{trivlist} \item{\sc {Proof.}}}
  {\hfill $\square$ \end{trivlist}}

\long\def\symbolfootnote[#1]#2{\begingroup%
\def\thefootnote{\fnsymbol{footnote}}\footnote[#1]{#2}\endgroup}

\title{Deformation theory of representations of prop(erad)s}

\author{Sergei Merkulov, Bruno Vallette}
\begin{document}

\sloppy


\begin{abstract}

We study the deformation theory of morphisms of properads and
props thereby extending  Quillen's
deformation theory for commutative rings to a non-linear framework. The associated chain
complex is endowed with an $L_\infty$-algebra structure. Its Maurer-Cartan elements correspond to deformed
structures, which allows us to give a geometric interpretation of
these results.

To do so, we endow the category of prop(erad)s with a model category
structure. We provide a complete study of models for prop(erad)s.
A new effective method to make minimal models explicit, that
extends the Koszul duality theory, is introduced and the associated
notion is called \emph{homotopy Koszul}.

As a corollary, we obtain the (co)homology theories of (al)gebras
over a prop(erad) and of homotopy (al)gebras as well. Their
underlying chain complex is endowed  an $L_\infty$-algebra structure in
general and a Lie algebra structure only in the Koszul case. In
particular, we make the deformation complex of morphisms
from the properad of associative bialgebras explicit . For any minimal model
of this properad, the boundary map of this chain complex is shown
to be the one defined by Gerstenhaber and Schack. As a corollary,
this paper provides a complete proof of the existence of an
$L_\infty$-algebra structure on the Gerstenhaber-Schack bicomplex
associated to the deformations of associative bialgebras.

\end{abstract}
\maketitle

\section*{Introduction}
The theory of props and properads, which generalizes the theory of
operads, provides us with a universal language to describe many
algebraic, topological and differential geometric structures. Our
main purpose in this paper is to establish a new and surprisingly
strong link between the theory of prop(erad)s and the theory of
$L_\infty$-algebras.

\bip

We introduce several families of $L_\infty$-algebras canonically associated with
 prop(erad)s, moreover, we develop new methods which explicitly compute the associated
$L_\infty$-brackets in
terms of prop(erad)ic compositions and differentials. Many important dg Lie algebras in algebra and geometry
(such as Hochschild, Poisson, Schouten, Fr\"olicher-Nijenhuis and many others) are proven to be of this prop(erad)ic
origin.

\bip

The $L_\infty$-algebras we construct in our paper out of dg
prop(erad)s encode many important properties of the latter. The
most important one controls the deformation theory of morphisms of
prop(erad)s and, in particular, the deformation theory of
(al)gebras over prop(erad)s. Applications of our results to the
deformation theory of many algebraic and geometric structures
becomes therefore another major theme of our paper.

\bip

On the technical side, the story develops (roughly speaking) as
follows: first we associate canonically to a pair, $(\cF(V), \p)$
and $(\Qo, d)$, consisting of a differential graded (dg, for
short) quasi-free prop(erad) $\cF(V)$ on a $\bS$-bimodule $V$ and
an arbitrary dg prop(erad)  $\Qo$, a structure of
$L_\infty$-algebra on the (shifted) graded vector space,
 $s^{-1}\Hm(V, \Qo)$, of morphisms
of $\bS$-bimodules; then we prove the Maurer-Cartan elements of
this $L_\infty$-algebra are in {\em one-to-one}\, correspondence with the
set of all dg morphisms,
$$
\{(\cF(V), \p)\lon (\Qo, d)\},
$$ of dg prop(erad)s. This canonical $L_\infty$-algebra is used then to define, for any
particular morphism $\ga: (\cF(V), \p)\rar (\Qo, d)$, another
{twisted} $L_\infty$-algebra which controls deformation theory of
the morphism $\ga$. In the special case when $(\Qo, d)$ is the
endomorphism prop(erad), $(\End_X, d_X)$, of some dg vector space
$X$, our theory gives $L_\infty$-algebras which control
deformation theory of many classical algebraic and geometric
structures on $X$, for example, associative algebra structure, Lie
algebra structure, commutative algebra structure, Lie bialgebra
structure, associative bialgebra structure, formal Poisson
structure, Nijenhuis structure etc. As the case of associative
bialgebras has never been rigorously treated in the literature
before, we discuss this example in full details; we prove, in
particular, that the first term of the canonical
$L_\infty$-structure controlling deformation theory of bialgebras
is precisely the Gerstenhaber-Schack differential.\\

We derive and study the deformation complex and its
$L_\infty$-structure from several different perspectives. One of
them can be viewed as a nontrivial generalization to the case of
prop(erad)s of Van der Laan's approach \cite{VanderLaan02} to the
deformation theory of algebras over operads, while others are
completely new and provide us with, perhaps, a conceptual
explanation of the observed (long ago) phenomenon that deformation
theories are controlled by dg Lie and, more generally, $L_\infty$
structures.   \\

First, we define the deformation complex of a morphism of
prop(erad)s $\Po \to \Qo$ generalizing Quillen's definition of the
deformation complex of a morphism of commutative rings. When
$(\F(\Co), \p)$ is a quasi-free resolution of $\Po$, we prove that
this chain complex is isomorphic, up to a shift of degree, to the
space of morphisms of $\Sy$-bimodule $\Hom^\Sy_\bullet(\Co, \Qo)$,
where $\Co$($\simeq s^{-1}V$) is a homotopy coprop(erad), that is
the dual notion of prop(erad) with relations up to homotopy. Since
$\Qo$ is a (strict) prop(erad), we prove that the space
$\Hom^\Sy_\bullet(\Co, \Qo)$ has a rich algebraic structure,
namely it is a homotopy non-symmetric prop(erad), that is a
prop(erad) without the action of the symmetric groups and with
relations up to homotopy. From this structure, we extract a
canonical $L_\infty$-structure on $\Hom^\Sy_\bullet(\Co,
\Qo)\simeq s^{-1}\Hm(V, \Qo)$. We also obtain  higher operations
with $m+n$ inputs acting on $\Hom^\Sy_\bullet(\Co, \Qo)$ which are
 important in applications. In the case of, for example,
 the non-symmetric operad, $\mathcal{A}ss$, of associative algebras
 the deformation complex is the Hochschild cochain complex
of an associative algebra, and the higher homotopy operations are
shown to be non-symmetric braces operations which
 play a fundamental role in the proof
of Deligne's conjecture
(see
\cite{Tamarkin98, Voronov00, KontsevichSoibelman00, McClureSmith02, BergerFresse04, Kaufmann07}).\\

Recall that M. Markl proved in \cite{Markl04Linfini} a first
interesting partial result, that is for a given minimal model
$(\F(\Co), \p)$  a prop(erad) $\Po$ concentrated in degree $0$,
there exists a $L_\infty$-structure on the space of derivations
from $\F(\Co)$ to $\End_X$, where $X$ is a $\Po$-(al)gebra. By
using a different conceptual method based on the notions of
homotopy (co)prop(erad)s and convolution prop(erad)s, we prove
here that for any representation $\Qo$ of any prop(erad) $\Po$,
there exists a homotopy prop(erad) structure on the space of
derivations from any quasi-free resolution of $\Po$ to $\Qo$.
Moreover this construction
is shown to be functorial, that is does not depend on the model
chosen. From this
we derive functorially the general $L_\infty$-structure.\\

Another approach of deriving the deformation complex and its $L_\infty$-structure
is based on a canonical enlargement of the category of dg prop(erad)s via an extension of the notion
of {\em morphism}; the set of morphisms, $\Mor_\bZ(\cP_1,\cP_2)$, in this enlarged category is
identified with a certain {\em dg affine scheme}\, naturally associated with both $\cP_1$ and $\cP_2$;
moreover, when the dg prop(erad) $\cP_1$ is quasi-free,
the dg affine scheme $\Mor(\cP_1,\cP_2)$ is proven to be a {\em smooth}\, dg manifold for any $\cP_2$ and hence
gives canonically rise to a $L_\infty$-structure.\\

The third major theme of our work is the theory of models and
minimal models. To make explicit the deformation complex, we need
models, that is quasi-free resolutions of prop(erad)s. We extend
the bar and cobar construction to prop(erad)s and show that the
bar-cobar construction provides a canonical cofibrant resolution
of a prop(erad). Since this construction is not very convenient to
work with because it is too big, we would rather use minimal
models. We give a complete account to the theory of minimal models
for prop(erad)s. We prove that minimal models for prop(erad)s are not in general
determined by resolutions of their genus $0$ parts, namely
dioperads, giving thereby a negative  answer to a question
raised by M. Markl and A.A. Voronov \cite{MarklVoronov03}, that is
we prove that the free functor from dioperads to prop(erad)s is
not exact. We provide an explicit example of a Koszul dioperad
which does not induce the prop(erad)ic resolution of the
associated
prop(erad).\\

A properad is Koszul if and only if it admits a quadratic model.
In this case, Koszul duality theory of properad \cite{Vallette07}
provides an effective method to compute this special minimal
model. Unfortunately, not all properads are Koszul. For instance,
the properad coding associative bialgebras is not. We include this
example in a new notion, called \emph{homotopy Koszul}. A homotopy
Koszul properad is shown to have a minimal model that can be
explicitly computed. Its space of generators is equal to the
Koszul dual of a quadratic properad associated to it. And the
differential is made explicit by use of the (dual) formulae of J.
Gran\aa ker \cite{Granaker06} of transfer of homotopy coproperad
structure, that is by perturbing the differential. We apply this
notion to show that morphisms of homotopy $\Po$-algebras are in
one-to-one correspondence with
Maurer-Cartan elements of a convolution $L_\infty$-algebra.\\

In the appendix, we endow the category of dg prop(erad)s with a
 model category structure which is used throughout the text.\\

The paper is organized as follows. In \S 1 we remind key facts
about properads and props and we define the notion of
\emph{non-symmetric prop(erad)}. In \S 2 we introduce and study
the convolution prop(erad) canonically associated with a pair,
$(\cC,\Po)$, consisting of an arbitrary coprop(erad) $\cC$ and an
arbitrary prop(erad) $\Po$; our main results is the construction of a Lie algebra structure on this
convolution properad, as well as higher operations. In \S 3 we discuss bar
and cobar constructions for (co)prop(erad)s. We introduce the
notion of \emph{twisting morphism (cochain)} for prop(erads) and
prove Theorem~\ref{BarCobarResolution} on bar-cobar resolutions
extending thereby earlier results of \cite{Vallette03} from
weight-graded dg properads to arbitrary dg properads. In \S 4 we
recall to the notion and properties of homotopy properads which
were  first introduced in \cite{Granaker06} and we define the
notions of \emph{homotopy (co)prop(erad)}. We apply these notions
to convolution prop(erad)s. In \S 5, we give a complete study of
minimal models for properads. In \S 6 we remind geometric
interpretation of $L_\infty$-algebras, and then use this geometric
language to prove Theorem~\ref{mor} which associates to pair,
$(\cF(V), \p)$ and $(\cP, d)$, consisting of quasi-free prop(erad)
$(\cF(V), \p)$ and an arbitrary dg prop(erad) $(\Po, d)$, a
structure of $L_\infty$-algebra on the (shifted) graded vector
space,
 $s^{-1}\Hm(V, \cP)$; we also show in \S 6 full details behind the construction
of the above mentioned  enlargement of the category of dg prop(erad)s.
In \S 7, we define the deformation complex
 following Quillen's methods and identify it with
 $s^{-1}\Hm(V, \cP)$ in Theorem~\ref{boundary=Q}.
 We show next how this canonical $L_\infty$-algebra gives rise to twisted
 $L_\infty$-algebras
which control deformation theories of particular morphisms $\ga:
(\cF(V), \p)\rar (\Qo, d)$. In \S 8, we illustrate this general
construction with several examples from algebra and geometry. We
make explicit the deformation complex of representation of the
properad of associative bialgebras and prove that it corresponds
with the one defined by Gerstenhaber-Schack. In the Appendix, we
show that the category of dg prop(erad)s is a cofibrantly
generated model category.

 \tableofcontents

In this paper, we will always work over a field $\KK$ of
characteristic $0$.

\section{(Co)Properads, (Co)Props and their non-symmetric versions}\label{Conventions}

In this section, we recall briefly the definitions of (co)properad
and (co)prop. For the reader who does not know what a properad or what a prop is, we
strongly advise to read the first sections of \cite{Vallette03}
before reading the current article since we will make use of the notions everywhere in the sequel. Generalizing the
notion of non-symmetric operads to prop(erad), we introduce the
notions of \emph{non-symmetric properad} and \emph{non-symmetric
prop}.

\subsection{$\Sy$-bimodules, graphs, composition
products}\label{S-bimodule}

A \emph{(dg) $\Sy$-bimodule} is a collection $\cP=\{\Po(m, n)
\}_{m, \, n\in \NN}$ of dg modules over the symmetric groups $\Sy_n$ on
the right and $\Sy_m$ on the left. These two actions are supposed
to commute. In the sequel, we will mainly consider \emph{reduced}
$\Sy$-bimodules, that is $\Sy$-bimodules $\Po$ such that $\Po(m,
n)=0$ when $n=0$ or $m=0$. We use the homological convention, that
is the degree of differentials is $-1$. An $\Sy$-bimodule $\cP$ is
\emph{augmented} when it naturally splits as $\cP= \overline{\cP}
\oplus I$ where $I=\{I(m,n\}$ is an $\Sy$-bimodule with all components  $I(m,n)$ vanishing except
 for $I(1,1)$ which equals $\KK$.
We denote the module of morphisms of $\Sy$-bimodules by
$\Hom(\cP,{\mathcal Q})$ and the module of equivariant morphisms,
with respect
to the action of the symmetric groups, by $\Hom^\Sy(\cP,{\mathcal Q})$.\\

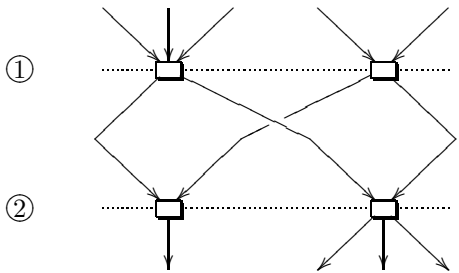
\begin{figure}[h]
$$ \xymatrix@R=17pt@C=20pt{
 & \ar[dr] & \ar[d] &\ar[dl] &\ar[dr] & &\ar[dl] \\
*+[o][F-]{1} & \ar@{..}[r]& *+[F-,]{\ } \ar@{-}[dl]
\ar@{-}[drr] \ar@{..}[rrr]& & & *+[F-,]{\ }
\ar@{-}[dr]   \ar@{-}[dll] |(0.75) \hole  \ar@{..}[r]& \\
 & *=0{} \ar[dr] & &*=0{}\ar[dl] &*=0{}\ar[dr] & &*=0{} \ar[dl]\\
*+[o][F-]{2}& \ar@{..}[r] & *+[F-,]{\ } \ar[d]
\ar@{..}[rrr]& & & *+[F-,]{\ } \ar[dl] \ar[d]
\ar[dr]
\ar@{..}[r] & \\
& & & & & & } $$ \caption{Example of a 2-level graph.}
\label{2-levelled graph}
\end{figure}

A graph is given by two sets, the set $V$ of vertices and the set $E$ of edges, and relations among which say when an edge is attached to one or two vertices (see \cite{GetzlerKapranov98} (2.5)). The egdes of the graph considered in the sequel will always be directed by a global flow (\emph{directed graphs}). The edges can be divided into two parts: the ones with one vertex at each end, called \emph{internal edges}, and the ones with just one vertex on one end, called \emph{external edges}. The \emph{genus} of a graph is the first Betti number of the underlying topological space of a graph. A \emph{2-levelled directed graphs} is a directed graph such that the vertices are divided into two parts, the ones belonging to a bottom level and the ones belonging to a top level. In the category of $\Sy$-bimodule, we define  two \emph{composition} products,
$\boxtimes$ based on the composition of operations indexing the
vertices of a 2-levelled directed graphs, and $\boxtimes_c$ based on the
composition of operations indexing the vertices of a 2-levelled
directed connected graph (see Figure~\ref{2-levelled graph} for an
example).
Let $\G$ be such a graph with $N$ internal edges between
vertices of the two levels. This set of edges between vertices of
the first level and vertices of the second level induces a
permutation of $\Sy_N$.

Let $\Po$ and $\Qo$ be two $\Sy$-bimodules, their composition
product is given by the explicit formula
$$\mathcal{P}\boxtimes \mathcal{Q}(m,\, n) :=
 \bigoplus_{N\in \mathbb{N}^*} \left( \bigoplus_{\ol,\, \ok,\, \oj,\, \oi} \KK[\mathbb{S}_m]
 \otimes_{\mathbb{S}_\ol}
\mathcal{P}(\ol,\, \ok)\otimes_{\mathbb{S}_{\ok}} \KK[\Sy_N]
\otimes_{\mathbb{S}_\oj} \mathcal{Q}(\oj,\, \oi)
\otimes_{\mathbb{S}_\oi} \KK[\mathbb{S}_n]
\right)_{\Sy_b^{\textrm{op}}\times \Sy_a} ,$$ where the second
direct sum runs over the $b$-tuples $\ol$, $\ok$ and the
$a$-tuples $\oj$, $\oi$ such that  $|\ol|=m$, $|\ok|=|\oj|=N$,
$|\oi|=n$ and where the coinvariants correspond to the following
action of $\Sy_b^{\textrm{op}}\times \Sy_a$ :
\begin{eqnarray*}
&&\theta \otimes p_1\otimes \cdots \otimes p_b \otimes \sigma
\otimes q_1 \otimes \cdots \otimes q_a \otimes\omega
  \sim \\
&&\theta \,\tau^{-1}_\ol \otimes p_{\tau^{-1}(1)}\otimes \cdots
\otimes p_{\tau^{-1}(b)} \otimes
 \tau_\ok\,
 \sigma \, \nu_\oj \otimes q_{\nu(1)} \otimes \cdots \otimes q_{\nu(a)}
 \otimes \nu^{-1}_{\oi} \,
 \omega,
\end{eqnarray*}
for $\theta \in \Sy_m$, $\omega \in \Sy_n$, $\sigma \in \Sy_N$ and
for $\tau \in \mathbb{S}_b$ with $\tau_\ok$ the corresponding
block permutation, $\nu \in \mathbb{S}_a$ and $\nu_\oj$ the
corresponding block permutation. This product is associative but
has no unit. To fix this issue, we restrict to connected graphs.

The permutations of $\Sy_N$ associated to connected graphs are
called \emph{connected} (for more details see Section $1.3$ of
\cite{Vallette03}). We denote the set of connected permutations by
$\Sy^c$. We define the \emph{connected composition product} by the
following formula
$$\mathcal{P}\boxtimes_c \mathcal{Q}(m,\, n) :=
 \bigoplus_{N\in \mathbb{N}^*} \left( \bigoplus_{\ol,\, \ok,\, \oj,\, \oi} \KK[\mathbb{S}_m]
 \otimes_{\mathbb{S}_\ol}
\mathcal{P}(\ol,\, \ok)\otimes_{\mathbb{S}_{\ok}} \KK[\Sc_\kj]
\otimes_{\mathbb{S}_\oj} \mathcal{Q}(\oj,\, \oi)
\otimes_{\mathbb{S}_\oi} \KK[\mathbb{S}_n]
\right)_{\Sy_b^{\textrm{op}}\times \Sy_a}.$$ The unit $\II$ for
this monoidal product is given by $$\left\{
\begin{array}{l}
I(1,\, 1):=\KK, \quad \textrm{and}   \\
I(m,\,n):=0 \quad \textrm{otherwise}.
\end{array} \right.$$
We denote by $(\Sy\textrm{-biMod}, \boxtimes_c,\, I)$ this
monoidal category.\\

We define the \emph{concatenation product} of two bimodules $\Po$
and $\Qo$ by
$$ \Po\otimes \Qo (m,\,n) := \bigoplus_{m'+m''=m \atop n'+n''=n}
\KK[\Sy_{m'+m''}]\otimes_{\Sy_{m'}\times\Sy_{m''}} \Po(m',\,n')
\otimes_\KK \Qo(m'',\, n'') \otimes_{\Sy_{n'}\times
\Sy_{n''}}\KK[\Sy_{n'+n''}].$$
This product corresponds to taking
the (horizontal) tensor product of the elements of $\Po$ and $\Qo$
(see Figure $3$ of \cite{Vallette03} for an example). It is
symmetric, associative and unital. On the contrary to the two
previous products, it is linear on the left and on the right. We
denote by $\mathcal{S}_\otimes(\Po)$ the free symmetric monoid
generated by an $\Sy$-bimodule $\Po$ for the concatenation product
(and $\bar{\mathcal{S}}_\otimes(\Po)$ its augmentation ideal).
There is a natural embedding $\Po \bc_c \Qo \mono \Po \bc \Qo$.
And we obtain the composition product from the connected
composition product by concatenation, that is
$\bar{\mathcal{S}}_\otimes(\Po \boxtimes_c \Qo) \cong \Po
\boxtimes \Qo$. (From this relation, we can see that $I\bc
\Po=\bar{\mathcal{S}}_{\otimes}(\Po)$ and not $\Po$.)

\subsection{Properad}

A \emph{properad} is a monoid in the monoidal category
$(\Sy\textrm{-biMod}, \boxtimes_c,\ I)$. We denote the set of
morphisms of properads by $\Mor(\Po,\Qo)$. A properad $\Po$ is
\emph{augmented} if there exists a morphism of properads
$\varepsilon \, : \, \Po \to I$. We denote by $\oPo$ the kernel of
the augmentation $\varepsilon$ and call it the \emph{augmentation
ideal}. When $(\Po, \mu, \eta,\, \varepsilon)$ is an augmented
properad, $\Po$ is canonically isomorphic to $\II \oplus \oPo$. We
denote by $\big( I\oplus \underbrace{\oPo}_{r} \big) \boxtimes_c
\big( I\oplus \underbrace{\oPo}_{s} \big)$ the sub-$\Sy$-bimodule
of $\Po \boxtimes_c \Po$ generated by compositions of $s$
non-trivial elements of $\Po$ on the first level with $r$
non-trivial elements of $\Po$ on the second level. The
corresponding restriction of the composition product $\mu$ on this
sub-$\Sy$-bimodule is denoted $\mu_{(r,\,s)}$. The bilinear part
of $\Po\boxtimes_c \Po$ is the $\Sy$-bimodule $\big( I\oplus
\underbrace{\oPo}_{1} \big) \bc_c \big( I\oplus
\underbrace{\oPo}_{1} \big)$. It corresponds to the compositions
of only $2$ non-trivial operations of $\Po$. We denote it by
$\Po\bbc \Po$. The composition of two elements $p_1$ and $p_2$ of
$\oPo$ is written $p_1\bbc p_2$ to lighten the notations. The
restriction $\mu_{(1,\, 1)}$ of the composition product $\mu$ of a
properad $\Po$ on $\Po\bbc \Po$ is called the \emph{partial
product}.

A properad is called \emph{reduced} when the underlying
$\Sy$-bimodule is reduced, that is when  $\Po(m, n)=0$ for $n=0$
or $m=0$.

\subsection{Connected coproperad}
Dually, we defined the notion of \emph{coproperad}, which is a
comonoid in $(\Sy\textrm{-biMod}, \boxtimes_c)$. Recall that the
partial coproduct $\Delta_{(1,\, 1)}$ of a coproperad $\Co$ is the
projection of the coproduct $\Delta$ on $\Co \bbc \Co:= \big(
\II\oplus \underbrace{\Co}_{1} \big) \bc_c \big( \II \oplus
\underbrace{\Co}_{1} \big)$. More generally, one can define the
$(r, s)$-part of the coproduct by the projection of the image of
$\Delta$ on $\big( \II\oplus \underbrace{\Co}_{r} \big) \bc_c
\big(
\II\oplus \underbrace{\Co}_{s} \big)$.\\

Since the dual of the notion of coproduct is the notion of
product, we have to be careful with coproperad. For instance, the
target space of a morphism of coproperads is a direct sum of
modules and not a product. (The same problem appears at the level
of algebras and coalgebras). We generalize here the notion of
\emph{connected} coalgebra, which is the dual notion of Artin
rings,  introduced by D. Quillen in \cite{QuillenRational}
Appendix B, Section 3 (see also J.-L. Loday and M. Ronco
\cite{LodayRonco06} Section 1).

Let $(\Co,\, \Delta,\, \varepsilon,\, u)$ be an coaugmented (dg)
coproperad. Denote by $\oC:=\Ker(\Co \xrightarrow{\varepsilon} I)$
its \emph{augmentation}. We have $\Co = \oC \oplus I$ and
$\Delta(I)=I \boxtimes_c I$. For $X\in \oC$, denote by $\oD$ the
non-primitive part of the coproduct, that is
$\Delta(X)=I\boxtimes_c X + X \boxtimes_c I + \oD(X)$. The
\emph{coradical filtration} of $\Co$ is defined inductively as
follows
\begin{eqnarray*}
&& F_0:=\KK\,  I \\
&& F_r:= \{ X \in \Co \,  | \, \oD(X)\in F_{r-1} \boxtimes_c
F_{r-1}\}.
\end{eqnarray*}
\draftnote{$\Delta$ ou $\Delta_{(1,1)}$ ??} An augmented
coproperad is \emph{connected} if the coradical filtration is
exhaustive $\Co = \bigcup_{r\ge 0} F_r$. This condition implies
that $\Co$ is \emph{conilpotent} which means that for every $X\in
\Co$, there is an integer $n$ such that $\oD^n(X)=0$. This
assumption is always required to construct morphisms or
coderivations between coproperads (see next paragraph and
Lemma~\ref{uniquecoDer} for instance).\\

For the same reason, we will sometimes work with the invariant
version of the composition product denoted $\Po \bc^\Sy_c \Qo$
when working with coproperads. It is defined by the same formula
as for $\bc_c$ but where we consider the invariant elements under
the actions of the symmetric groups instead of the coinvariants
(see Lemma~\ref{convolution} for instance). When we want to
emphasize the difference between invariants and coinvariants, we
use the notations $\bc^\Sy$ and $\bc_{\Sy}$. Otherwise, we use
only $\bc $ since the two are isomorphic in characteristic $0$.

\subsection{Free properad and cofree connected coproperad}
\label{Free properad} Recall from \cite{Vallette04} the
construction of the free properad. Let $V$ be an $\Sy$-bimodule.
Denote by $V^+:=V\oplus I$ its augmentation and by
$V_n:=(V^+)^{\boxtimes_c n}$ the $n$-fold ``tensor" power of
$V^+$. This last module can be thought of as $n$-levelled graphs
with vertices indexed by $V$ and $I$. We define on $V_n$ an
equivalence relation $\sim$ by identifying two graphs when one is
obtained from the other by moving an operation from a level to an
upper or lower level. (Note that this permutation of the place of
the operations induces signs). We consider the quotient
$\widetilde{V}_n:=V_n/\sim$ by this relation. Finally, the free
properad $\F(V)$ is given by a particular colimit of the
$\widetilde{V}_n$. The dg $\Sy$-bimodule $\F(V)$ is generated by
graphs without levels with vertices indexed by elements of $V$. We
denote such graphs by $\G(v_1, \ldots, v_n)$, with $v_1, \ldots,
v_n \in V$. Since $\G(v_1, \ldots, v_n)$ represents an equivalence
class of levelled graphs, we can chose, up to signs, an order for
the vertices. (Any graph $\G$ with $n$ vertices is the quotient by
the relation $\sim$ of a graph with $n$ levels and only one
non-trivial vertex on each level). The composition product of
$\F(V)$ is given by the grafting. It is naturally graded by the
number of vertices. This grading is called the \emph{weight}. The
part of weight $n$ is denoted by
$\F(V)^{(n)}$.\\

Since we are working over a field of characteristic $0$, the
cofree connected coproperad on a dg $\Sy$-bimodule $V$ has the
same underlying space as the free properad, that is
$\F^c(V)=\F(V)$. The coproduct is given by pruning the graphs into
two parts. This coproperad verifies the universal property only
among connected coproperads (see Proposition $2.7$ of
\cite{Vallette03})


\subsection{Props}

We would like to define the notion of \emph{prop} as a monoid in
the category of $\Sy$-bimodules with the composition product
$\bc$. Since this last one has no unit and is not a monoidal
product, strictly speaking, we have to make this definition
explicit.

\begin{dei}[Prop]
A \emph{prop} $\Po$ is an $\Sy$-bimodule endowed with two maps
$\Po \bc \Po \xrightarrow{\mu} \Po$ and $I \xrightarrow{\eta} \Po$
such that the first is associative and the second one verifies
$$\xymatrix{ I\bc_c \Po\ \ar@{>->}[r] \ar[drr]^{\sim}  & I \bc \Po\ \ar@{>->}[r]^{\eta\bc \Po} &  \Po \bc \Po\ar[d]^{\mu} &
\ar@{>->}[l]_{\Po \bc \eta}\  \Po \bc I  & \ar@{>->}[l]\ \Po\bc_c I \ar[dll]_{\sim} \\
& & \Po. & & }$$
\end{dei}

This definition is equivalent to  the original definition of Adams and
MacLane \cite{AdamsBook78, MacLane65}. The original definition consists
of two coherent bilinear products, the vertical and horizontal
compositions of operations. The definition of the composition
product given here includes these two previous compositions at the
same time. The partial product $\Po \bbc \Po
\xrightarrow{\mu_{(1,1)}} \Po$ composes two operations. If they
are connected by at least one edge, this composition is the
vertical composition, otherwise this composition can be seen as
the horizontal composition of operations. This presentation will
allow us later to define the bar construction, resolutions and
minimal models for props.

It is straightforward to extend the results of the preceding
subsections to props. There exists notions of augmented props,
free prop, coprop and connected cofree coprop. We refer the reader
to \cite{Vallette03} Section $2$ for a complete treatment.


\subsection{(Co)triple interpretation}

The free prop(erad) functor induces a triple $\F \, :
\Sy-\textrm{biMod} \to \Sy-\textrm{biMod}$ such that an algebra
over it is a prop(erad) (see D. Borisov and Y.I. Manin
\cite{BorisovManin06}). When $(\Po, \mu)$ is a prop(erad), we will
denote by $\widetilde{\mu}_\Po \, : \, \F(\Po)^{(\ge 2)}\to \Po$
the associated map. Dually, the notion of coprop(erad) is
equivalent to the notion of coalgebra over the cotriple $\F^c \,
:\Sy-\textrm{biMod} \to \Sy-\textrm{biMod}$. When $(\Co, \Delta)$
is a coprop(erad), we will denote by $\widetilde{\Delta}_\Co \, :
\, \Co \to \F^c(\Co)^{(\ge 2)}$ the associated map.

\subsection{Non-symmetric prop(erad)}

In the sequel, we will have to work with algebraic structures endowed with operations having no symmetries.
One can model them with properads but the action of the symmetric group gives no real information.
Therefore, we introduction the notion of \emph{non-symmetric properad} which will be enough.
Since this structure is the direct generalization of
the notion of non-symmetric operad, we call it \emph{non-symmetric
properad}. All the definitions and propositions of this section
can be generalized directly to props. For simplicity, we only make
them explicit in the case of properads.

\begin{dei}\label{N-bimodule}
A \emph{(dg) $\NN$-bimodule} is a collection $\{\mathsf{P}(m, n)
\}_{m, \, n\in \NN^*}$ of dg modules.
\end{dei}

\begin{dei}[Non-symmetric connected composition product]
Let $\mathsf{P}$ and $\mathsf{Q}$ be two $\NN$-bimodules, we
define their \emph{non-symmetric connected composition product} by
the following formula
$$\mathsf{P}\boxtimes_c \mathsf{Q}(m,\, n) :=
 \bigoplus_{N\in \mathbb{N}^*} \left( \bigoplus_{\ol,\, \ok,\, \oj,\, \oi}
\mathsf{P}(\ol,\, \ok)\otimes \KK[\Sc_\kj] \otimes
\mathsf{Q}(\oj,\, \oi) \right)_{\Sy_b^{\textrm{op}}\times \Sy_a}
,$$ where the second direct sum runs over the $b$-tuples $\ol$,
$\ok$ and the $a$-tuples $\oj$, $\oi$ such that  $|\ol|=m$,
$|\ok|=|\oj|=N$, $|\oi|=n$ and where the coinvariants correspond
to the following action of $\Sy_b^{\textrm{op}}\times \Sy_a$ :
\begin{eqnarray*}
&&p_1\otimes \cdots \otimes p_b \otimes \sigma \otimes q_1 \otimes
\cdots \otimes q_a
  \sim p_{\tau^{-1}(1)}\otimes \cdots \otimes p_{\tau^{-1}(b)} \otimes
 \tau_\ok\,
 \sigma \, \nu_\oj \otimes q_{\nu(1)} \otimes \cdots \otimes q_{\nu(a)}
 ,
\end{eqnarray*}
for $\sigma \in \Sc_\kj$ and for $\tau \in \mathbb{S}_b$ with
$\tau_\ok$ the corresponding block permutation, $\nu \in
\mathbb{S}_a$ and $\nu_\oj$ the corresponding block permutation.
Since the context is obvious, we still denote it by $\boxtimes_c$.
\end{dei}

The definition of the monoidal product for $\Sy$-bimodule is based
on 2-levelled graphs with leaves, inputs and outputs labelled by
integers. This definition is based on non-labelled 2-levelled
graphs. We define the \emph{non-symmetric composition product}
$\bc$ by the same formula with the set of all permutations of
$\Sy_N$ instead of connected permutations.

\begin{pro}
The category $(\NN\textrm{-biMod}, \boxtimes_c,\, I)$ of
$\NN$-bimodules with the product $\boxtimes_c$ and the unit $I$ is
a monoidal category.
\end{pro}

\begin{proo}
The proof is similar to the one for $\Sy$-bimodules (see
\cite{Vallette03} Proposition $1.6$).
\end{proo}

\begin{dei}[Non-symmetric properad]
A \emph{non-symmetric properad} $(\mathsf{P}, \mu, \eta)$ is a
monoid in the monoidal category $(\NN\textrm{-biMod},
\boxtimes_c,\, I)$.
\end{dei}

\begin{Ex}
A non-symmetric properad $\mathsf{P}$ concentrated in arity $(1,
n)$, with $n\ge 1$, is the same as a non-symmetric operad.
\end{Ex}


\subsection{Representations of prop(erad)s, gebras}

Let $\Po$ and $\Qo$ be two prop(erad)s. A morphism $\Po
\xrightarrow{\Phi} \Qo$ of $\Sy$-bimodules is a \emph{morphism of
prop(erad)s} if it commutes with the products and the units of
$\Po$ and $\Qo$. In this case, we say that $\Qo$ is a
\emph{representation} of $\Po$.\\

We will be mainly interested in representations of the following
form.  Let $X$ be a dg vector space. We consider the
$\Sy$-bimodule $\End_X$ defined by
$\End_X(m,n):=\Hom_\KK(X^{\otimes n}, X^{\otimes m})$. The natural
composition of maps provides this $\Sy$-bimodule with a structure
of prop and properad. It is called the
\emph{endomorphism prop(erad)} of the space $X$.\\

Props and properads are meant to model the operations acting on
types of algebras or bialgebras in a generalized sense. When $\Po$
is a prop(erad), we call \emph{$\Po$-gebra} a dg vector space $X$
with a morphism of prop(erad)s $\Po \to \End_X$, that is a
representation of $\Po$ of the form $\End_X$. When $\Po$ is an
operad, a $\Po$-gebra is an algebra over $\Po$. To code operations
with multiple inputs and multiple outputs acting on an algebraic
structure, we cannot use operads  anymore and we need to use
prop(erad)s. The categories of (involutive) Lie bialgebras and (involutive) Frobenius bialgebras  are categories of gebras over a properad
(see Section~\ref{examples}). The categories of (classical) associative
bialgebras and infinitesimal Hopf algebras (see \cite{Aguiar00}) are governed by non-symmetric properads.
In these cases, the associated prop is
freely obtained from a properad. Therefore, the prop does not
model more relations than the properad and the two categories of
gebras over the prop and the properad are equal. For more details,
see the beginning of Section~\ref{ComparisonModels}.


\section{Convolution prop(erad)}\label{Convolution prop(erad)}

When $A$ is an associative algebra and $C$ a coassociative
coalgebra, the space of morphisms $\Hom_\KK(C, \, A)$ from $C$ to
$A$ is naturally an associative algebra with the convolution
product. We generalize this property to prop(erad)s, that is the
space of morphisms  $\Hom(\Co, \, \Po)$
from a coprop(erad) $\Co$ and a prop(erad) $\Po$ is a prop(erad).
From this rich structure, we get general operations, the main one being the \emph{intrinsic} Lie bracket used to
study the deformation theory of algebraic structures later in
\ref{Deformation theory}.

\subsection{Convolution prop(erad)}

For two $\Sy$-bimodules $M=\{M(m, n)\}_{m,n}$ and
$N=\{N(m,n)\}_{m,n}$, we denote by $\Hom(M, N)$ the collection $\{
\Hom_\KK\big(M(m,n), N(m,n) \big)\}_{m,n}$ of morphisms of
$\KK$-modules. It is an $\Sy$-bimodule with the action by
conjugation, that is
$$(\sigma. f. \tau)(x) := \sigma.(f(\sigma^{-1}.x.\tau^{-1})).\tau, $$
for $\sigma \in \Sy_m$, $\tau \in \Sy_n$ and $f \in\Hom\big(M,
N\big)(m,n)$. An invariant element under this action is an
equivariant map from $M$ to $N$, that is
$\Hom(M, N)^\Sy=\Hom^\Sy(M, N) $.   \\

When $\Co$ is a coassociative coalgebra and $\Po$ is an
associative algebra, $\Hom(\Co, \Po)$ is an associative algebra
 known as the  \emph{convolution algebra}. When $\Co$ is a cooperad and $\Po$ is an operad, $\Hom(\Co, \Po)$
is an operad called the \emph{convolution operad} by C. Berger and
I. Moerdijk in \cite{BergerMoerdijk03} Section $1$. We extend this
construction to properads and props.

\begin{lem}\label{convolution}
Let $\Co$ be a coprop(erad) and $\Po$ be a prop(erad). The space
of morphisms $\Hom(\Co, \Po)=\Po^\Co$ is a prop(erad).
\end{lem}

\begin{proo}
We use the notations of Section~\ref{S-bimodule} (see also Section
1.2 of \cite{Vallette03}). We define an associative and unital map
$\mu_{\Po^\Co} \, : \, \Po^\Co \boxtimes_\Sy \Po^\Co \to \Po^\Co$
as follows. Let $\G^2(f_1, \ldots, f_r ;\, g_1, \ldots , g_s)\in
\Po^\Co\boxtimes \Po^\Co (m,n)$ be a 2-levelled graph whose
vertices of the first level are labelled by $f_1, \ldots, f_r$ and
whose vertices of the second level are labelled by $g_1, \ldots,
g_s$. The image of $\G^2(f_1, \ldots, f_r ;\, g_1, \ldots , g_s)$
under $\mu_{\Po^\Co}$ is the composite
$$\Co \xrightarrow{\Delta_\Co} \Co \boxtimes^\Sy \Co \mono \Co \boxtimes \Co
\xrightarrow{\widetilde{\G}^2(f_1, \ldots, f_r ;\, g_1, \ldots ,
g_s)} \Po \boxtimes \Po \epi \Po \boxtimes_\Sy \Po
\xrightarrow{\mu_\Po} \Po,$$ where $\widetilde{\G}^2(f_1, \ldots,
f_r ;\, g_1, \ldots , g_s)$ applies $f_i$ on the element of $\Co$
indexing the $i^{\textrm{th}}$ vertex of the first level and $g_j$
on the element of $\Co$ indexing the  $j^{\textrm{th}}$ vertex of
the second level of an element of $\Co \boxtimes \Co$. Since the
action of the symmetric groups on $\Po^\Co$ is defined by
conjugation and since the image of the coproduct lives in the
space of invariants, this map factors through the coinvariants,
that is $\Po^\Co \boxtimes_\Sy \Po^\Co \to \Po^\Co$.

The unit is given by the map $\Co \xrightarrow{\varepsilon} I
\xrightarrow{\eta} \Po$. The associativity of $\mu_{\Po^\Co}$
comes directly from the coassociativity of $\Delta_{\Co}$ and the
associativity of $\mu_\Po$.
\end{proo}

\begin{dei}
The properad $\Hom(\Co, \Po)$ is called the \emph{convolution
prop(erad)} and is denoted by $\Po^\Co$.
\end{dei}

Assume now that $(\Co, d_\Co)$ is a dg coprop(erad) and $(\Po,
d_\Po)$ is a dg prop(erad). The \emph{derivative} of a graded
linear map $f$ from $\Co$ to $\Po$ is defined as follows
$$D(f):= d_\Po \circ f - (-1)^{|f|}f \circ d_\Co .$$

A $0$-cycle for this differential is a morphism of chain
complexes, that is it commutes with the differentials. In Section
\ref{Derivative as dg affine scheme}, we give a geometric
interpretation of this derivative. The derivative is a derivation
for the product of the prop(erad) $\Hom(\Co, \Po)$ that verifies
$D^2=0$. We sum up these relations in the following proposition. The same result holds
in the non-symmetric case.

\begin{pro}\label{dg convolution prop(erad)}
When $(\Co, d_\Co)$ is a dg coprop(erad) and $(\Po,  d_\Po)$ is a
dg prop(erad), $(\Hom(\Co,\Po), D)$ is a dg prop(erad).

When $(\Co, d_\Co)$ is a dg non-symmetric coprop(erad) and $(\Po,  d_\Po)$ is a
dg non-symmetric prop(erad), $(\Hom(\Co,\Po), D)$ is a dg non-symmetric prop(erad).
\end{pro}


\subsection{Lie-admissible products and Lie brackets associated to a properad}

In \cite{KapranovManin01}, the authors proved that the total
space $\oplus_n \Po(n)$, as well as the space of coinvariants $\oplus_n \Po(n)_{\Sy_n}$, of an operad is endowed with a natural Lie
bracket. This Lie bracket is the anti-symmetrization of the preLie
product $p\circ q = \sum_i p \circ_i q$ defined by the sum on all
possible ways of composing two operations $p$ and $q$. Notice that this result was implicitly stated by Gerstenhaber  in  \cite{Gerstenhaber63}. We
generalize this results to properads.\\

For any pair of elements, $\mu$ and $\nu$, in a (non-symmetric)
properad $\Po$, we denote by $\mu\circ\nu$ the sum of all the
possible compositions of $\mu$ by $\nu$ along any $2$-levelled
graph with two vertices in $\Po$. For  another element $\eta$ in
$\Po$, the composition $(\mu \oc \nu )\oc \eta$ is spanned by
graphs  of the form

$$\vcenter
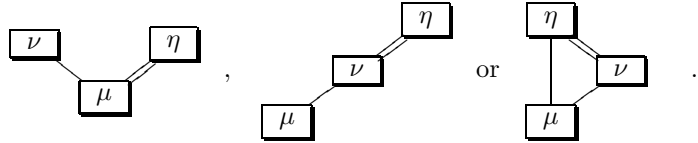
{\xymatrix@R=8pt@C=8pt{*+[F-,]{\ \nu \ } \ar@{-}[dr] &&
*+[F-,]{\ \eta \ } \ar@{-}@<0.5ex>[dl] \ar@{-}[dl] \\
& *+[F-,]{\ \mu \ } & }} \quad , \quad
\vcenter{\xymatrix@R=8pt@C=8pt{ && *+[F-,]{\ \eta \ } \ar@{-}[dl]  \ar@{-}@<0.5ex>[dl]\\
 & *+[F-,]{\ \nu \ } \ar@{-}[dl]  \\
 *+[F-,]{\ \mu \ } & & }} \quad \textrm{or} \quad
\vcenter{\xymatrix@R=8pt@C=8pt{*+[F-,]{\ \eta \ } \ar@{-}[dr]  \ar@{-}@<-0.5ex>[dr] \ar@{-}[dd] && \\
 & *+[F-,]{\ \nu \ } \ar@{-}[dl]  \\
 *+[F-,]{\ \mu \ } & & }}\ . $$
Let us denote by $\mu \circ (\nu,\, \eta)$  the summand spanned by graphs the first type.

In the same way,  $\mu \oc (\nu \oc \eta)$ is spanned by graphs  of
the form
$$\vcenter
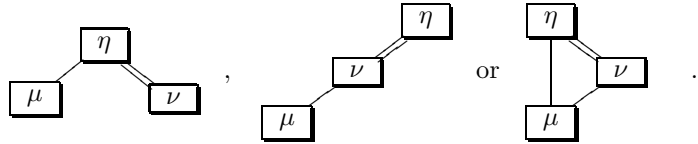
{\xymatrix@R=8pt@C=8pt{& *+[F-,]{\ \eta \ } &  \\
*+[F-,]{\ \mu \ } \ar@{-}[ur] && *+[F-,]{\ \nu \ }
\ar@{-}@<0.5ex>[ul] \ar@{-}[ul]}} \quad , \quad
\vcenter{\xymatrix@R=8pt@C=8pt{ && *+[F-,]{\ \eta \ } \ar@{-}[dl]  \ar@{-}@<0.5ex>[dl]\\
 & *+[F-,]{\ \nu \ } \ar@{-}[dl]  \\
 *+[F-,]{\ \mu \ } & & }} \quad \textrm{or} \quad
\vcenter{\xymatrix@R=8pt@C=8pt{*+[F-,]{\ \eta \ } \ar@{-}[dr]  \ar@{-}@<-0.5ex>[dr] \ar@{-}[dd] && \\
 & *+[F-,]{\ \nu \ } \ar@{-}[dl]  \\
 *+[F-,]{\ \mu \ } & & }} \ . $$
and we denote by $(\mu,\, \nu)\circ \eta$ the summand of  $\mu \oc (\nu \oc \eta)$
spanned by graphs of the first (from the left) type. With
these notations, we have in $\Po$ the following formula
$$(\mu \circ \nu) \circ \eta - \mu \circ (\nu \circ \eta )=
\mu \circ (\nu, \, \eta) -(\mu,\, \nu)\circ\eta.$$

When $\Po=A$ is concentrated in arity $(1,1)$, it is an
associative algebra. In this case, the product $\circ$ is the
associative product of $A$. When $\Po$ is an operad, the operation
$(\mu,\, \nu)\circ\eta$ vanishes and the product $\mu \circ \nu$
is right symmetric, that is $(\mu \circ \nu) \circ \eta - \mu
\circ (\nu \circ \eta )=(\mu \circ \eta) \circ \nu - \mu \circ
(\eta \circ \nu)$. Such a product is called \emph{preLie}. In the
general case of properads, this product verifies a weaker relation
called $\emph{Lie-admissible}$ because its anti-symmetrized
bracket verifies the Jacobi identity. Denote by $\As(\mu,\, \nu,\,
\eta):=(\mu \oc \nu)\oc \eta -\mu\oc (\nu\oc \eta)$ the associator
of $\mu$, $\nu$ and $\eta$.

\begin{dei}[Lie-admissible algebra]
A graded vector space $A$ with a binary product $\oc$ is called a
(graded) \emph{Lie-admissible} algebra if one has $\displaystyle
\sum_{\sigma \in \Sy_3}\sgn(\sigma) \As(- ,\, -,\, -)^\sigma=0$,
where, for instance, $\As(-,\, -,\, -)^{(23)}$ applied to $\mu$,
$\nu$ and $\eta$ is equal to $(-1)^{|\nu||\eta|}((\mu\oc \eta) \oc
\nu - \mu\oc (\eta\oc \nu))$. A \emph{differential graded Lie-admissible algebra} (or dg Lie-admissible algebra for short) is a dg module $(A, d_A)$ endowed with a Lie-admissible product $\oc$ such that the $d_A$ is a derivation.
\end{dei}

\begin{pro}\label{properad->Lie-admissible}
Let $\Po$ be a dg properad or a non-symmetric dg properad, the
space $\bigoplus_{m,n} \Po(m,n)$, endowed with the product $\oc$,
is a dg Lie-admissible algebra.
\end{pro}

\begin{proo}
Let $H=\{id,\, (23)\}$ and $K=\{id,\, (12) \}$ be two subgroups of
$\Sy_3$. We have
\begin{eqnarray*}
\sum_{\sigma \in \Sy_3} \sgn(\sigma) \As(-,\,-,\, -)^\sigma &=&
\sum_{\sigma \in \Sy_3} \sgn(\sigma) \big( (-\circ\, (-\circ\,
-))^\sigma -  ((-\circ \, -)
\circ\,-)^\sigma\big) \\
&=& \sum_{\tau H \in \Sy_3\backslash H} \sgn(\tau)
\underbrace{\left(  (-\circ\,  ( -,\, -))^\tau -  (-  \circ\,(
-,\,
-))^{\tau (23)} \right)}_{=0} - \\
&& \sum_{\omega K \in \Sy_3\backslash K} \sgn(\omega)
\underbrace{\left( (( -, \, -)\circ\, -)^\omega - ((-,\, -)
\circ\, -)^{\omega (12)} \right)}_{=0} \\
&=& 0.
\end{eqnarray*}
\end{proo}

Actually on the direct sum $\bigoplus_{m,n} \Po(m,n)$ of the
components of a properad, there are higher operations with $r+s$
inputs which turns it into a ``non-differential
$B_\infty$-algebra''. We refer to the next section for more details.\\

For a prop $\Po$, we still define the product $p\circ q$ on
$\bigoplus_{m,n} \Po(m,n)$ by all the possible ways of composing
the operations $p$ and $q$, that is all vertical composites and
the horizontal one.

\begin{pro}\label{prop->Associative algebra}
\footnote{This result was mentioned to the second author by M.M.
Kapranov (long time ago).} Let $\Po$ be a dg prop or a
non-symmetric dg prop, the space $\bigoplus_{m,n} \Po(m,n)$,
endowed with the product $\oc$, is a dg associative algebra.
\end{pro}

\begin{proo}
We denote by $p\, \circ_v q$ the sum of all vertical (connected)
composites of $p$ and $q$ and by $p \circ_h q$ the horizontal
composite. We continue to use the notation $p\circ_v (q, r)$ to
represent the composite of an operation $p$ connected to two
operations $q$ and $r$ above. We have (in degree $0$)
\begin{eqnarray*}
 && (p\circ q)\circ r = (p \circ_v q + p \circ_h q)\circ r= \\ &&
p \circ_v q \circ_v r + p \circ_v (q, r) + (p \circ_v q) \circ_h r
+ (p \circ_v r) \circ_h q +p \circ_h (q \circ_v r) + (p, q)
\circ_v r + p \circ_h q \circ_h r,
\end{eqnarray*}
and
\begin{eqnarray*}
 && p\circ (q\circ r) = p\circ (q\circ_v r + q\circ_h r ) = \\ &&
p \circ_v q \circ_v r + (p, q)\circ_v r +  p \circ_h (q \circ_v r)
+ (p \circ_v q) \circ_h r + q \circ_h (p \circ_v r) + p \circ_v
(q, r) + p \circ_h q \circ_h r.
\end{eqnarray*}
Since the horizontal product is commutative, $(p \circ_v r)
\circ_h q$ is equal to $q \circ_h (p \circ_v r)$, which finally
implies $(p\circ q)\circ r = p\circ (q\circ r) $.
\end{proo}

These structures pass to coinvariants $\bigoplus \Po_\Sy:=\bigoplus_{m,n} \Po(m,n)_{\Sy_m^{\textrm{op}}\times \Sy_n}$ as follows.

\begin{pro}\label{OnCoinvariants}
Let $\Po$ be a dg properad (respectively dg prop), the dg Lie-admissible (associative) product $\circ$ on
$\bigoplus \Po$ induces a dg Lie-admissible (associative) product on the space of coinvariants
 $\bigoplus \Po_\Sy$.
\end{pro}

\begin{proo}
It is enough to prove that the space $C:=\lbrace p-\tau.p.\nu\, ;\, p\in   \Po(m,n),\, \tau \in \Sy_m,\, \nu \in\Sy_n   \rbrace $ is a two-sided ideal for the Lie-admissible product $\circ$. Let us denote $p\circ q$ by
$\sum_{\sigma} \mu(p, \, \sigma,\, q)$, where $\mu$ is the composition map of the properad $\Po$ and where $\sigma$ runs thought connected permutations. For any $\tau \in \Sy_m$, we have
$$(p-\tau.p)\circ q=\sum_\sigma \big(\mu(p, \, \sigma,\, q) -\mu(\tau.p, \, \sigma,\, q)\big)=
\sum_\sigma \big(\mu(p, \, \sigma,\, q) -\tau_\sigma.\mu(p, \, \sigma,\, q)\big)\in C
, $$
where $\tau_\sigma$ is a permutation which depends on $\sigma$. For any $\nu \in \Sy_n$, we have
\begin{eqnarray*}
(p-p.\nu)\circ q&=&\sum_\sigma \mu(p, \, \sigma,\, q) -\sum_\sigma\mu(p, \, \nu.\sigma,\, q)
=\sum_\sigma \mu(p, \, \sigma,\, q) -\sum_\sigma\mu(p, \, \sigma',\, q).\nu_\sigma \\
&=&\sum_\sigma \big(\mu(p, \, \sigma,\, q)- \mu(p, \, \sigma,\, q).\nu_{\sigma''}\big)\in C,\end{eqnarray*}
since the connected permutations $\sigma'$ obtained runs thought the same set of connected permutations as $\sigma$.
Therefore, $C$ is a right ideal. The same arguments prove that $C$ is a left ideal.
\end{proo}

In the sequel, we will have to work with the space of invariants $\bigoplus \Po^\Sy:=\bigoplus_{m,n} \Po(m,n)^{\Sy_m^{\textrm{op}}\times \Sy_n}$ ,
and not coinvariants, of a properad. Since we work over a field of
characteristic zero, both are canonically isomorphic. Let $V$ be a
vector space with an action of a finite group $G$. The subspace of
invariants is defined by $V^G:=\{v\in V \, |\,  v.g=v \ ,\,
\forall g\in G \}$ and the quotient space of coinvariants is
defined by $V_G:=V/<v-v.g ,\, \forall (v,g) \in V\times G > $. The
map from $V^G$ to $V_G$ is the composite of the inclusion $V^G
\mono V$ followed by the projection $V \epi V_G$. The inverse map
$V_G \to V^G $ is given by $\displaystyle [v]\mapsto
\frac{1}{|G|}\sum_{g\in G } v.g$, where $[v]$ denotes the class of
$v$ in $V_G$.

\begin{cor}\label{OnInvariants}
Let $\Po$ be a dg properad (respectively dg prop), its total space of invariant elements
$\bigoplus \Po^\Sy$ is a dg Lie-admissible algebra (dg associative algebra).
\end{cor}


The Lie-admissible relation of a product $\oc$ is equivalent to
the Jacobi identity $[[-,\, -],\, -]+[[-,\,
-],\, -]^{(123)} +[[-,\,
-],\, -]^{(132)}=0$ for its induced bracket $[\mu,\, \nu]:=\mu\oc
\nu-(-1)^{|\mu||\nu|}\nu\oc \mu$.

\begin{thm}\label{LieBracket}
Let $\Po$ be a dg properad (respectively dg prop), its total space $\bigoplus \Po$, the total space of coinvariant elements $\bigoplus \Po_\Sy$ and the total space of invariant elements $\bigoplus \Po^\Sy$ are dg Lie algebras.
\end{thm}

The first of this statement is also true for non-symmetric dg prop(erad)s.

\subsection{LR-algebra associated to a properad}\label{LR-algebra}
On the total space of a properad, we have constructed a binary product $\circ$ in the previous section. We now define
more general operations with multiple inputs.

\begin{dei}[LR-operations]
Let $(\Po,\, \mu)$ be a properad and  $p_1, \ldots ,\,  p_r$ and $q_1,\ldots ,\, q_s$ be elements of
$\Po$. Their \emph{LR-operation} $\{p_1, \ldots ,\,
p_r\}\{q_1,\ldots ,\, q_s\}$ is defined by
$$\sum_{\sigma}\mu(p_1,\ldots ,\,  p_r;\, \sigma;\, q_1,\ldots ,\, q_s),$$
where $\sigma$ runs through connected permutations.
\end{dei}

In order words, the LR-product is the sum over all possible ways to compose the elements of $\Po$.

These operations are obviously graded symmetric with respect to
Koszul-Quillen sign convention, that is
$$\{p_1, \ldots, p_r\}\{q_1, \ldots, q_s\}=\varepsilon(\sigma, p_1, \ldots, p_r).\varepsilon(\tau, q_1, \ldots, q_s)
\{p_{\sigma(1)}, \ldots, p_{\sigma(r)}\} \{q_{\tau(1)}, \ldots,
q_{\tau(s)}\},$$ for $\sigma\in \Sy_r $ and $\tau\in \Sy_s  $. The element $\varepsilon(\sigma, p_1,
\ldots, p_r)\in\lbrace +1, -1\rbrace$ stands for the
Koszul-Quillen signs induced by the permutations of the graded
elements $p_1, \ldots, p_r$ under $\sigma$.
Notice that the Lie-admissible product is equal to $p\circ q :=
\{p\}\{q\}$. By convention, we set $\{\, \} \{\, \}=0$, $\{\, \}
\{q\}=q$, $\{p \} \{\, \}=p$ and $\{\, \} \{q_1,\ldots, q_s\}=0$
for $s>1$, $\{p_1, \ldots, p_r \} \{\, \}=0$ for $r>1$. The name
\emph{LR-operations} stands for Left-Right operations as well as
for Loday-Ronco operations since such operations are used in
\cite{LodayRonco06} to extend Cartier-Milnor-Moore Theorem to
non-cocommutative Hopf algebras.

\begin{pro}
The LR-operations satisfy the relations of a ``non-differential
$B_\infty$-algebra'', that is, for all $o_1, \ldots, o_r, p_1,
\ldots, p_s, q_1,\ldots, q_t$ in $\Po$.
\begin{eqnarray*} \sum_{\Theta} \varepsilon  \big\lbrace
\{o_{1},\ldots ,o_{i_1}\} \{ p_{1}, \ldots , p_{j_1}\}, \ldots ,
\{o_{i_1+\cdots+i_{a-1}+1},\ldots ,o_{r}\} \{
p_{j_1+\cdots+j_{a-1}+1}, \ldots , p_{r}\} \big\rbrace \big\lbrace
q_1, \ldots, q_t\big\rbrace \\ = \sum_{\Theta'} \varepsilon'
\big\lbrace o_1,\ldots, o_s \big\rbrace \big\lbrace
\{p_{1},\ldots, p_{k_1}\} \{ q_{1}, \ldots , q_{l_1}\}, \ldots ,
\{p_{k_1+\cdots+k_{b-1}+1},\ldots p_{s}\} \{
q_{l_1+\cdots+l_{b-1}+1}, \ldots , q_{t}\} \big\rbrace,
\end{eqnarray*}
where $\Theta$ runs over $1 \leqslant a \leqslant Max(r,s)$, $i_1,
\ldots, i_a \geqslant 0$ such that $i_1+ \cdots +i_a=r$,  $j_1,
\ldots, j_a \geqslant 0$ such that $j_1+ \cdots +j_a=s$ and where
$\Theta'$ runs over $1 \leqslant b \leqslant Max(s,t)$, $k_1,
\ldots, k_b \geqslant 0$ such that $k_1+ \cdots +k_b=s$,  $l_1,
\ldots, l_b \geqslant 0$ such that $l_1+ \cdots +l_b=t$. The sign
$\varepsilon$ comes from the permutations of the $o$ and the $p$
and the sign $\varepsilon'$ comes from the permutations of the $p$
and the $q$.
\end{pro}

\begin{proo}
It is a direct translation to LR-operations of the associativity of the operad $\Po$. See also,  Example 1.7 (d) of
\cite{LodayRonco06} and Lemma~\ref{convolution}.
\end{proo}

Therefore, the total space $\oplus \Po$ of a properad $\Po$, with the LR-operations, forms a
``non-differential $B_\infty$'', structure that we call a
\emph{LR-algebra}.
The same result also holds for non-symmetric prop(erad)s.

\begin{pro}\label{LRalgebra}
The $\Po$ be a dg prop(erad), its total space $\oplus \Po$, the total space of coinvariants elements $\oplus \Po_\Sy$ and the total space of invariants elements $\oplus \Po^\Sy$ form a LR-algebra.
\end{pro}

\begin{proo}\label{LR on CoIn}
The structure of LR-algebra of $\oplus \Po$ factors through the coinvariant elements $\oplus \Po_\Sy$ by the same arguments as in Proposition~\ref{OnCoinvariants}. Since the space of coinvariant and invariant elements are isomorphic, we can transfer this structure to invariant elements.
\end{proo}

\subsection{Lie-admissible bracket and LR-algebra of a convolution
properad}\label{operations on convolution properad} Since $\Hom(\Co,\Po)$ is an properad, it has a Lie-admissible bracket and more generally it enjoys a structure of LR-algebra by the preceding sections. We make these structures explicit here. We will use them later on
in our study of deformation theory (see Section~\ref{Deformation
theory}).

\begin{dei}[Convolution product]
Let $f$ and $g$ be two elements of $\Hom(\Co,\, \Po)$. Their
\emph{convolution product} $f \cp g$ is defined by the following
composite
$$\xymatrix@C=40pt{{\Co} \ar[r]^(0.45){\Delta_{(1,\, 1)}} &
{\Co \bbc \Co} \ar[r]^(0.45){f\bbc g} & {\Po \bbc \Po}
\ar[r]^(0.55){\mu} & {\Po.}}$$
\end{dei}

Since the partial coproduct of a coproperad (or a cooperad) is not
coassociative in general, the convolution product is not
associative.

\begin{pro}\label{Hom Lie-admissible}
Let $\Po$ be a dg prop(erad) and $\Co$ be a dg coprop(erad). The
convolution product $\star$ on $\oplus \Hom(\Co,\, \Po)$ is equal to
the  product $\circ$ associated to the  convolution dg
prop(erad). In the case of dg (co)properads, it is dg
Lie-admissible and for dg (co)props, it is dg associative.

This convolution product is stable on the space of invariant elements $\oplus  \Hom^\Sy(\Co,\, \Po)$ with respect to the action of the symmetric groups.
\end{pro}

\begin{proo}
The image of the map $\Delta_{(1,1)}$ is a sum over all possible
$2$-levelled graphs with two vertices indexed by some elements of
$\Co$. Therefore, the map $\star$ is equal to the sum of all
possible compositions of $f$ and $g$.

Saying that $f$ and $g$ are invariant elements in $\Hom(\Co,\, \Po)$ means that they are equivariant maps. Since
the composition map $\mu$ of $\Po$ and the partial coproduct $\Delta_{(1,\,1)}$ are also equivariant maps, we have
\begin{eqnarray*}
(\sigma. f\star g. \tau)(c) &=& \sigma.(f\star g(\sigma^{-1}.c.\tau^{-1})).\tau=
\sigma.(\mu\circ (f\otimes g)\circ\Delta_{(1,\,1)} (\sigma^{-1}.c.\tau^{-1})).\tau\\
&=&\sigma.\sigma^{-1}.( f\star g)(c).\tau^{-1}. \tau=f\star g(c).
\end{eqnarray*}
\end{proo}

Using the projections $\Delta_{(r,\, s)}$ of the coproduct of $\Co$, we
make explicit the LR-operations with $r$ and $s$ inputs of $\Hom(\Co,\, \Po)$ as follows.

\begin{pro}
Let $f_1, \ldots ,\,  f_r$ and $g_1,\ldots ,\, g_s$ be elements of
$\Hom(\Co,\, \Po)$. Their \emph{LR-operation} $\{f_1, \ldots ,\,
f_r\}\{g_1,\ldots ,\, g_s\}$ is equal to
\begin{eqnarray*}
 \Co \xrightarrow{\Delta_{(r,\, s)}} &{\big( \II\oplus
\underbrace{\Co}_{r} \big) \bc^\Sy \big( \II\oplus
\underbrace{\Co}_{s}
\big)}&\mono   \\
&\underbrace{\Co\otimes\cdots \otimes\Co}_{r} \bc
\underbrace{\Co\otimes\cdots \otimes\Co}_{s}& \xrightarrow{{\{f_1,
\ldots, f_r\} \bc \{g_1,\ldots, g_r\}}} {\Po \bc \Po} \epi {\Po \bc_\Sy \Po}
\xrightarrow{\mu}\Po,
\end{eqnarray*}
where $\{f_1, \ldots, f_r\}=\sum_{\sigma \in \Sy_r}
\varepsilon(\sigma, f_1, \ldots, f_r)f_{\sigma(1)} \otimes  \cdots
\otimes f_{\sigma(r)}$. The element $\varepsilon(\sigma, f_1,
\ldots, f_r)\in\lbrace +1, -1\rbrace$ stands for the
Koszul-Quillen signs induced by the permutations of the graded
elements $f_1, \ldots, f_r$ under $\sigma$. This means that we
apply $\{f_1, \ldots, f_r\}$ and $\{g_1, \ldots, g_s\}$ everywhere
we can.
\end{pro}

\begin{proo}
The proof is similar to the previous one.
\end{proo}

\begin{thm}\label{Convolution LR and Lie algebra}
Let $\Co$ be a dg coprop(erad) and $\Po$ be a dg prop(erad), the space $\oplus  \Hom(\Co, \, \Po)$ is a dg LR-algebra and thus a dg Lie algebra, structures that are stable on the space of equivariant maps $\oplus  \Hom^\Sy(\Co,\, \Po)$.
\end{thm}

\begin{proo}
Since the $\Delta_{(r,\, s)}$ and $\mu$ are equivariant maps, the LR-operations are stable on the space of equivariant maps $\Hom^\Sy(\Co,\, \Po)$ by their explicit form given in the previous proposition.
\end{proo}

\begin{Rq}
In the case of the convolution properad, we do not have to transfer the structure of LR-algebra or Lie algebra from $\Hom(\Co, \, \Po)$ to $\Hom^\Sy(\Co,\, \Po)$ through the coinvariant-invariant isomorphism. These structures on directly stable on the space of invariant elements.
\end{Rq}

When $\Co=C$ is a coassociative coalgebra and $\Po=A$ an
associative algebra, the product $\bc$ is equal to $\otimes$ and
is bilinear. In this case, the partial coproduct of $C$ is equal
to the coproduct of $C$ and is coassociative. (All the
$\Delta_{(r,\, s)}$ are null for $r>1$ or $s>1$). In this case,
the product $\star$ is the classical convolution product on
$\Hom(C,\, A)$, which is associative.\\

When $\Co$ is a cooperad and $\Po$ is an operad we have
$\Delta_{(r,\, s)}=0$ for $r>1$ as the operations  $\{f_1, \ldots
,\, f_r\}\{g_1,\ldots ,\, g_s\}$ are null unless $r=1$. The
remaining operations $\{f \}\{g_1,\ldots ,\, g_s\}$ are graded
symmetric \emph{brace} operations coming from the brace-type
relations verified by the operadic product (see  \cite{GuinOudom04, LadaMarkl05}). Remark
that when $\Co$ is a non-symmetric cooperad and $\Po$ a
non-symmetric operad, we can define non-symmetric braces on
$\Hom(\Co,\Po)$ without the sum over all permutations. In this
case, we find the classical non-symmetric braces of
\cite{Gerstenhaber63}, see also
\cite{GerstenhaberVoronov95,Vallette06Preprint}. The convolution
product verifies the relation  $(f\cp g)\cp h -f\cp (g\cp
h)=\{f\}\{g,\, h\}$. Therefore, in the operadic case, the (graded)
symmetry of the brace products implies that the associator $(f\cp
g)\cp h-f\cp (g\cp h)$ is symmetric in $g$ and $h$. Hence the
convolution product $\cp $ on $\Hom(\Co,\, \Po)$ is a graded
preLie product. For an interpretation of the LR-operations
 (or braces operations) on cohomology theories, we refer the reader to Section~\ref{deformation theory}.


\section{Bar and cobar construction}

In this section, we recall the definitions of the bar and cobar
constructions for (co)properads and extend it to (co)props. We
prove adjunction between these two constructions using the notion
of \emph{twisting morphism}, that is Maurer-Cartan elements in the
convolution prop(erad). Finally, we show that the bar-cobar
construction provides us with a canonical cofibrant resolution.

\subsection{Infinitesimal bimodule over a prop(erad)}\label{Infinitesimal bimodule}

The notion of bimodule $M$ over a prop(erad) $\Po$, defined in a
categorical way, is given by two compatible actions, left $\Po\bc
M \to M$ and right $M \bc \Po \to M$. For our purposes, we need a
\emph{linearized or infinitesimal} version of bimodules. Such a
phenomenon cannot be seen on the level of associative algebras
since the monoidal product $\otimes$ defining them is bilinear.

The $\Sy$-bimodule $(M\oplus N)\bc O$ is the sum over $2$-levelled graphs with vertices on the top level labelled by elements of $O$ and with vertices on the bottom level labelled by elements of $M$ or $N$. We denote by $(\underbrace{M}_{r}\oplus N)\bc O$ the sub-$\Sy$-module of $(M\oplus N)\bc O$ with exactly $r$ elements of $M$ on the bottom level.

\begin{dei}[Infinitesimal bimodule]
Let $(\Po, \mu)$ be a prop(erad). An \emph{infinitesimal
$\Po$-bimodule} is an $\Sy$-bimodule $M$ endowed with two action
maps of degree $0$
$$\lambda \ : \ \Po \bc (\Po \oplus \underbrace{M}_{1}) \to M
\quad \textrm{and} \quad \rho \ : \ (\Po \oplus
\underbrace{M}_{1})\bc \Po \to M,$$ such that the following
diagrams commute \begin{itemize} \item Compatibility between the
left action $\lambda$ and the composition product $\mu$ of $\Po$~:
$$\xymatrix@C=50pt@R=30pt{\Po \bc  \Po \bc (\Po \oplus
\underbrace{M}_{1}) \ar[r]^{\Po \bc (\lambda + \mu)} \ar[d]^{\mu
\bc (\Po \oplus M)} &
\ar[d]^{\lambda} \Po \bc  (\Po \oplus \underbrace{M}_{1})  \\
\Po \bc (\Po \oplus \underbrace{M}_{1}) \ar[r]^{\lambda} & M,}$$

\item Compatibility between the right action $\rho$ and the
composition product $\mu$ of $\Po$~:
$$ \xymatrix@C=50pt@R=30pt{(\Po
\oplus \underbrace{M}_{1}) \bc \Po \bc  \Po
 \ar[r]^{(\rho + \mu) \bc \Po  } \ar[d]^{ (\Po \oplus
M) \bc \mu} & \ar[d]^{\rho} (\Po \oplus \underbrace{M}_{1})\bc
\Po    \\
(\Po \oplus \underbrace{M}_{1})\bc \Po  \ar[r]^{\rho} & M,}$$

\item Compatibility between the left and the right action~:

$$\xymatrix@C=50pt@R=30pt{ \Po \bc (\Po \oplus \underbrace{M}_{1})
\bc \Po \ar[r]^{(\lambda+\mu)\bc \Po} \ar[d]^{\Po \bc (\rho +
\mu)} & (\Po \oplus \underbrace{M}_{1})\bc \Po  \ar[d]^{\rho} \\
\Po \bc (\Po \oplus \underbrace{M}_{1}) \ar[r]^{\lambda}  & M.}
$$
\end{itemize}
The notation  $\Po \bc \Po \bc (\Po \oplus \underbrace{M}_{1})$
stands for the sub-$\Sy$-bimodule of $\Po \bc \Po \bc (\Po \oplus
{M})$ with only one $M$ on the upper level. It is represented by
linear combinations of 3-levelled graphes whose vertices are
indexed by elements of $\Po$ and just one of $M$ on the first
level. The other $\Sy$-bimodules with just one element coming from
$M$ are denoted in the same way, $ \Po \bc (\Po \oplus
\underbrace{M}_{1})\bc \Po$ has one element of $M$ on the second
level and $(\Po \oplus \underbrace{M}_{1}) \bc \Po \bc \Po$ has
one element of $M$ on the third level.
\end{dei}

One purpose of this notion is to define the notion of \emph{abelian or
infinitesimal} extension of a prop(erad) $\Po$. It is defined by a prop(erad) structure on $\Po\oplus M$, when $M$ is an infinitesimal bimodule over $\Po$ (see Section~\ref{Def a la Quillen}
Lemma~\ref{infinitesimal extension}). Another important property
is that, for any sub-$\Sy$-bimodule $J$ of $\Po$, the ideal
generated by $J$ in
$\Po$ is equal to the free infinitesimal $\Po$-bimodule on $J$.\\

Since the prop(erad) $\Po$ has a unit, it is equivalent to have
two actions $\lambda \,: \, \Po \bc_{(1,1)} M \to M$ and $\rho\,
:\, M \bc_{(1,1)} \Po \to M$ that are compatible with the
composition product of prop(erad) $\Po$. Notice that the category
of infinitesimal bimodules over a
prop(erad) forms an abelian category. \\

\begin{Ex}
Any morphism of prop(erad)s $f\, : \, \Po \to \Qo$ defines an
infinitesimal $\Po$-bimodule structure on $\Qo$~:
$$\Po \bc_{(1,1)} \Qo \xrightarrow{f\bc \Qo}\Qo \bc_{(1,1)} \Qo
\xrightarrow{\mu_\Qo} \Qo \quad \textrm{and} \quad \Qo \bc_{(1,1)}
\Po \xrightarrow{\Qo\bc f}\Qo \bc_{(1,1)} \Qo
\xrightarrow{\mu_\Qo} \Qo.$$
\end{Ex}

\subsection{(Co)Derivations}

Let $(\Po, \mu)$ be a dg prop(erad) and $(M, \lambda, \rho)$ be an
infinitesimal $\Po$-bimodule

\begin{dei}[Derivation]
A homogenous morphism $\partial\, :\, \Po \to M$ is a
\emph{homogenous derivation}  if
$$\partial\circ \mu_{(1,\,
1)}(-,-)=\rho(\partial(-), -) + \lambda(-,
\partial(-)).$$
 This formula, applied to elements $p_1 \bbc
p_2$ of $\Po\bbc \Po$, where $p_1$ and $p_2$ are homogenous
elements of $\Po$, gives
\begin{eqnarray*}
\partial\circ\mu(p_1\bbc p_2) &=&
\rho\big(\partial(p_1)\bbc p_2 \big) +
(-1)^{|\partial||p_1|}\lambda \big( p_1\bbc
\partial(p_2)\big).
\end{eqnarray*}
A \emph{derivation} is a sum of homogenous derivations. The set of
homogenous derivations of degree $n$ is denoted by $\Der^n(\Po,\,
M)$ and the set of derivations is denoted $\Der^\bullet(\Po,\, M)$
\end{dei}

\begin{Ex}
The differential of a dg prop(erad) $\Po$ is a derivation of
degree $-1$, that is an element of $\Der^{-1}(\Po,\, \Po)$.
\end{Ex}

In this section, we only consider derivations $\Der(\Po,\, \Qo)$,
where the  infinitesimal $\Po$-bimodule structure on $\Qo$ is
given by a morphism of prop(erad)s $\Po \to \Qo$. In the rest of
the text, we need the following lemma which makes explicit the  derivations on a free prop(erad). For a prop(erad) $(\Qo,
\mu_\Qo)$, any graph $\G$ of $\F(\Qo)^{(n)}$ represents a class
$\overline{\G}$ of levelled graphs of $\Qo^{\boxtimes n}$. We
recall that there is a morphism $\widetilde{\mu}_\Qo \, : \,
\F(\Qo) \to \Qo$, the counit of adjunction, equal to
$\widetilde{\mu}_\Qo(\G):=\mu_\Qo^{\circ (n-1)}(\overline{\G})$.
The morphism $\widetilde{\mu}_\Qo$ is the only morphism of
prop(erad)s extending the map $\Qo \xrightarrow{\Id} \Qo$.

\begin{lem}
\label{uniqueDer} Let $\rho \, :\, \F(V) \to \Qo$ be a morphism of
\draftnote{augmented dg ?} prop(erad)s of degree $0$. Every
derivation from the free dg prop(erad) $\F(V)$ to $\Qo$ is
characterized by its restriction on $V$, that is there is a
canonical one-to-one correspondence $\Der^n_\rho(\F(V),\,
\Qo)\cong \Hom_n^\Sy(V,\, \Qo)$.

For every morphism of dg $\Sy$-bimodules $\theta \, : \, V \to
\Qo$, we denote the unique derivation which extends $\theta$ by
$\partial_\theta$. The image of an element $\G(v_1, \ldots ,v_n)$
of $\F(V)^{(n)}$ under $\partial_\theta$ is
$$\partial_\theta(\G(v_1, \ldots ,v_n))=\sum_{i=1}^n (-1)^{|\theta|.(|v_1|+\cdots + |v_{i-1}|)}
\widetilde{\mu}_\Qo\big(\G(\rho(v_1), \ldots ,\rho(v_{i-1}),
\theta(v_i), \rho(v_{i+1}), \ldots , \rho(v_n)) \big) .$$
\end{lem}

\begin{proo}
Let us denote by $\theta$ the restriction of the derivation
$\partial$ on $V$, that is $\theta=\partial_V \, :\, V \to
\overline{\Qo}$. From $\theta$, we can construct the whole
derivation $\partial$ by induction on the weight $n$ of the free
prop(erad) $\F(V)$ as follows.

For $n=1$, we have $\partial_\theta^1=\theta \, : \, V \to \Qo$.
Suppose now that $\partial_\theta^n \, : \, \F(V)^{(n)} \to \Qo$
is constructed and it is given by the formula of the Lemma. Any
simple element of $\F(V)^{(n+1)}$ represented by a graph with
$n+1$ vertices indexed by elements of $V$ is the concatenation of
a graph with $n$ vertices with an extra vertex from the top or the
bottom. In the last case, $\partial_\theta^{n+1}$ is given the
commutative diagram
$$\xymatrix@C=70pt{\F(V)^{(n+1)} \ar[r]^{\partial_\theta^{n+1}}   & \Qo \\
V \boxtimes_{(1,1)} \F(V)^{(n)} \ar[u]_{\mu_{\F(V)}} \ar[r]^{\rho
\boxtimes \partial_\theta^{n}+ \partial_\theta^{n}  \boxtimes \rho
} & \Qo\boxtimes_{(1,1)} \Qo \ar[u]_{\mu_\Qo}.}$$ The other case
is dual. It is easy to check that the formula is still true for
elements of $\F(V)^{(n+1)}$, that is graphs with $n+1$ vertices.
Finally, since $\rho$ is a morphism of prop(erad)s,
$\partial_\theta$ is well defined and is a derivation.
\end{proo}

\begin{Ex}
A differential $\partial$ on a free prop(erad) $\F(V)$ is a
derivation of $\Der^{-1}_\Id(\F(V),\, \F(V))$ such that
$\partial^2=0$.
\end{Ex}

\begin{dei}[quasi-free prop(erad)]
A dg prop(erad) $(\F(V),\, \partial)$ such that the underlying
prop(erad) is free is called a \emph{quasi-free} prop(erad).
\end{dei}

Notice that in a quasi-free prop(erad), the differential is not
freely generated and is a derivation of the form given above.\\

Dually, let $(\Co, \Delta^\Co)$ and $(\Do, \Delta^\Do)$ be two
coaugmented dg coprop(erad)s and let $\rho \, : \, \Co \to \Do$ be
a morphism of coaugmented dg coprop(erad)s of degree $0$. One can
define the dual notion of infinitesimal comodule over a
coprop(erad) and general coderivations. Since we only need
coderivations between two coprop(erad)s, we do not go into such
details here.

\begin{dei}[Coderivation]
A homogeneous morphism $d \, : \, \Co \to \Do$ is a
\emph{homogeneous coderivation} of $\rho$ if the following diagram
is commutative
$$\xymatrix@C=60pt{\Co \ar[r]^{d} \ar[d]^{\Delta^\Co_{(1,1)}}& \Do \ar[d]^{\Delta^\Do_{(1,1)}} \\
\Co \boxtimes \Co \ar[r]^{d\boxtimes \rho + \rho \boxtimes d}   &
\Do \boxtimes \Do.}$$ A \emph{coderivation} is a sum of homogenous
coderivations. The space of coderivations is denoted by
$\textrm{Coder}^\bullet_\rho(\Co, \, \mathcal{D})$.
\end{dei}

\begin{Ex}
The differential of a dg coprop(erad) $\Co$ is a coderivation of
degree $-1$.
\end{Ex}

\begin{Rq}
For a cooperad $\mathcal{D}$, we can define a more general notion
of coderivation form a $\mathcal{D}$-cobimodule to $\mathcal{D}$
by a similar formula. The definition given here is a particular
case. Since $\rho \, : \, \Co \to \Do$ is a morphism of
coprop(erad)s, it provides $\Co$ with a natural structure of
$\mathcal{D}$-cobimodule.
\end{Rq}

As explained in the first section, the dual statement of
Lemma~\ref{uniqueDer} holds only for connected coprop(erad)s.

\begin{lem}\label{uniquecoDer}
Let $\Co$ be a connected coprop(erad) and let $\rho \, :\, \Co \to
\F^c(W)$ be a morphism of augmented coprop(erad)s. Every
coderivation from $\Co$ to the cofree connected coprop(erad)
$\F^c(W)$ is characterized by its projection on $W$, that is there
is a canonical one-to-one correspondence
$\textrm{Coder}^n_\rho(\Co,\, \F^c(W))\cong \Hom_n^\Sy(\oC,\, W)$.
\end{lem}

\begin{proo}
The proof is similar to the one of Lemma~\ref{uniqueDer} and goes
by induction on $r$, where $F_r$ stands for the coradical
filtration of $\Co$. The assumption that the coprop(erad) $\Co$ is
connected ensures that the image of an element $X$ of $F_r$ under
$d$ lives in $\bigoplus_{n\leq r} \F^c(W)^{(n)}$. Therefore,
$d(X)$ is finite and $d$ is well defined.
\end{proo}

We denote by $d_\nu$ the unique coderivation which extends a map
$\nu \, : \, \oC \to W$.

\begin{Ex}
A differential $d$ on a cofree coprop(erad) $\F^c(W)$ is a
coderivation of $\Der^{-1}_\Id(\F^c(W),\, \F^c(W))$ such that
$d^2=0$. By the preceding lemma, it is characterized by the
composite $\F^c(W)\xrightarrow{d} \F^c(W) \epi W$. Its explicit
formula can be found in Lemma~\ref{Form of the coderivation}.
\end{Ex}

\begin{dei}[quasi-cofree coprop(erad)]
 A dg coprop(erad) $(\F^c(W),\, d)$ such that the
underlying coprop(erad) is connected cofree is called a
\emph{quasi-cofree} coprop(erad).
\end{dei}

\subsection{(De)Suspension}\label{suspension}

The homological \emph{suspension} of a dg $\Sy$-bimodule $M$ is
denoted by $sM:=\KK s \otimes M$ with $|s|=1$, that is
$(sM)_i\cong M_{i-1}$. Dually, the homological \emph{desuspension}
of $M$ is denoted by $s^{-1}M := \KK s^{-1}\otimes M$ with
$|s^{-1}|=-1$, that is $(s^{-1}M)_i\cong M_{i+1}$.\\

Let $(\Po,\, d)$ be an augmented dg $\Sy$-bimodule, that is
$\Po=\oPo \oplus I$. A map of augmented $\Sy$-bimodules $\mu \, :
\, \F^c(\oPo) \to \Po$ consists of a family of morphisms of dg
$\Sy$-bimodules $\mu_n \, :\, \F^c(\Po)^{(n)}\to \Po$ for each
integer $n\ge 1$. (For $n=0$, the map $\mu$ is the identity $I \to
I$.) There is a one-to-one correspondence between maps
$\{\F^c(\oPo)\to \Po\}$ and maps $\{\F^c(s\oPo)\to s\Po\}$. To
each map $\mu\, : \, \F^c(\oPo)\to \Po$, we associate the map
$s\mu\, : \, \F^c(s\oPo)\to s\Po$ defined as follows for $n\ge 1$,
$$(s{\mu})_n \,: \, \F^c(s\oPo)^{(n)} \xrightarrow{\tau_n} s^n \F^c(\oPo)^{(n)}
\xrightarrow{s^{-(n-1)}}  s \F^c(\oPo)^{(n)} \xrightarrow{s\otimes
\mu_n} s\Po,$$
where the map $\tau_n$ moves the place of the suspension elements from the vertices outside the graph.
Since it involves permutations between suspensions $s$ and
elements of $\Po$, the map $\tau_n$ yields signs by Koszul-Quillen
rule. Using the fact that an element of $\F^c(\oPo)$ is an
equivalent class of graphs with levels (see \ref{Free properad}),
one can make these signs explicit. The exact formula relating
$(s\mu)$ to $\mu$ is
$$\mu\big(\G(p_1,\ldots,p_n)\big)=(-1)^{\varepsilon(p_1,\ldots,p_n)}s^{-1}(s\mu)\big(\G(s p_1,\ldots,s p_n)\big),$$
where
$\varepsilon(p_1,\ldots,p_n)=(n-1)|p_1|+(n-2)|p_2|+\cdots+|p_{n-1}|$.

The degrees of $\mu$ and $s\mu$ are related by the formula
$|(s\mu)_n|=|\mu_n| -(n-1)$. Therefore, the degree of $\mu_n$ is
$n-2$ if and only if the degree of $(s\mu)_n$ is
$-1$.\\

Dually, for any map of augmented $\Sy$-bimodules $\delta \, : \,
\Co \to \F(\oC)$, we denote by $\delta_n$ the composite $\Co
\xrightarrow{\delta} \F(\oC)  \epi \F(\oC)^{(n)}$. There is a
one-to-one correspondence between maps $\{\Co
 \to \F(\oC) \}$ and maps $\{ s^{-1}\Co \to \F(s^{-1}\oC)  \}$. To each map
$\delta\, : \, \Co
 \to \F(\oC)$, we associate the map
$s^{-1}\delta\, : \, s^{-1}\Co \to \F^c(s^{-1}\oC)$ defined as
follows, for $n\ge 1$,
$$(s^{-1}{\delta})_n \,: \,  s^{-1} \Co \xrightarrow{s^{-(n-1)}\otimes \, \delta_n}
s^{-n}\F(\oC)^{(n)} \xrightarrow{\sigma_n}
\F(s^{(-1)}\oC)^{(n)}.$$

We have  $|(s^{-1}\delta)_n|=|\delta_n| -(n-1)$. The degree of
$\delta_n$ is $n-2$ if and only if the degree of
$(s^{-1}\delta)_n$ is $-1$.

\subsection{Twisting morphism}\label{Twisting morphism}
We generalize the notion of \emph{twisting morphism} (or twisting
cochains) of associative algebras (see E. Brown \cite{Brown59} and J.C. \cite{Moore70}) to
prop(erad)s.\\

Let $\Co$ be a dg coprop(erad) and $\Po$ be a dg prop(erad). We
proved in Theorem~\ref{Convolution LR and Lie algebra} that $\Hom^\Sy(\Co,
\Po)$ is a dg Lie-admissible algebra with the convolution product.

\begin{dei}
A morphism $\Co \xrightarrow{\alpha} \Po$, of degree $-1$, is
called a \emph{twisting morphism} if it is a solution of the
\emph{Maurer-Cartan} equation
$$D(\alpha) + \alpha \star \alpha=0.$$
\end{dei}

Denote by $\textrm{Tw}(\Co, \Po)$ the set of twisting morphisms
 in $\Hom^\Sy(\Co, \Po)$, that is Maurer-Cartan elements in the convolution prop(erad).
Since twisting morphisms have degree $-1$, it is equivalent for
them to be solution of the classical Maurer-Cartan equation in the
associated dg Lie algebra, that is $D(\alpha) +
\frac{1}{2}\lbrack \alpha, \alpha\rbrack =0$.\\

When $\Po$ is augmented and $\Co$ coaugmented, we will consider
either a twisting morphism between $\Co$ and $\Po$, which sends
$I$ to $0$, or the associated morphism which sends $I$ to $I$ and
$\oC$ to $\oPo$.\\

The following constructions show that the bifunctor $\textrm{Tw}(-, -)$
can be represented on the left
and on the right.

\subsection{Bar construction}

We recall from \cite{Vallette03} Section $4$, the definition of
the \emph{bar construction} for properads and extend it to props.
It is  a functor
$$B \ : \ \{\textrm{aug. dg prop(erad)s}    \} \longrightarrow
\{\textrm{coaug. dg coprop(erad)s} \}.  $$

Let $(\Po,\, \mu,\, \eta,\, \epsilon)$ be an augmented prop(erad).
Denote by $\oPo$ its augmentation ideal $\Ker (\Po
\xrightarrow{\epsilon} I)$. The prop(erad) $\Po$ is naturally
isomorphic to $\Po=I \oplus \oPo$. The bar construction $B(\Po)$
of $\Po$ is a dg coprop(erad) whose underlying space is the cofree
coprop(erad) $\F^c(s\oPo)$ on the suspension of $\oPo$.

The partial product of $\Po$ induces a map of augmented
$\Sy$-bimodules defined by the composite
$$\mu_2 \ : \ \overline{\F}^c(\oPo) \epi \F^c(\oPo)^{(2)} \cong \oPo \boxtimes_{(1,1)}
\oPo \xrightarrow{\mu_{(1,1)}} \oPo.$$ We have seen in the
previous section that $\mu_2$ induces a map $s\mu_2$. Consider the
map $\KK s \otimes \KK s \xrightarrow{\Pi_s} \KK s$ of degree $-1$
defined by $\Pi_s(s\otimes s):=s$. The map $s\mu_2$ is equal to
the composite
\begin{eqnarray*}
s\mu_2 \ &:& \ \overline{\F}^c(s\oPo) \epi \F^c(s\oPo)^{(2)}\cong
(\KK s \otimes \oPo) \boxtimes_{(1,1)} (\KK s \otimes \oPo)\\ &&
\xrightarrow{\Id \otimes \tau \otimes \Id} (\KK s \otimes \KK s)
\otimes (\oPo \boxtimes_{(1,1)} \oPo) \xrightarrow{\Pi_s\otimes
\mu_{(1,1)}} \KK s \otimes \oPo.
\end{eqnarray*}

Since $\F^c(s\oPo)$ is a cofree connected coprop(erad), by
Lemma~\ref{uniquecoDer} there exists a unique coderivation
$d_2:=d_{s\mu_2} \, : \, \F^c(s\oPo)\to \F^c(s\oPo)$ which extends
$s\mu_2$. When $(\Po,\ d_\Po)$ is an augmented dg prop(erad), the
differential $d_\Po$ on $\Po$ induces an internal differential
$d_1$ on $\F^c(s\oPo)$. The total complex of this bicomplex is the
\emph{bar construction} $ B(\Po,\, d_\Po) := (\F^c(s\oPo),\,
d=d_1+d_2)$ of the augmented dg prop(erad) $(\Po,\, d_\Po)$.\\

Notice that the relation $d^2=0$ can be understood conceptually
from the Lie-admissible relations verified by the partial product
of the prop(erad) $\Po$.

\subsection{Cobar construction}

Dually, the \emph{cobar construction} (\cite{Vallette03} Section
$4$)  for coprop(eard)s is a functor
$$\Omega \ : \ \{\textrm{coaug. dg coprop(erad)s}    \} \longrightarrow
\{\textrm{aug. dg prop(erad)s} \}.  $$

Let $(\Co,\, \Delta,\, \varepsilon,\, u)$ be a coaugmented
coprop(erad). Denote by $\oC$ its augmentation $\Ker (\Co
\xrightarrow{\varepsilon} I)$. In this case, $\Co$ splits
naturally as $\Co=I \oplus \oC$. The cobar construction
$\Omega(\oC)$ of $\oC$ is a dg prop(erad) whose underlying space
is the free prop(erad) $\F(s^{-1}\oC)$ on the desuspension of
$\oC$.

The partial coproduct of $\Co$ induces a natural map of augmented
$\Sy$-bimodules defined by
$$\Delta_2 \ : \ \oC \xrightarrow{\Delta_{(1,1)}}\oC \boxtimes_{(1,1)} \oC \cong \F(\oC)^{(2)}
\mono \overline{\F}(\oC). $$

This map gives a map $s^{-1}\Delta_2 \, : \, s^{-1}\oC \to
\F(s^{-1}\oC)$. Consider $\KK s^{-1}$ equipped with the diagonal
map $\KK s^{-1} \xrightarrow{\Delta_s} \KK s^{-1} \otimes \KK
s^{-1}$ of degree $-1$ defined by the formula
$\Delta_s(s^{-1}):=s^{-1}\otimes s^{-1}$. The map $s^{-1}\Delta_2$
is equal to

\begin{eqnarray*}
s^{-1}\Delta_2  &:&  \KK s^{-1}\otimes \oC
\xrightarrow{\Delta_s\otimes \Delta_{(1,1)}}
 \KK s^{-1}\otimes \KK s^{-1}\otimes \oC\boxtimes_{(1,1)} \oC  \xrightarrow{\Id\otimes  \tau \otimes \Id }\\
&& (\KK s^{-1}\otimes\Co) \boxtimes_{(1,1)} (\KK s^{-1}\otimes
\Co) \cong \F(s^{-1} \oC)^{(2)}\mono \F(s^{-1} \oC).
\end{eqnarray*}

Since $\F(s^{-1}\oC)$ is a free prop(erad), by
Lemma~\ref{uniqueDer} there exists a unique derivation
$\partial_2:=\partial_{s^{-1}\Delta_2} \, : \, \F(s^{-1} \oC) \to
\F(s^{-1} \oC)$ which extends $s^{-1}\Delta_2$. When $(\Co,\
d_\Co)$ is an augmented dg coprop(erad), the differential $d_\Co$
on $\Co$ induces an internal differential $\partial_1$ on
$\F(s^{-1} \oC)$. The total complex of this bicomplex is the
\emph{cobar construction} $ \Omega(\Co,\, d_\Co) :=
(\F(s^{-1}\oC),\,
\partial=\partial_1+\partial_2)$ of the augmented dg coprop(erad) $(\Co,\,
d_\Co)$.

\subsection{Bar-Cobar Adjunction}

As for derivations, a morphism of prop(erad)s is characterized by
the image of the indecomposable elements. We recall this fact and
the dual statement in the following lemma.

\begin{lem}\label{unique morphism}
Let $V$ be an $\Sy$-bimodule and let $\Qo$ be a prop(erad), there
is a canonical one-to-one correspondence $\mathrm{Mor}_{\rm
prop(erad)s}\big( \F(V), \Qo\big) \cong \Hom^\Sy(V, \Qo)$.

Dually, let $W$ be an $\Sy$-bimodule and let $\Co$ be a
coprop(erad), there is a canonical one-to-one correspondence
$\mathrm{Mor}_{\rm coprop(erad)s} \big( \Co, \F^c(W) \big) \cong
\Hom^\Sy(\Co, W)$.
\end{lem}

Let $(\Co, d_\Co)$ be a dg coprop(erad) and $(\Po,  d_\Po)$ be a
dg prop(erad). We will apply this result to the bar and the cobar
construction of $\Po$ and $\Co$ respectively, that is we want to
describe the space of morphisms of \textbf{dg}-prop(erad)s
$\Mor_{\textrm{dg prop(erad)s}}\left(\Omega(\Co),\, \Po\right)$
for instance. By the preceding lemma, this space is isomorphic to
the space of morphisms of $\Sy$-bimodules $\Hom^\Sy_0(s^{-1} \Co,
\Po)$ of degree $0$ whose unique extension commutes with the
differentials. Therefore, this space of morphisms is the subspace
of $\Hom^\Sy_{-1}(\Co, \Po)$ whose elements satisfy a certain
relation, which is exactly the Maurer-Cartan equation.

\begin{pro}\label{Bar-cobar adjunction} For every augmented dg prop(erad) $\Po$ and every
coaugmented dg coprop(erad) $\Co$, there are canonical one-to-one
correspondences
$$\mathrm{Mor}_{\textrm{dg prop(erad)s}}\left(\Omega(\Co),\,
\Po\right)\cong \mathrm{Tw}(\Co, \Po)\cong
\mathrm{Mor}_{\textrm{dg coprop(erad)s}}\left(\Co,\,
\B(\Po)\right).$$
\end{pro}

\begin{proo}
Since $\Omega(\Co)=\F(s^{-1}(\oC))$, by Lemma~\ref{unique
morphism} every morphism $\varphi$ of $\Sy$-bimodules in
$\Hom^\Sy_0(s^{-1} \Co, \Po)$ extends to a unique morphism of
prop(erad)s between $\Omega(\Co)$ and $\Po$. The latter one
commutes with the differentials if and only if the following
diagram commutes
$$\xymatrix@R=30pt@C=30pt{s^{-1}\oC \ar[r]^{\varphi} \ar[d]^{\partial} & \Po \ar[dr]^{d_\Po}&  \\
\F(s^{-1}\oC)^{(\leqslant 2)} \ar[r]^(0.6){\F(\varphi)} & \F(\Po)
\ar[r]^{\widetilde{\mu}^\Po} & \Po. }$$ For an element $c\in \oC$,
we use Sweedler's notation to denote the image of $c$ under
$\Delta_2$, that is $\Delta_2(c)=\sum c' \boxtimes_{(1,1)} c''$.
The diagram above corresponds to the relation
$$d_\Po\circ \varphi (s^{-1} c)= \varphi \circ \partial_1(s^{-1}
c) + \mu^\Po \circ (\varphi \boxtimes_{(1,1)} \varphi) \circ
s^{-1}\Delta_2 (s^{-1} c). $$ Denote by $\alpha$ the desuspension
of $\varphi$, that is $\alpha(c)=-\varphi({s^{-1}}c)$. Since
$\partial_1(s^{-1}c)=-s^{-1}\partial_\Co(c)$, the relation becomes
$$-d_\Po\circ \alpha (c)= \alpha \circ \partial_\Co(
c) + \mu^\Po \circ (\alpha \boxtimes_{(1,1)} \alpha) \circ
\Delta_2 (c), $$ which is the Maurer-Cartan equation.
\end{proo}

Therefore, the bar and cobar constructions form a pair of adjoint
functors
$$\Omega \ : \ \{\textrm{coaug. dg coprop(erad)s}\} \rightleftharpoons
\{\textrm{aug. dg prop(erad)s}\}\ : \ B.$$

If we apply the isomorphisms of Proposition~\ref{Bar-cobar
adjunction} to $\Co=B(\Po)$, the morphism associated to the
identity on $B(\Po)$ is the \emph{counit} of the adjunction
$\epsilon \, : \, \Omega(B(\Po))\to \Po$. In this case, we get a
universal twisting morphism $B(\Po) \to \Po$.

The morphism associated to the identity of $\Omega(\Co)$ when
$\Po=\Omega(\Co)$  is the \emph{counit} of the adjunction $\Co \to
B(\Omega(\Co))$. In this case, we get a universal twisting
morphism $\Co \to \Omega(\Co)$.

\begin{pro}
Any twisting morphism $\alpha \, :\, \Co \to \Po$ factors through
$B(\Po) \to \Po$ and $\Co \to \Omega(\Co)$.
$$\xymatrix@R=10pt@C=10pt{ & \Omega(\Co)\ar[dr] & \\
\Co\ar[ur] \ar[dr]\ar[rr]^{\alpha}&& \Po \\
&B(\Po)\ar[ur] &} $$
\end{pro}

\begin{proo}
It is a corollary of Proposition~\ref{Bar-cobar adjunction}.
\end{proo}

\subsection{Props vs properads}

The main difference for (co)bar construction between props and
properads lies on the type of graphs and compositions. The
underlying module of the bar construction of a prop $\Po$ is
spanned by not necessarily connected graphs whose vertices are
labelled with elements of $\Po$. The boundary map is the unique
coderivation which extends the partial product. It is given explicitly
by the sum of the compositions of pair of vertices that are
either adjacent (see Section~\ref{admissible subgraph}) or belong
to two different connected graphs. Whereas for a properad, the
underlying module is spanned by connected labelled graphs and the
boundary map just composes adjacent pairs of operations.

\subsection{Bar-cobar resolution}

In \cite{Vallette03} Theorem $5.8$, we proved that the unit of
adjunction $\epsilon$ is a canonical resolution in the weight
graded case. We extend this result to any dg properad here.

\begin{thm}\label{BarCobarResolution}
For every augmented dg properad $\Po$, the bar-cobar construction
is a resolution of $\Po$
$$\epsilon\ : \ \Omega(B(\Po))\xrightarrow{\backsimeq}\Po.$$
\end{thm}

\begin{proo}
 The bar-cobar construction
of $\Po$ is the chain complex defined on the underlying
$\Sy$-bimodule $\F\big(s^{-1} \overline{\F}^c(s \bar{\Po}) \big)$.
The differential $d$ is the sum of three terms
$d=\partial_2+d_2+d_{\Po}$, where $d_\Po$ is induced by the
differential on $\Po$, $d_2$ is induced by the differential of the
bar construction $B(\Po)$ and $\partial_2$ is the unique
derivation which extends the partial coproduct of $\F^c(s
\bar{\Po})$.

Define the filtration $F_s:=\bigoplus_{r\le s} \F\big(s^{-1}
\overline{\F}^c(s \bar{\Po}) \big)_r$, where $r$ denotes the total
number of elements of $\bar{\Po}$. Let $E^\bullet_{st}$ be the
associated spectral sequence.

This filtration is bounded below and exhaustive. Therefore, we can
apply the classical convergence theorem for spectral sequences
(see \cite{Weibel} \draftnote{indiquer le chapitre}) and prove
that $E^\bullet$ converges to the homology of the bar-cobar
construction.

We have that $E^0_{st}=\F_{s+t}\big(s^{-1} \F^c(s \bar{\Po})
\big)_s$, where $s+t$ is the total homological degree. From
$d_2(F_s)\subset F_{s-1}$, $d_\Po(F_s)\subset F_s$ and
$\partial_2(F_s)\subset F_s$, we get that $d^0=\partial_2+d_\Po$.
The problem is now reduced to the computation of the homology of
the cobar construction of the dg cofree connected coproperad
$\F^c(s \bar{\Po})$ on the dg $\Sy$-bimodule $s \bar{\Po}$. This
complex is equal to the bar-cobar construction of the weight
graded properad ($\Po$, $\mu'$), where $\Po^{(0)}=I$ and
$\Po^{(1)}=\bar{\Po}$, such that the composition $\mu'$ is null.
We conclude using Theorem~{5.8} of \cite{Vallette03}.
\end{proo}

\begin{pro}
The bar-cobar resolution provides a canonical cofibrant resolution
to any non-negatively graded dg properad.
\end{pro}

We refer the reader to the Appendix~\ref{model category} for the
model category structure on dg prop(erads)

\begin{proo}
The bar-cobar resolution is quasi-free. We conclude by
Corollary~\ref{quasi-free -> cofibrant}.
\end{proo}

\section{Homotopy (Co)prop(erad)s}

An associative algebra is a vector space endowed with a binary
product that verifies the strict associative relation. J. Stasheff
defined in \cite{Stasheff63} a lax version of this notion. It is
the notion of an associative algebra up to homotopy or (strong)
homotopy algebra. Such an algebra is a vector space equipped with
a binary product that is associative only up to an infinite
sequence of homotopies. In this section, we recall the
generalization of this notion, that is the notion of
\emph{(strong) homotopy properad} due to J. Gran\aa ker
\cite{Granaker06}. We extend it to props and we also define in
details the dual notion of \emph{(strong) homotopy coprop(erad)},
 which will be essential to deal with minimal models in the next section .
 The notions of \emph{homotopy non-symmetric (co)properad} and
\emph{homotopy non-symmetric (co)prop} are obtained in the same
way.

\subsection{Definitions}

Following the same ideas as for algebras (associative or Lie, for
instance), we define the notion of \emph{homotopy (co)prop(erad)}
via (co)derivations and (co)free (co)prop(erad)s.

\begin{dei}[Homotopy prop(erad)]
A structure of \emph{homotopy prop(erad)} on an augmented dg
$\Sy$-bimodule $(\Po, d_\Po)$ is a coderivation $d$ of degree $-1$
on $\F^c(s\oPo)$ such that $d^2=0$.
\end{dei}

A structure of homotopy prop(erad) is equivalent to a structure of
quasi-cofree coprop(erad) on $s\oPo$. We call the latter the
\emph{(generalized) bar construction of $\Po$} and we still denote
it by $B(\Po)$. Since $\F^c(s\oPo)$ is a cofree connected
coprop(erad), by Lemma~\ref{uniquecoDer} the coderivation $d$ is
characterized by the composite
$$s\mu \ : \ \F^c(s\oPo) \xrightarrow{d} \F^c(s\oPo) \epi s\Po,$$
that is $d=d_{s\mu}$. The map $s\mu$ of degree $-1$ is equivalent
to a unique map $\mu \, : \, \F^c(\oPo) \to \Po$, such that $\mu_n
\, : \, \F^c(\oPo)^{(n)}\to \Po$ has degree $n-2$. The condition
$d^2=0$ written with the $\{ \mu_n\}_n$ is made explicit in
Proposition~\ref{Properadinfty}.

\begin{Ex}
A dg prop(erad) is a homotopy prop(erad) such that every map
$\mu_n=0$ for $n\ge3$. In this case, $(\F^c(s\oPo), d)$ is the bar
construction of $\Po$.
\end{Ex}

We define the notion of \emph{homotopy coprop(erad)} by a direct
dualization of the previous arguments.

\begin{dei}[Homotopy coprop(erad)]
A structure of \emph{homotopy coprop(erad)} on an augmented dg
$\Sy$-bimodule $(\Co,\, d_\Co)$ is a derivation $\partial$ of
degree $-1$ on $\F(s^{-1}\oC)$ such that $\partial^2=0$.
\draftnote{idem}
\end{dei}

A structure of homotopy coprop(erad) is equivalent to a structure
of quasi-free prop(erad) on $s^{-1}\oC$. We call the latter the
\emph{(generalized) cobar construction of $\Co$} and we still
denote it by $\Omega(\Co)$. By Lemma~\ref{uniqueDer}, the
derivation $\partial$ is characterized by its restriction on
$s^{-1}\oC$
$$s^{-1}\Delta \ : \ s^{-1}\oC \mono \oF(s^{-1}\oC) \xrightarrow{\partial} \oF(s^{-1}\oC),$$
that is $\partial=\partial_{s^{-1}\Delta}$. The map $s^{-1}\Delta$
of degree $-1$ is equivalent a map $\Delta \, :\, \Co \to
\F(\oC)$, such that the component $\Delta_n \, : \, \Co \to
\F(\oC)^{(n)}$ has degree $n-2$. The condition $\partial^2=0$ is
equivalent to relations for the $\{ \Delta_n \}_n$ that we make
explicit in Proposition~\ref{Coproperadinfty}.

\begin{Ex}
A dg coprop(erad) is a homotopy coprop(erad) such that every map
$\Delta_n=0$ for $n\ge3$. In this case, $(\F(s^{-1}\oC),
\partial)$ is the cobar construction of $\Co$.
\end{Ex}

When $\Po$ is concentrated in arity $(1, 1)$, the definition of a
homotopy properad on $\Po$ is exactly the same than the definition
of an strong homotopy algebra given by J. Stasheff in
\cite{Stasheff63}. Dually, when $\Co$ is concentrated in arity
$(1,\, 1)$, we get the notion of strong homotopy coassociative
algebra.

When $\Po$ is concentrated in arity $(1, n)$ for $n\ge 1$, we have
the notion of \emph{strong homotopy operad} (see
\cite{VanderLaan02}). The dual notion gives the definition of a
\emph{strong homotopy cooperad}.

\begin{Rq}
By abstract nonsense, the notion of homotopy prop(erad) should
also come from Koszul duality for colored operads (see
\cite{VanderLaan03}). There exists a colored operad
whose``algebras" are (partial) prop(erad)s. Such a colored operad
is quadratic (the associativity relation of the partial product of
a prop(erad) is an equation between compositions of two elements.)
It should be a Koszul colored operad. An ``algebra" over the
Koszul resolution of this colored operad is exactly a homotopy
prop(erad).
\end{Rq}

\subsection{Admissible subgraph}\label{admissible subgraph}

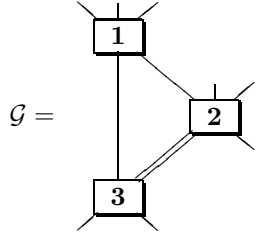
\begin{figure}[h]
$$\G=\vcenter{\xymatrix@R=6pt@C=6pt{
& & & & \\
&*+[F-,]{\ \textbf{1} \ } \ar@{-}[ur] \ar@{-}[u] \ar@{-}[ul] \ar@{-}[dddd]  \ar@{-}[ddrr] & & & \\
&&&&\\
& & & *+[F-,]{\ \textbf{2} \ } \ar@{-}[u] \ar@{-}[dr] \ar@{=}[ddll] \ar@{-}[ur] & \\
&&&& \\
& *+[F-,]{\ \textbf{3} \ } \ar@{-}[dr]\ar@{-}[dl]& & & \\
& & & & }}  $$ \caption{The vertices $1$, $2$ and $2$, $3$ are
adjacent. The vertices $1$ and $3$ are not adjacent.}
\label{AdjacentVertices}
\end{figure}

Let $\mathcal{G}$ be a connected graph directed by a flow and
denote by $\mathcal{V}$ its set of vertices. We define a partial
order on $\mathcal{V}$ by the following covering relation  : $i
\prec j$ if $i$ is below $j$ according to the flow and if there is
no vertex between them. In this case, we say that $i$ and $j$ are
\emph{adjacent} (see also \cite{Vallette03} p. 34). Examples of
adjacent and non-adjacent vertices can be found in
Figure~\ref{AdjacentVertices}.

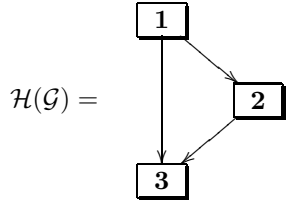
\begin{figure}[h]
$$\mathcal{H(G)}=\vcenter{\xymatrix@R=6pt@C=6pt{
&*+[F-,]{\ \textbf{1} \ } \ar[dddd]  \ar[ddrr] & & & \\
&&&&\\
& & & *+[F-,]{\ \textbf{2} \ }  \ar[ddll] & \\
&&&& \\
& *+[F-,]{\ \textbf{3} \ } & & & }}  $$ \caption{The Hasse
diagramm associated to the graph of Figure~\ref{AdjacentVertices}}
\label{HasseDiagramm}
\end{figure}

Denote this poset by $\Pi_{\mathcal{G}}$ and consider its Hasse
diagram $\mathcal{H(G)}$, that is the diagram composed by the
elements of the poset with one edge between two of them, when they
are related by a covering relation. See Figure~\ref{HasseDiagramm}
for an example.

Actually, $\mathcal{H(G)}$ is obtained from $\G$ be removing the
external edges and by replacing several edges between two vertices
by only one edge. Since $\mathcal{G}$ is connected and directed by
a flow, the Hasse diagram $\mathcal{H(G)}$ has the same
properties. A \emph{convex subset} \draftnote{definition coming
from combinatorics ref ?}$\mathcal{V'}$ of $\mathcal{V}$ is a set
of vertices of $\mathcal{G}$ such that for every pair $i\leq j$ in
$\mathcal{V'}$ the interval $[i,j]$ of $\Pi_{\G}$ is included in
$\mathcal{V'}$. If $\G$ is a connected graph of genus $0$, the set
of vertices of any connected subgraph of $\G$ is convex. This
property does not hold any more for connected graphs of higher
genus.

\begin{lem}
Let $\G$ be a connected directed graph without oriented loops and
let $\G'$ be a connected subgraph of $\G$. The set of vertices of
$\G'$ is convex if and only if the contraction of $\G'$ inside of
$\G$ gives a graph without oriented loops.
\end{lem}

\begin{figure}[h]
$$\G'=\vcenter{\xymatrix@R=6pt@C=6pt{
&&&&\\
& & & *+[F-,]{\ \textbf{2} \ } \ar@{-}[u] \ar@{-}[dr] \ar@{=}[ddll] \ar@{-}[ur] & \\
&&&& \\
& *+[F-,]{\ \textbf{3} \ } \ar@{-}[u] \ar@{-}[dr]\ar@{-}[dl]& & & \\
& & & & }} \quad \quad ;\quad \quad
\vcenter{\xymatrix@R=6pt@C=6pt{
& & & & \\
&*+[F-,]{\ \textbf{1} \ } \ar@{-}[ur] \ar@{-}[u] \ar@{-}[ul] \ar@{-}[dddd]  \ar@{-}[dr] & & & \\
&&&&\\
& & \vdots &  & \\
&&&& \\
& *+[F-,]{\ \textbf{3} \ } \ar@{=}[ur]  \ar@{-}[dr]\ar@{-}[dl]& & & \\
& & & & }}
  $$ \caption{Example of a admissible subgraph $\G'$ of $\G$ and an example of a non-admissible subgraph of $\G$.}
\label{AdmissibleSubgraph}
\end{figure}
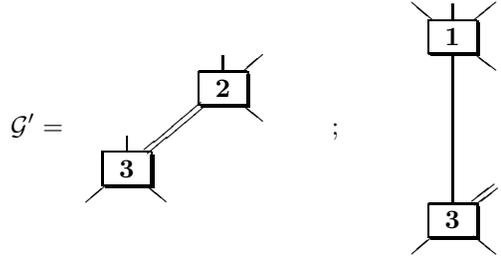

A connected subgraph $\G'$ with this property is called
\emph{admissible} in \cite{Granaker06}. We denote by $\G/\G'$ the
graph obtained by the contraction of $\G'$ inside $\G$. See
Figure~\ref{AdmissibleSubgraph} for an example of a admissible
subgraph and an example of a non-admissible subgraph of $\G$. By
extension, an admissible subgraph of a non-necessarily connected
graph is a union of admissible subgraphs (eventually empty) of
each connected component.

\subsection{Interpretation in terms of graphs}
\label{graph interpretation}

Let $\mu \, : \, \F^c(\oPo) \to \Po$ be a morphism of augmented dg
$\Sy$-bimodules. We denote by $\mu\big( \G(p_1, \ldots, p_n)\big)$
the image of an element $\G(p_1, \ldots, p_n)$ of
$\F^c(\oPo)^{(n)}$ under $\mu$. Let $\G'$ be a admissible subgraph
of $\G$ with $k$ vertices. Denote by $\G/\mu\G'(p_1,\ldots,p_n)$
the element of $\F^c(\oPo)^{(n-k+1)}$ obtained by composing
$\G'(p_{i_1},\ldots ,p_{i_k})$ in $\G(p_1, \ldots, p_n)$ under
$\mu$. When the $p_i$ and $\mu$ are not of degree zero, this
composition induces natural signs that we make explicit in the
sequel. Let start with a representative element of a class of
graph $\G(p_1, \ldots, p_n)$ whose vertices are indexed by
elements $p_i$, that is to say we have chosen an order between the
$p_i$ (see Section~\ref{Free properad}). The vertices of $\G'$ are
indexed by elements $p_{i_1},\ldots ,p_{i_k}$. We denote by
$J=(i_1,\ldots,i_k)$ the associated ordered subset of
$[n]=\{1,\ldots,n \}$ and $p_J=p_{i_1},\ldots ,p_{i_k}$. Since
$\G'$ is an admissible subgraph, its set of vertices forms a
convex subset of the set of vertices of $\G$ (or a disjoint union
of convex subsets if $\G$ is not connected). Therefore, it is
possible to change the order of the vertices of $\G$ such that the
vertices of $\G'$ are next to each others. That is there exists
two ordered subsets $I_1$ and $I_2$ of $[n]$ such that the
underlying subsets $I_1$, $I_2$ and $J$ without order form a
partition of $[n]$ and such that $\G(p_1, \ldots,
p_n)=(-1)^{\varepsilon_1}\G(P_{I_1}, P_J, P_{I_2})$. The sign
$(-1)^{\varepsilon_1}$ is given by the Koszul-Quillen sign rule
from the permutation of the $p_i$. Now we can apply $\mu$  to get
$$\G/\mu\G'(p_1,\ldots,p_n)=(-1)^{\varepsilon_1+\varepsilon_2}\G/\G'\big(
 P_{I_1}, \mu(\G'(P_J)), P_{I_2} \big), $$
where $\varepsilon_2=| P_{I_1}|.|\mu|$. It is an easy exercise to
prove that this definition of the signs does not depend on the
different choices.

\begin{lem}\label{Form of the coderivation}
Let $\nu$ be a map $\F^c(W) \to W$ of degree $-1$. The unique
coderivation $d_\nu \in \CoDer^{-1}_\Id(\F^c(W),\, \F^c(W))$ which
extends $\nu$ is given by
$$d_\nu\big(\G(w_1,\ldots, w_n)\big)=\sum_{\G'\subset \G}   \G/\nu\G'(w_1,\ldots,w_n),$$
where the sum runs over admissible subgraphs $\G'$ of $\G$.
\end{lem}

\begin{proo}
This formula defines a coderivation. Since the composite of
$d_\nu$ with the projection on $W$ is equal to $\nu$, we conclude
by the uniqueness property of coderivations of
Lemma~\ref{uniquecoDer}.
\end{proo}

\begin{pro}\label{Properadinfty}
A map $\mu \, : \, \F^c(\oPo) \to \Po$ defines a structure of
homotopy prop(erad) on the augmented dg $\Sy$-bimodule $\Po$ if
and only if, for every $\G(p_1,\ldots, p_n)$ in $\F^c(\oPo)$, we
have
$$\sum_{\G'\subset \G} (-1)^{\varepsilon(\G', p_1, \ldots, p_n)} \mu\big( \G/\mu\G'(p_1,\ldots,p_n)\big)=0,$$
where the sum runs over admissible subgraphs $\G'$ of $\G$.
\end{pro}

\begin{proo}
By definition, $\mu$ induces a structure of homotopy prop(erad) if
and only if  ${d}^2_{s\mu}=0$. This last condition holds if and
only if the composite $\textrm{proj}_{s\Po} \circ
{d}^2_{s\mu}=(s\mu) \circ d_{s\mu}$ is zero, where
$\textrm{proj}_{s\Po}$ is the projection on $s\Po$. From
Lemma~\ref{Form of the coderivation}, this is equivalent to
$$\sum_{\G'\subset \G}
(s\mu)\big(\G/(s\mu)\G'(sp_1,\ldots,sp_n)\big)=0,$$ where the sum
runs over admissible subgraphs $\G'$ of $\G$. Recall from
Section~\ref{suspension} that the signs between $(s\mu)$ and $\mu$
are
$$\mu\big(\G(p_1,\ldots,p_n)\big)=(-1)^{\varepsilon(p_1,\ldots,p_n)}s^{-1}(s\mu)\big(\G(s p_1,\ldots,s p_n)\big),$$
where
$\varepsilon(p_1,\ldots,p_n)=(n-1)|p_1|+(n-2)|p_2|+\cdots+|p_{n-1}|$.
Therefore, $\mu$ induces a structure of homotopy prop(erad) if and
only if
$$\sum_{\G'\subset \G} (-1)^{\varepsilon(\G', p_1, \ldots, p_n)} \mu\big( \G/\mu\G'(p_1,\ldots,p_n)\big)=0,$$
where $(-1)^{\varepsilon(\G', p_1, \ldots, p_n)}$ is product of
the sign coming the composition with $s\mu$ and the sign coming
from the formula between $\mu$ and $s\mu$.
\end{proo}

\begin{Rq}
In the case of associative algebras, the graphs involved are
ladders (branches, directed graphs just one incoming edge and one outgoing edge for each vertex) and we recover exactly the original definition
of J. Stasheff \cite{Stasheff63}.
\end{Rq}

Dually, we have the following characterization of homotopy
coprop(erad)s. Let $\G$ be a graph whose
$\textrm{i}^{\textrm{th}}$ vertex has $n$ inputs and $m$ outputs.
For every graph $\G'$ with $n$ inputs and $m$ outputs, denote by
$\G \circ_i \G'$ the graph obtained by inserting $\G'$ in $\G$ at
the place of the $\textrm{i}^{\textrm{th}}$ vertex.

\begin{pro}\label{Coproperadinfty}
A map $\Delta \, :\, \Co \to \F(\oC)$ defines a structure of
homotopy coprop(erad) on the augmented dg $\Sy$-bimodule $\Co$ if
and only if, for every $c\in \oC$, we have
$$\sum (-1)^{\rho(\G^2_i, c_1, \ldots, c_l)}  \G^1 \circ_i \G^2_i (c_1, \ldots, c_{i-1},
c'_1, \ldots, c'_k, c_{i+1}, \ldots, c_l)=0,$$ where the sum runs
over elements $\G^1(c_1, \ldots, c_l)$ and $\G^2_i(c'_1, \ldots ,
c'_k)$ such that $\Delta(c)=\sum \G^1(c_1, \ldots, c_l)$ and
$\Delta(c_i)=\sum \G^2_i(c'_1, \ldots , c'_k)$.
\end{pro}

\begin{proo}
By definition, $\Delta$ induces a structure of homotopy
coprop(erad) if and only if $\partial_{s^{-1}\Delta}^2=0$. Since
$\partial_{s^{-1}\Delta}$ is a derivation,
$\partial_{s^{-1}\Delta}^2=0$ is equivalent to
$\partial_{s^{-1}\Delta}\circ ({s^{-1}\Delta}) (s^{-1} c)=0$, for
every $c\in \oC$. Denote $(s^{-1}\Delta) (s^{-1}c)=\sum
\G^1(s^{-1} c_1, \ldots, s^{-1} c_l)$ and $(s^{-1} \Delta)
(s^{-1}c_i)=\sum \G^2_i(s^{-1}c'_1, \ldots, s^{-1}c'_k)$. By the
explicit formula for $\partial_{s^{-1}\Delta}$ given in
Lemma~\ref{uniqueDer} applied to $\rho=\Id_{\F(s^{-1}\oC)}$, we
have
\begin{eqnarray*}
\partial_{s^{-1}\Delta}\circ ({s^{-1}\Delta}) (s^{-1}
c) & =& \partial_{s^{-1}\Delta}\big(\sum  \G^1(s^{-1}c_1, \ldots,
s^{-1}c_l) \big)\\
&=& \sum   \G^1 \circ_i \G^2_i (s^{-1} c_1, \ldots, s^{-1}
c_{i-1}, s^{-1} c'_1, \ldots, s^{-1}c'_k, s^{-1} c_{i+1}, \ldots,
s^{-1}c_l) =0
\end{eqnarray*}
We get back to the map $\Delta$ with the formula
$$\Delta(c)=(-1)^{\varepsilon(c_1,\ldots, c_l)}\sum \G^1(c_1, \ldots, c_l) , $$
where
$\varepsilon(c_1,\ldots,c_l)=(l-1)|c_1|+(l-2)|c_2|+\cdots+|c_{l-1}|$.
We conclude as in proof of Proposition~\ref{Properadinfty}.
\end{proo}

\subsection{Homotopy non-symmetric (co)prop(erad)}\label{non-symmetric properad up to
homotopy}

It is straightforward to generalize the two previous subsections
to non-symmetric (co)prop(erad)s. One has just to consider
non-labelled graphs instead of graphs with leaves, inputs and
outputs labelled by integers. Therefore, there is a bar and a
cobar construction between non-symmetric dg prop(erad)s and
non-symmetric dg coprop(erad)s. The notion that will be used in
the sequel is the notion of \emph{homotopy non-symmetric
prop(erad)}. It is defined by a coderivation on the non-symmetric
cofree (connected) coprop(erad). Equivalently, we can describe it
in terms of non-labelled graphs like in
Proposition~\ref{Properadinfty}. The chain complex defining the
cohomology of a gebra over a non-symmetric prop(erad) has always such a
structure (see
Section~\ref{Deformation theory}).

\subsection{Homotopy properads and associated homotopy Lie
algebras}\label{oplusP1} It was proven in \cite{KapranovManin01}
that for any operad, $\Po=\{\Po(n)\}$, the vector space $\bigoplus_n
\Po(n)$ has  a natural structure of Lie algebra which descends to
the space of coinvariants $\bigoplus_n \Po(n)_{\bS_n}$, which is isomorphic to the space of
invariants $\bigoplus_n \Po(n)^{\bS_n}$. In
\cite{VanderLaan02} this result was generalized to homotopy
operads and the associated $L_\infty$-algebras. In this section,
we further extend the results of
\cite{KapranovManin01,VanderLaan02} from homotopy operads to
arbitrary homotopy prop(erad)s:
$\Po=\{\Po(m,n)\}$.\\

Recall that a structure of $L_\infty$-algebra on $\fg$ is given by
a square-zero coderivation on $\mathcal{S}^c(s\fg)$, where
$\mathcal{S}^c(s\fg)$ stands for the cofree cocommutative
coalgebra on the suspension of $\fg$. Hence, such a structure is
completely characterized by the image of the coderivation on
$s\fg$, $\mathcal{S}^c(s\fg)\to s\fg$. Equivalently, an
$L_\infty$-algebra is an algebra over the minimal (Koszul)
resolution of the operad $\Li$. We refer the reader to
Section~\ref{dg manifold, dg affine scheme} for more details on
$L_\infty$-algebras.\\

Let $\Po$ be an $\Sy$-bimodule. We denote by $\oplus \Po$ the
direct sum of all the components of $\Po$, that is $\bigoplus_{m,
n} \Po(m,n)$. We consider the map $\Theta \, : \,
\mathcal{S}^c(\oplus \Po)\to \F^c(\Po)$ defined by
$\Theta(p_1\odot \cdots \odot p_n):=\sum \G(p_1, \ldots, p_n)$,
where the sum runs over the classes of graphs under the action of
the automorphism group of the graph. This sum is finite and since a graph is a
quotient of a levelled graph (see Section~\ref{Free properad}),
the signs are well defined.

\begin{thm}\label{Prop->Homo Lie}
Let $\Po$ be a homotopy properad, the direct sum $\oplus \Po$ of
its components has an induced $L_\infty$-structure.
\end{thm}

\begin{proo}
We define the partial cotriple coproduct of a cofree coprop(erad)
by the composite~:
$$\Delta' \ : \ \F^c(V) \xrightarrow{\widetilde{\Delta}} \F^c(\F^c(V))
\epi \F^c(V, \underbrace{\F^c(V)}_1),$$ where $\F^c(V,
\underbrace{\F^c(V)}_1)$ represents graphs indexed by elements of
$V$ and one element of $\F^c(V)$. Similarly, we define the partial
cotriple coproduct of the cofree cocommutative coalgebra by :
$$ \delta' \ : \ \mathcal{S}^c(V) \xrightarrow{\widetilde{\delta}} \mathcal{S}^c(\mathcal{S}^c(V))
\epi \mathcal{S}^c(V, \underbrace{\mathcal{S}^c(V)}_1).$$

Let $s\mu \, : \F^c(s\oPo) \to s \oPo$ be a map of degree $-1$
defining a homotopy properad structure on $\Po$, that is the
following composite
$$\F^c(s\oPo) \xrightarrow{\Delta'} \F^c(s\oPo, \underbrace{\F^c(s\oPo)}_{1})
\xrightarrow{\F^c(s\oPo, s\mu)} \F^c(s\oPo) \xrightarrow{s\mu}
s\oPo$$ is zero. A map $l \, : \, \mathcal{S}^c(s\fg)\to s\fg$
induces a square-zero coderivation on $\mathcal{S}^c(s\fg)$ means
that the following composite is equal to zero
$$\mathcal{S}^c(s\fg) \xrightarrow{\delta'}
\mathcal{S}^c(s\fg, \underbrace{\mathcal{S}^c(s\fg)}_1)
\xrightarrow{\mathcal{S}^c(s\fg, l)} \mathcal{S}^c(s\fg)
\xrightarrow{l}  s\fg.$$ We define the induced $L_\infty$
structure by
$$l \ : \ \mathcal{S}^c(s(\oplus \oPo))\xrightarrow{\Theta} \F^c(s\oPo) \xrightarrow{s\mu} s\Po.$$
The relation of the $L_\infty$ structure for $l$ lifts to the
relation of the homotopy prop(erad) by the following commutative
diagram~:
$$\xymatrix@C=35pt{\mathcal{S}^c(s(\oplus \oPo)) \ar[d]^{\Theta} \ar[r]^(0.4){\delta'} &
\mathcal{S}^c(s(\oplus \oPo),\underbrace{\mathcal{S}^c(s(\oplus
\oPo))}_{1})
\ar[r]^(0.6){\mathcal{S}^c(s\oPo, l)}& \mathcal{S}^c(s(\oplus \oPo)) \ar[r]^{l}& s(\oplus \oPo)\\
 \F^c(s\oPo) \ar[r]^(0.45){\Delta'}& \F^c(s\oPo, \underbrace{\F^c(s\oPo)}_{1}) \ar[r]^(0.55){\F^c(s\oPo, s\mu)}&
 \F^c(s\oPo)\ar[ur]^{s\mu} \ , & }$$ which
concludes the proof.
\end{proo}

When $\Po$ is a (strict) prop(erad), the induced structure is the
(strict) Lie algebra coming from the anti-symmetrization of the
Lie-admissible algebra of Proposition~\ref{properad->Lie-admissible}.
Theorem~\ref{Prop->Homo Lie} generalizes the well-know fact that a homotopy (associative)
algebra is a homotopy Lie algebra by anti-symmetrization of the structure maps.
\\

The same statement holds for the space of coinvariants elements and the space of invariant elements.

\begin{thm}\label{Prop->Homo Lie On Co-In}
Let $\Po$ be a homotopy properad, the total space of coinvariant elements $\oplus \Po_\Sy$ and the
the total space of invariant elements $\oplus \Po^\Sy$  have an induced $L_\infty$-structure.
\end{thm}

\begin{proo}
We apply the same arguments as in the proof of Proposition~\ref{OnCoinvariants}.
\end{proo}

We prove below that the maps $\Po \to \oplus \Po$ and $\Po \to \oplus \Po^\Sy$ are functors for
the category of the homotopy prop(erad)s to one of  homotopy Lie algebras (see
Proposition~\ref{homotopy morphisms prop->Lie infinite}). The same
result holds for non-symmetric homotopy properads as well.

\subsection{Homotopy convolution prop(erad)}

In this section, we extend the definition of the convolution
prop(erad) to the homotopy case.

\begin{thm}\label{Convolutionproperadinfty}
When $(\Co,\Delta)$ is a (non-symmetric) homotopy coprop(erad) and $(\Po,\mu)$ is
a  (non-symmetric) prop(erad), the convolution prop(erad) $\Po^\Co=\Hom(\Co,\Po)$
is a homotopy (non-symmetric) prop(erad).

The same result holds when $\Co$ is a (non-symmetric) coprop(erad)
and $\Po$ a homotopy (non-symmetric) prop(erad).
\end{thm}

\begin{proo}
To an element $\G(f_1,\ldots,f_n)$ of $\F^c(\oPo^{\oC})^{(n)}$, we
consider the map $\widetilde{\G}(f_1,\ldots,f_n) \, : \,
\F(\oC)^{(n)} \to \F^c(\oPo)^{(n)}$ defined by
$\G'(c_1,\ldots,c_n)\mapsto (-1)^{\xi}\G(f_1(c_1),\ldots,
f_n(c_n))$ if $\G'\cong\G$ and $0$ otherwise, where
$\xi=\sum_{i=2}^n |f_i|(|c_1|+\cdots +|c_{i-1}|)$. We define maps
$\mu_n \, :\, \F^c(\oPo^{\oC})^{(n)} \to \Po^\Co$ by the formula
$$\mu_n\big(\G(f_1,\ldots,f_n)\big):=   \widetilde{\mu}_\Po
  \circ \widetilde{\G}(f_1,\ldots,f_n)  \circ \Delta_n.$$
The degree of $\Delta_n$ is $n-2$ and the degree of
 $\widetilde{\mu}_\Po$ is
zero. Therefore, the degree of $\mu_n$ is $n-2$.

The map $\mu$ verifies the relation of
Proposition~\ref{Properadinfty}
\begin{eqnarray*}
\sum_{\G'\subset \G} \pm\, \mu \big( \G/\mu\G'(f_1,\ldots,f_n)
\big) &=& \sum \pm\, \widetilde{\mu}_\Po\circ
\widetilde{\G/\G'}(f_1, \ldots, \mu_k(\G'(f_{i_1}, \ldots,
f_{i_k})), \ldots, f_n )\circ \Delta_{l}\\
&=&\sum \pm\, \widetilde{\mu}_\Po\circ \widetilde{\G/\G'}(f_1,
\ldots, \widetilde{\mu}_\Po \circ \widetilde{\G'}(f_{i_1}, \ldots,
f_{i_k})\circ \delta_k , \ldots, f_n )\circ \Delta_{l},
\end{eqnarray*}
where the sum runs over admissible subgraphs $\G'$ of $\G$. We
denote by $k$ the number of vertices of $\G'$ and $l=n-k+1$.We use
the generic notation $i$ for the new vertex of $\G/\G'$ obtained
after contracting $\G'$. For every element $c\in \oC$, we denote
by $\Delta(c)=\sum \G^1(c_1, \ldots,c_{l})$ and $\Delta(c_i)=\sum
\G^2_i(c'_1, \ldots,c'_k)$. The associativity of the product of
$\Po$ gives
\begin{eqnarray*}
&&\sum_{\G'\subset \G} (-1)^{\varepsilon(\G', f_1, \ldots, f_n)}\,
\mu \big( \G/\mu\G'(f_1,\ldots,f_n) \big)(c)=\\
&& \widetilde{\mu}_\Po \circ \widetilde{\G}(f_1,\ldots,f_n)\circ
\big( \sum (-1)^{\rho(\G^2_i, c_1,\ldots,c_l)}\, \G^1\circ_i
\G^2_i (c_1, \ldots, c'_1, \ldots, c'_k, \ldots, c_l)\big).
\end{eqnarray*}
Since $(\Co, \Delta)$ is a homotopy coprop(erad), the last term
vanishes by Proposition~\ref{Coproperadinfty}.

The same statement in the non-symmetric case is proven in the same
way and the dual statement also.
\end{proo}

\begin{Rq}
In the particular case when $\Co$ is a homotopy coalgebra and
$\Po$ an associative algebra, $\Hom(\Co,\Po)$ is a homotopy
algebra. In the same way, when $\Co$ is a homotopy operad and
$\Po$ an operad, $\Hom(\Co,\Po)$ is an homotopy operad (see Lemma
$5.10$ of \cite{VanderLaan02}).
\end{Rq}
\begin{thm}\label{conv Lie infty}
When $(\Co, \Delta)$ is a homotopy coprop(erad) and $(\Po, \mu)$
is a prop(erad) (or when $(\Co, \Delta)$ is a coprop(erad) and
$(\Po, \mu)$ is a homotopy prop(erad)), the total space of the convolution prop(erad)
$\Po^\Co=\Hom(\Co,\Po)$ is a  homotopy Lie algebras.

The total subspace  $\Hom^\Sy(\Co,\Po)$ of invariant elements is a sub-$L_\infty$-algebra.
\end{thm}

\begin{proo}
The first part is a direct corollary of Theorem\ref{Convolutionproperadinfty} and Theorem\ref{Prop->Homo Lie}.
Since the structure maps of this $\L_\infty$-algebra are composite of equivariant maps ($\Delta_n$, $\widetilde{\mu}_\Po$), the induce an $L_\infty$-algebra structure on the total space of $\Hom^\Sy(\Co,\Po)$. (This is similar to the one used in the proof of Proposition~\ref{Hom Lie-admissible}).
\end{proo}

In the latter case, the $L_\infty$-`operations' or homotopies are
explicitly given by the following formula. The image of $f_1,
\ldots, f_n \in \Hom^\Sy(\Co,\Po)$ under $l_n$, for $n>1$, is
given by
$$l_n(f_1,\ldots,f_n)=\sum_{\sigma \in \Sy_n}
(-1)^{\textrm{sgn}(\sigma, f_1, \ldots, f_n)}\widetilde{\mu}_\Po
\circ (f_{\sigma(1)}\otimes  \cdots \otimes f_{\sigma(n)}) \circ
\Delta_n,$$ where $(-1)^{\textrm{sgn}(\sigma, f_1, \ldots, f_n)}$
is the Koszul-Quillen sign appearing after permutating the $f_i$
with $\sigma$. The first `operation' $l_1$ is the differential, that is
$l_1(f):=D(f)=d_\Po\circ f - (-1)^{|f|}f \circ d_\Co$.\\

In this homotopy Lie algebra, the generalized Maurer-Cartan
equation is well defined since the formal infinite sum
$\displaystyle Q(\alpha):=\sum_{n\ge 1} \frac{1}{n!}\, l_n(\alpha,
\ldots, \alpha)$ is in fact equal to the composite
$D+\widetilde{\mu}_\Po \circ \F(\alpha) \circ \Delta$ in
$\Hom(\Co,\Po)$, when $\Co$ is a homotopy coprop(erad) and to
$D+{\mu} \circ \F(\alpha) \circ \widetilde{\Delta}_\Co$ when $\Po$
is a homotopy prop(erad). (See \ref{filtered} for the general
definition of filtered $L_\infty$-algebras).

\begin{dei}
Let $(\Co,\Delta)$ be a homotopy coprop(erad) and $(\Po,\mu)$ be a
 prop(erad) (or $(\Co,\Delta)$ a coprop(erad) and $(\Po,\mu)$ a homotopy
 properad). A morphism $\Co \xrightarrow{\alpha} \Po$, of degree
 $-1$, is called a \emph{twisting morphism} if it is a solution of
 the (generalized) Maurer-Cartan equation
 $$ Q(\alpha):=\sum_{n\ge 1} \frac{1}{n!}\, l_n(\alpha,
\ldots, \alpha)=0,$$ in the homotopy Lie algebra
$\Hom^\Sy(\Co,\Po)$. We denote this set by $\textrm{Tw}(\Co,
\Po)$.
\end{dei}

We can represent the bifunctor $\textrm{Tw}(-,\, -)$ in the same
as in the strict case (see Proposition~\ref{Bar-cobar
adjunction}).

\begin{pro}\label{Twisting homotopy morphism}
Let $(\Co,\Delta)$ be a homotopy coprop(erad) and $(\Po,\mu)$ be a
 prop(erad). There is a natural bijection
 $$\mathrm{Mor}_{\textrm{dg prop(erad)s}}\left(\Omega(\Co),\,
\Po\right)\cong \mathrm{Tw}(\Co, \Po) .$$

Let $(\Co,\Delta)$ be a coprop(erad) and $(\Po,\mu)$ be a homotopy
 prop(erad). There is a natural bijection
 $$\mathrm{Tw}(\Co, \Po)\cong
\mathrm{Mor}_{\textrm{dg coprop(erad)s}}\left(\Co,\,
\B(\Po)\right) .$$
\end{pro}

\begin{proo}
The proof is a direct generalization of the proof of
Proposition~\ref{Bar-cobar adjunction}.
\end{proo}

\subsection{Morphisms of homotopy (co)prop(erad)s}
\label{morphism of homotopy prop}

In this section, we recall the notion of morphism between two
homotopy properads due to \cite{Granaker06}. We extend it to homotopy (co)props and make them
explicit in terms of Maurer-Cartan elements in some convolution $L_\infty$-algebra.\\

Since a homotopy properad is equivalent to its associated
(generalized) bar construction, the notion of \emph{morphism of
homotopy properads} (or \emph{weak morphism}) is defined as
follows.

\begin{dei}\cite{Granaker06}
Let $\Po_1$ and $\Po_2$ be two homotopy prop(erad)s. A morphism
between $\Po_1$ and $\Po_2$ is a morphism of dg coprop(erad)s
between their bar constructions~: $B(\Po_1) \to B(\Po_2)$.
\end{dei}

A morphism of dg coprop(erad)s $\Phi\, :\,  B(\Po_1)=\F^c(s
\bar{\Po_1}) \to B(\Po_2)=\F^c(s \bar{\Po_2})$ is characterized by
its image on $s \bar{\Po_2}$. We denote by $s^{-1} \varphi\, :\,
B(\Po_1) \to  \bar{\Po_2}$ the composite of $\Phi$ with the
projection on $s \bar{\Po_2}$ followed by the desuspension. Notice
that the degree of $s^{-1} \varphi$ is $-1$. By
Proposition~\ref{Twisting homotopy morphism}, $\Phi$ is a morphism
of dg coprop(erad)s if and only if $s^{-1}\varphi$ is a
Maurer-Cartan element in $\Hom^\Sy(B(\Po_1), \Po_2)$, that is
$$Q(s^{-1} \varphi)= \sum_{n\ge 1} \frac{1}{n!}\, l_n(s^{-1}\varphi, \ldots,
s^{-1}\varphi)=D(s^{-1} \varphi) + \mu_{\Po_2} \circ \F^c(s^{-1}
\varphi) \circ \widetilde{\Delta}=0,$$ where $\widetilde{\Delta}$
is the coproduct map $B(\Po_1)=\F^c(s \bar{\Po}_1) \to
\F^c(\F^c(s\bar{\Po}_1))$.

\begin{pro}
A morphism of $\Sy$-bimodules $\varphi\, :\, B(\Po_1) \to s
\bar{\Po_2}$ induces a morphism of homotopy properads between
$\Po_1$ and $\Po_2$ if and only if $s^{-1} \varphi$ is a
Maurer-Cartan element in the $L_\infty$-algebra
$\Hom^\Sy(B(\Po_1), \Po_2)$, that is $Q(s^{-1} \varphi)=0$.
\end{pro}

Like in Section~\ref{graph interpretation}, we make explicit the
above definition in terms of graphs.

\begin{pro}
A map $s^{-1} \varphi \, : \, B(\Po_1) \to \bar{\Po}_2$ is a
morphism of homotopy prop(erad)s if and only if, for every class
of graphs $\G$ under the action of the automorphism group, the
following relation holds

$$\sum s\mu^{\Po_2}_k \big(\G/\varphi\G_1\sqcup \ldots \sqcup \varphi \G_k\big) =
\sum \varphi \big(\G/(s\mu^{\Po_1})\G'\big), $$ where the first
sum runs over all partition of the graph $\G$ into admissible
subgraphs $\G_1 \sqcup  \ldots \sqcup \G_k$ and where the second
sum runs over all admissible subgraph $\G'$ of $\G$. Once again,
the signs are induced by Koszul-Quillen rule, when apply to
elements $sp_1, \ldots, s p_n$, such that $n$ is the number of
vertices of $\G$.
\end{pro}

\begin{proo}
The map $s^{-1} \varphi \, :\,  B(\Po_1) \to \bar{\Po}_2$ induces
a unique morphism of coprop(erad)s $\Phi \, : \, B(\Po_1) \to
B(\Po_2)$ which commutes with the differentials if and only if the
above relation is verified. (The left hand term is the projection
on $\bar{\Po}_2$ of the composite $d_{B(\Po_2)}   \circ \Phi$ and
the right hand term is the projection on the same space of the
composite $\Phi \circ d_{B(\Po_1)}$, that is $\varphi \circ
d_{B(\Po_1)}$.)

\end{proo}

When applied to $A_\infty$-algebras, the underlying graphs are
ladders and this proposition gives the classical notion of weak
morphisms, that is morphisms between $A_\infty$-algebras. \\

Dually, we define the notion of morphism between homotopy
coprop(erad)s.

\begin{dei}
Let $\Co_1$ and $\Co_2$ be two homotopy prop(erad)s. A morphism
between $\Co_1$ and $\Co_2$ is a morphism of dg prop(erad)s
between their cobar constructions~: $\Omega(\Co_1) \to
\Omega(\Co_2)$.
\end{dei}

A morphism of dg prop(erad)s $\Psi\, :\,  \Omega(\Co_1)=\F(s^{-1}
\bar{\Co_1}) \to \Omega(\Co_2)=\F(s^{-1} \bar{\Co_2})$ is
characterized by the image of $s^{-1} \bar{\Co_1}$. We denote by
$s^{-1} \psi\, :\, \bar{\Co_1}\to \Omega(\Co_2)$ the desuspension
of the restriction of $\Psi$ on $s^{-1} \bar{\Co_1}$. By
Proposition~\ref{Twisting homotopy morphism}, $\Psi$ is a morphism
of dg prop(erad)s if and only if $s^{-1} \psi$ is a twisting
morphism in $\Hom^\Sy(\bar{\Co}_1,\Omega(\Co_2))$, that is
$$Q(s^{-1}\psi)= \sum_{n\ge 1} \frac{1}{n!}\, l_n(s^{-1}\psi,
\ldots, s^{-1}\psi)=D(s^{-1}\psi)+\widetilde{\mu} \circ
\F(s^{-1}\psi) \circ {\Delta}_{\Co_1}=0,$$ where $\widetilde{\mu}$
is the composition map $\F(\Omega(\Co_2))=\F(\F(s^{-1}
\bar{\Co}_2)) \to \F(s^{-1} \bar{\Co}_2)=\Omega(\Co_2)$.

\begin{pro}
A morphism of $\Sy$-bimodules $\psi\, :\, s^{-1} \bar{\Co_1}\to
\Omega(\Co_2)$ induces a morphism of homotopy coproperads between
$\Co_1$ and $\Co_2$ if and only if $s^{-1}\psi$ is a Maurer-Cartan
element in the $L_\infty$-algebra $\Hom^\Sy(\Co_1,\Omega(\Co_2))$,
that is $Q(s^{-1}\psi)=0$.
\end{pro}



We now prove that the convolution prop(erad) is a construction
functorial with respect to the first argument.

\begin{thm}\label{functoriality convolution properad}
Let $\Psi$ be a morphism of homotopy coprop(erad)s between $\Co_1$
and $\Co_2$. Let $\Po$ be a prop(erad). There exists a natural
morphism of homotopy prop(erad)s between $\Hom(\Co_2, \Po)$ and
$\Hom(\Co_1, \Po)$  induced by $\Psi$.
\end{thm}
The same statement holds in the non-symmetric case.
\begin{proo}
Let $\Psi$ denote the morphism of dg prop(erads) $\Omega(\Co_1)
\to \Omega(\Co_1)$ and $s^{-1}\psi$ the induced twisting morphism
$\bar{\Co_1}\to \Omega(\Co_2)$, that is $Q(s^{-1}\psi)=0$. We
define the morphism of coprop(erad)s $\Phi \,: \, B(\Hom(\Co_2,
\Po)) \to B(\Hom(\Co_1, \Po))$ by its  image $\varphi$ on $s\,
\overline{\Hom}(\Co_1, \Po)=\Hom(s^{-1}\bar{\Co}_1, \bar{\Po})$ as
follows. Let $\G(f_1, \ldots,\, f_n)\in B(\Hom(\Co_2,
\Po))=\F^c(s\, \overline{\Hom}(\Co_1,
\Po))=\F^c(\Hom(s^{-1}\bar{\Co}_1, \bar{\Po}))$. The image of
$\G(f_1, \ldots,\, f_n)$ under $\varphi$ is equal to the composite
$$\varphi(\G(f_1, \ldots,\, f_n)) \ : \ s^{-1}\bar{\Co}_1
\xrightarrow{\psi}
\F(s^{-1}\bar{\Co_2})\xrightarrow{\widetilde{\G}(f_1, \ldots,\,
f_n)} \F(\oPo)\xrightarrow{\widetilde{\mu}_\Po}\Po.$$ It remains
to prove that $s^{-1} \varphi$ is a twisting element in
$\Hom(B(\Hom(\Co_2, \Po)), \Hom(\Co_1, \Po)) $ , that is $Q(s^{-1}
\varphi)=0$. By the definition of $Q$ in this homotopy prop(erad)
and by the `associativity' of $\widetilde{\mu}_\Po$, $Q(s^{-1}
\varphi)(\G(f_1, \ldots, f_n))$ is equal to the composite
$$ \bar{\Co}_1 \xrightarrow{\Delta_{\Co_1}}
\F(\bar{\Co_1}) \xrightarrow{\F(s^{-1}\psi)}
\F(\F(s^{-1}\bar{\Co_2})) \xrightarrow{\widetilde{\mu}}
\F(s^{-1}\bar{\Co_2}) \xrightarrow{\widetilde{\G}(f_1, \ldots,\,
f_n)} \F(\oPo)\xrightarrow{\widetilde{\mu}_\Po}\Po,$$ where
$\widetilde{\mu}$ is the `triple' map associated to the free
prop(erad) $\F(s^{-1}\bar{\Co_2})$. Therefore
$Q(s^{-1}\varphi)(\G(f_1, \ldots, f_n))=\widetilde{\mu}_\Po \circ
\widetilde{\G}(f_1, \ldots,\, f_n) \circ Q(s^{-1}\psi)$ which
vanishes since $Q(s^{-1} \psi)=0$.
\end{proo}

The dual statement is also true and can be proved in the same way.
It will appear in a future work of the second author in relation with the
 the transfer of algebraic structures up to homotopy through a deformation-retract (homological perturbation lemma).

\begin{pro}\label{homotopy morphisms prop->Lie infinite}
The constructions given in Theorem~\ref{Prop->Homo
Lie} and Theorem~\ref{Prop->Homo Lie On Co-In} provide us with three functors,
$$
\mbox{\sf Category of homotopy properads} \lon \mbox{\sf Category
of homotopy Lie algebras}.
$$
\end{pro}

\begin{proo}
Let $\Phi \, : B(\Po_1) \to B(\Po_2)$ be a morphism of
coprop(erad)s defining a morphism of homotopy prop(erad)s between
$\Po_1$ and $\Po_2$. The associated projection $\varphi$ verifies
$Q(s^{-1}\varphi)=0$, that is
$$\F^c(s\oPo_1) \xrightarrow{\widetilde{\Delta}} \F^c(\F^c(s\oPo_1)) \xrightarrow{\F^c(s^{-1}\varphi)}
\F^c(\oPo_2) \xrightarrow{\mu_{\Po_2}} \Po_2$$ equals $0$. We
define the following map
$$f \ : \ \mathcal{S}^c(s(\oplus \oPo_1)) \xrightarrow{\Theta} \F^c(s \oPo_1)
\xrightarrow{\varphi} s(\oplus \oPo_2).$$ The map $f$ is a
morphism of $L_\infty$-algebras. Its desuspension $s^{-1}f$
verifies the Maurer-Cartan equation in the $L_\infty$-algebra
$\Hom(\mathcal{S}^c(s(\oplus \oPo_1)), \oplus \oPo_2)$ (see
\cite{Dolgushev07}). The Maurer-Cartan equation for $s^{-1}f$
lifts to the Maurer-Cartan equation for $s^{-1}\varphi$ via
$\Theta$, that is the following diagram is commutative
$$\xymatrix@C=35pt{\mathcal{S}^c(s(\oplus \oPo_1)) \ar[d]^{\Theta} \ar[r]^(0.4){\widetilde{\delta}} &
\mathcal{S}^c(\mathcal{S}^c(s(\oplus \oPo_1)))
\ar[r]^(0.6){\mathcal{S}^c(s^{-1}f)}& \mathcal{S}^c(\oplus \oPo_2) \ar[r]^{l_{\oplus \Po_2}}& \oplus \oPo_2\\
 \F^c(s\oPo_1) \ar[r]^(0.45){\widetilde{\Delta}}& \F^c(\F^c(s\oPo_1)) \ar[r]^(0.55){\F^c(s^{-1}\varphi)}&
 \F^c(\oPo_2)\ar[ur]^{\mu_{\Po_2}} \ , & }$$ which concludes the proof.
\end{proo}

\begin{cor}
\label{functoriality convolution L-infinite} Let $\Psi$ be a
morphism of homotopy coprop(erad)s between $\Co_1$ and $\Co_2$.
Let $\Po$ be a prop(erad). There exists a natural morphism of
$L_\infty$-algebras between $\Hom(\Co_2, \Po)$ and $\Hom(\Co_1,
\Po)$ induced by $\Psi$. Its restriction to $\Hom^\Sy(\Co_2, \Po)$ gives a
 natural morphism of
$L_\infty$-algebras between $\Hom^\Sy(\Co_2, \Po)$ and $\Hom^\Sy(\Co_1,
\Po)$.
\end{cor}

\begin{proo}
The first part is a direct corollary of Theorem~\ref{functoriality convolution
properad} and Proposition~\ref{homotopy morphisms prop->Lie
infinite}. Since these constructions are composite of equivariant maps, they are stable on the space
of invariant elements $\Hom^\Sy(\Co_2, \Po)$ and $\Hom^\Sy(\Co_1,
\Po)$.
\end{proo}


\section{Models}

In this section, we recall the definitions of \emph{minimal} and
\emph{quadratic model} for properads and we formally extend them
to props. Recall that a model is a quasi-free resolution. Our
viewpoint here is to classify properads according to the form of
their minimal model, when it exists. For instance, a properad is
\emph{Koszul} if and only if it admits a quadratic model. To
clarify the genus of some resolutions, we introduce the notion of
\emph{contractible prop(erad)s}. Such properads have genus $0$
quadratic models.

\subsection{Minimal models}

Recall that a \emph{quasi-free} prop(erad) is a (dg) prop(erad)
whose underlying $\Sy$-bimodule, that is forgetting the differential map, is a free prop(erad) $\F(M)$. It
is not necessarily a free dg prop(erad) since the differential
$\partial$ may not be freely generated by the differential of $M$.

\begin{dei}[Model]
Let $\Po$ be a prop(erad). A \emph{model} of $\Po$ is a quasi-free
prop(erad) $(\F(M),\, \partial)$ equipped with a quasi-isomorphism
$\F(M) \xrightarrow{\sim} \Po$.
\end{dei}

Theorem~\ref{BarCobarResolution} proves that every augmented
prop(erad) has a canonical model given by the bar-cobar
construction. Some prop(erad)s admit more simple models. The
differential $\partial$ of a quasi-free prop(erad) $\F(M)$ is by
definition a derivation. Lemma~\ref{uniqueDer} shows that it is
characterized by its restriction $\partial_M \, :\,  M \to \F(M)$
on $M$.

\begin{dei}[Decomposable differential]
The differential $\partial$ of a quasi-free prop(erad) is called
\emph{decomposable} if the image of its restriction to $M$,
$\partial_M \, : \, M \to \F(M)$, is composed by decomposable
elements, that is $\textrm{Im}(\partial_M)\subset\bigoplus_{n\ge
2} \F(M)^{(n)}.$
\end{dei}

\begin{dei}[Minimal model]
A model $(\F(M),\, \partial)$ is called \emph{minimal} if its
differential $\partial$ is decomposable.
\end{dei}


\subsection{Form of minimal models}

From Theorem~\ref{BarCobarResolution}, we know that every
augmented (dg) properad admits a resolution of the form
$\Omega(B(P))$. A natural way to get a minimal model from this
would be to consider the homology of the bar construction, try to
endow it with a structure of homotopy coproperad and then take the
generalized cobar construction of it. In this section, we prove
that when minimal models exist, they are of this form.

\begin{pro}
\label{HomologyBarQuasiFree} Let $(\F(M), \, \partial)$ be a
quasi-free properad with a decomposable differential generated by
a non-negatively graded $\Sy$-module $M$. Then the homology of the
bar construction $B(\F(M))$ of $(\F(M), \,
\partial)$ is equal to the suspension of $M$.

\end{pro}

\begin{proo}
The bar construction of the dg-properad $\Po:=\F(M)$ is defined by
the underlying $\Sy$-bimodule $B(\Po):=\F^c(s \bar{\Po})=\F^c(s \,
\bar{\F}(M))$. The differential $d$ is the sum of two terms
$d_0+\tilde{\partial} $. The component $\tilde{\partial}$ comes
from $\partial$ and $d_0$ is the unique coderivation which
extends the partial product of $\F(M)$.

Consider the filtration $F_s:= \bigoplus_{r\le s} \F^c\big(s
\bar{\F}(M)\big)_r$, where $r$ is the sum of the degrees of the
elements of $M$. Let's denote by $E^\bullet_{st}$ the associated
spectral sequence.

Since the chain complex $M$ is bounded below, this filtration is
bounded below $F_{-1}=0$. It is obviously exhaustive, therefore
the classical theorem of convergence of spectral sequences shows
that $E^\bullet$ converges to the homology of $B(\F(M))$.

We have $\tilde{\partial}(F_s) \subset F_{s-1}$ and
$d_0(F_s)\subset F_s$. Hence, the first term of the spectral
sequence is $E^0_{st}=\F^c_{s+t}\big(s \bar{\F}(M)\big)_s$, where
$s+t$ is the total homological degree, and $d^0=d_0$. We have
reduced the problem to computing the homology of the bar
construction of the free properad on $M$, which is equal to
$\Sigma M$ by Corollary~$5.10$ of \cite{Vallette03} (where we
choose to put each element of $M$ in weight $1$).
\end{proo}

The next proposition shows that, when a minimal model of a
properad $\Po$ exists, it is necessarily given by a quasi-free
properad on the homology of the bar construction of $\Po$.

\begin{thm}\label{from of minimal model}
Let $\Po$ be an augmented dg properad and let $(\F(M),\,
\partial)$ be a minimal model of $\Po$. The $\Sy$-bimodule $sM$ is
isomorphic to the homology of the bar construction of $\Po$.
\end{thm}

\begin{proo}
In \cite{Vallette03}, we proved in Proposition 4.9 that the bar
construction preserves quasi-isomorphisms. Therefore, the bar
construction of $\F(M)$ is quasi-isomorphic to the bar
construction of $\Po$. We conclude by
Proposition~\ref{HomologyBarQuasiFree}.
\end{proo}

We denote by $\Po^{\ac}:=H_\bullet(B(\Po))$ the homology of the
bar construction of $\Po$. When $(\F(s^{-1}\Po^{\ac}),\,
\partial)$ is a minimal model of $\Po$, the derivation
${\partial}$ is equivalent to a structure of homotopy coproperad
on $\Po^{\ac}$ such that $\delta_1=0$. That is
$(\F(s^{-1}\Po^{\ac}),\,
\partial)$ is the generalized cobar construction
$\Omega(\Po^{\ac})$ of the homotopy coproperad $\Po^{\ac}$. As a
conclusion, we have that following corollary which gives the form
of minimal models.

\begin{cor}\label{form minimal model}
A minimal model of an augmented dg-properad $\Po$ is always the
 cobar construction $\Omega(\Po^{\ac})$ on the homology
of $B(\Po)$ endowed with a structure of homotopy coproperad.
\end{cor}

In the sequel, we will only consider props freely generated by a
properad, in the sense of the horizontal (concatenation) product.
The minimal model of such props is given by the generalized cobar
construction of the associated homotopy coproperad, viewed as a
homotopy coprop. And the result of the preceding lemma still
holds.

\subsection{Quadratic models and Koszul duality theory}
In general, it is a difficult problem to find the minimal model of
a prop(erad). One can first compute the homology of the bar
construction and then provide a structure of homotopy coproperad
on it, that is with higher homotopy cooperations. For some weight
graded properads, there exist simple minimal models which are
given by the Koszul duality theory. These properads are called
Koszul.

\begin{dei}[Quadratic differential]
The differential $\partial$ of a quasi-free prop(erad) is called
\emph{quadratic} if the image of $\partial_M \, : \,  M \to \F(M)$
is in $\F(M)^{(2)}.$
\end{dei}

\begin{dei}[Quadratic model]
A model $(\F(M),\, \partial)$ is called \emph{quadratic} if its
differential $\partial$ is quadratic.
\end{dei}

When $\Po$ is a weight graded properad, its bar construction
splits as a direct sum of finite chain complexes indexed by the
weight (cf. \cite{Vallette03} Section~$7.1.1$). In this case, we
can talk about top dimensional homology groups.

\begin{thm}\label{Koszul=quadratic model}
Let $\Po$ be a weight graded properad concentrated in homological
degree $0$. The following assertions are equivalent.
\begin{enumerate}
\item The homology of $B(\Po)$ is concentrated in top dimension.
\item The $\Sy$-bimodule $\Po^{\ac}$ is a strict coproperad. \item
The properad $\Po$ admits a quadratic model~: $\Omega(\Po^{\ac})
\xrightarrow{\sim} \Po$.
\end{enumerate}
\end{thm}

\begin{proo}
$(1)\Rightarrow (2)$ is given by Proposition 7.2 of
\cite{Vallette03}.

$(2)\Rightarrow (3)$ is given by Theorem~$5.9$ of
\cite{Vallette03}. When $\Po^{\ac}$ has a structure of strict
coproperad, its cobar construction is a resolution of $\Po$ and
the differential of it is quadratic.

$(3)\Rightarrow (1)$ Since $\Po$ is isomorphic to
$\F(M_0)/(\partial(M_1))$, which $\partial$ quadratic, this
presentation is quadratic. Define an extra weight on $M$ by the
formula $\omega(M_n):=n+1$. With this weight, the
quasi-isomorphism $\F(M) \xrightarrow{\rho} \Po$ is a morphism of
weight graded dg-properads. The induced morphism $B(\rho)$ on the
bar construction preserves this grading. Therefore we have
$H_n(B(\Po)^{(n)})=H_n(B(\F(M))^{(n)})=(sM)_n$ and the homology of
the bar construction of $\Po$ is concentrated in top dimension.
\end{proo}

In this case, the properad $\Po$ is called a \emph{Koszul}
properad. The coproperad $\Po^{\ac}$ is its \emph{Koszul dual} and
$\Po$ has a quadratic model which is the cobar construction on
$\Po^{\ac}$. In other words, a properad is Koszul when its bar
construction is \emph{formal}, that is when $B(\Po)$ is
quasi-isomorphic to its homology $\Po^{\ac}$ as a dg-coproperad.
This case is simple and particularly efficient. When
$\Po=F(V)/(R)$ has a quadratic presentation with a finite
dimensional space of generators $V$, then the linear dual (twisted
by the signature representation) of the coproperad $\Po^{\ac}$ is
a properad equal, up to suspension, to
$\Po^{!}=F(V^\vee)/(\R^\perp)$ where $V^\vee$ is the linear dual
of $V$ twisted by the signature representation. This relation
provides a concrete method to compute the minimal model of Koszul
properads. The next step is to be able to prove that it is Koszul.
Koszul duality theory provides a smaller chain complex $\Po^{\ac}
\bc \Po$ which is acyclic if and only if the properad $\Po$ is
Koszul. Therefore, there are simple methods to show that a
properad is Koszul. When a properad is defined by two Koszul
properads with a distributive law, Proposition~$8.4$ of
\cite{Vallette03} shows that it is Koszul. In the operadic case,
there are basically two other efficient methods. First if the
homology of the free $\Po$-algebra is acyclic then the operad
$\Po$ is Koszul (see Proposition~$5.3.5$ of \cite{Fresse04}).
Finally, when the operad is set theoretic, we can use the
associated poset to prove that it is Koszul (see
\cite{Vallette06}).

\subsection{Homotopy Koszul properads}
 \label{SectionHoKoszul}
 If a properad is Koszul, then we have clearly cut means to construct its minimal model.
However, the ordinary notion of Koszulness does not cover many
important examples. For example, the properad of associative
bialgebras is not Koszul since it is not quadratic and any Koszul
properad has a quadratic presentation by Corollary~$7.5$ of
\cite{Vallette03}. So we are left in such cases with no concrete
methods of proving that a particular properad $\Po$ admits a
minimal model,  and, if so, constructing it explicitly.  It is
already a highly non-trivial problem in general to find the set of
generators for a  minimal model, not speaking about the
differential. In this section we extend the notion of Koszulness
in such a way that some of the above
problems become effectively solvable.

\begin{dei}
Let $\Po= \cF( V )/ (\Ro)$ be a properad generated by an
$\Sy$-bimodule $V=\{V(m,n)\}$ concentrated in degree zero, and
with  an ideal generated by $\Ro\subset \cF( V )^{(\geq 2)}$. Let
$\pi_k: \cF( V) \rar
 \cF( V )^{(k)}$ be the natural projection, and let us set,
 $$
 \Ro_k:= \pi_k(\Ro), \ \ \ \mbox{for}\ k=2,3,\ldots.
 $$
Let us also denote by $\Po^{(\geq k)}$ the image of $ \cF( V
)^{(\geq k)}$ under the natural epimorphism  $\cF( V ) \epi \Po$.

The properad $\Po$ is called {\em homotopy Koszul}\, if
\begin{itemize}
\item[(i)] the quadratic properad $\Po_2:= \cF( V )/ (\Ro_2)$ is
Koszul, \item[(ii)] $\Po$ and $\Po_2$ are isomorphic as
$\Sy$-bimodules, \item[(iii)] there is an extra grading on the
properad $\Po=\oplus_{\lambda}\Po(\lambda)$, with $\Po(\lambda)$
being a collection of finite-dimensional $\Sy$-bimodules.
\end{itemize}
\end{dei}

In practice the conditions (i)-(iii) above are often not hard to
check (see examples below). As an extra grading one can use, for
example, the path grading of a free properad introduced by
Kontsevich and studied in \cite{MarklVoronov03}. The main
motivation behind the definition is the following.

\begin{thm}\label{homotopykoszul}
If a properad $\Po$ is homotopy Koszul, then it admits a minimal
model of the form $(\cF (s\bar{\Po}_2^{\ac}), \delta)$, where
$\Po_2^{\ac}$ is the coproperad Koszul dual to $\Po_2$.
\end{thm}

\begin{proo} Consider the bounded above increasing filtration $F_{-p}\Po:=\Po^{(\geq p)}$ of the properad $\Po$.
As $F_{-p}\Po\cap P(\lambda)$ are finite-dimensional vector spaces, the spectral sequences associated with this filtration
(see below) have good convergence properties.
Since $\Po$ is isomorphic to $\Po_2$ as an $\Sy$-bimodule,
the associated graded properad,
$$
\bigoplus_{p\geq 0} \frac{\Po^{(\geq p)}}{\Po^{(\geq p+1)}},
$$
is isomorphic to $\Po_2$ as a properad. Then we have,

\sip

{\bf Claim 1}. {\em The homologies of the bar constructions,
$B(\Po)$ and $B(\Po_2)$, are isomorphic as $\Sy$-bimodules, i.e.\
$H_\bullet(B(\Po))\simeq \Po_2^{\ac}$ as $\Sy$-bimodules}.

\sip

\noindent Indeed, the filtration $F_{-p}\Po:=\Po^{(\geq p)}$ induces an associated
filtration of the complex $B(\Po)$ (as differential in $B(\Po)$ is built from compositions
in $\Po$ which respect the filtration $F_{-p}\Po$). By the above observation,
 the $0$th term, $E^0$, of the associated spectral sequence, $\{E^r, d^r\}$,
is exactly the complex $B(\Po_2)$,
$E^0_{pq}=B(\Po_2)^{(-p)}_{p+q}$ and $d^0=d_{B(\Po_2)}$. As
$\Po_2$ is Koszul, $E^1=H_\bullet(B(\Po_2))$ is exactly the Koszul
dual coproperad $\Po_2^{\ac}$, that is $E^1_{pq}=0$ for $q\neq
-2p$ and
$E^1_{pq}=H_{-p}\big(B(\Po_2)^{(-p)}\big)=(\Po_2^{\ac})^{(-p)}$
when $q=-2p$. The induced differentials, $d^r$ for $r\geq 1$, are
zero because of the homological degree 0 assumption on $\Po$. Thus
the spectral sequence $\{E^r, d^r\}$ degenerates at the first
term. The extra grading on the properad $\Po$ induces an extra
grading $\lambda$ on $B(\Po)$ which makes $F_{-p}(B(\Po))\cap
B(\Po)(\lambda)$ into a bounded filtration of $B(\Po)(\lambda)$.
Hence it converges to $H_\bullet(B(\Po))(\lambda)$ by the
Classical Convergence Theorem~$5.5.1$ of \cite{Weibel}, thereby
proving Claim 1.

\sip

Choosing a homological splitting of the complex $B(\Po)$,
$$\xymatrix{  H_\bullet(B(\Po)) \ar@<1ex>[r]^(0.5){i} & *{\
B(\Po)\quad \ \  } \ar@(dr,ur)[]_h \ar@<1ex>[l]^(0.5){p}   },$$
one can use dual transfer formulae of \cite{Granaker06} for
homotopy coproperads to induce on the $\Sy$-bimodule
$H_\bullet(B(\Po))\simeq \Po_2^{\ac}$ the associated  strongly
homotopy coproperad structure, that is a
 differential, $\delta$, in the free properad
 \footnote{In fact the Gran{\aa}ker formulae provide us in general
 with a  differential $\delta$ in a {\em completed}\, (with respect to the number of vertices) free properad:
 there is no guarantee
 that such  $\delta$ applied to a generator is a {\em finite}\ sum of terms but we can only be sure that $\delta$ is
 continuous  with respect to the topology induced by the number of vertices  filtration. However, our assumption
 on existence of an extra gradation in $\Po$ implies that $\delta$ is well-defined in the ordinary category of properads:
 it is  {\em finite}\,  on every generator so that
  $(\cF( s^{-1}H_\bullet(\bar{B}(\Po))), \delta)$ makes sense
  without completion.

 \sip

 It is important to notice that had we chosen to work with topological properads (with topology induced by the number
 of vertices or genus filtrations), the condition (iii) in the definition of homotopy Koszulness can be safely omitted
 --- Theorem~\ref{homotopykoszul} stays true in the category of (completed) topological properads because all the spectral sequences
 we used in the proof stay convergent by classical Complete Convergence Theorem 5.5.10 (see p.139 in \cite{Weibel}).
 As an example of the deformation quantization prop \cite{Merkulov04}  shows, working with topological prop(erad)s is unavoidable
 in application of the theory of prop(erad)s to geometry and mathematical physics.
 }
 $\cF\big(s^{-1}H_\bullet(\bar{B}(\Po))\big)=
  \Omega\big(H_\bullet({B}(\Po))\big)$ generated by $H_\bullet(B(\Po))$.
 In general, this differential
 is {\em not}\, quadratic, i.e.\ the induced homotopy coproperad structure on $H_\bullet(B(\Po))$ is {\em not}\, equal
 to the coproperad structure on $\Po_2^{\ac}$. Moreover, the
 chosen homological splitting provides us canonically with a morphism
 of homotopy coproperads which extends $i$,
 $$
 H_\bullet(B(\Po)) \lon B(\Po),
 $$
i.e.\ with a morphism of dg properads,
$$
\phi: (\cF( s^{-1}H_\bullet(\bar{B}(\Po))), \delta) \lon \Omega
(B(\Po)).
$$
As $\Omega (B(\Po))\stackrel{\simeq}{\lon} \Po$ is a resolution
of
$\Po$ by Theorem~\ref{BarCobarResolution}, the required
Theorem~\ref{homotopykoszul} follows immediately from the
following

\sip

{\bf Claim 2}.  {\em  Under the assumption on the properad $\Po$ the morphism $\phi$ is a quasi-isomorphism}.

\sip

Indeed, the introduced above filtration of the bar construction,
$B(\Po)$, induces a filtration, $F_{-p}H_\bullet(B(\Po))$, of its
homology  with the associated  graded coproperad being exactly
$\Po_2^{\ac}$. This filtration of $H_\bullet(B(\Po))$ induces in
turn a filtration of the complex $(\cF(
s^{-1}H_\bullet(\bar{B}(\Po))), \delta)$. The $0$th term of the
associated spectral sequence is precisely the minimal model,
$(\cF( s^{-1}\bar{\Po}_2^{\ac}), \delta)$, of the properad
$\Po_2$. As the latter is Koszul by assumption, its homology is
equal to $\Po_2$. By homological degree assumption on $\Po$, the
induced differential on the next term of the spectral sequence
vanishes so that it degenerates.  The extra grading assumption on
$\Po$ implies that this spectral sequences converges to the
homology $(\cF( s^{-1}H_\bullet(\bar{B}(\Po))), \delta)$ which is
equal, therefore, as an $\Sy$-bimodule to $\Po_2\simeq \Po$. This
fact completes the proof of Claim 2 and hence of the Theorem.
\end{proo}

\bip

The operad $\Po_2$ is Koszul means that  the differential of the
minimal model $(\Omega(\Po^{\ac}_2), \delta_2)$ is quadratic, that
is $\delta_2 \, :\, s^{-1}\bar{\Po}_2^{\ac} \to
\F(s^{-1}\bar{\Po}_2^{\ac})^{(2)}$. Since the transfer of homotopy
coproperad structures does not change the map $\Delta_2$ defining
the homotopy coproperad structure on $H_\bullet(B(\Po))$ but just
add extra terms $\Delta_n$, for $n\geq 3$, the final differential
$\delta$ defining a minimal model of $\Po$ is equal to $\delta_2$
plus extra terms $\delta_n$ for $n\geq 3$ such that $\delta_n \,
:\, s^{-1}\bar{\Po}_2^{\ac} \to \F(s^{-1}\bar{\Po}_2^{\ac})^{(n)}
$, that is to say, $\delta$ is a perturbation of $\delta_2$. \\

The coproperad $\Po_2^{\ac}$ is computable by Koszul duality
theory. Therefore the above Theorem gives an immediate estimate of
the set of generators for a minimal model of a homotopy Koszul
properad. Moreover, the differential in this quasi-free model can
in principle be computed via ordinary  homotopy transfer formulae.

\sip

The class of properads which are homotopy Koszul but not  Koszul
is non-empty and contains an important example of the properad,
$\ab$, of (co)associative bialgebras which  can be defined as a
quotient,
$$
\ab:= {\cF(V)}/(\Ro)
$$
of the free properad, $\cF(V)$, generated by the $\bS$-bimodule
$V=\{V(m,n)\}$,
\[
V(m,n):=\left\{ \Ba{rr} \bK[\bS_2]\ot \bK[\bS_1]
\equiv\mbox{span}\left\langle
\begin{xy}
 <0mm,-0.55mm>*{};<0mm,-2.5mm>*{}**@{-},
 <0.5mm,0.5mm>*{};<2.2mm,2.2mm>*{}**@{-},
 <-0.48mm,0.48mm>*{};<-2.2mm,2.2mm>*{}**@{-},
 <0mm,0mm>*{\bullet};<0mm,0mm>*{}**@{},
 <0mm,-0.55mm>*{};<0mm,-3.8mm>*{_1}**@{},
 <0.5mm,0.5mm>*{};<2.7mm,2.8mm>*{^2}**@{},
 <-0.48mm,0.48mm>*{};<-2.7mm,2.8mm>*{^1}**@{},
 \end{xy}
\, ,\,
\begin{xy}
 <0mm,-0.55mm>*{};<0mm,-2.5mm>*{}**@{-},
 <0.5mm,0.5mm>*{};<2.2mm,2.2mm>*{}**@{-},
 <-0.48mm,0.48mm>*{};<-2.2mm,2.2mm>*{}**@{-},
 <0mm,0mm>*{\bullet};<0mm,0mm>*{}**@{},
 <0mm,-0.55mm>*{};<0mm,-3.8mm>*{_1}**@{},
 <0.5mm,0.5mm>*{};<2.7mm,2.8mm>*{^1}**@{},
 <-0.48mm,0.48mm>*{};<-2.7mm,2.8mm>*{^2}**@{},
 \end{xy}
   \right\rangle  & \mbox{if}\ m=2, n=1,\vspace{3mm}\\
\bK[\bS_1]\ot \bK[\bS_2]\equiv \mbox{span}\left\langle
\begin{xy}
 <0mm,0.66mm>*{};<0mm,3mm>*{}**@{-},
 <0.39mm,-0.39mm>*{};<2.2mm,-2.2mm>*{}**@{-},
 <-0.35mm,-0.35mm>*{};<-2.2mm,-2.2mm>*{}**@{-},
 <0mm,0mm>*{\bullet};<0mm,0mm>*{}**@{},
   <0mm,0.66mm>*{};<0mm,3.4mm>*{^1}**@{},
   <0.39mm,-0.39mm>*{};<2.9mm,-4mm>*{^2}**@{},
   <-0.35mm,-0.35mm>*{};<-2.8mm,-4mm>*{^1}**@{},
\end{xy}
\, ,\,
\begin{xy}
 <0mm,0.66mm>*{};<0mm,3mm>*{}**@{-},
 <0.39mm,-0.39mm>*{};<2.2mm,-2.2mm>*{}**@{-},
 <-0.35mm,-0.35mm>*{};<-2.2mm,-2.2mm>*{}**@{-},
 <0mm,0mm>*{\bullet};<0mm,0mm>*{}**@{},
   <0mm,0.66mm>*{};<0mm,3.4mm>*{^1}**@{},
   <0.39mm,-0.39mm>*{};<2.9mm,-4mm>*{^1}**@{},
   <-0.35mm,-0.35mm>*{};<-2.8mm,-4mm>*{^2}**@{},
\end{xy}
\right\rangle
\ & \mbox{if}\ m=1, n=2, \vspace{3mm}\\
0 & \mbox{otherwise}, \Ea \right.
\]
representing a binary product and a binary coproduct without
symmetries, modulo the ideal generated by relations
\[
\Ro:\ \ \
\begin{xy}
 <0mm,0mm>*{\bullet};<0mm,0mm>*{}**@{},
 <0mm,-0.49mm>*{};<0mm,-3.0mm>*{}**@{-},
 <0.49mm,0.49mm>*{};<1.9mm,1.9mm>*{}**@{-},
 <-0.5mm,0.5mm>*{};<-1.9mm,1.9mm>*{}**@{-},
 <-2.3mm,2.3mm>*{\bullet};<-2.3mm,2.3mm>*{}**@{},
 <-1.8mm,2.8mm>*{};<0mm,4.9mm>*{}**@{-},
 <-2.8mm,2.9mm>*{};<-4.6mm,4.9mm>*{}**@{-},
   <0.49mm,0.49mm>*{};<2.7mm,2.3mm>*{^3}**@{},
   <-1.8mm,2.8mm>*{};<0.4mm,5.3mm>*{^2}**@{},
   <-2.8mm,2.9mm>*{};<-5.1mm,5.3mm>*{^1}**@{},
 \end{xy}
\ - \
\begin{xy}
 <0mm,0mm>*{\bullet};<0mm,0mm>*{}**@{},
 <0mm,-0.49mm>*{};<0mm,-3.0mm>*{}**@{-},
 <0.49mm,0.49mm>*{};<1.9mm,1.9mm>*{}**@{-},
 <-0.5mm,0.5mm>*{};<-1.9mm,1.9mm>*{}**@{-},
 <2.3mm,2.3mm>*{\bullet};<-2.3mm,2.3mm>*{}**@{},
 <1.8mm,2.8mm>*{};<0mm,4.9mm>*{}**@{-},
 <2.8mm,2.9mm>*{};<4.6mm,4.9mm>*{}**@{-},
   <0.49mm,0.49mm>*{};<-2.7mm,2.3mm>*{^1}**@{},
   <-1.8mm,2.8mm>*{};<0mm,5.3mm>*{^2}**@{},
   <-2.8mm,2.9mm>*{};<5.1mm,5.3mm>*{^3}**@{},
 \end{xy}
\ \ \ \  , \ \ \ \
 \begin{xy}
 <0mm,0mm>*{\bullet};<0mm,0mm>*{}**@{},
 <0mm,0.69mm>*{};<0mm,3.0mm>*{}**@{-},
 <0.39mm,-0.39mm>*{};<2.4mm,-2.4mm>*{}**@{-},
 <-0.35mm,-0.35mm>*{};<-1.9mm,-1.9mm>*{}**@{-},
 <-2.4mm,-2.4mm>*{\bullet};<-2.4mm,-2.4mm>*{}**@{},
 <-2.0mm,-2.8mm>*{};<0mm,-4.9mm>*{}**@{-},
 <-2.8mm,-2.9mm>*{};<-4.7mm,-4.9mm>*{}**@{-},
    <0.39mm,-0.39mm>*{};<3.3mm,-4.0mm>*{^3}**@{},
    <-2.0mm,-2.8mm>*{};<0.5mm,-6.7mm>*{^2}**@{},
    <-2.8mm,-2.9mm>*{};<-5.2mm,-6.7mm>*{^1}**@{},
 \end{xy}
\ - \
 \begin{xy}
 <0mm,0mm>*{\bullet};<0mm,0mm>*{}**@{},
 <0mm,0.69mm>*{};<0mm,3.0mm>*{}**@{-},
 <0.39mm,-0.39mm>*{};<2.4mm,-2.4mm>*{}**@{-},
 <-0.35mm,-0.35mm>*{};<-1.9mm,-1.9mm>*{}**@{-},
 <2.4mm,-2.4mm>*{\bullet};<-2.4mm,-2.4mm>*{}**@{},
 <2.0mm,-2.8mm>*{};<0mm,-4.9mm>*{}**@{-},
 <2.8mm,-2.9mm>*{};<4.7mm,-4.9mm>*{}**@{-},
    <0.39mm,-0.39mm>*{};<-3mm,-4.0mm>*{^1}**@{},
    <-2.0mm,-2.8mm>*{};<0mm,-6.7mm>*{^2}**@{},
    <-2.8mm,-2.9mm>*{};<5.2mm,-6.7mm>*{^3}**@{},
 \end{xy}
\ \ \ \ , \ \ \ \
 \begin{xy}
 <0mm,2.47mm>*{};<0mm,-0.5mm>*{}**@{-},
 <0.5mm,3.5mm>*{};<2.2mm,5.2mm>*{}**@{-},
 <-0.48mm,3.48mm>*{};<-2.2mm,5.2mm>*{}**@{-},
 <0mm,3mm>*{\bullet};<0mm,3mm>*{}**@{},
  <0mm,-0.8mm>*{\bullet};<0mm,-0.8mm>*{}**@{},
<0mm,-0.8mm>*{};<-2.2mm,-3.5mm>*{}**@{-},
 <0mm,-0.8mm>*{};<2.2mm,-3.5mm>*{}**@{-},
     <0.5mm,3.5mm>*{};<2.8mm,5.7mm>*{^2}**@{},
     <-0.48mm,3.48mm>*{};<-2.8mm,5.7mm>*{^1}**@{},
   <0mm,-0.8mm>*{};<-2.7mm,-5.2mm>*{^1}**@{},
   <0mm,-0.8mm>*{};<2.7mm,-5.2mm>*{^2}**@{},
\end{xy}
\ - \
\begin{xy}
 <0mm,0mm>*{\bullet};<0mm,0mm>*{}**@{},
 <0mm,-0.49mm>*{};<0mm,-3.0mm>*{}**@{-},
 <-0.5mm,0.5mm>*{};<-3mm,2mm>*{}**@{-},
 <-3mm,2mm>*{};<0mm,4mm>*{}**@{-},
 <0mm,4mm>*{\bullet};<-2.3mm,2.3mm>*{}**@{},
 <0mm,4mm>*{};<0mm,7.4mm>*{}**@{-},
<0mm,0mm>*{};<2.2mm,1.5mm>*{}**@{-},
 <6mm,0mm>*{\bullet};<0mm,0mm>*{}**@{},
 <6mm,4mm>*{};<3.8mm,2.5mm>*{}**@{-},
 <6mm,4mm>*{};<6mm,7.4mm>*{}**@{-},
 <6mm,4mm>*{\bullet};<-2.3mm,2.3mm>*{}**@{},
 <0mm,4mm>*{};<6mm,0mm>*{}**@{-},
<6mm,4mm>*{};<9mm,2mm>*{}**@{-}, <6mm,0mm>*{};<9mm,2mm>*{}**@{-},
<6mm,0mm>*{};<6mm,-3mm>*{}**@{-},
   <-1.8mm,2.8mm>*{};<0mm,7.8mm>*{^1}**@{},
   <-2.8mm,2.9mm>*{};<0mm,-4.3mm>*{_1}**@{},
<-1.8mm,2.8mm>*{};<6mm,7.8mm>*{^2}**@{},
   <-2.8mm,2.9mm>*{};<6mm,-4.3mm>*{_2}**@{},
 \end{xy}.
\]

These relations stand respectively for the associativity of the
product, the coassociativity of the coproduct and the relation
between them, that is the coproduct is a morphism of algebras or
equivalently the product is a morphism of coalgebras. As the ideal
contains 4-vertex graphs, the properad $\ab$ is not quadratic.
Hence $\ab$ can not be Koszul in the ordinary sense. However, we
have  the following

\begin{pro}
The properad $\ab$ is homotopy Koszul.
\end{pro}

\begin{proo}
(i) The properad $\ab_2$ is Koszul as it is generated by the
bimodule $V$ with the relations,
\[
\begin{xy}
 <0mm,0mm>*{\bullet};<0mm,0mm>*{}**@{},
 <0mm,-0.49mm>*{};<0mm,-3.0mm>*{}**@{-},
 <0.49mm,0.49mm>*{};<1.9mm,1.9mm>*{}**@{-},
 <-0.5mm,0.5mm>*{};<-1.9mm,1.9mm>*{}**@{-},
 <-2.3mm,2.3mm>*{\bullet};<-2.3mm,2.3mm>*{}**@{},
 <-1.8mm,2.8mm>*{};<0mm,4.9mm>*{}**@{-},
 <-2.8mm,2.9mm>*{};<-4.6mm,4.9mm>*{}**@{-},
   <0.49mm,0.49mm>*{};<2.7mm,2.3mm>*{^3}**@{},
   <-1.8mm,2.8mm>*{};<0.4mm,5.3mm>*{^2}**@{},
   <-2.8mm,2.9mm>*{};<-5.1mm,5.3mm>*{^1}**@{},
 \end{xy}
\ - \
\begin{xy}
 <0mm,0mm>*{\bullet};<0mm,0mm>*{}**@{},
 <0mm,-0.49mm>*{};<0mm,-3.0mm>*{}**@{-},
 <0.49mm,0.49mm>*{};<1.9mm,1.9mm>*{}**@{-},
 <-0.5mm,0.5mm>*{};<-1.9mm,1.9mm>*{}**@{-},
 <2.3mm,2.3mm>*{\bullet};<-2.3mm,2.3mm>*{}**@{},
 <1.8mm,2.8mm>*{};<0mm,4.9mm>*{}**@{-},
 <2.8mm,2.9mm>*{};<4.6mm,4.9mm>*{}**@{-},
   <0.49mm,0.49mm>*{};<-2.7mm,2.3mm>*{^1}**@{},
   <-1.8mm,2.8mm>*{};<0mm,5.3mm>*{^2}**@{},
   <-2.8mm,2.9mm>*{};<5.1mm,5.3mm>*{^3}**@{},
 \end{xy}
\ \ \ \  , \ \ \ \
 \begin{xy}
 <0mm,0mm>*{\bullet};<0mm,0mm>*{}**@{},
 <0mm,0.69mm>*{};<0mm,3.0mm>*{}**@{-},
 <0.39mm,-0.39mm>*{};<2.4mm,-2.4mm>*{}**@{-},
 <-0.35mm,-0.35mm>*{};<-1.9mm,-1.9mm>*{}**@{-},
 <-2.4mm,-2.4mm>*{\bullet};<-2.4mm,-2.4mm>*{}**@{},
 <-2.0mm,-2.8mm>*{};<0mm,-4.9mm>*{}**@{-},
 <-2.8mm,-2.9mm>*{};<-4.7mm,-4.9mm>*{}**@{-},
    <0.39mm,-0.39mm>*{};<3.3mm,-4.0mm>*{^3}**@{},
    <-2.0mm,-2.8mm>*{};<0.5mm,-6.7mm>*{^2}**@{},
    <-2.8mm,-2.9mm>*{};<-5.2mm,-6.7mm>*{^1}**@{},
 \end{xy}
\ - \
 \begin{xy}
 <0mm,0mm>*{\bullet};<0mm,0mm>*{}**@{},
 <0mm,0.69mm>*{};<0mm,3.0mm>*{}**@{-},
 <0.39mm,-0.39mm>*{};<2.4mm,-2.4mm>*{}**@{-},
 <-0.35mm,-0.35mm>*{};<-1.9mm,-1.9mm>*{}**@{-},
 <2.4mm,-2.4mm>*{\bullet};<-2.4mm,-2.4mm>*{}**@{},
 <2.0mm,-2.8mm>*{};<0mm,-4.9mm>*{}**@{-},
 <2.8mm,-2.9mm>*{};<4.7mm,-4.9mm>*{}**@{-},
    <0.39mm,-0.39mm>*{};<-3mm,-4.0mm>*{^1}**@{},
    <-2.0mm,-2.8mm>*{};<0mm,-6.7mm>*{^2}**@{},
    <-2.8mm,-2.9mm>*{};<5.2mm,-6.7mm>*{^3}**@{},
 \end{xy}
\ \ \ \ , \ \ \ \
 \begin{xy}
 <0mm,2.47mm>*{};<0mm,-0.5mm>*{}**@{-},
 <0.5mm,3.5mm>*{};<2.2mm,5.2mm>*{}**@{-},
 <-0.48mm,3.48mm>*{};<-2.2mm,5.2mm>*{}**@{-},
 <0mm,3mm>*{\bullet};<0mm,3mm>*{}**@{},
  <0mm,-0.8mm>*{\bullet};<0mm,-0.8mm>*{}**@{},
<0mm,-0.8mm>*{};<-2.2mm,-3.5mm>*{}**@{-},
 <0mm,-0.8mm>*{};<2.2mm,-3.5mm>*{}**@{-},
     <0.5mm,3.5mm>*{};<2.8mm,5.7mm>*{^2}**@{},
     <-0.48mm,3.48mm>*{};<-2.8mm,5.7mm>*{^1}**@{},
   <0mm,-0.8mm>*{};<-2.7mm,-5.2mm>*{^1}**@{},
   <0mm,-0.8mm>*{};<2.7mm,-5.2mm>*{^2}**@{},
\end{xy},
\]
which verify the Distributive Law (see Section~\ref{Distributuive law} and Proposition~$8.5$ of
\cite{Vallette03}).

 (ii) The $\Sy$-bimodule
isomorphism $\ab\simeq \ab_2$ was established in
\cite{EnriquezEtingof03}.

(iii) The ideal generated by $\Ro$ preserves the path grading (see
\cite{MarklVoronov03} for its definition and main properties) of
the free properad $\cF(V)$ and hence induces an associated
filtration on $\ab$ which satisfies the last condition in the
definition of a homotopy Koszulness properad.
\end{proo}

\begin{cor}(cf.\ \cite{Markl06})\label{generators}
The properad $\ab$ admits a minimal resolution, $\cF( \Co)$,
generated by the $\Sy$-bimodule
 $\Co=\{\Co(m,n)\}_{m,n\geq 1, m+n\geq 3}$, with
\[
\Co(m,n):= s^{m+n-3}\bK[\bS_m]\ot
\bK[\bS_n]=\mbox{span}\left\langle
\begin{xy}
 <0mm,0mm>*{\bullet};<0mm,0mm>*{}**@{},
 <0mm,0mm>*{};<-8mm,5mm>*{}**@{-},
 <0mm,0mm>*{};<-4.5mm,5mm>*{}**@{-},
 <0mm,0mm>*{};<-1mm,5mm>*{\ldots}**@{},
 <0mm,0mm>*{};<4.5mm,5mm>*{}**@{-},
 <0mm,0mm>*{};<8mm,5mm>*{}**@{-},
   <0mm,0mm>*{};<-8.5mm,5.5mm>*{^1}**@{},
   <0mm,0mm>*{};<-5mm,5.5mm>*{^2}**@{},
   <0mm,0mm>*{};<4.5mm,5.5mm>*{^{m\hspace{-0.5mm}-\hspace{-0.5mm}1}}**@{},
   <0mm,0mm>*{};<9.0mm,5.5mm>*{^m}**@{},
 <0mm,0mm>*{};<-8mm,-5mm>*{}**@{-},
 <0mm,0mm>*{};<-4.5mm,-5mm>*{}**@{-},
 <0mm,0mm>*{};<-1mm,-5mm>*{\ldots}**@{},
 <0mm,0mm>*{};<4.5mm,-5mm>*{}**@{-},
 <0mm,0mm>*{};<8mm,-5mm>*{}**@{-},
   <0mm,0mm>*{};<-8.5mm,-6.9mm>*{^1}**@{},
   <0mm,0mm>*{};<-5mm,-6.9mm>*{^2}**@{},
   <0mm,0mm>*{};<4.5mm,-6.9mm>*{^{n\hspace{-0.5mm}-\hspace{-0.5mm}1}}**@{},
   <0mm,0mm>*{};<9.0mm,-6.9mm>*{^n}**@{},
 \end{xy}
\right\rangle .
\]
\end{cor}
\begin{proo}
The Koszul dual properad of $\ab_2$ is the properad generated by a
binary product and a binary coproduct which are associative and
coassociative. All the composites with the product and the
coproduct vanish except $\vcenter{\xymatrix@M=0pt@R=6pt@C=6pt{\ar@{-}[dr] &  &\ar@{-}[dl]  \\  &\ar@{-}[d] & \\
& \ar@{-}[dl]\ar@{-}[dr]& \\ & & }}$. The only non-vanishing
element of this properad are obtained by composing first some
products and then coproducts. We conclude that
 $\ab_2^{\ac}(m,n)= s^{m-2}\bK[\bS_m]\ot
s^{n-2}\bK[\bS_n]$  for $m,n\geq 1, m+n\geq 3$ and zero otherwise. Then Theorem~\ref{homotopykoszul}
implies the claim.
\end{proo}

We refer the reader to Section~\ref{Homotopy Morphism=MC} for
another application of the notion of homotopy Koszulness.

\subsection{Models for associative algebras, non-symmetric operads, operads, properads, props}
\label{ComparisonModels}

There are several different notions of algebraic objects in the
literature that are used to model the operations acting on some
algebraic category. We briefly recall them in the following
table.\\

\begin{tabular}{|p{2.2cm}|c|c|c|c|c|}
\hline \textbf{Operations}  &
$\vcenter{\xymatrix@M=0pt@R=6pt@C=6pt{&& \\  & \ar@{-}[d]&  \\
& \bullet \ar@{-}[d] & \\ & & \\ &&}}$ &
$\vcenter{\xymatrix@M=0pt@R=6pt@C=6pt{&& \\ \ar@{-}[dr]  & \ar@{-}[d]& \ar@{-}[dl]  \\
& \bullet \ar@{-}[d] & \\ & & \\ & &}  } \atop \textrm{no
symmetry}$ &
$\vcenter{\xymatrix@M=0pt@R=6pt@C=6pt{&& \\ \ar@{-}[dr]  & \ar@{-}[d]& \ar@{-}[dl]  \\
& \bullet \ar@{-}[d] & \\ & & \\ &&}} $
&  $\vcenter{\xymatrix@M=0pt@R=6pt@C=6pt{&& \\ \ar@{-}[dr]  & \ar@{-}[d]& \ar@{-}[dl]  \\
& \bullet \ar@{-}[dl] \ar@{-}[dr] & \\ & & \\ &&}}$ &
 $\vcenter{\xymatrix@M=0pt@R=6pt@C=6pt{&& \\ \ar@{-}[dr]  & \ar@{-}[d]& \ar@{-}[dl]  \\
& \bullet \ar@{-}[dl] \ar@{-}[dr] & \\ & & \\ &&}}$ \\
\hline \textbf{Composition}  &  $\vcenter{\xymatrix@M=0pt@R=6pt@C=6pt{&& \\  & \ar@{-}[d]&  \\
\ar@{..}[r]& \bullet \ar@{-}[d]\ar@{..}[r]  & \\\ar@{..}[r] & \bullet \ar@{-}[d] \ar@{..}[r]& \\ & & \\
&&}}$ & $\vcenter{\xymatrix@M=0pt@R=6pt@C=6pt{&&&&&& \\
\ar@{-}[dr] & \ar@{-}[d]& \ar@{-}[dl] & \ar@{-}[d] & \ar@{-}[dr]
& & \ar@{-}[dl] \\
\ar@{..}[r] & \bullet \ar@{-}[drr] \ar@{..}[rr] & &
\bullet \ar@{-}[d] \ar@{..}[rr] & &\bullet   \ar@{-}[dll] \ar@{..}[r] &  \\
\ar@{..}[rrr] & & & \bullet \ar@{-}[d] \ar@{..}[rrr]& &  &   \\
& & & & & &   } } \atop \textrm{planar}$ &
$\vcenter{\xymatrix@M=0pt@R=6pt@C=6pt{&&&&&& \\ \ar@{-}[dr]  &
\ar@{-}[d]& \ar@{-}[dl] & \ar@{-}[d] & \ar@{-}[dr]
& & \ar@{-}[dl] \\
\ar@{..}[r] & \bullet \ar@{-}[drr] \ar@{..}[rr] & &
\bullet \ar@{-}[d] \ar@{..}[rr] & &\bullet   \ar@{-}[dll] \ar@{..}[r] &  \\
\ar@{..}[rrr] & & & \bullet \ar@{-}[d] \ar@{..}[rrr]& &  &   \\
& & & & & &   } } \atop \textrm{non-planar}$  &
$\vcenter{\xymatrix@M=0pt@R=6pt@C=6pt{&&&& \\ \ar@{-}[dr] & & \ar@{-}[dl] & \ar@{-}[d] & \\
\ar@{..}[r]&\bullet \ar@{-}[ddrr]   \ar@{-}@/_/[dd] \ar@{..}[rr]
&&\bullet \ar@{-}[ddll]\ar@{-}@/^/[dd] \ar@{..}[r] &
\\ &&&&
\\ \ar@{..}[r] &\bullet \ar@{-}[d] \ar@{..}[rr]&&\bullet \ar@{-}[dl] \ar@{-}[d] \ar@{-}[dr] \ar@{..}[r] &
\\&&&&
\\&&&& }}$ & $\vcenter{\xymatrix@M=0pt@R=6pt@C=6pt{
&&&& &&\\
 \ar@{-}[dr] & & \ar@{-}[dl] & \ar@{-}[d] & \ar@{-}[dr]&\ar@{-}[d]&\ar@{-}[dl]\\
\ar@{..}[r]&\bullet \ar@{-}[ddrr]   \ar@{-}@/_/[dd] \ar@{..}[rr]
&&\bullet \ar@{-}[ddll]\ar@{-}@/^/[dd] \ar@{..}[rr] &&\bullet
\ar@{-}@/_/[dd]
\ar@{-}@/^/[dd] \ar@{..}[r] &\\
&&&&&&\\
 \ar@{..}[r] &\bullet \ar@{-}[d] \ar@{..}[rr]&&\bullet
\ar@{-}[dl] \ar@{-}[d] \ar@{-}[dr] \ar@{..}[rr] &&\bullet \ar@{-}[d] \ar@{..}[r]&\\
&&&&&&\\
&&&& &&}}$ \\
\hline \textbf{Monoidal} \textbf{category} & $(\textrm{Vect},
\otimes)$ & $(\textrm{gVect}, \circ)$ & $(\Sy\textrm{-Mod},
\circ)$ & $(\Sy\textrm{-biMod}, \boxtimes_c)$ &
$(\Sy\textrm{-biMod}, \boxtimes)$ \\
\hline \textbf{Monoid}  &  ${\textrm{Associative} \atop
\textrm{algebras}}$ & ${\textrm{Non-symmetric} \atop
\textrm{operads}}$ & $\textrm{Operads}$   &
$\textrm{Properads}$ & $\textrm{Props}$  \\
\hline \textbf{Modules}  & $\textrm{Modules}$  &
${\textrm{Non-symmetric} \atop \textrm{algebras}}$  &
$\textrm{Algebras}$ & $\textrm{(Bial)gebras}$  &  $\textrm{(Bial)gebras}$ \\
\hline \textbf{Free monoid}  & ${\textrm{Ladders} \atop
\textrm{(Tensor module)}}$ &  ${\textrm{Planar} \atop
\textrm{trees}}$  & ${\textrm{Trees}}$ & ${\textrm{Connected} \atop \textrm{graphs}}$ &  ${\textrm{Graphs}}$ \\
\hline
\end{tabular}

\vskip8pt
 To each pair of such objects, there is a forgetful
functor and a left adjoint~:
$$\xymatrix{ \textrm{Associative algebras} \ar@<1ex>[r]  & \ar@<1ex>[l] \textrm{Non-symmetric operads}
 \ar@<1ex>[r] & \ar@<1ex>[l] \textrm{Operads} \ar@<1ex>[r] & \ar@<1ex>[l] \textrm{Properads}
  \ar@<1ex>[r] & \ar@<1ex>[l] \textrm{Props}.  }$$
Let us make them explicit.

\begin{itemize}
\item To any prop $\Po$, the associated properad
$\mathcal{U}^{\textrm{props}}_{\textrm{properads}}(\Po)$ is given
by the same underlying $\Sy$-bimodule where we  only consider
vertical compositions of operations based on connected graphs.
That is we forget the horizontal composition. Its left adjoint
$\F^{\textrm{props}}_{\textrm{properads}}(\Po)$ is given by the
free symmetric tensor on $\Po$ for the horizontal tensor product.
(This functor was introduced in \cite{Vallette03} at Section~$1$
where it is denoted by $\mathcal{S}$.) In other words, we freely
generate horizontal compositions from a properad to get a prop.\\

\item The operad obtained from a properad $\Po$ is the
$\Sy$-module
$\UU^{\textrm{properads}}_{\textrm{operads}}(\mathcal{P})(n):=\mathcal{P}(1,n)$
equipped with the restriction to one rooted trees composition. Its
left adjoint functor is
$\F^{\textrm{properads}}_{\textrm{operads}}(\Po)(m,n):=\Po(n)$ for
$m=1$ and $0$ for $m>1$.

\item For any operad $\Po$, we consider the non-symmetric operad
$\UU^{\textrm{operads}}_{\textrm{non-symm. operads}}(\Po)=\Po$
where we forget the action of the symmetric group. The left
adjoint is given by $$\F^{\textrm{operads}}_{\textrm{non-symm.
operads}}(\Po)(n)=\Po(n)\otimes \KK[\Sy_n].$$ (see M. Aguiar and
M. Livernet \cite{AguiarLivernet05})

\item The pair of adjoint functors between associative algebras
and non-symmetric operads is defined in the same then the pair of
functors between operads and properads. In one way, we just
consider the unital operation (arity $(1)$) of a non-symmetric
operad. In this other way, for an associative algebra we define a
non-symmetric operad concentrated in arity $(1)$.
\end{itemize}

\begin{pro}
All these functors are exact, that is the image of a
quasi-isomorphism is a quasi-isomorphism.
\end{pro}

\begin{proo}
It is trivial for the forgetful functors and for the functors
$F^{\textrm{operads}}_{\textrm{operads}}$ and
$F^{\textrm{non-symm. operads}}_{\textrm{ass. algebras}}$ because
the underlying dg-module does not change. Since the functor
$F^{\textrm{operads}}_{\textrm{non-symm. operads}}$ is given by
tensoring $\Sy_n$-modules with the flat $\KK$ module $\KK[\Sy_n]$
(the characteristic of $\KK$ is $0$), it is exact. Over a field of
characteristic $0$, the functor
$F^{\textrm{props}}_{\textrm{properads}}$ is also exact .
\end{proo}

This proposition justifies the following philosophy. To study the
deformation theory of elements of an algebraic category, that is a
class of gebras (modules, algebras, bialgebras), one should first
model this category using the simplest possible object of the
previous table. For instance, associative, diassociative,
dendriform algebras \cite{LodayFrabettiChapotonGoichot01} are
encoded each time by a non-symmetric operad. Commutative, Lie,
PreLie, Gerstenhaber, Poisson algebras are modelled by operads.
Lie bialgebras, infinitesimal Hopf algebras \cite{Aguiar00},
(associative) bialgebras (see \ref{Associative bialgebras}) are
representations of properads. Non-unital infinitesimal Hopf
algebras, Semi Hopf algebras, Lily bialgebras \cite{Loday06} can
only be represented by a prop.

Then to study the deformation theory of this algebraic category,
that is to define the stable notion up to homotopy (see
\ref{Homotopy gebra}) or the deformation complex (see
\ref{Deformation theory}), one has to find a cofibrant resolution
(bar-cobar, minimal model for instance) of the related operad,
properad or prop $\Po$. This resolution contains all the necessary
data since a resolution for the induced prop is ``freely''
obtained by the free exact functor.

\subsection{Models generated by genus $0$ differentials}

Let $\mathcal{A}$ be a category of gebras defined by some products
and some coproducts with relations that can be written as linear
combinations of connected graphs of genus $0$, for example Lie
bialgebras, Frobenius bialgebras, Infinitesimal bialgebras (see
\cite{Gan03, Vallette03}). In this case, the class of gebras can
be faithfully modelled with a
smaller algebraic object called a dioperad \cite{Gan03}. \\

A dioperad is a properad with only compositions of operations
based on genus $0$ connected graphs. Hence, there is a natural
forgetful functor from properads to dioperads. To any properad
$\Po$, the associated dioperad
$\UU^{\textrm{properads}}_{\textrm{dioperads}}(\Po)$ has the same
underlying $\Sy$-bimodule and we only consider vertical
compositions of operations based on connected graphs of genus $0$.
Let us denote by $\Box$ the restriction of $\bc$ to genus $0$
graphs. With this notation, a dioperad is a monoid $(\mathcal{D},
\mu_\mathcal{D})$ in the monoidal category $(\Sy\textrm{-biMod},
\Box)$. From now on, let us denote the genus in exponent. For
instance, $\F^0$ will denote the free dioperad functor
$\F^{\textrm{dioperads}}_{\Sy\textrm{-biMod}}$ and $\F$ will
simply denote the free properad functor
$\F^{\textrm{properads}}_{\Sy\textrm{-biMod}}$.

\begin{pro}\label{Free:Dioperad->Properad}
The left adjoint of the forgetful functor
$\UU^{\textrm{properads}}_{\textrm{dioperads}}(\Po) \, :\,
 \mathrm{Properads}  \to \mathrm{Dioperads}$
is given by
$$\F(\mathcal{D})/I, $$
where $I$ is the (properadic) ideal generated by the image under
$\mu_\mathcal{D}-\Id$ of $\F^0(\mathcal{D})^{(2)}$, that is the
connected graphs of genus $0$ with two vertices.

In other words, this construction is the quotient of the free
properad on $\mathcal{D}$, consider as an $\Sy$-bimodule, by the
(dioperadic) composition of any pair of adjacent vertices with
only one edge in between.
\end{pro}

Notice that this construction is the same as the universal
enveloping algebra of a Lie algebra. Therefore, we will often call
it \emph{the universal enveloping properad of a dioperad} and
$\F^{\textrm{properads}}_{dioperads}$ the \emph{universal
enveloping functor}.

\begin{proo}
The proof is the same as the proof of the universal property of
the universal enveloping algebra of a Lie algebra. Hence it is
left to the reader.
\end{proo}

A direct corollary gives that the universal enveloping properad of
a dioperad defined by generators and relations is a properad given
the same generators and relations.

\begin{cor}\label{diop->prope V,R}
Let $\mathcal{D}$ be a dioperad defined by generators and
relations : $\mathcal{D}=\F^0(V)/(R)$, where $(R)$ is the
(dioperadic) ideal generated by $R$. The universal enveloping
properad is equal to
$$\F^{\textrm{properads}}_{dioperad}(\mathcal{D})=\F(V)/(R), $$
where $(R)$ is the (properadic) ideal generated by $R$.
\end{cor}

Even if an algebraic category $\mathcal{A}$ can be modelled by a
dioperad, the induced cofibrant resolution of this dioperad does
not contains all the data necessary for the study of deformation
theory of $\mathcal{A}$ because the universal enveloping functor
$\F^{\textrm{properads}}_{\textrm{dioperads}}$ is not exact as the
following counter-example shows.\\

Let $\IBi$ be the properad which models infinitesimal bialgebras (see \cite{Vallette03} Section~$2.9$). We consider its Koszul dual properad without the relation $\vcenter{\xymatrix@M=0pt@R=5pt@C=5pt{& \ar@{-}[d] &   \\
&\ar@{-}[dl] \ar@{-}[dr] & \\
\ar@{-}[dr]& & \ar@{-}[dl] \\
& & \\ & \ar@{-}[u] & }}=0$. Let us denote it by $\NCFrob$ because it models some kind of non-commutative  Frobenius bialgebras. A $\NCFrob$-bialgebra is a vector space $X$ equipped with a binary associative
  product $\mu \, :\,  X\otimes X \to X$ and a binary coassociative
coproduct $\Delta \, : \, X\to X\otimes X $ such that $\Delta$ is
a morphism of bimodules. This means
$$\Delta \circ \mu =(\Id \otimes \mu) \circ (\Delta \otimes \Id)=
(\mu \otimes \Id) \circ (\Id \otimes \Delta).$$  The graphical
picture of all the relations is the following~:

\begin{eqnarray*}
V \ &:& \ \vcenter{\xymatrix@M=0pt@R=6pt@C=6pt{\ar@{-}[dr] &
&\ar@{-}[dl]  \\  &\ar@{-}[d] & \\
& & \\ & & }} \oplus
\vcenter{\xymatrix@M=0pt@R=6pt@C=6pt{  & \ar@{-}[d]&  \\
&\ar@{-}[dl] \ar@{-}[dr]& \\ & &\\&&}}     \\
 R \ &:& \ \vcenter{\xymatrix@M=0pt@R=6pt@C=6pt{ & & & & \\
\ar@{-}[ddrr] &
&\ar@{-}[dl] & & \ar@{-}[ddll]  \\
&  & &  & \\
& &\ar@{-}[d] & & \\
& & }} = \vcenter{\xymatrix@M=0pt@R=6pt@C=6pt{& & & & \\
\ar@{-}[dr] &
&\ar@{-}[dr] & & \ar@{-}[dl]  \\
& \ar@{-}[dr] & &\ar@{-}[dl]  & \\
& &\ar@{-}[d] & & \\
& & }} \quad , \quad
\vcenter{\xymatrix@M=0pt@R=6pt@C=6pt{ & & & & \\
 & &\ar@{-}[d] & & \\
& & \ar@{-}[dl]  \ar@{-}[dr] & & \\
& \ar@{-}[dl]  \ar@{-}[dr] & &\ar@{-}[dr] & \\
& & & &}} =
\vcenter{\xymatrix@M=0pt@R=6pt@C=6pt{ & & & & \\
 & &\ar@{-}[d] & & \\
& & \ar@{-}[dl]  \ar@{-}[dr] & & \\
& \ar@{-}[dl]   & &\ar@{-}[dl] \ar@{-}[dr] & \\
& & & & }}  \quad , \quad
\vcenter{\xymatrix@M=0pt@R=6pt@C=6pt{\ar@{-}[dr] &  &\ar@{-}[dl]  \\  &\ar@{-}[d] & \\
& \ar@{-}[dl]\ar@{-}[dr]& \\ & & }} =
 \vcenter{\xymatrix@M=0pt@R=6pt@C=6pt{  & \ar@{-}[d] &  &\ar@{-}[d]  \\  & \ar@{-}[dl]
 \ar@{-}[dr] &  & \ar@{-}[dl]\\
\ar@{-}[d]&  & \ar@{-}[d] & \\ & & & }} =
\vcenter{\xymatrix@M=0pt@R=6pt@C=6pt{ \ar@{-}[d] & &\ar@{-}[d] &
\\ \ar@{-}[dr] &
& \ar@{-}[dl] \ar@{-}[dr]& \\
& \ar@{-}[d]&  & \ar@{-}[d] \\ & & & }} \ \\
&& \vcenter{\xymatrix@R=1pt@C=1pt{ &\ar@{-}[d] & & &\ar@{-}[d] \\
 & *=0{}\ar@{-}[dl]\ar@{-}[dr] & & & *=0{}\ar@{-}[d] \\
 *=0{}\ar@{-}[ddrr] & & *=0{}\ar@{-}[ddll] | \hole  & & *=0{}\ar@{-}[d] \\
 & & & & *=0{}\ar@{-}[d] \\
 *=0{}\ar@{-}[d] & & *=0{}\ar@{-}[dr] & & *=0{}\ar@{-}[dl] \\
 *=0{}\ar@{-}[d] & & & *=0{}\ar@{-}[d] & \\
 & & & & }}
=
\vcenter{\xymatrix@R=1pt@C=1pt{ *=0{}\ar@{-}[dddd]  & & & *=0{}\ar@{-}[d]  & \\
 & & & *=0{}\ar@{-}[dr] \ar@{-}[dl]  & \\
 & & *=0{}\ar@{-}[ddrr]  & & *=0{}\ar@{-}[ddll] | \hole  \\
 & & & & \\
 *=0{}\ar@{-}[dr] & & *=0{}\ar@{-}[dl] & & *=0{}\ar@{-}[d]  \\
 &  *=0{}\ar@{-}[d]& & & *=0{}\ar@{-}[d]  \\
 & & & &  }}
=
\vcenter{\xymatrix@R=1pt@C=1pt{ *=0{}\ar@{-}[dd] &  &  &  *=0{}\ar@{-}[d] &   \\
 &  &  &  *=0{}\ar@{-}[dr]\ar@{-}[dl] &   \\
 *=0{}\ar@{-}[ddrr] &  & *=0{}\ar@{-}[ddll] | \hole  &  &   *=0{}\ar@{-}[dddd] \\
 &  &  &  &   \\
 *=0{}\ar@{-}[dr] &  &  *=0{}\ar@{-}[dl] &  &   \\
 &  *=0{}\ar@{-}[d] &  &  &   \\
 &  &  &  &   }}
=
\vcenter{\xymatrix@R=1pt@C=1pt{ &  *=0{}\ar@{-}[d] &  & &   *=0{}\ar@{-}[dd] \\
 &  *=0{}\ar@{-}[dr]\ar@{-}[dl] &  &  &   \\
 *=0{}\ar@{-}[dddd] &  & *=0{}\ar@{-}[ddrr]  &  &  *=0{}\ar@{-}[ddll] | \hole  \\
 &  &  &  &   \\
 &  &  *=0{}\ar@{-}[dr] &  & *=0{}\ar@{-}[dl]   \\
 &  &  &  *=0{}\ar@{-}[d] &   \\
 &  &  &  &   }}=0.
 \end{eqnarray*}

Since the relations are linear combinations of connected graphs of
genus $0$, this category is faithfully modelled by the dioperad
$\NCFrob^0=\F^0(V)/(R)$. The exponent $0$ stands for the restriction
to graphs of genus $0$. It was proved in \cite{Gan03} that
$\NCFrob^{0}$ is a Koszul dioperad, since its Koszul dual dioperad $\IBi^0$ is Koszul by means of
distributive laws. That is the dioperad $\NCFrob^{0}$ admits a quadratic
dioperadic (genus $0$) model $(\F^0(\Co),
\partial^0) \xrightarrow{\sim} \NCFrob^{0}$, where $\Co$ is the codioperad
$s^{-1}{\IBi^{0}}^{\vee}$.
  (Notice that there is no
direct proof of this fact.) The differential $\partial^0$ splits
each element of $\Co$ into two vertices with only
one edge in between. \\

Consider now the properad
$\NCFrob=\F(V)/(R)$, which is the image under the universal
enveloping functor $\F^{\textrm{properads}}_\textrm{dioperads}$ of
$\Frob^{0}$ by Corollary~\ref{diop->prope V,R}. The image of the
chain complex $(\F^0(\Co),
\partial^0)$ under the functor
$\F^{\textrm{properads}}_{\textrm{dioperads}}$ is the quasi-free
properad on $\Co$ with the differential $\partial^0$, that is the
cobar construction of $\Co$, where this later is considered as a
coproperad. The homology of this chain complex is not
concentrated in degree $0$. \\

We build a cycle based on graphs of genus $2$ from the following
picture~:

$$\xymatrix@R=40pt@C=80pt{ \mu\circ \Delta\circ  \mu\circ \Delta
\ar[d]^{\mu\circ R_{rm}\circ \Delta} \ar[r]^{-\mu\circ R_{lm}\circ
\Delta}    &  \mu \circ (\mu \ot
\Id) \circ (\Id \ot \Delta )\circ \Delta  \ar[d]^{\quad \mu \circ (\mu \ot \Id )\circ R_c} \\
\ar[r]^{-R_a \circ (\Delta \ot \Id) \circ \Delta} \mu \circ (\Id
\ot \mu) \circ (\Delta \ot \Id )\circ \Delta &  \mu \circ (\mu \ot
\Id) \circ (\Delta \ot \Id )\circ \Delta,}$$

where $R_{rm}$ stands for the ``right module'' relation $\mu \circ
\Delta \to (\Id \otimes \mu) \circ (\Delta \otimes \Id)$, $R_{lm}$
for the ``left module'' relation $\mu \circ \Delta \to (\mu
\otimes \Id) \circ (\Id \otimes \Delta)$, $R_a$ the associativity
relation $\mu \circ (\mu \ot \Id) \to  \mu \circ (\Id \ot \mu)$
and $R_c$ the coassociativity relation $(\Delta \ot \Id)\circ
\Delta \to  (\Id \ot \Delta)\circ \Delta $.

The graphical picture is as follows~:

\begin{eqnarray*}
 R_{rm} \ :
\vcenter{\xymatrix@M=0pt@R=6pt@C=6pt{\ar@{-}[dr] &  &\ar@{-}[dl]  \\  &\ar@{-}[d] & \\
& \ar@{-}[dl]\ar@{-}[dr]& \\ & & }}
 \to
 \vcenter{\xymatrix@M=0pt@R=6pt@C=6pt{  & \ar@{-}[d] &  &\ar@{-}[d]  \\  & \ar@{-}[dl]
 \ar@{-}[dr] &  & \ar@{-}[dl]\\
\ar@{-}[d]&  & \ar@{-}[d] & \\ & & & }}  \quad &,& \quad R_{lm} \
:
\vcenter{\xymatrix@M=0pt@R=6pt@C=6pt{\ar@{-}[dr] &  &\ar@{-}[dl]  \\  &\ar@{-}[d] & \\
& \ar@{-}[dl]\ar@{-}[dr]& \\ & & }}
 \to\vcenter{\xymatrix@M=0pt@R=6pt@C=6pt{ \ar@{-}[d] &
&\ar@{-}[d] &
\\ \ar@{-}[dr] &
& \ar@{-}[dl] \ar@{-}[dr]& \\
& \ar@{-}[d]&  & \ar@{-}[d] \\ & & & }}\quad , \\
R_a \ :\vcenter{\xymatrix@M=0pt@R=6pt@C=6pt{ & & & & \\
\ar@{-}[ddrr] &
&\ar@{-}[dl] & & \ar@{-}[ddll]  \\
&  & &  & \\
& &\ar@{-}[d] & & \\
& & }} \to \vcenter{\xymatrix@M=0pt@R=6pt@C=6pt{& & & & \\
\ar@{-}[dr] &
&\ar@{-}[dr] & & \ar@{-}[dl]  \\
& \ar@{-}[dr] & &\ar@{-}[dl]  & \\
& &\ar@{-}[d] & & \\
& & }} \quad  &,& \quad R_c \ :
\vcenter{\xymatrix@M=0pt@R=6pt@C=6pt{ & & & & \\
 & &\ar@{-}[d] & & \\
& & \ar@{-}[dl]  \ar@{-}[dr] & & \\
& \ar@{-}[dl]  \ar@{-}[dr] & &\ar@{-}[dr] & \\
& & & &}} \to
\vcenter{\xymatrix@M=0pt@R=6pt@C=6pt{ & & & & \\
 & &\ar@{-}[d] & & \\
& & \ar@{-}[dl]  \ar@{-}[dr] & & \\
& \ar@{-}[dl]   & &\ar@{-}[dl] \ar@{-}[dr] & \\
& & & & }}
\end{eqnarray*}

Then, the cycle is based upon the following picture~:

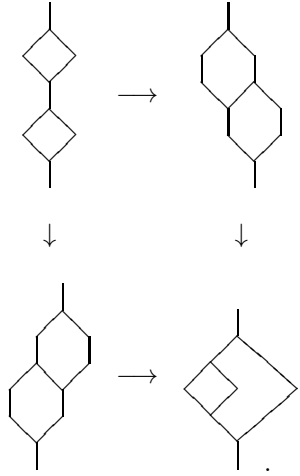
$$\begin{array}{ccc} \vcenter{\xymatrix@M=0pt@R=10pt@C=10pt{
&\ar@{-}[d]& \\ & \ar@{-}[dl]\ar@{-}[dr]&
\\ \ar@{-}[dr] &  &\ar@{-}[dl]  \\  &\ar@{-}[d] & \\
& \ar@{-}[dl]\ar@{-}[dr]& \\  \ar@{-}[dr]& &\ar@{-}[dl] \\ &
\ar@{-}[d]&\\ && }} &\longrightarrow &
\vcenter{\xymatrix@M=0pt@R=10pt@C=10pt{&\ar@{-}[d]&&\\
&\ar@{-}[dr]\ar@{-}[dl]&&\\
\ar@{-}[d]&&\ar@{-}[d]&\\
\ar@{-}[dr]&&\ar@{-}[dl]\ar@{-}[dr]&\\
&\ar@{-}[d]&&\ar@{-}[d]\\
&\ar@{-}[dr]&&\ar@{-}[dl]\\
&&\ar@{-}[d]&\\
&&&}}
\\
&& \\
\downarrow & \quad & \downarrow \\
&& \\
\vcenter{\xymatrix@M=0pt@R=10pt@C=10pt{ && \ar@{-}[d]& \\
&&\ar@{-}[dr]\ar@{-}[dl]& \\
&\ar@{-}[d]&&\ar@{-}[d]\\
&\ar@{-}[dr]\ar@{-}[dl]&&\ar@{-}[dl]\\
\ar@{-}[d]&&\ar@{-}[d]& \\
\ar@{-}[dr]&&\ar@{-}[dl]&\\
&\ar@{-}[d]&&\\
&&&}} &\longrightarrow & \vcenter{\xymatrix@M=0pt@R=10pt@C=10pt{ & & & & \\
 & &\ar@{-}[d] & & \\
& & \ar@{-}[dl]  \ar@{-}[dr] & & \\
& \ar@{-}[dl]  \ar@{-}[dr] & &\ar@{-}[dr] & \\
\ar@{-}[ddrr] &
&\ar@{-}[dl] & & \ar@{-}[ddll]  \\
&  & &  & \\
& &\ar@{-}[d] & & \\
& &&.& }}
\end{array}$$

We denote with the same notation the corresponding homotopies,
that is elements of $\Co$~:
\begin{eqnarray*}
\partial_0(R_{rm})=\mu \circ
\Delta - (\Id \otimes \mu) \circ (\Delta \otimes \Id) \quad &,&
\quad \partial_0(R_{lm})=\mu \circ \Delta - (\mu \otimes \Id)
\circ (\Id \otimes \Delta), \\
\partial_0(R_a)= \mu \circ (\mu \ot \Id) -  \mu \circ (\Id \ot \mu)  \quad &,&
\quad \partial_0(R_c)= (\Delta \ot \Id)\circ \Delta- (\Id \ot
\Delta)\circ \Delta.
\end{eqnarray*}

The previous picture proves that
$$\xi := \mu\circ R_{rm}\circ  \Delta -\mu\circ R_{lm}\circ \Delta
-R_a \circ (\Delta \ot \Id)\circ \Delta +\mu \circ (\mu \ot \Id
)\circ R_c
$$ is a cycle in $(\F(\Co),
\partial^0)$, that is $\partial^0(\xi)=0$.

\begin{lem}\label{CounterExample}
The cycle $\xi$ is not a boundary under $\partial^0$.
\end{lem}

\begin{proo}
The degree of $\xi$ is $1$. Suppose that there exists an element
$\zeta$ of degree $2$ such that $\partial^0(\zeta)=\xi$. This
element belongs to
$$\zeta \in \F(\Co_0\oplus \underbrace{\Co_1}_{(2)})
\oplus \F(\Co_0\oplus \underbrace{\Co_2}_{(1)}).$$ Let us denote
by $\zeta=\zeta_1+\zeta_2$ each component. The image under the
quadratic differential $\partial^0$ of any element of
$\F(\underbrace{\Co_0}_{(k)} \oplus \underbrace{\Co_2}_{(1)})$ is
an element of $\F(\underbrace{\Co_0}_{(k+1)} \oplus
\underbrace{\Co_1}_{(1)})$. And since the genus of the
differential $\partial^0$ is $0$, $\zeta_2$ is in
$\F^2(\underbrace{\Co_0}_{(1)} \oplus \underbrace{\Co_2}_{(1)})$,
that is the part of genus $2$ of $\F(\Co_0\oplus \Co_2)$. The
$\Sy$-bimodule $\Co_0$ is equal to
$V=\vcenter{\xymatrix@M=0pt@R=6pt@C=6pt{\ar@{-}[dr] &
&\ar@{-}[dl]  \\  &\ar@{-}[d] & \\
& & \\ & & }} \oplus
\vcenter{\xymatrix@M=0pt@R=6pt@C=6pt{  & \ar@{-}[d]&  \\
&\ar@{-}[dl] \ar@{-}[dr]& \\ & &\\&&}}$, that is binary. Hence
$\F(\underbrace{\Co_0}_{(1)} \oplus \underbrace{\Co_2}_{(1)}) $ is
concentrated in genus $0$ and $1$, which proves $\zeta_2=0$.

Since the image of $\F(\underbrace{\Co_0}_{(k)} \oplus
\underbrace{\Co_1}_{(2)})$ under $\partial^0$ is in
$\F(\underbrace{\Co_0}_{(k+2)} \oplus \underbrace{\Co_1}_{(1)})$,
$\zeta_1$ must belong to $\F(\Co_1)^{(2)}$. More precisely,
$\zeta_1$ is a element of  $\F^2(\Co_1)^{(2)}$ because the
differential $\partial^0$ preserves the genus. The $\Sy$-bimodule
$\Co_1$ is generated by the four elements $R_{rm}\in \Co(2,2)$,
$R_{lm}\in \Co(2,2)$, $R_{a}\in \Co(1,3)$ and $R_{c}\in \Co(3,1)$.
The only way to get an element of genus $2$ is to graft one
element from $\Co(1,3)$ to an element from $\Co(3,1)$. Finally
$\zeta$ is linear combination of $R_c \circ \sigma \circ Ra$, with
$\sigma \in \Sy_3$. And in this case, $\partial^0(\zeta)$ cannot
contain elements like $\mu\circ R_{rm}\circ  \Delta -\mu\circ
R_{lm}\circ \Delta$ whence the contradiction.
\end{proo}

This counter-example answers a question raised by
\cite{MarklVoronov03}, that is the functor
$\F^{\textrm{props}}_{\textrm{dioperads}}$ is not exact.

\begin{thm}
The universal enveloping functor
$\F^{\textrm{properads}}_{\textrm{dioperads}}$ is not exact.
\end{thm}

 For this reason, we are reluctant to include
dioperads in the preceding table. It is not enough in general to
find a resolution of the genus $0$ part of a properad to generate
a complete resolution of it. Nevertheless, it is sometimes the
case. We have emphasized the class of properads that admits a
quadratic model, that is Koszul properad. We do the same thing
with properads for which there exists a model with a genus $0$
differential.

\begin{dei}[Contractible properad]
We call \emph{contractible properad} any properad $\Po$ that
admits a model $(\F(\Co), \partial^0)\xrightarrow{\sim} \Po$ with
$\partial^0|_{\Co} \, : \Co \to \F^0(\Co)$, that is the part of
genus $0$ of the free properad on $\Co$.
\end{dei}

It is equivalent to ask that $\Co$ is a homotopy coproperad with
structure maps $\delta_n \, : \, \Co \to \F^0(\Co)^{(n)}$ with
image of genus $0$. In other words, $\Co$ is a \emph{homotopy
codioperad}.

\begin{pro}
Let $\Po=\F(V)/(R)$ be a properad defined by genus $0$ relations,
$R\subset \F^0(V)$. The properad $\Po$ is a contractible properad
if and only if the associated dioperad $\mathcal{D}:=\F^0(V)/(R)$
admits a quasi-free (dioperadic) resolution $(\F^0(\Co),
\partial^0) \xrightarrow{\sim} \mathcal{D}$, which is a
quasi-isomorphism preserved by the universal enveloping functor
$\F^{\textrm{properads}}_{\textrm{dioperads}}$.
\end{pro}

\begin{proo}
If $\Po$ is contractible, we denote by $(\F(\Co),
\partial^0) \xrightarrow{\sim} \Po$ its genus $0$ differential
model. Since $\partial^0$ preserves the genus, the chain complex
$(\F(\Co), \partial^0)$ is equal to the direct sum of
sub-complexes $\bigoplus_{g\ge 0} (\F^g(\Co), \partial^0)$. Hence,
the genus $0$ chain complex is a a resolution of $\mathcal{D}$.
And by Corollary~\ref{diop->prope V,R} the image under the
universal enveloping functor
$\F^{\textrm{properads}}_{\textrm{dioperads}}$ of the
quasi-isomorphism $(\F^0(\Co), \partial^0) \xrightarrow{\sim}
\mathcal{D}$ is the resolution $(\F(\Co), \partial^0)
\xrightarrow{\sim} \Po$. The other way is trivial.
\end{proo}

A Koszul contractible properad $\Po$ is a properad with a minimal
model $(\F(\Co), \partial^0)\xrightarrow{\sim} \Po$ whose
differential $\partial^0$ is quadratic and genus $0$. It is
equivalent to say that $\Co$ is a codioperad. If a properad
$\Po=\F(V)/(R)$ with
 genus $0$ relations is contractible Koszul, then the associated
 dioperad $\mathcal{D}=\F^0(V)/(R)$ is Koszul in the sense of
 \cite{Gan03}. But it is not true that any Koszul dioperad is a
 Koszul contractible properad as the example of $\NCFrob$ shows.
 Lemma~\ref{CounterExample} shows that it is not contractible.
Moreover
 we shall see below that it is not Koszul as a properad either.

\begin{pro}\label{Koszul contractible properad}
Let $\Po=\F(V)/(R)$ be a Koszul properad defined by a finite
dimensional $\Sy$-bimodule $V$ and by genus $0$ relations,
$R\subset \F^0(V)$. If the Koszul dual properad of $\Po$ is equal,
as an $\Sy$-bimodule, to the Koszul dual dioperad of the
associated dioperad $\mathcal{D}:=\F^0(V)/(R)$ then the properad
$\Po$ is  contractible.
\end{pro}

\begin{proo}
In this case, the Koszul dual coproperad $\Po^{\ac}={\Po^!}^\vee$
is equal to the Koszul dual dioperad
$\mathcal{D}{\ac}={\mathcal{D}^!}^\vee$. Hence  the image of the
partial coproduct $\Delta_{(1,1)} \, : \,    \Po^{\ac} \to
\Po^{\ac} \bc \Po^{\ac}$ is actually in $\Po^{\ac} \Box \Po^{\ac}$
which is the part of genus $0$ of $\Po^{\ac} \bc \Po^{\ac}$.
\end{proo}

The Koszul dual properad is equal to the Koszul dual dioperad if
and only if the part of genus $>0$ of $\Po^!$ vanished, that is
$\F^g(V^\vee)/(R^\perp)=0$ for $g>0$. Proposition~\ref{Koszul
contractible properad} allows us to give examples of Koszul
contractible properads. One way to prove that a properad is Koszul
is by means of \emph{distributive laws} (see Proposition~$8.4$ of
\cite{Vallette03}). Let $\Po$ be a quadratic properad of the form
$\Po=\F(V,\, W)/(R\oplus D \oplus S)$, where $R\subset \F^{(2)}(V)$, $S\subset
\F^{(2)}(W)$ and
 where $$D\subset (I\oplus \underbrace{W}_1)\boxtimes_c
(I\oplus \underbrace{V}_1)  \bigoplus (I\oplus
\underbrace{V}_1)\boxtimes_c (I\oplus \underbrace{W}_1).$$
The two pairs of $\Sy$-bimodules $(V,\, R )$ and $(W,\, S)$
generate two properads denoted $\mathcal{A}:=\F(V)/(R)$ and
$\mathcal{B}:=\F(W)/(S)$.

\begin{dei}[Distributive law]\label{Distributuive law}
Let $\lambda$ be a morphism of $\Sy$-bimodules
$$\lambda \ : \   (I\oplus \underbrace{W}_1)\boxtimes_c
(I\oplus \underbrace{V}_1)  \to (I\oplus
\underbrace{V}_1)\boxtimes_c (I\oplus \underbrace{W}_1).$$
such that the $\Sy$-bimodule $D$ is defined by the image of
$$(id,\, -\lambda) :\  (I\oplus \underbrace{W}_1)\boxtimes_c
(I\oplus \underbrace{V}_1) \to  (I\oplus
\underbrace{W}_1)\boxtimes_c (I\oplus \underbrace{V}_1)  \bigoplus
(I\oplus \underbrace{V}_1)\boxtimes_c (I\oplus
\underbrace{W}_1).$$
We call $\lambda$ a \emph{distributive law}
and denote $D$ by $D_\lambda$ if the two following morphisms are injective
$$ \left\lbrace
\begin{array}{c}
 \underbrace{\mathcal{A}}_1 \boxtimes_c \underbrace{\mathcal{B}}_2 \to \Po  \\
\underbrace{\mathcal{A}}_2 \boxtimes_c \underbrace{\mathcal{B}}_1 \to \Po.
\end{array}
\right.$$
\end{dei}

The last condition must be seen as a coherence axiom, which  ensures that the natural
morphism $\mathcal{A} \boxtimes \mathcal{B} \to \Po$ is injective. In this case, Proposition~$8.4$ of \cite{Vallette07} states that $\Po$ is Koszul if $\mathcal{A}$ and $\mathcal{B}$ are Koszul. A properad is called \emph{binary} if it is
generated by binary products and coproducts.

\begin{pro}
Let $\Do=\F^0(V)/(R)$ be a binary Koszul dioperad defined by a
distributive law such that $V$ is finite dimensional. Then the
associated properad $\Po:=\F(V)/(R)$ is Koszul and contractible.
\end{pro}

\begin{proo}
If a binary dioperad $\Do$ defined by a distributive law verifies
the hypotheses of Proposition~$5.9$ of \cite{Gan03}, then the
associated properad $\Po$ is also defined by distributive law and
verifies the hypotheses od Proposition~$8.4$ of \cite{Vallette03}.
In this case, the Koszul dual coproperad, given by
Proposition~$8.2$ of \cite{Vallette03} has a genus $0$ coproduct.
\end{proo}

\begin{cor}
The properads $\BLi$ of Lie bialgebras and $\IBi$ of infinitesimal
Hopf algebras are Koszul contractible.
\end{cor}

In this case, the Koszul dual (co)dioperad provides the good space
of ``homotopies'' for the resolution of the properad. Therefore, it
gives the proper notion of homotopy $\Po$-gebra (see \ref{Homotopy
gebra}). An example of
this fact for $\BLi$ can be found in Section~\ref{Poisson structure}in \cite{Gan05,Merkulov06}.\\

\begin{Rq}
Dually, in this case, the products of operations based on strictly
positive genus graphs of the Koszul dual properad always vanish.
If $g$ denotes the genus of the underlying graph, it means that
any such product is equivalent to products based on graphs with
$g$ simple loops $\vcenter{\xymatrix@M=0pt@R=5pt@C=5pt{& \ar@{-}[d] &   \\
&\ar@{-}[dl]
\ar@{-}[dr] & \\
\ar@{-}[dr]& & \ar@{-}[dl]\\ &\ar@{-}[d] & \\ & & }}$, using the
relations of the products and the relations of coproducts.
Therefore, it is zero because of the relation
$\vcenter{\xymatrix@M=0pt@R=5pt@C=5pt{& \ar@{-}[d] &   \\
&\ar@{-}[dl] \ar@{-}[dr] & \\
\ar@{-}[dr]& & \ar@{-}[dl] \\
& & \\ & \ar@{-}[u] & }}=0$ in the Koszul dual properad. This
statement is a non-trivial result about the coherence of the
relations of a properad.
\end{Rq}

To any binary properad $\Po$, we associate a properad $\Po_\Diamond$ which codes $\Po$-gebras satisfying the extra loop
relation $\vcenter{\xymatrix@M=0pt@R=5pt@C=5pt{& \ar@{-}[d] &   \\
&\ar@{-}[dl] \ar@{-}[dr] & \\
\ar@{-}[dr]& & \ar@{-}[dl] \\
& & \\ & \ar@{-}[u] & }}=0$. Since the properad $\BLi$ is Koszul, its Koszul dual properad
$\Frob_\Diamond$ is also Koszul by Koszul duality theory. This
means that $\Frob_\Diamond$ has a quadratic model. Since the
properad $\BLi$ has non trivial higher genus compositions, this
model is not contractible, that is the boundary map creates higher
genus graphs. The example $\Frob_\Diamond$ provides an example of
a Koszul non-contractible properad. (We do not know how to prove this result without the help
of Koszul duality for properads). \\

Let $\Co$ denote the Koszul dual coproperad of $\NCFrob$, that is
$\Co=s^{-1}{\IBi_\Diamond}^{\vee}$. Recall that a properad $\Po$
is Koszul if and only if the cobar construction of the Koszul dual
coproperad $\Omega(\Po^{\ac})=(\F(\Po^{\ac}),
\partial)$ is a resolution of $\Po$. This statement is not true for $\NCFrob$.
The cycle $\xi$ given above induces a non-trivial element in
homology.

\begin{lem}
The cycle $\xi$ is not a boundary under $\partial$.
\end{lem}

\begin{proo}
We use the same notations as in Lemma~\ref{CounterExample} but
applied to $\partial$ instead of $\partial^0$. The space $\Co_1$
is generated by the elements $R_{lm}$, $R_{rm}$, $R_a$, $R_c$ and some $R_i$ for $i=1, \ldots, 4$.
For the same reason, $\zeta_1$ must be an an
element of $\F(\Co_1)^{(2)}$. Since the image under $\partial$ of
any element of $\Co_1$ is a graph with two adjacent vertices
indexed by $\vcenter{\xymatrix@M=0pt@R=6pt@C=6pt{\ar@{-}[dr] &
&\ar@{-}[dl]  \\  &\ar@{-}[d] & \\
& & \\ & & }}$ or $
\vcenter{\xymatrix@M=0pt@R=6pt@C=6pt{  & \ar@{-}[d]&  \\
&\ar@{-}[dl] \ar@{-}[dr]& \\ & &\\&&}}$, the element $\mu \circ
R_{lm} \circ \Delta$ cannot belong to $\partial(\zeta_1)$. Hence
$\mu \circ R_{lm} \circ \Delta$ must be an element of
$\partial(\zeta_2)$. Since $\partial$ is quadratic, there exists
an element $S$ in $\Co_2$ such that $\partial(S)=\mu \circ R_{lm}
+ \cdots$ or $\partial(S)= R_{lm}\circ \Delta + \cdots$. Such an
$S$ has to be an element of either $\IBi_\Diamond^\vee(1,2)^{(3)}$
 or $\IBi_\Diamond^\vee(2,1)^{(3)}$. Consider the first case, the second one being symmetrical. The only element in 
 $\IBi_\Diamond^\vee(1,2)^{(3)}$ whose partial coproduct includes  $\mu \circ R_{lm}$ is the dual of the composite of 
 $$\vcenter{\xymatrix@M=0pt@R=6pt@C=6pt{\ar@{-}[dr] &
&\ar@{-}[dl]  \\  &\ar@{-}[d] & \\
& & \\ & & }} \boxtimes \left( \  \vcenter{\xymatrix@M=0pt@R=6pt@C=6pt{\ar@{-}[dr] &  &\ar@{-}[dl]  \\  &\ar@{-}[d] & \\
& \ar@{-}[dl]\ar@{-}[dr]& \\ & & }}
 -\vcenter{\xymatrix@M=0pt@R=6pt@C=6pt{ \ar@{-}[d] &
&\ar@{-}[d] &
\\ \ar@{-}[dr] &
& \ar@{-}[dl] \ar@{-}[dr]& \\
& \ar@{-}[d]&  & \ar@{-}[d] \\ & & & }} \  \right)$$
 in  $\IBi_\Diamond(1,2)^{(3)}$. The associativity relation and the loop relation in $\IBi_\Diamond$ show that this composite is equal to zero, which concludes the proof.  
\end{proo}

\begin{thm}
The properad $\NCFrob$ of non-commutative Frobenius bialgebras and the properad
$\IBi_\Diamond$ of involutive infinitesimal bialgebras are not Koszul.
\end{thm}

We hope that this helps to clarify the general picture of models
for prop(erad)s.


\section{Homotopy $\Po$-gebra}

In this section, we define the notion of $\Po$-gebra up to homotopy or \emph{homotopy $\Po$-gebra}.
We make explicit structures of homotopy $\Po$-gebras in terms of  Maurer-Cartan elements. We also define and make explicit morphisms of homotopy $\Po$-algebras, when $\Po$ is an operad, in terms of Maurer-Cartan elements in an $L_\infty$-algebra. This last part uses the notion of homotopy Koszul (colored) operads defined in the previous section.

\subsection{$\Po$-gebra, $\Po_{(n)}$-gebra and homotopy
$\Po$-gebra}\label{Homotopy gebra}

Let $\Po$ be a dg prop(erad) and $\Omega(\Co)$ be a model of
$\Po$.

\begin{dei}[Homotopy $\Po$-gebra]
A structure of \emph{homotopy $\Po$-gebra} on a dg module $X$ is a
morphism of dg prop(erad)s :  $\Omega(\Co) \to \End_X$.
\end{dei}

Any $\Po$-gebra is a homotopy $\Po$-gebra of particular type. In
this case, the morphism of dg-properads factors through $\Po$,
that is $\Omega(\Co) \xrightarrow{\sim} \Po \to \End_X$. For the
Koszul operads $\mathcal{A}ss$, $\Com$, $\Li$, this notion
coincides with homotopy associative, commutative, Lie algebras.
For the properads $\mathcal{B}i\mathcal{L}ie$ and
$\mathcal{A}ss\mathcal{B}i$, we get the notions of homotopy Lie
bialgebras and homotopy bialgebras. Since
$\mathcal{B}i\mathcal{L}ie$ is contractible, the explicit
definition given in \cite{Gan05,Merkulov06} coincides
with this one.\\

Theorem~\ref{mor} shows that a structure of homotopy $\Po$-gebra
on $X$ is equivalent to a morphism of $\Sy$-bimodules in
$s^{-1}\Hom^\Sy_0(\Co,\End_X)$ which is a
 Maurer-Cartan element in the $L_\infty$-convolution algebra
 $\Hom^\Sy(\Co,\End_X)$.

\begin{thm}\label{P infinit gebra =MC}
A $\Po$-gebra structure on $X$ is equivalent to a Maurer-Cartan
element in $\Hom^\Sy(\Co, \End_X)$.
\end{thm}

This notion is well defined
 and independent of the choice of a model. By
 Theorem~\ref{independence of homotopy convolution}, if
 $\Omega(\Co_1)$ and $\Omega(\Co_2)$ are two models of
 $\Po$, then the convolution $L_\infty$-algebras are
 quasi-isomorphic, which induces a bijection between the set of
 Maurer-Cartan elements. \\

We can discuss the form of the solutions of the Maurer-Cartan
equation. It gives the following definition.

\begin{dei}[$\Po_{(n)}$-gebra]
A dg module $X$ endowed with a Maurer-Cartan element $\gamma$ in
$\Hom^\Sy(\Co,\End_X)$ such that $\gamma(c)=0$ for every $c\in
\Co_{k>n}$ is called a \emph{$\Po_{(n)}$-gebra}.
\end{dei}

This notion is the direct generalization of the notion of
$A_{(n)}$-algebra of Stasheff \cite{Stasheff63} or
$L_{(n)}$-algebras. A $\Po_{(n)}$-gebra is a homotopy $\Po$-gebra
with strict relations from degree $n$.

\subsection{Morphisms of homotopy $\Po$-algebras as Maurer-Cartan elements}
\label{Homotopy Morphism=MC} Another application of the notion of
homotopy Koszul can be found in the study of morphisms between
homotopy $\Po$-algebras. A \emph{colored properad} is an operad
such that the inputs and outputs are labelled by an extra
labelling and such that the composition is coherent with respect
to this extra labelling. That is if the 'colors' (labelling) do
not match, the composition of operations vanishes. It is proven in
\cite{VanderLaan03} how to extend Koszul duality of operads to
colored operads. It is
straightforward to generalize Theorem~\ref{homotopykoszul} to this case.\\

Let $\Po=\F(V)/(R)$ be a Koszul operad. One can define the
2-colored operad $\Po_{\bullet \to \bullet}$ by $\Po=\F(V_1\oplus
V_2\oplus f)/(R_1\oplus R_2\oplus R_{\bullet \to \bullet})$, where
$V_1$ and $R_1$ (resp. $V_2$ and $R_2$) are copies of $V$ and $R$
with inputs and outputs labelled by the color $1$ (resp. $2$), $f$
is a generator of arity $(1,1)$ which goes from $1$ to $2$ and
 $R_{\bullet \to \bullet}$ is generated by $v\circ f^{\ot n} -
f \circ v$ for any element $v \in V(n)$ (see
\cite{Markl04Homotopy} for more details). The purpose of this
definition lies in the following result. A structure of
$\Po_{\bullet \to \bullet}$-algebra is the data of two
$\Po$-algebras with a morphism of $\Po$-algebras between them.

\begin{lem}
When $\Po$ is Koszul generated by a finite dimensional
$\Sy$-module $V$ such that $V(1)=0$, the $2$-colored operad
$\Po_{\bullet\to \bullet}$ is homotopy Koszul.
\end{lem}

\begin{proo}
(i) The operad ${(\Po_{\bullet \to \bullet}})_2$ is equal to
$\F(V_1\oplus V_2\oplus f)/(R_1\oplus R_2\oplus R_\bullet)$, where
$R_\bullet=f \circ V_1$. Hence, it is equal to $({\Po_{\bullet \to
\bullet}})_2\cong\Po_1 \oplus \Po_2 \oplus \Po\circ (I\oplus
\underbrace{f}_{\geqslant 1})$. Its Koszul dual is equal to
$({\Po_{\bullet \to \bullet}})^{\ac}_2=\Po_1^{\ac} \oplus
\Po_2^{\ac} \oplus s(f \circ \Po^{\ac})$. Therefore,
$(\F(s^{-1}\bar{\Po}_1^{\ac} \oplus s^{-1}\bar{\Po}_2^{\ac} \oplus
f \circ \bar{\Po}^{\ac})), \delta_2)$ is a quadratic model of
$({\Po_{\bullet\to\bullet}})_2$, because $\delta_2$
is equal to $3$ copies of the Koszul resolution of $\Po$. \\

(ii) Since $\Po_{\bullet \to \bullet}\cong \Po_1 \oplus \Po_2
\oplus \Po\circ (I\oplus \underbrace{f}_{\geqslant 1})$, it is
equal to $({\Po_{\bullet\to\bullet}})_2$.

(iii) Since $V$ is finite dimensional and $V(1)=0$, the filtration
with the number of leaves gives a suitable filtration.
\end{proo}

In this case, the minimal model of $\Po_{\bullet \to \bullet}$ is
given by  $(\F(s^{-1}\bar{\Po}_1^{\ac} \oplus
s^{-1}\bar{\Po}_2^{\ac} \oplus f \circ \bar{\Po}^{\ac})), \delta)$
by Theorem~\ref{homotopykoszul}.

\begin{pro}
An algebra over the model of $\Po_{\bullet\to\bullet}$ is the data
of two homotopy $\Po$-algebras with a homotopy (or weak) morphism
between them.
\end{pro}

\begin{proo}
A morphism of $2$-colored operads
$(\F(s^{-1}{\bar{\Po}^{\ac}_{\bullet \to \bullet  2}}), \delta)
\to \End_{X,Y}$ defines a homotopy $\Po$-algebra structure on $X$
and $Y$. The component on $\{\Hom(X^{\ot n}, Y) \}_{n \geqslant
1}$ is equivalent to a morphism of dg $\Po^{\ac}$-coalgebras
$\Po^{\ac}(X)\to \Po^{\ac}(Y)$, that is between the bar
constructions of $X$ and $Y$.
\end{proo}

\begin{thm}
Morphisms of homotopy $\Po$-algebras between $X$ and $Y$ are in
one-to-one  correspondence with  Maurer-Cartan elements in the
$L_\infty$-algebra $({\Hom^\Sy(\Po^{\ac}_{\bullet\to \bullet}}_2,
\End_{X,Y}), \delta) $
\end{thm}

Notice that this result was already proved by hands in
\cite{Dolgushev07}
in the case of homotopy Lie algebras.\\

Finally, a structure of homotopy $\Po$-algebra on $X$ is a
Maurer-Cartan element in the strict Lie algebra
$\Hom^\Sy(\Po^{\ac}, \End_X) $, whereas a morphism of homotopy
$\Po$-algebras between $X$ and $Y$ is a (generalized)
Maurer-Cartan element in the homotopy Lie algebra
$\Hom^\Sy(\Po^{\ac}, \End_{X,Y})$. The conceptual explanation of
this phenomenon is the following. In the first case, we have a
quadratic model of the Koszul operad $\Po$ and the second case, we
use a non-quadratic model of the homotopy Koszul $2$-colored
operad $\Po_{\bullet \to \bullet}$.


\sip

\section{$L_\infty$-algebras, dg manifolds, dg affine schemes and morphisms of prop(erad)s}

\sip

\subsection{$L_\infty$-algebras, dg manifolds and dg affine schemes.}
\label{dg manifold, dg affine scheme} Structure of a
$L_\infty$-algebra on a $\bZ$-graded vector space $\fg$ is, by
definition, a degree $-1$ coderivation, $Q: \odot^{\geq 1} s\fg\
\rar \odot^{\geq 1} s\fg$, of the free cocommutative coalgebra
without counit,
$$
\odot^{\geq 1}s\fg:=\bigoplus_{n\geq 1} \odot^n (s\fg)\ \subset \
\odot^\bullet s\fg:=\bigoplus_{n\geq 0} \odot^n (s\fg),
$$
 which satisfies the condition $Q^2=0$. It is often very helpful to use geometric intuition and language when working
 with $L_\infty$-algebras. Let us view the vector space $s\fg$ as a formal graded manifold (so that a choice of a basis
 in $\fg$ provides us with natural smooth coordinates on $s\fg$).  If $\fg$  is finite-dimensional, then the structure
 ring, $\f_{s\fg}$, of formal smooth functions on the formal manifold $s\fg$ is equal to the completed graded commutative algebra  $\widehat{\odot^\bullet} (s\fg)^*:= \prod_{n\geq 0} \odot^n (s\fg)^*$ which is precisely the dual of the coalgebra $\odot^\bullet s\fg$. This dualization sends the augmentation, $\odot^{\geq 1}s\fg$,
 of the latter into the ideal, $I:=\prod_{n\geq 1} \odot^n (s\fg)^*$, of
the distinguished point $0\in {s\fg}$, while the
coderivation $Q$  into as a degree $-1$ derivation of $\f_{s\fg}$, i.e.\ into a formal vector field (denoted by the
same letter $Q$)
 on the manifold ${s\fg}$
which vanishes at the distinguished point (as $QI\subset I)$ and
satisfies the condition $[Q,Q]=2Q^2=0$. Such vector fields are
often called {\em homological}.

\sip

In this geometric picture of $L_\infty$-algebra structures on
$\fg$, the subclass of dg Lie algebra structures gets represented
by at most quadratic homological vector fields $Q$, that is \texttt{$Q((s\fg)^*)\subset (s\fg)^*\oplus \odot^2 (s\fg)^*$}
Such a vector field has a well-defined value
at an arbitrary point $s\ga\in s\fg$, not only at the distinguished point $0\in s\fg$, i.e.\ it defines
a smooth homological vector field on $s\fg$ viewed as an ordinary (non-formal) graded manifold.
 Given a
particular dg Lie algebra $(\fg, d, [\ ,\ ])$, the associated
homological vector field $Q$ on $s{\fg}$ has the value at a point $s\ga \in s\fg$ given explicitly
by
\Beq
\label{Qu}
 Q(\ga):= d\ga + \frac{1}{2}[\ga,\ga],
\Eeq
where we used a canonical identification of the tangent space, ${\mathcal T}_{\ga}$, at $s\ga\in s\fg$ with $\fg$. One checks,
\begin{eqnarray*}
 Q^2(\ga)
&=& Q\left(d\ga + \frac{1}{2}[\ga,\ga]\right)\\
&=& - d(Q(\ga)) + \left[Q(\ga),\ga\right]\\
&=&  -d\left(d\ga + \frac{1}{2}[\ga,\ga]\right) +\left[d\ga +
\frac{1}{2}[\ga,\ga], \ga\right]\\
&=& 0.
\end{eqnarray*}
Notice that the zero locus of $Q$ is the set of Maurer-Cartan
elements in $\fg$.

\bip

A serious deficiency of the above geometric interpretation
of $L_\infty$-algebras is the necessity to work with the dual
objects, $(\f_{s\fg}, Q)$, which make sense only for finite
dimensional $\fg$. So we follow a suggestion of Kontsevich
\cite{Kontsevich03} and understand from now on a {\em dg (smooth formal)
manifold} as a pair, $(\odot^{\geq 1}X, Q)$, consisting of a
cofree cocommutative algebra on a $\bZ$-graded vector space $X$
together with a degree $-1$ codifferential $Q$. Note that the dual
of $\odot^{\geq 1}X$ is a well defined graded commutative algebra
(without assumption on finite-dimensionality of $X$) and that dual
of $Q$ is a well-defined derivation of the latter. We identify
from now on $Q$ with its dual and call it a {\em homological}\,
vector field on the dg manifold\footnote{A warning about shift of
grading: according to our definitions,
 a homological vector field on a graded vector space $X$ is the same as a
 $L_\infty$-structure on $s^{-1}X$.}
 $X$. This abuse of terminology is very helpful as it permits us to employ geometric intuition
and
 use
simple formulae of type (\ref{Qu}) to
define (in a mathematically rigorous way!)
codifferentials $Q$ on
$\odot^{\geq 1}X$.
Such codifferentials, $Q: \odot^{\geq 1}X\rar
\odot^{\geq 1}X$, are completely determined by the associated
compositions,
$$
Q_{proj}: \odot^{\geq 1}X \stackrel{Q}{\lon} \odot^{\geq 1}X
\stackrel{proj}{\lon} X.
$$
The restriction of $Q_{proj}$ to $\odot^{n}X\subset \odot^{\geq
1}X$ is denoted by $Q^{(n)}$, $n\geq 1$.

\bip

Since we work with dual notions (coalgebras, coderivations), we will need the notion of \emph{coideal}, which is the categorical dual to the notion of ideal. Hence, a coideal $I$ of a coalgebra $C$ is defined to be a quotient of $C$ such that the kernel of the associated projection $C\epi I$ is a subcoalgebra of $C$. For a complete study of this notion, we refer the reader to Appendix B ``Categorical Algebra'' of \cite{Vallette06Preprint}. This notion should not be confused with the notion of coideal used in Hopf algebra theory. Since a Hopf is an algebra and a coalgebra at the same time, a coideal in that sense is a submodule such the induced quotient carries again a bialgebra structure. \\

If $I$ is a  coideal of the coalgebra $\odot^{\geq 1}X$, we denote
the associated sub-coalgebra of $\odot^{\geq 1}X$ by $(\f_I:= I
\setminus \odot^{\geq 1}X, Q)$. The latter is  defined by the push-out diagram
in the category of coalgebras,
$$
\xymatrix {\f_I \ar[r]^{} \ar[d]  &\odot^{\geq 1}X\ar[d]\\
0 \ar[r]^{}  & I.}
$$
If the coideal $I$ is preserved by
$Q$ (i.e.\ admits a codifferential such that the right vertical arrow is a morphism of dg coalgebras), then the data $(\f_I, Q)$ is
naturally a differential graded coalgebra which we often call
a {\em dg affine scheme}
(cf.\ \cite{Barannikov06}). 
The coideal may not, in general, be homogeneous so the ``weight"
gradation, $\bigoplus_n \odot^n X$, may not survive in $\f_I$. A
generic dg affine scheme by no means corresponds
 to
a $L_\infty$-algebra but, as we shall see below,  some interesting
examples (with non-trivial and non-homogeneous coideals) do.

\sip

A {\em morphism}\, of dg affine schemes is, by definition,
a morphism of the associated dg coalgebras,
$(I_1\setminus \odot^{\geq 1}X_1, Q_1)\rar
(I_2\setminus \odot^{\geq 1}X_2, Q_2)$.

\subsection{Another geometric model for a $L_\infty$-structure} One can interpret a $L_\infty$-structure
on a  graded vector space $\fg$ as a linear total degree 1 polyvector field on the dual vector
space $\fg^*$ viewed as a graded affine manifold. Note that there is no need to employ the degree shifting
functors $s$ and $s^{-1}$ in this approach.

Indeed, let $(\wedge^\bullet {\mathcal T}_{\fg^*}, [\ ,\ ]_S)$ be the Schouten Lie algebra
of polynomial polyvector fields on the affine manifold $\fg^*$. A
  generic total degree 1
 polynomial polyvector field, $\nu= \{\nu_n\in \wedge^n {\mathcal T}_{\fg^*}\}_{n\geq 0}$, can be identified
 with a collection of its Taylor components with respect to affine coordinates on $\fg^*$, i.e. with a
 collection of linear maps,
 $$
 \nu_{m,n}: \odot^m \fg^* \lon \wedge^n \fg^*, m\geq 0,n\geq 0,
 $$
 of degree $n-2$. If $\nu$ is a linear polyvector field and lies in the Lie subalgebra $\wedge^{\geq 1} {\mathcal T}_{\fg^*}$,
 then only the Taylor components $\{\nu_{1,n}\}_{n\geq 1}$ can be non-zero. Their duals, $\nu_n:=(\nu_{1,n})^*$, is a collection
 of linear maps, $\nu_n: \wedge^n\fg \rar \fg$, $n\geq 1$, of degree $2-n$. It is easy to check the following
\begin{pro}
The data $\{\nu_n\}_{n\geq 1}$ defines a structure of $L_\infty$-algebra on $\fg$ if and only if the linear
polyvector field $\nu$ on $\fg^*$ satisfies the equation $[\nu, \nu]_S=0$.
\end{pro}
\begin{cor}
There is a one-to-one correspondence between structures of $L_\infty$-algebra on a finite-dimensional vector
space $\fg$ and linear degree one polyvector fields, $\nu\in \wedge^{\geq 1} {\mathcal T}_{\fg^*}$, satisfying the
equation $[\nu,\nu]_S=0$.
\end{cor}

Kontsevich's formality morphism \cite{Kontsevich03}, ${\mathcal F}$, associates to an arbitrary Maurer-Cartan element in the Schouten
Lie algebra $\nu\in \wedge^{\bullet} {\mathcal T}_{\fg^*}$ a Maurer-Cartan element, ${\mathcal F}(\nu)$, in the Hochschild
dg Lie algebra, $\oplus_{n\geq 0} \Hom_{poly}(\f_{\fg^*}^{\ot n}, \f_{\fg^*})[[\hbar]]$, of polydifferential operators on
the graded commutative algebra, $\f_{\fg^*}:=\odot^\bullet\fg$, of smooth functions on the affine manifold $\fg^*$.
If $\nu$ is a linear polyvector field on $\fg^*$ satisfying the equation $[\nu,\nu]_S=0$, then one can set to zero all contributions
to the formality morphism ${\mathcal F}$ coming from graphs with closed directed paths (wheels) \cite{Shoikhet99} and the
resulting element ${\mathcal F}_{no-wheels}(\nu)\in \oplus_{n\geq 0} \Hom_{poly}(\f_{\fg^*}^{\ot n}, \f_{\fg^*})[[\hbar]]$ is still
Maurer-Cartan. It is easy to check that
${\mathcal F}_{no-wheels}(\nu)$ has no summand with weight $n=0$ and hence defines an $A_\infty$-structure on $\odot^\bullet \fg$
which also makes sense for $\hbar=1$. Moreover, as ${\mathcal F}_{no-wheels}(\nu)$ involves no wheels
 (and hence no associated
traces
of linear maps), this $A_\infty$-structure makes sense for arbitrary (not necessarily finite-dimensional)
$L_\infty$-algebra.
\begin{dei}
\label{universal}
Let $\{\nu_n: \wedge^n\fg \rar \fg\}_{n\geq 1}$ be  a
 $L_\infty$-structure on a graded vector space $\fg$. The $A_\infty$-structure,
${\mathcal F}_{no-wheels}(\nu)$,
on $\odot^\bullet \fg$ obtained via Kontsevich's ``no-wheels" quantization of the associated linear polyvector field $\nu$
is called the {\em universal enveloping algebra}\, of the $L_\infty$-algebra.
\end{dei}

{In a recent interesting paper \cite{Baranovsky07} Baranovsky also
defined a universal enveloping for a $L_\infty$-algebra $\fg$ as a
certain $A_\infty$-structure on the space $\odot^\bullet \fg$. In
his approach the $A_\infty$-structure is constructed with the help
of the homological perturbation and the natural homotopy transfer
of the canonical dg associative algebra structure on the cobar
construction on the dg coalgebra $\odot^\bullet s\fg$.}

\subsection{Maurer-Cartan elements in a filtered $L_\infty$-algebra}\label{filtered}
 A $L_\infty$-algebra $(\fg, Q=\{Q^{(n)}\}_{n\geq 1})$ is called
{\em filtered}\, if $\fg$ admits a non-negative decreasing
Hausdorff filtration,
$$
\fg_0=\fg \supseteq \fg_1\supseteq\ldots\supseteq \fg_i \supseteq
\ldots,
$$
such that the linear map $Q^{(n)}:\odot^n(s\fg)\rar s\fg$ takes values
in $s\fg_n$ for all $n\geq n_0$ and some $n_0\in \bN$. In this
case $Q$ extends naturally to a coderivation of the cocommutative
coalgebra, ${\odot^{\geq 1}} s\hat{\fg}$, with $\hat{\fg}$ being
the completion of $\fg$ with respect to the topology induced by
the filtration, and the equation,
$$
Q(\sum_{n\geq 1}\frac{1}{n!}\ga^{\odot n})=0,
$$
for a degree zero element $\ga\in s\hat{\fg}$ (i.e.\ for a degree
$-1$ element in $\hat{\fg}$) makes sense. Its solutions are called
{\em (generalized) Maurer-Cartan elements}\, (or, shortly, {\em MC
elements}) in $(\fg, Q)$. Geometrically, an MC element is a degree
$-1$ element in $\hat{\fg}$ at which the homological vector field
$Q$ vanishes. From now on we do not distinguish between $\fg$ and
its completion $\hat{\fg}$.

\sip

To every MC element $\ga$ in a filtered $L_\infty$-algebra $(\fg,
Q)$ there corresponds, by Theorem 2.6.1 in \cite{Merkulov00},    a
twisted $L_\infty$-algebra, $(\fg, Q_\ga)$, with
$$
Q_\ga(\al):= Q(\sum_{n\geq 0}\frac{1}{n!}\ga^{\odot n}\odot \al)
$$
for an arbitrary $\al\in \odot^{\geq 1} s\fg$. The geometric meaning of this twisted
$L_\infty$-structure is simple \cite{Merkulov00}: if a homological
vector field $Q$ vanishes at a degree 0 point $\ga\in s\fg$, then
applying to $Q$ a formal diffeomorphism, $\phi_\ga$, which is a
translation sending $\ga$ into the origin $0$ (and which is
nothing but the unit shift, $e^{\ad \ga}$, along the formal
integral lines of the constant vector field $-\ga$) will give us a
new formal vector field, $Q_\ga:=d\phi_\ga( Q)$, which is {\em
homological} and {\em vanishes} at the distinguished point; thus
$Q_\ga$ defines a $L_\infty$ structure on the underlying space
$\fg$. In fact, we can apply this ``translation diffeomorphism"
trick to arbitrary (i.e.\ not necessarily
 MC) elements $\ga$ of degree 0 in $s\fg$
and get {\em homological}\, vector fields, $Q_\ga:=d\phi_\ga( Q)$,
which do not vanish at $0$ and hence define generalized $L_\infty$
structures on $\fg$ with ``zero term" $Q_\ga^{(0)}\neq 0$.



\subsection{Extended morphisms of dg props as a dg affine scheme.} \label{Derivative as dg affine scheme}
Let $(\Po, \p_\Po)$ and $(\mathcal{E}, \p_\mathcal{E})$ be dg
prop(erad)s
with differentials
$\p_\Po$ and $\p_\mathcal{E}$ of degree $-1$. Let
$\Hm(\Po,\mathcal{E})$ denote the graded vector space of all
possible morphisms $\Po\rar \mathcal{E}$ in category of
$\bZ$-graded $\bS$-bimodules, and let $\Mor(\Po,\mathcal{E})$
denote the set  of all possible morphisms $\Po\rar \mathcal{E}$ in
category of prop(erad)s, (note that we do {\em not}\,
 assume that elements of $\Hm(\Po,\mathcal{E})$  or $\Mor(\Po,\mathcal{E})$  respect differentials). It is clear that
\[
\Mor(\Po,\mathcal{E})=\left\{\ga\in \Hom_\bZ(\Po,\mathcal{E})\,
|\, \ga\circ\mu_\Po\left(\Po \bbc \Po\right)=
\mu_\mathcal{E}\left(\ga(\Po) \bbc \ga(\Po)\right)\ \mbox{and}\
|\ga|=0\right\}.
\]
We need a $\bZ$-graded extension of this set,
\Beq\label{extmor}
\Mor_\bZ(\Po,\mathcal{E}):=\left\{\ga\in
\Hom_\bZ(\Po,\mathcal{E})\, |\, \ga\circ\mu_\Po\left(\Po \bbc
\Po\right)= \mu_\mathcal{E}\left(\ga(\Po) \bbc
\ga(\Po)\right)\right\},
\Eeq
which we define by the same algebraic equations but dropping the
assumption on the degree and  homogeneity
 of $\ga$.

\begin{lem}\label{def-Q}
The vector space  $\Hm(\Po,\mathcal{E})$ is naturally a dg
manifold.
\end{lem}

\begin{proo} We define a degree $-1$ coderivation of the free cocommutative coalgebra,
$\odot^{\geq 1}\Hm(\Po,\mathcal{E})$ by setting (in the dual
picture, cf.\ \S~\ref{dg manifold, dg affine scheme})
\Beq\label{Q}
Q(\ga) := \p_\mathcal{E} \circ \ga -
(-1)^{\ga}\ga\circ \p_\Po
\Eeq
for an arbitrary $\ga\in
\Hm(\Po,\mathcal{E})$.
As
\Beqrn
Q^2(\ga) &=& Q(\p_\mathcal{E} \circ \ga - (-1)^{\ga}\ga\circ \p_\Po)\\
&=& -  \p_\mathcal{E} \circ Q(\ga) - (-1)^{\ga}Q(\ga)\circ \p_\Po)\\
&=& -  (-1)^{\ga} \p_\mathcal{E} \circ \ga\circ \p_\Po  + (-1)^{\ga}\p_\mathcal{E}\circ \ga \circ \p_\Po\\
&=& 0,
\Eeqrn
$Q$ is a linear homological field on
$\Hm(\Po,\mathcal{E})$. (By the way, the zero locus of $Q$
is a linear subspace of $ \Hm(\Po,\mathcal{E})$ describing
morphisms of {\em complexes}.)
\end{proo}

\begin{pro}\label{scheme}
The set $\Mor_\bZ(\Po,\mathcal{E})$ is naturally a dg affine
scheme.
\end{pro}

\begin{proo}
Let $I$ be the coideal in $\odot^{\geq
1}\Hm(\Po,\mathcal{E})$ cogenerated by the algebraic
relations,
\Beq\label{I}
\ga\circ\mu_\Po\left(\Po \bbc \Po\right) -
\mu_\mathcal{E}\left(\ga(\Po) \bbc \ga(\Po)\right),
\Eeq
on the ``variable" $\ga\in \Hm(\Po,\mathcal{E})$. The
sub-coalgebra,
$$
\f_{{\mathrm M\mathrm o \mathrm r}_\bZ(\Po,\mathcal{E})}:=
I\setminus \odot^{\geq 1}\Hm(\Po,\mathcal{E}),
$$
of $\odot^{\geq 1}\Hm(\Po,\mathcal{E})$ makes the set
$\Mor_\bZ(\Po,\mathcal{E})$ into a $\bZ$-graded affine scheme.
Next we show that the homological vector field $Q$ defined in
Lemma~\ref{def-Q} is tangent to $\Mor_\bZ(\Po,\mathcal{E})$.
Indeed, identifying $Q$ and $I$ with their duals (as in
subsection~\ref{dg manifold, dg affine scheme} and the proof of
Lemma~\ref{def-Q}), we have
\Beqrn Q\left(\ga\circ\mu_\Po\left(\Po
\bbc \Po\right) - \mu_\mathcal{E}\left(\ga(\Po) \bbc
\ga(\Po)\right)\right)\hspace{-2mm}
&=& \hspace{-2mm} Q(\ga)\circ \mu_\Po\left(\Po \bbc
\Po\right) - \mu_\mathcal{E}\left(Q(\ga)(\Po) \bbc
\ga(\Po)\right)\\
\hspace{-2mm} && - (-1)^{|\ga|} \mu_\mathcal{E}\left(\ga(\Po)
\bbc Q(\ga)(\Po)\right).
\Eeqrn
Consistency of $\p_\Po$ and
$\p_\mathcal{E}$ with $\mu_\Po$ and, respectively,
$\mu_\mathcal{E}$ implies,
\Beqrn
 Q(\ga)\circ \mu_\Po\left(\Po \bbc \Po\right) &=& \p_\mathcal{E}\circ \ga \circ\mu_\Po\left(\Po \bbc \Po\right)
- (-1)^{\ga}\ga\circ \p_\Po\circ  \mu_\Po\left(\Po \bbc \Po\right) \\
&=& \p_\mathcal{E}\circ \ga\circ \mu_\Po\left(\Po \bbc \Po\right)
- (-1)^{\ga}\ga\circ  \mu_\Po\left(\p_\Po(\Po) \bbc \Po\right)\\
&& -
(-1)^{\ga}\ga\circ  \mu_\Po\left(\Po \bbc \p_\Po(\Po)\right) \\
&=_{\bmod I}& \p_\mathcal{E}\circ \mu_\mathcal{E}\left(\ga(\Po)
\bbc \ga(\Po)\right) -
(-1)^{\ga}\mu_\mathcal{E}\left(\ga\circ \p_\Po(\Po) \bbc \ga(\Po)\right)\\
&& -\mu_\mathcal{E}\left(\ga(\Po) \bbc  \ga\circ \p_\Po(\Po)\right)\\
&=_{\bmod I}& \mu_\mathcal{E}\left(Q(\ga)(\Po) \bbc
\ga(\Po)\right) + (-1)^{|\ga|} \mu_\mathcal{E}\left(\ga(\Po) \bbc
Q(\ga)(\Po)\right). \Eeqrn Thus $Q(I)\subset I$, and hence $Q$
gives rise to a degree $-1$ codifferential on the coalgebra
$\f_{{\mathrm M\mathrm o \mathrm r}_\bZ(\Po,\mathcal{E})}$ proving
the claim.
\end{proo}


In the following theorem, we study the properties of the convolution $L_\infty$-algebra defined in Theorem~\ref{conv Lie infty}.

\begin{thm}\label{mor}
Let $(\Po=\cF(s^{-1} \Co), \p_\cP)$ be a quasi-free prop(erad)
 generated by an $\bS$-bimodule
$s^{-1} \Co$ (so that $\Co$ is a homotopy coprop(erad)), and let
$(\mathcal{E},\p_\mathcal{E})$ be an arbitrary dg prop(erad). Then
\begin{itemize}
\item[(i)] The graded vector space,
$\Hm(\Co,\mathcal{E})$, is canonically a
$L_\infty$-algebra;
\item[(ii)] The canonical $L_\infty$-structure
in (i) is filtered and its
 $MC$ elements are morphisms,  $(\Po, \p_\Po)\rar (\mathcal{E},\p_\mathcal{E}) $, of\, {\em dg} prop(erad)s;
\item[(iii)] if $\p_\Po(s^{-1}\Co)\subset \cF(s^{-1}\Co)^{(\leq
2)}$, where $\cF(s^{-1}\Co)^{(\leq 2)}$
 is the subspace of $\cF(s^{-1}\Co)$ spanned by decorated graphs with at most two vertices, then
 $\Hm(\Co,\mathcal{E})$ is canonically a dg Lie algebra.
\end{itemize}
\end{thm}
\begin{proo} (i) If $\Po$ is free as a prop(erad), then extended morphisms from $\Po$ to $\mathcal Q$ are uniquely determined
by their values on the generators $s^{-1}\Co$ so that
$\f_{{\mathrm M\mathrm o \mathrm r}_\bZ(\Po,\mathcal{E})} =
\odot^{\geq 1}\Hm(s^{-1}\Co,\mathcal{E})$ and the claim
follows from the definition of $L_\infty$-structure in \S~\ref{dg
manifold, dg affine scheme}.
\bip

(ii) The canonical $L_\infty$ structure on $\Hm(\Co,\mathcal{E})$
is given by the restriction of the homological vector field
(\ref{Q}) on $s\Hm(\Po, \mathcal{E})$ to the subspace
$\Hm(\Co,\mathcal{E})$. This field is a formal power series in
coordinates on $\Hm(\Co,\mathcal{E})$ and its part, $Q^{(n)}$,
corresponding to monomials of (polynomial) degree $n$ is given
precisely by
\Beq\label{Q_n}
Q^{(1)}(\ga):=\p_\mathcal{E} \circ \ga -
(-1)^{\ga}\ga\circ \p_\Po^{(1)} \quad \textrm{and}\quad
Q^{(n)}(\ga)
:=
 - (-1)^{\ga}\ga\circ \p_\Po^{(n)} \ \textrm{for}\ n>1, \Eeq
where $\p_\Po^{(n)}$ is the composition\footnote{Note that for any
differential $\p_\Po$ in a free properad $\cF(s^{-1}\Co)$
the induced map $\p_\Po^{(1)}: s^{-1}\Co \rar s^{-1}\Co$ is also a differential.}
$$
\p_\Po^{(n)}: s^{-1}\Co \stackrel{\p_\Po}{\lon}
\cF(s^{-1}\Co)\stackrel{proj}{\lon} \cF(s^{-1}\Co)^{(n)}.
$$
Note that the first summand on the r.h.s.\ of (\ref{Q})
contributes only to $Q^{(1)}$.

\sip

Define an exhaustive increasing filtration on the $\bS$-bimodule
${ \Co}$ by
$$
{ \Co}_0=0,\ \ {\Co}_i:=s\bigcap_{n\geq i} \Ker \p_\Po^{(n)}\
\mbox{for}\ i\geq 1,
$$
and the associated decreasing filtration on
$\Hm(\Co,\mathcal{E})$ by
$$
\Hm(\Co,\mathcal{E})_i:=\{ \ga \in
\Hm(\Co,\mathcal{E})\ |\ \ga(v)=0 \ \forall v\in {
\Co}_i\}, \ i\geq 0.
$$
Then, for all $n \geq 2$ and any $f_1, \ldots, f_n\in
\Hm(\Co,\mathcal{E})$, equality (\ref{Q_n}) implies that the value
of the map $Q^{(n)}(f_1, \ldots, f_n)\in s^{-1}\Hm(\Co,\mathcal{E})$
on arbitrary element of $\Co_n\subset \Ker \p_{\Po}^{(n)}$  is
equal to zero, i.e.
$$
Q^{(n)}(f_1, \ldots, f_n) \in \Hm(\Co,\mathcal{E})_n.
$$
Which in turn implies  the claim that the canonical $L_\infty$
structure on $\Hm(\Co,\mathcal{E})$ is filtered with
respect to the constructed filtration. The claim about MC elements
follows immediately from the definition (\ref{Q}) of the
homological vector field.

\bip

(iii) As $ \p_\Po^{(n)}=0$ for $n>2$ we conclude using formula
(\ref{Q_n}) that $Q^{(n)}=0$ for all $n>2$.
\end{proo}

\bip

A special case of the above Theorem when both $\Po$ and
$\mathcal{E}$ are operads was proven earlier by van der Laan
\cite{VanderLaan02} using different ideas.

\sip

The main point of our proof of Theorem~\ref{mor} is an observation
that, for a free prop(erad) $\Po=\cF(s^{-1} \Co)$, the set,
$\Mor_{\bZ}(\Po,\cE)$, of extended morphisms from $\Po$ to an
arbitrary prop(erad) $\cE$, i.e.\ the set of solutions of
equation~(\ref{extmor}), can be canonically identified with the
graded {\em vector}\, space $s\Hom_\bullet^\Sy(\cC,\cE)$. This
simple fact makes the dg affine scheme   $\Mor_{\bZ}(\Po,\cE)$
into a dg {\em smooth}\, manifold and hence provides us with a
canonical $L_\infty$-structure on  $\Hom_\bullet^\Sy(\Co,\cE)$. A
similar phenomenon occurs for the set of extended morphisms,
${\mathrm C\mathrm o \Mor}_{\bZ}(\cE^c, \Po^c)$, from an arbitrary
coprop(erad) $\cE^c$ to a cofree coprop(erad) $\Po^c:= \cF^c(s\cE)$,
and hence the arguments very similar to the ones used in the proof
of Theorem~\ref{mor} (and which we leave to the reader as an
exercise) establish the following,

\begin{thm}\label{comor}
Let $(\cP^c=\cF^c(s \Co), d_\cP)$ be a quasi-free coprop(erad),
that is $\Co$ is a homotopy prop(erad),  and let $(\cE^c,d_\mathcal{E})$
be an arbitrary dg coprop(erad). Then
\begin{itemize}
\item[(i)] The graded vector space,
$\Hm(\cE^c,\Co)$, is canonically a
$L_\infty$-algebra;
\item[(ii)] The canonical $L_\infty$-structure
in (i) is filtered and its
 $MC$ elements are morphisms,  $(\cE^c, d_\cE)\rar (\cP^c,d_\cP)$, of\, {\em dg} coprop(erad)s;
\item[(iii)] if $d_\Po$ is quadratic, that is $\Co$ is a usual
prop(erad), then
 $\Hm(\cE^c,\Co)$ is canonically a dg Lie algebra.
\end{itemize}
\end{thm}

\sip

For finite-dimensional $\Qo$ and $\Co$  Theorems~\ref{mor} and \ref{comor} are, of course,
equivalent to each other.

\sip

A morphism of $L_\infty$-algebras, $
(\fg_1, \cQ_1) \rar (\fg_2, \cQ_2)$, is, by definition \cite{Kontsevich99},
a morphism,  $\lambda: (\odot^\bullet(s \fg_1), \cQ_1) \rar
(\odot^\bullet(s \fg_2), \cQ_2)$, of the associated dg coalgebras.
It is called a {\em quasi-isomorphism}\,
if the composition,
$$
s\fg_1 \stackrel{i}{\lon} \odot^\bullet (s\fg_1) \stackrel{\lambda}{\lon}
 \odot^\bullet (s\fg_2)  \stackrel{p}{\lon} s\fg_2
$$
induces  an isomorphism, $H(s\fg_1, \cQ_1^{(1)})\rar H(s\fg_2, \cQ_2^{(1)})$,
of the associated homology groups with respect to the {\em linear}\, (in cogenerators) parts
of the codifferentials. Here $i$ is a natural inclusion and $p$ a natural projection.

\sip


By analogy,
a map $\phi: (\cF(s^{-1} \Co_1), \p_1) \rar (\cF(s^{-1} \Co_2), \p_2)$ of quasi-free properads
is called
a {\em tangent}\, quasi-isomorphism if the composition
 $$
s^{-1} \Co_1 \stackrel{i}{\lon} \cF(s^{-1} \Co_1) \stackrel{\phi}{\lon}
 \cF(s^{-1} \Co_2)  \stackrel{p}{\lon} s^{-1} \Co_2
$$
induces  an isomorphism of cohomology groups,
 $H(s^{-1} \Co_2, \p_2^{(1)})\rar H(s^{-1} \Co_2, \p_2^{(1)})$.

\sip

If we assume that properads $\cF(s^{-1} \Co_1)$ and $ (\cF(s^{-1} \Co_2)$
are completed  by the number of vertices (see \S 5.4) and that their differentials are bounded,
$$
\p_i(s^{-1} \Co_i) \subset \cF(s^{-1} \Co_i)^{(\leq n_i)}\ \ \mbox{for some}\ n_i\in {\mathbb N},\ i=1,2,
$$

then it is not
 hard to show (using filtrations by the number of vertices as in \S 5.4
 and the classical Comparison Theorem of spectral sequences) that any
 continuous tangent
quasi-isomorphism  $\phi: (\cF(s^{-1} \Co_1), \p_1) \rar (\cF(s^{-1} \Co_2), \p_2)$
 is actually a quasi-isomorphism in the ordinary sense.


%
\bip

\begin{thm}\label{morqis}
(i) Let $(\Po_1:=\cF(s^{-1} \Co_1), \p_1)$ and $(\Po_2:=\cF(s^{-1}
\Co_2), \p_2)$  be  quasi-free prop(erad)s,
$(\mathcal{E},\p_\mathcal{E})$ a dg prop(erad), and
$(\Hm(\Co_1,\mathcal{E}), Q_1)$ and $(\Hm(\Co_1,\mathcal{E}),
Q_2)$ the associated  $L_\infty$-algebras.  Then any morphism,
$$
\phi: (\Po_1, \p_1) \lon (\Po_2, \p_2),
$$
of dg prop(erad)s induces canonically an associated morphism,
$$
\phi_{ind}: (\Hm(\Co_2,\mathcal{E}), Q_2) \lon
(\Hm(\Co_1,\mathcal{E}), Q_1)
$$
of $L_\infty$-algebras. Moreover, if $\phi$ is a tangent quasi-isomorphism
of dg prop(erad)s, then $\phi_{ind}$ is a quasi-isomorphism of
$L_\infty$-algebras.

\sip

(ii) Let $(\Po:=\cF(s^{-1} \Co), \p)$  be   a quasi-free prop(erad),
$(\mathcal{E}_1,\p_{\mathcal{E}_1})$ and $(\mathcal{E}_2,\p_{\mathcal{E}_2})$ arbitrary
dg prop(erad)s, and
$(\Hm(\Co,\mathcal{E}_1), Q_1)$ and $(\Hm(\Co,\mathcal{E}_2),
Q_2)$ the associated  $L_\infty$-algebras.  Then any morphism,
$$
\psi: (\cE_1, \p_{\cE_1}) \lon (\cE_2, \p_{\cE_2}),
$$
of dg prop(erad)s induces canonically an associated morphism,
$$
\psi_{ind}: (\Hm(\Co,\mathcal{E}_1), Q_1) \lon
(\Hm(\Co,\mathcal{E}_2), Q_2)
$$
of $L_\infty$-algebras. Moreover, if $\psi$ is a quasi-isomorphism
of dg prop(erad)s, then $\psi_{ind}$ is a quasi-isomorphism of
$L_\infty$-algebras.

\end{thm}
\begin{proo} (i) The map $\phi$ induces  a degree 0 linear map,
$$
\Ba{ccc}
 \Hm(\Po_2, {\mathcal E}) & \lon &\Hm(\Po_1, {\mathcal E})\\
 \ga & \lon & \ga\circ \phi
 \Ea
 $$
Using  definition (\ref{Q}) of the codifferentials $Q_1$ and $Q_2$, and the fact that
 $\phi$ respects differentials $\p_1$ and $\p_2$, we obtain,
for any $\ga\in\Hm(\Po_2, {\mathcal E})$,
 \Beqrn
Q_1(\ga\circ \phi)&=& \p_{\mathcal E}\circ \ga\circ \phi -
(-1)^{\ga} \ga\circ\phi\circ \p_1\\
&=& \p_{\mathcal E}\circ \ga\circ \phi -
(-1)^{\ga} \ga\circ \p_2\circ\phi \\
&=& Q_2(\ga)\circ\phi,
 \Eeqrn
and hence conclude that $\phi$ induces a morphism of dg
coalgebras,
$$
\phi_{ind}: (\odot^{\geq 1} \Hm(\Po_2, {\mathcal E}), Q_2) \lon
(\odot^{\geq 1} \Hm(\Po_1, {\mathcal E}), Q_1).
$$
As
$$
\phi\circ\mu_{\Po_1}\left(\Po_1 \bbc \Po_1\right)=
\mu_{\Po_2}\left(\phi(\Po_1) \bbc \phi(\Po_1)\right)\subset \mu_{\Po_2}\left(\Po_2 \bbc \Po_2\right)
$$
we have
$$
\ga\circ\phi\circ\mu_{\Po_1}\left(\Po_1 \bbc \Po_1\right) -
\mu_\mathcal{E}\left(\ga\circ\phi(\Po_1) \bbc \ga\circ\phi(\Po_1)\right)\subset
\hspace{50mm}
$$
$$
\hspace{60mm}
\subset
\ga\circ\mu_{\Po_2}\left(\Po_2 \bbc \Po_2\right) -
\mu_\mathcal{E}\left(\ga(\Po_2) \bbc \ga(\Po_1)\right).
$$
Thus the map $\phi_{ind}$ sends cogenerators  (\ref{I}) of the coideal $I_2$
in $\odot^{\geq 1}\Hm(\Po_2, {\mathcal E})$
into cogenerators of the coideal $I_1$ in
 $\odot^{\geq 1}\Hm(\Po_1, {\mathcal E})$,
and hence gives rise to a morphism of dg coalgebras,
$$
\phi_{ind}: \left(\f_{{\mathrm Mor}_\bZ(\Po_2, {\mathcal E})}, Q_2\right) \lon
\left(\f_{{\mathrm Mor}_\bZ(\Po_1, {\mathcal E})}, Q_1\right),
$$
i.e.\ to a morphism of dg affine schemes,
$\phi_{ind}: \left(\Mor_\bZ(\cP_2,\mathcal E), Q_2\right)\rar
 \left(\Mor_\bZ(\cP_1,\mathcal E), Q_1\right)$.\\

If the dg prop(erad)s $\Po_1$ and $\Po_2$ are quasi-free, then
 the above morphism of
dg affine schemes is the same as a morphism of smooth dg manifolds, i.e.\ a morphism,
$$
\phi_{ind}: (\Hm(\Co_2,\mathcal{E}), Q_2) \lon
(\Hm(\Co_1,\mathcal{E}), Q_1),
$$
of $L_\infty$-algebras. The last statement of Theorem~\ref{morqis} follows
immediately from the formulae
(\ref{Q_n}) for $n=1$ and the K\"unneth formula completing the proof of Claim (i).

\sip

Claim (ii) is much easier than Claim (i): it follows
directly from the formulae
(\ref{Q_n}) for $n=1$.
\end{proo}

\sip


\subsection{Proofs via local coordinate computations} Calculations in local coordinates is a powerful and useful
tool in differential geometry. In this section we show a new proofs of Lemma~\ref{def-Q},
Proposition~\ref{scheme} and Theorem~\ref{mor} by explicitly describing all the notions and  constructions
of \S\ref{Derivative as dg affine scheme}
in local coordinates and justifying thereby the geometric language we used in that section.
 For simplicity, we show the proofs only for the case when
$(\cP, \p_P)$ and $(\cE, \p_\cE)$  are dg associative
algebras, that is, dg properads concentrated in biarity $(1,1)$
(a generalization
to arbitrary dg (prop)erads is straightforward); moreover, to simplify Koszul signs in the formulae
below we also assume
 that  both $\cP$ and  $\cE$ are free modules
over a graded commutative ring, $R=\oplus_{i\in \bZ} R^i$, with  degree 0 generators $\{e_a\}_{a\in I}$
and, respectively, $\{e_\al\}_{\al\in J}$. Then multiplications and differentials in $\cP$ and $\cE$
have the following coordinate representations,
$$
e_a\cdot e_b= \sum_{c\in I} \mu_{ab}^c e_c, \ \ \ e_\al\cdot a_\be= \sum_{\ga\in J} \mu_{\al\be}^\ga e_\ga,
$$
$$
\p_{\cP}e_a=\sum_{b\in I} D_a^b e_b, \ \ \  \p_{\cE}e_\al=\sum_{\be\in J} D_\al^\be e_\be,
$$
for some coefficients $\mu_{ab}^c, \mu_{\al\be}^\ga\in R^0$ and $D_a^b, D_\al^\be\in R^{-1}$.
Equations $\p_\cP^2=\p_\cE^2=0$ as well as equations for compatibility
of differentials with products are given in coordinates as follows,
\Beq\label{D-squraed}
\sum_{b\in I}D_a^bD_b^c=0, \ \ \ \ \ \ \sum_{\be\in J}D_\al^\be D_\be^\ga=0,
\Eeq
\Beq\label{D-Leibniz}
D_a^m\mu_{mb}^c + D_b^m\mu_{am}^c=\mu_{ab}^mD_m^c, \ \ \ \ \ \ \
D_\al^\nu\mu_{\nu \be}^\ga +  D_\be^\nu\mu_{\al\nu}^\ga=\mu_{\al\be}^\nu D_\nu^\ga.
\Eeq

A generic homogeneous map of graded vector space, $\ga: \cP\rar \cE$,
of degree $i\in \bZ$ is unique;y determined by its values on the generators,
$$
\ga(e_a)=\sum_{\al\in J} \ga_{a (i)}^{\al} e_\al,
$$
for some coefficients $\ga_{a (i)}^{\al} \in R^i$. We shall understand these coefficients as coordinates
on the flat manifold $\Hom_\bullet^\bS(\cP,\cE)$.

\sip

Consider now a completed free graded commutative algebra, $R[[ \ga_{a (i)}^{\al}]]$,
generated by formal variables $\ga_{a (i)}^{\al}$
to which we assign degree $i$. This algebra  is precisely the algebra of smooth functions,
$\f_{\Hom_\bullet^\bS(\cP,\cE)}$,
on the manifold $\Hom_\bullet^\bS(\cP,\cE)$. Let us consider a  degree
$-1$ vector field
(that is, a derivation of $\f_{\Hom_\bullet^\bS(\cP,\cE)}$),
$$
Q= \left(\sum_{\al,\be,a, i}D_\be^\al \ga_{a (i)}^\be - \sum_{a,b,\al,i}(-1)^i
D_a^b \ga_{b (i)}^{\al}\right)
\frac{\p}{\p \ga_{a (i)}^{\al}},
$$
on  $\Hom_\bullet^\bS(\cP,\cE)$.
In view of (\ref{D-squraed}), we have
\Beqrn
[Q,Q]&=& 2\left(\sum_{\al,\be,\ga, a, i}-D_\be^\al D^\be_\ga \ga_{a (i)}^\ga - \sum_{\al,\be,a,b, i}(-1)^{i}
D_a^b D_\be^\al\ga_{b (i)}^{\be}
\right)\frac{\p}{\p \ga_{a (i)}^{\al}}\\
&& + 2\left(-\sum_{\al,a,b,c, i}(-1)^{i}D_\be^\al D_a^b \ga_{b (i)}^{\be}
-
\sum_{\al,\be,\ga, a, i}D_a^b D^c_b \ga_{c (i)}^\al
\right)\frac{\p}{\p \ga_{a (i)}^{\al}}\\
&=& 0,
\Eeqrn
proving thereby Lemma~\ref{def-Q} which claims that $(\Hom_\bullet^\bS(\cP,\cE), Q)$ is a dg manifold.

\sip

The space of extended morphisms of associative $R$-algebras $\Mor_\bZ(\Po,\mathcal{E})$ is, by its definition,
a (singular, in general) subspace of the manifold
$\Hom_\bullet^\bS(\cP,\cE)$ given explicitly by the following equations,
$$
\Mor_\bZ(\Po,\mathcal{E}):=\left\{\ga_{a (i)}^{\al}\in \Hom_\bullet^\bS(\cP,\cE):
 \sum_{c\in I}\mu_{ab}^c \ga_{c (i)}^{\al} - \sum_{\be,\ga\in J \atop j+k=i} \mu_{\be\ga}^\al
\ga_{a (j)}^{\be}\ga_{b (k)}^{\ga}=0\right\}.
$$
Let $I$ be an ideal in $\f_{\Hom_\bullet^\bS(\cP,\cE)}=\bK[[ \ga_{a (i)}^{\al}]]$ generated by the functions
$\{ \sum_{c\in I}\mu_{ab}^c \ga_{c (i)}^{\al} - \sum_{\be,\ga\in J\atop j+k=i} \mu_{\be\ga}^\al
\ga_{a (j)}^{\be}\ga_{b (k)}^{\ga}\}$. Then the structure sheaf, $\f_{\Mor_\bZ(\Po,\mathcal{E})}$, of the scheme
$\Mor_\bZ(\Po,\mathcal{E})$ is given, by definition, by the quotient algebra $\f_{\Hom_\bullet^\bS(\cP,\cE)}/I$
(which, in general, is not freely generated, i.e. is not smooth). We claim that the vector field $Q$
on $\Hom_\bullet^\bS(\cP,\cE)$ ia tangent to the subspace $\Mor_\bZ(\Po,\mathcal{E})$. Indeed, in view of (\ref{D-Leibniz}),
we have
\Beqrn
Q\left( \sum_{c\in I}\mu_{ab}^c \ga_{c (i)}^{\al} - \sum_{\be,\ga\in J \atop j+k=i} \mu_{\be\ga}^\al
\ga_{a (j)}^{\be}\ga_{b (k)}^{\ga}\right)
&=&\sum_{\sigma\in J} D_\sigma^\al\left( \sum_{c\in I}\mu_{ab}^c \ga_{c (i)}^{\sigma}
- \sum_{\be\ga\in J\atop j+k=i} \mu_{\be\ga}^\sigma
\ga_{a (j)}^{\be}\ga_{b (k)}^{\ga}\right)\\
&& -(-1)^i
\sum_{e\in E} D_a^m\left( \sum_{c\in I}\mu_{mb}^c \ga_{c (i)}^{\sigma}
- \sum_{\be\ga\in J\atop j+k=i} \mu_{\be\ga}^\sigma
\ga_{m (j)}^{\be}\ga_{b (k)}^{\ga}\right)\\
&& -(-1)^i
\sum_{e\in E} D_b^m\left( \sum_{c\in I}\mu_{am}^c \ga_{c (i)}^{\sigma}
- \sum_{\be\ga\in J\atop j+k=i} \mu_{\be\ga}^\sigma
\ga_{a (j)}^{\be}\ga_{m (k)}^{\ga}\right)\\
\Eeqrn
Thus $Q(I)\subset I$ so that $Q$ makes $\f_{\Mor_\bZ(\Po,\mathcal{E})}$ into a {\em differential}\, graded
algebra proving thereby Proposition~\ref{scheme}.

To prove Theorem~\ref{mor} we have to assume from now on that $\cP$ is a free algebra, $\ot^\bullet V$, generated by
some free $R$-module $V$. Let $\{e_A\}_{A\in K}$ stand for a set of generators of $V$ so that the basis $\{e_a\}$ we used above
can be identified with the following set,
$$
\{e_a\}_{a\in I} =\{e_A, \ e_{A_1A_2}:= e_{A_1}\ot e_{A_2},\ e_{A_1A_2A_3}=  e_{A_1}\ot e_{A_2}\ot e_{A_3},\
\ldots\}_{A_\bullet
\in K}.
$$
The differential $\p_\cP$ is now completely determined by its values on the generators $\{e_A\}$,
$$
\p_\cP e_A=\sum_{k\geq1}\sum_{A_1,\ldots,A_k\in K} D_A^{A_1\ldots A_k} e_{A_1\ldots A_k},
$$
for some coefficients  $D_A^{A_1\ldots A_k}\in R^{-1}$. On the other hand, the $R$-algebra of smooth formal functions on the manifold
$\Hom_\bullet^\bS(\cP,\cE)$ gets the following explicit representation,
$$
\f_{\Hom_\bullet^\bS(\cP,\cE)} = R[[\ga_{A(i)}^\al,\ \ga_{A_1A_2(i)}^\al,\  \ga_{A_1A_2A_3(i)}^\al, \ldots]].
$$
The key point is that the system of equations defining the subspace
$\Mor_\bZ(\Po,\mathcal{E})\subset \Hom_\bullet^\bS(\cP,\cE)$ can now be easily solved,
\Beqrn
\ga_{A_1A_2(i)}^\al &=& \sum_{\be_1,\be_2\in J\atop i_1+i_2=i} \mu_{\be_1\be_2}^\al \ga_{A_1(i_1)}^{\be_1}
\ga_{A_2(i_2)}^{\be_2}, \\
\ga_{A_1A_2A_3(i)}^\al &=& \sum_{\be_\bullet,\ga \in J\atop i_1+ i_2 +i_3=i}
\mu_{\be_1\ga}^\al\mu_{\be_2\be_3}^{\ga}
 \ga_{A_1(i_1)}^{\be_1}  \ga_{A_2(i_2)}^{\be_2}\ga_{A_3(i_3)}^{\be_3},
,\\
\ldots &&
\Eeqrn
in terms of the
independent variables $\ga_{A(i)}^\al$. Thus $\Mor_\bZ(\Po,\mathcal{E})$ is itself a
smooth formal manifold
with the structure sheaf $\f_{\Mor_\bZ(\Po,\mathcal{E})}\simeq R[[\ga_{A(i)}^\al]]$.
The vector field $Q$ on the manifold $\Hom_\bullet^\bS(\cP,\cE)$
 restricts to a smooth degree -1 homological vector field on
 the subspace $\Mor_\bZ(\Po,\mathcal{E})\subset \Hom_\bullet^\bS(\cP,\cE)$
which is given explicitly as follows,
$$
Q|_{\Mor_\bZ(\Po,\mathcal{E})} =
\left(\sum_{\al,\be,a, i}D_\be^\al \ga_{A (i)}^\be - \sum_{A, A_\bullet,\al,i}(-1)^i
D_A^{A_1 \ldots A_k} \ga_{A_1\ldots A_k (i)}^{\al}\right)
\frac{\p}{\p \ga_{A (i)}^{\al}}
$$
where, for $k\geq 2$,
$$
\ga_{A_1A_2\ldots A_k(i)}^\al = \sum_{\be_\bullet,\ga_\bullet \in J\atop i_1+\ldots +i_k=i}
\mu_{\be_1\ga_1}^\al\mu_{\be_2\ga_2}^{\ga_1}\ldots \mu_{\be_{k-1}\be_k}^{\ga_{k-2}}
 \ga_{A_1(i_1)}^{\be_1}  \ga_{A_2(i_2)}^{\be_2}\ldots  \ga_{A_k(i_k)}^{\be_k}.
$$
Thus $\Mor_\bZ(\Po,\mathcal{E})=\Hom_\bZ(V, \cE)$ is canonically a dg manifold, i.e.\
$s\Hom_\bZ(V, \cE)$ is canonically a $L_\infty$-algebra, and Theorem~\ref{mor}(i) is proved.
Theorem~\ref{mor}(ii) follows from the above explicit expression for the homological vector field
$Q|_{\Mor_\bZ(\Po,\mathcal{E})}$ as its zero set is precisely the set of morphisms
$\cP\rar \cE$ which commute with the differentials. Finally, if the differential $\p_\cP$ is at most
quadratic in generators, then
 $D_A^{A_1\ldots A_k}=0$ for $k\geq 3$ and hence
$Q|_{\Mor_\bZ(\Po,\mathcal{E})}$ is evidently at most quadratic homological vector field so that
Theorem~\ref{mor}(iii) is also done.

In a similar purely geometric way one can prove anew Theorem~\ref{morqis}.  We leave the
details as an exercise to the interested reader.

\subsection{Enlarged category of dg prop(erad)s} For any dg prop(erad)s $(\cP_1,\p_1)$,
 $(\cP_2,\p_2)$ and  $(\cP_3,\p_3)$, the natural composition map,
$$
\Ba{ccc}
 \Hm(\Po_2,\Po_3)\ot \Hm(\Po_1,\Po_2) & \lon &\Hm(\Po_1, \Po_3)\\
 \ga_2\ot \ga_1 & \lon & \ga_2\circ \ga_1
 \Ea
$$
respects the relations (\ref{I}) and hence induces a map of coalgebras (cf.\ the proof
of Theorem~\ref{morqis}(i)),
$$
\circ: \f_{{\mathrm Mor}_\bZ(\Po_2, \Po_3)}\ot \f_{{\mathrm Mor}_\bZ(\Po_1, \Po_2)} \lon
\f_{{\mathrm Mor}_\bZ(\Po_1, {\mathcal Q})}.
$$

\bip

\begin{pro}
The map $\circ$ respects the codifferentials (\ref{Q}), i.e.\ induces a morphism of dg affine schemes,
$$
\left(\Mor_\bZ(\cP_2,\cP_3), Q_{23}\right)\times \left(\Mor_\bZ(\cP_1,\cP_2), Q_{12}\right)
\lon \left(\Mor_\bZ(\cP_1,\cP_3), Q_{13}\right).
$$
\end{pro}
\sip

\begin{proo} We have, for any $\ga_1\in  \Hm(\Po_1,\Po_2)$ and  $\ga_2\in  \Hm(\Po_2,\Po_3)$,
\Beqrn
Q_{13}(\ga_2\circ \ga_1)&\stackrel{{\mathrm b \mathrm y}(\ref{Q})}{=}& \p_3\circ \ga_2\circ \ga_1 - (-1)^{\ga_1+\ga_2}
\ga_2\circ \ga_1\circ \p_1\\
&=&\p_3\circ \ga_2\circ \ga_1 - (-1)^{\ga_2}\ga_2\circ \p_2\circ \ga_1 \\
&& + (-1)^{\ga_2}\circ \p_2\circ \ga_1 - (-1)^{\ga_1+\ga_2}
\ga_2\circ \ga_1\circ \p_1\\
&=& Q_{23}(\ga_2)\circ \ga_1 + (-1)^{\ga_2}\ga_2\circ Q_{12}(\ga_1).
\Eeqrn
\end{proo}

As the composition $\circ$ is obviously associative, we end up with the following canonical enlargement
of the category of dg prop(erad)s.

\begin{cor}
The data,
\Beqrn
 \mbox{\sf Objects} &:=&  \mbox{ dg prop(erad)s}, \\
\mbox{\sf Hom}(\cP_1,\cP_2) &:=& \mbox{the dg affine}\
 {\mbox scheme}\ \left(\Mor_\bZ(\cP_1,\cP_2), Q_{12}\right)
\Eeqrn
is a category. Moreover, the composition,
$$
\circ:\mbox{\sf Hom}(\cP_2,\cP_3)\times
\mbox{\sf Hom}(\cP_1,\cP_2) \lon \mbox{\sf Hom}(\cP_1,\cP_3)
$$
is a morphism of dg affine schemes.
\end{cor}

Note that if $\cP_1$ is quasi-free then, by Theorem~\ref{mor},  $\mbox{\sf Hom}(\cP_1,\cP_2)$ is precisely
the filtered $L_\infty$-algebra whose Maurer-Cartan elements are ordinary morphisms of dg prop(erad)s
from $\cP_1$ to $\cP_2$. Note also that if $\phi: \cP_1\rar \cP_2$ is an ordinary morphism of quasi-free
 dg prop(erad)s, then
 the composition map,
$$
\circ:\mbox{\sf Hom}(\cP_2,\cP_3)\times
\phi \lon \mbox{\sf Hom}(\cP_1,\cP_3)
$$
is precisely the $L_\infty$-morphism of Theorem~\ref{morqis}(i).


\subsection{Families of natural $L_\infty$-structures on $\oplus \Po$}
\label{oplusP2}
 It was shown in Section~\ref{oplusP1} that for any homotopy properad $\Po$
 the associated direct sum $\oplus \Po := \bigoplus_{m,n} \Po(m,n)$
 has a natural structure of $L_\infty$-algebra which encodes all
 possible compositions in $\Po$. In this section we show a new proof of this result which is independent
 of Section~\ref{oplusP1} and the earlier works \cite{KapranovManin01,VanderLaan02} which treated the special case
 of operads. The present approach is based on certain universal properties of the properad of Frobenius algebras (and its
 non-commutative versions) and Theorem~\ref{mor}; it provides  a conceptual explanation
of the phenomenon in terms of convolution properads.

\begin{thm}\label{main}
Let $\cP=\{P(m,n)\}$ be a homotopy prop(erad). Then
\begin{itemize}
\item[(i)] $\bigoplus_{m,n} P(m,n)$ is canonically a
$L_\infty$-algebra; \item[(ii)] $\bigoplus_{m,n} P(m,n)^{\bS_m}$
is canonically a $L_\infty$-algebra; \item[(iii)] $\bigoplus_{m,n}
P(m,n)^{\bS_n}$ is canonically a $L_\infty$-algebra; \item[(iv)]
$\bigoplus_{m,n} P(m,n)^{\bS_m\times \bS_n}$ is canonically a
$L_\infty$-algebra; \item[(v)] there is a natural commutative
diagram of $L_\infty$-morphisms,


$$
\xymatrix{
& \bigoplus_{m,n} P(m,n)^{\bS_m} \ar[dr] & \\
 \bigoplus_{m,n} P(m,n) \ar[ur] \ar[dr]  & &\bigoplus_{m,n} P(m,n)^{\bS_m\times \bS_n}  \\
& \bigoplus_{m,n} P(m,n)^{\bS_n} \ar[ur] & }
$$
\end{itemize}

Finally, if $\cP$ is a dg properad, then all the above data are dg
Lie algebras and morphisms of dg Lie algebras.
\end{thm}

\begin{proo} Recall that the prop(erad) of Frobenius algebras can be defined it as a quotient,
$$
\Frob:= {\cF\langle V\rangle}/(R)
$$
of the free prop(erad), $\cF(V)$, generated by the $\bS$-bimodule
$V=\{V(m,n)\}$,
\[
V(m,n):=\left\{
\Ba{rr}
\Id_2\ot \Id_1\equiv\mbox{span}\left\langle
\begin{xy}
 <0mm,-0.55mm>*{};<0mm,-2.5mm>*{}**@{-},
 <0.5mm,0.5mm>*{};<2.2mm,2.2mm>*{}**@{-},
 <-0.48mm,0.48mm>*{};<-2.2mm,2.2mm>*{}**@{-},
 <0mm,0mm>*{\bullet};<0mm,0mm>*{}**@{},
 <0mm,-0.55mm>*{};<0mm,-3.8mm>*{_1}**@{},
 <0.5mm,0.5mm>*{};<2.7mm,2.8mm>*{^2}**@{},
 <-0.48mm,0.48mm>*{};<-2.7mm,2.8mm>*{^1}**@{},
 \end{xy}
=
\begin{xy}
 <0mm,-0.55mm>*{};<0mm,-2.5mm>*{}**@{-},
 <0.5mm,0.5mm>*{};<2.2mm,2.2mm>*{}**@{-},
 <-0.48mm,0.48mm>*{};<-2.2mm,2.2mm>*{}**@{-},
 <0mm,0mm>*{\bullet};<0mm,0mm>*{}**@{},
 <0mm,-0.55mm>*{};<0mm,-3.8mm>*{_1}**@{},
 <0.5mm,0.5mm>*{};<2.7mm,2.8mm>*{^1}**@{},
 <-0.48mm,0.48mm>*{};<-2.7mm,2.8mm>*{^2}**@{},
 \end{xy}
   \right\rangle  & \mbox{if}\ m=2, n=1,\vspace{3mm}\\
\Id_1\ot \Id_2\equiv
\mbox{span}\left\langle
\begin{xy}
 <0mm,0.66mm>*{};<0mm,3mm>*{}**@{-},
 <0.39mm,-0.39mm>*{};<2.2mm,-2.2mm>*{}**@{-},
 <-0.35mm,-0.35mm>*{};<-2.2mm,-2.2mm>*{}**@{-},
 <0mm,0mm>*{\bullet};<0mm,0mm>*{}**@{},
   <0mm,0.66mm>*{};<0mm,3.4mm>*{^1}**@{},
   <0.39mm,-0.39mm>*{};<2.9mm,-4mm>*{^2}**@{},
   <-0.35mm,-0.35mm>*{};<-2.8mm,-4mm>*{^1}**@{},
\end{xy}=
\begin{xy}
 <0mm,0.66mm>*{};<0mm,3mm>*{}**@{-},
 <0.39mm,-0.39mm>*{};<2.2mm,-2.2mm>*{}**@{-},
 <-0.35mm,-0.35mm>*{};<-2.2mm,-2.2mm>*{}**@{-},
 <0mm,0mm>*{\bullet};<0mm,0mm>*{}**@{},
   <0mm,0.66mm>*{};<0mm,3.4mm>*{^1}**@{},
   <0.39mm,-0.39mm>*{};<2.9mm,-4mm>*{^1}**@{},
   <-0.35mm,-0.35mm>*{};<-2.8mm,-4mm>*{^2}**@{},
\end{xy}
\right\rangle
\ & \mbox{if}\ m=1, n=2, \vspace{3mm}\\
0 & \mbox{otherwise}
\Ea
\right.
\]
modulo the ideal generated by relations
$$
\begin{xy}
 <0mm,0mm>*{\bullet};<0mm,0mm>*{}**@{},
 <0mm,-0.49mm>*{};<0mm,-3.0mm>*{}**@{-},
 <0.49mm,0.49mm>*{};<1.9mm,1.9mm>*{}**@{-},
 <-0.5mm,0.5mm>*{};<-1.9mm,1.9mm>*{}**@{-},
 <-2.3mm,2.3mm>*{\bullet};<-2.3mm,2.3mm>*{}**@{},
 <-1.8mm,2.8mm>*{};<0mm,4.9mm>*{}**@{-},
 <-2.8mm,2.9mm>*{};<-4.6mm,4.9mm>*{}**@{-},
   <0.49mm,0.49mm>*{};<2.7mm,2.3mm>*{^3}**@{},
   <-1.8mm,2.8mm>*{};<0.4mm,5.3mm>*{^2}**@{},
   <-2.8mm,2.9mm>*{};<-5.1mm,5.3mm>*{^1}**@{},
 \end{xy}
\ - \
\begin{xy}
 <0mm,0mm>*{\bullet};<0mm,0mm>*{}**@{},
 <0mm,-0.49mm>*{};<0mm,-3.0mm>*{}**@{-},
 <0.49mm,0.49mm>*{};<1.9mm,1.9mm>*{}**@{-},
 <-0.5mm,0.5mm>*{};<-1.9mm,1.9mm>*{}**@{-},
 <2.3mm,2.3mm>*{\bullet};
 <2.3mm,2.3mm>*{};<4.3mm,4.9mm>*{}**@{-},
 <2.3mm,2.3mm>*{};<0mm,4.9mm>*{}**@{-},
   <0.49mm,0.49mm>*{};<-2.3mm,2.3mm>*{^1}**@{},
   <-1.8mm,2.8mm>*{};<-0.2mm,5.3mm>*{^2}**@{},
   <-2.8mm,2.9mm>*{};<5.0mm,5.3mm>*{^3}**@{},
 \end{xy}
\ \ \ , \ \ \
 \begin{xy}
 <0mm,0mm>*{\bullet};<0mm,0mm>*{}**@{},
 <0mm,0.69mm>*{};<0mm,3.0mm>*{}**@{-},
 <0.39mm,-0.39mm>*{};<2.4mm,-2.4mm>*{}**@{-},
 <-0.35mm,-0.35mm>*{};<-1.9mm,-1.9mm>*{}**@{-},
 <-2.4mm,-2.4mm>*{\bullet};<-2.4mm,-2.4mm>*{}**@{},
 <-2.0mm,-2.8mm>*{};<0mm,-4.9mm>*{}**@{-},
 <-2.8mm,-2.9mm>*{};<-4.7mm,-4.9mm>*{}**@{-},
    <0.39mm,-0.39mm>*{};<3.3mm,-4.0mm>*{^3}**@{},
    <-2.0mm,-2.8mm>*{};<0.5mm,-6.7mm>*{^2}**@{},
    <-2.8mm,-2.9mm>*{};<-5.2mm,-6.7mm>*{^1}**@{},
 \end{xy}
\ - \
 \begin{xy}
 <0mm,0mm>*{\bullet};<0mm,0mm>*{}**@{},
 <0mm,0.69mm>*{};<0mm,3.0mm>*{}**@{-},
 <0.39mm,-0.39mm>*{};<2.4mm,-2.4mm>*{}**@{-},
 <-0.35mm,-0.35mm>*{};<-1.9mm,-1.9mm>*{}**@{-},
 <2.3mm,-2.4mm>*{\bullet};
 <2.3mm,-2.8mm>*{};<4.6mm,-4.9mm>*{}**@{-},
 <2.3mm,-2.9mm>*{};<-0mm,-4.9mm>*{}**@{-},
    <0.39mm,-0.39mm>*{};<-2.4mm,-4.0mm>*{^1}**@{},
    <-2.0mm,-2.8mm>*{};<0.5mm,-6.7mm>*{^2}**@{},
    <-2.8mm,-2.9mm>*{};<4.7mm,-6.7mm>*{^3}**@{},
 \end{xy}
\ \ \ , \ \ \
 \begin{xy}
 <0mm,2.47mm>*{};<0mm,0.12mm>*{}**@{-},
 <0.5mm,3.5mm>*{};<2.2mm,5.2mm>*{}**@{-},
 <-0.48mm,3.48mm>*{};<-2.2mm,5.2mm>*{}**@{-},
 <0mm,3mm>*{\bullet};<0mm,3mm>*{}**@{},
  <0mm,-0.8mm>*{\bullet};<0mm,-0.8mm>*{}**@{},
<-0.39mm,-1.2mm>*{};<-2.2mm,-3.5mm>*{}**@{-},
 <0.39mm,-1.2mm>*{};<2.2mm,-3.5mm>*{}**@{-},
     <0.5mm,3.5mm>*{};<2.8mm,5.7mm>*{^2}**@{},
     <-0.48mm,3.48mm>*{};<-2.8mm,5.7mm>*{^1}**@{},
   <0mm,-0.8mm>*{};<-2.7mm,-5.2mm>*{^1}**@{},
   <0mm,-0.8mm>*{};<2.7mm,-5.2mm>*{^2}**@{},
\end{xy}
\  - \
\begin{xy}
 <0mm,-1.3mm>*{};<0mm,-3.5mm>*{}**@{-},
 <0.38mm,-0.2mm>*{};<2.0mm,2.0mm>*{}**@{-},
 <-0.38mm,-0.2mm>*{};<-2.2mm,2.2mm>*{}**@{-},
<0mm,-0.8mm>*{\bullet};<0mm,0.8mm>*{}**@{},
 <2.4mm,2.4mm>*{\bullet};<2.4mm,2.4mm>*{}**@{},
 <2.77mm,2.0mm>*{};<4.4mm,-0.8mm>*{}**@{-},
 <2.4mm,3mm>*{};<2.4mm,5.2mm>*{}**@{-},
     <0mm,-1.3mm>*{};<0mm,-5.3mm>*{^1}**@{},
     <2.5mm,2.3mm>*{};<5.1mm,-2.6mm>*{^2}**@{},
    <2.4mm,2.5mm>*{};<2.4mm,5.7mm>*{^2}**@{},
    <-0.38mm,-0.2mm>*{};<-2.8mm,2.5mm>*{^1}**@{},
 \end{xy} \ .
$$
Here $\Id_n$ stands for the trivial one dimensional representation
of the group $\Sy_n$. It is clear that $\Frob(m,n)=\Id_m\ot \Id_n$
and the compositions in $\Frob$ are determined by the canonical
isomorphism $\bK \ot \bK \rar \bK$ (thus $\Frob$ is a prop(erad)ic
analogue of $\Com$ in the theory of operads). The dual  space,
$\Frob^*$, is naturally a coprop(erad) \footnote{In fact,
$\Frob^*$ is a completed coproperad with respect to the topology
induced by the number of vertices. The formulae for the composite
coproduct in infinite. But since we 'dualize' it by considering
the convolution homotopy properad $\Hom^{\Sy}(\Frob^*, -)$ it does
not matter.}. Homotopy prop(erad) structure on $\Po$ is the same
as a differential, $d_\Po$, in the free coprop(erad) $\cF^c(
s\bar{\Po} )$. Theorem~\ref{comor}(i) applied to the coprop(erad)s
$\Frob^*$ and $(\cF^c( s\bar{\Po} ), d_\Po)$, asserts that the
vector space
$$
\Hm\left(\Frob^*, \Po\right)=
 \bigoplus_{m,n}(\Id_m\ot \Id_n) \ot_{\Sy_m\times \Sy_n} \Po(m,n)= \bigoplus_{m,n} \Po(m,n)^{\Sy_m\times \Sy_n}
$$
is canonically a $L_\infty$-algebra. Hence the claim (iv).

\sip

Let us next define a non-commutative analogue of $\Frob$ as a quotient,
$$
\Frob_+^+:= {\cF( V)}/(R)
$$
of the free prop(erad), $\cF( V )$, generated by the
$\bS$-bimodule $V=\{V(m,n)\}$,
\[
V(m,n):=\left\{
\Ba{rr}
\bK[\Sy_2]\ot \Id_1\equiv\mbox{span}\left\langle
\begin{xy}
 <0mm,-0.55mm>*{};<0mm,-2.5mm>*{}**@{-},
 <0.5mm,0.5mm>*{};<2.2mm,2.2mm>*{}**@{-},
 <-0.48mm,0.48mm>*{};<-2.2mm,2.2mm>*{}**@{-},
 <0mm,0mm>*{\bullet};<0mm,0mm>*{}**@{},
 <0mm,-0.55mm>*{};<0mm,-3.8mm>*{_1}**@{},
 <0.5mm,0.5mm>*{};<2.7mm,2.8mm>*{^2}**@{},
 <-0.48mm,0.48mm>*{};<-2.7mm,2.8mm>*{^1}**@{},
 \end{xy}
\ , \
\begin{xy}
 <0mm,-0.55mm>*{};<0mm,-2.5mm>*{}**@{-},
 <0.5mm,0.5mm>*{};<2.2mm,2.2mm>*{}**@{-},
 <-0.48mm,0.48mm>*{};<-2.2mm,2.2mm>*{}**@{-},
 <0mm,0mm>*{\bullet};<0mm,0mm>*{}**@{},
 <0mm,-0.55mm>*{};<0mm,-3.8mm>*{_1}**@{},
 <0.5mm,0.5mm>*{};<2.7mm,2.8mm>*{^1}**@{},
 <-0.48mm,0.48mm>*{};<-2.7mm,2.8mm>*{^2}**@{},
 \end{xy}
   \right\rangle  & \mbox{if}\ m=2, n=1,\vspace{3mm}\\
\Id_1\ot \bK[\Sy_2]\equiv
\mbox{span}\left\langle
\begin{xy}
 <0mm,0.66mm>*{};<0mm,3mm>*{}**@{-},
 <0.39mm,-0.39mm>*{};<2.2mm,-2.2mm>*{}**@{-},
 <-0.35mm,-0.35mm>*{};<-2.2mm,-2.2mm>*{}**@{-},
 <0mm,0mm>*{\bullet};<0mm,0mm>*{}**@{},
   <0mm,0.66mm>*{};<0mm,3.4mm>*{^1}**@{},
   <0.39mm,-0.39mm>*{};<2.9mm,-4mm>*{^2}**@{},
   <-0.35mm,-0.35mm>*{};<-2.8mm,-4mm>*{^1}**@{},
\end{xy}
\ , \
\begin{xy}
 <0mm,0.66mm>*{};<0mm,3mm>*{}**@{-},
 <0.39mm,-0.39mm>*{};<2.2mm,-2.2mm>*{}**@{-},
 <-0.35mm,-0.35mm>*{};<-2.2mm,-2.2mm>*{}**@{-},
 <0mm,0mm>*{\bullet};<0mm,0mm>*{}**@{},
   <0mm,0.66mm>*{};<0mm,3.4mm>*{^1}**@{},
   <0.39mm,-0.39mm>*{};<2.9mm,-4mm>*{^1}**@{},
   <-0.35mm,-0.35mm>*{};<-2.8mm,-4mm>*{^2}**@{},
\end{xy}
\right\rangle
\ & \mbox{if}\ m=1, n=2, \vspace{3mm}\\
0 & \mbox{otherwise}
\Ea
\right.
\]
modulo the ideal generated by relations
$$
\begin{xy}
 <0mm,0mm>*{\bullet};<0mm,0mm>*{}**@{},
 <0mm,-0.49mm>*{};<0mm,-3.0mm>*{}**@{-},
 <0.49mm,0.49mm>*{};<1.9mm,1.9mm>*{}**@{-},
 <-0.5mm,0.5mm>*{};<-1.9mm,1.9mm>*{}**@{-},
 <-2.3mm,2.3mm>*{\bullet};<-2.3mm,2.3mm>*{}**@{},
 <-1.8mm,2.8mm>*{};<0mm,4.9mm>*{}**@{-},
 <-2.8mm,2.9mm>*{};<-4.6mm,4.9mm>*{}**@{-},
   <0.49mm,0.49mm>*{};<2.7mm,2.3mm>*{^3}**@{},
   <-1.8mm,2.8mm>*{};<0.4mm,5.3mm>*{^2}**@{},
   <-2.8mm,2.9mm>*{};<-5.1mm,5.3mm>*{^1}**@{},
 \end{xy}
\ - \
\begin{xy}
 <0mm,0mm>*{\bullet};<0mm,0mm>*{}**@{},
 <0mm,-0.49mm>*{};<0mm,-3.0mm>*{}**@{-},
 <0.49mm,0.49mm>*{};<1.9mm,1.9mm>*{}**@{-},
 <-0.5mm,0.5mm>*{};<-1.9mm,1.9mm>*{}**@{-},
 <2.3mm,2.3mm>*{\bullet};
 <2.3mm,2.3mm>*{};<4.3mm,4.9mm>*{}**@{-},
 <2.3mm,2.3mm>*{};<0mm,4.9mm>*{}**@{-},
   <0.49mm,0.49mm>*{};<-2.3mm,2.3mm>*{^1}**@{},
   <-1.8mm,2.8mm>*{};<-0.2mm,5.3mm>*{^2}**@{},
   <-2.8mm,2.9mm>*{};<5.0mm,5.3mm>*{^3}**@{},
 \end{xy}
 =0
\ \ \ , \ \ \
 \begin{xy}
 <0mm,0mm>*{\bullet};<0mm,0mm>*{}**@{},
 <0mm,0.69mm>*{};<0mm,3.0mm>*{}**@{-},
 <0.39mm,-0.39mm>*{};<2.4mm,-2.4mm>*{}**@{-},
 <-0.35mm,-0.35mm>*{};<-1.9mm,-1.9mm>*{}**@{-},
 <-2.4mm,-2.4mm>*{\bullet};<-2.4mm,-2.4mm>*{}**@{},
 <-2.0mm,-2.8mm>*{};<0mm,-4.9mm>*{}**@{-},
 <-2.8mm,-2.9mm>*{};<-4.7mm,-4.9mm>*{}**@{-},
    <0.39mm,-0.39mm>*{};<3.3mm,-4.0mm>*{^3}**@{},
    <-2.0mm,-2.8mm>*{};<0.5mm,-6.7mm>*{^2}**@{},
    <-2.8mm,-2.9mm>*{};<-5.2mm,-6.7mm>*{^1}**@{},
 \end{xy}
\ - \
 \begin{xy}
 <0mm,0mm>*{\bullet};<0mm,0mm>*{}**@{},
 <0mm,0.69mm>*{};<0mm,3.0mm>*{}**@{-},
 <0.39mm,-0.39mm>*{};<2.4mm,-2.4mm>*{}**@{-},
 <-0.35mm,-0.35mm>*{};<-1.9mm,-1.9mm>*{}**@{-},
 <2.3mm,-2.4mm>*{\bullet};
 <2.3mm,-2.8mm>*{};<4.6mm,-4.9mm>*{}**@{-},
 <2.3mm,-2.9mm>*{};<-0mm,-4.9mm>*{}**@{-},
    <0.39mm,-0.39mm>*{};<-2.4mm,-4.0mm>*{^1}**@{},
    <-2.0mm,-2.8mm>*{};<0.5mm,-6.7mm>*{^2}**@{},
    <-2.8mm,-2.9mm>*{};<4.7mm,-6.7mm>*{^3}**@{},
 \end{xy}
 =0
\ \ \ , \ \ \
 \begin{xy}
 <0mm,2.47mm>*{};<0mm,0.12mm>*{}**@{-},
 <0.5mm,3.5mm>*{};<2.2mm,5.2mm>*{}**@{-},
 <-0.48mm,3.48mm>*{};<-2.2mm,5.2mm>*{}**@{-},
 <0mm,3mm>*{\bullet};<0mm,3mm>*{}**@{},
  <0mm,-0.8mm>*{\bullet};<0mm,-0.8mm>*{}**@{},
<-0.39mm,-1.2mm>*{};<-2.2mm,-3.5mm>*{}**@{-},
 <0.39mm,-1.2mm>*{};<2.2mm,-3.5mm>*{}**@{-},
     <0.5mm,3.5mm>*{};<2.8mm,5.7mm>*{^2}**@{},
     <-0.48mm,3.48mm>*{};<-2.8mm,5.7mm>*{^1}**@{},
   <0mm,-0.8mm>*{};<-2.7mm,-5.2mm>*{^1}**@{},
   <0mm,-0.8mm>*{};<2.7mm,-5.2mm>*{^2}**@{},
\end{xy}
\  - \
\begin{xy}
 <0mm,-1.3mm>*{};<0mm,-3.5mm>*{}**@{-},
 <0.38mm,-0.2mm>*{};<2.0mm,2.0mm>*{}**@{-},
 <-0.38mm,-0.2mm>*{};<-2.2mm,2.2mm>*{}**@{-},
<0mm,-0.8mm>*{\bullet};<0mm,0.8mm>*{}**@{},
 <2.4mm,2.4mm>*{\bullet};<2.4mm,2.4mm>*{}**@{},
 <2.77mm,2.0mm>*{};<4.4mm,-0.8mm>*{}**@{-},
 <2.4mm,3mm>*{};<2.4mm,5.2mm>*{}**@{-},
     <0mm,-1.3mm>*{};<0mm,-5.3mm>*{^1}**@{},
     <2.5mm,2.3mm>*{};<5.1mm,-2.6mm>*{^2}**@{},
    <2.4mm,2.5mm>*{};<2.4mm,5.7mm>*{^2}**@{},
    <-0.38mm,-0.2mm>*{};<-2.8mm,2.5mm>*{^1}**@{},
 \end{xy}
 =0 \ .
 $$
$$
 \begin{xy}
 <0mm,2.47mm>*{};<0mm,0.12mm>*{}**@{-},
 <0.5mm,3.5mm>*{};<2.2mm,5.2mm>*{}**@{-},
 <-0.48mm,3.48mm>*{};<-2.2mm,5.2mm>*{}**@{-},
 <0mm,3mm>*{\bullet};<0mm,3mm>*{}**@{},
  <0mm,-0.8mm>*{\bullet};<0mm,-0.8mm>*{}**@{},
<-0.39mm,-1.2mm>*{};<-2.2mm,-3.5mm>*{}**@{-},
 <0.39mm,-1.2mm>*{};<2.2mm,-3.5mm>*{}**@{-},
     <0.5mm,3.5mm>*{};<2.8mm,5.7mm>*{^2}**@{},
     <-0.48mm,3.48mm>*{};<-2.8mm,5.7mm>*{^1}**@{},
   <0mm,-0.8mm>*{};<-2.7mm,-5.2mm>*{^1}**@{},
   <0mm,-0.8mm>*{};<2.7mm,-5.2mm>*{^2}**@{},
\end{xy}
\  - \
\begin{xy}
 <0mm,-1.3mm>*{};<0mm,-3.5mm>*{}**@{-},
 <0.38mm,-0.2mm>*{};<2.2mm,2.2mm>*{}**@{-},
 <-0.38mm,-0.2mm>*{};<-2.2mm,2.2mm>*{}**@{-},
<0mm,-0.8mm>*{\bullet};<0mm,0.8mm>*{}**@{},
 <-2.4mm,2.4mm>*{\bullet};
 <-2.5mm,2.3mm>*{};<-4.4mm,-0.8mm>*{}**@{-},
 <-2.4mm,2.5mm>*{};<-2.4mm,5.2mm>*{}**@{-},
 <0mm,0mm>*{};<3mm,3.0mm>*{^2}**@{},
     <0mm,0mm>*{};<-2.2mm,5.7mm>*{^1}**@{},
   <0mm,0mm>*{};<-4.6mm,-3.1mm>*{^1}**@{},
   <0mm,-0.8mm>*{};<0.2mm,-5.4mm>*{^2}**@{},
    \end{xy} \ =\ 0.
\ \ \ , \ \ \
\begin{xy}
 <0mm,2.47mm>*{};<0mm,0.12mm>*{}**@{-},
 <0.5mm,3.5mm>*{};<2.2mm,5.2mm>*{}**@{-},
 <-0.48mm,3.48mm>*{};<-2.2mm,5.2mm>*{}**@{-},
 <0mm,3mm>*{\bullet};<0mm,3mm>*{}**@{},
  <0mm,-0.8mm>*{\bullet};<0mm,-0.8mm>*{}**@{},
<-0.39mm,-1.2mm>*{};<-2.2mm,-3.5mm>*{}**@{-},
 <0.39mm,-1.2mm>*{};<2.2mm,-3.5mm>*{}**@{-},
     <0.5mm,3.5mm>*{};<2.8mm,5.7mm>*{^2}**@{},
     <-0.48mm,3.48mm>*{};<-2.8mm,5.7mm>*{^1}**@{},
   <0mm,-0.8mm>*{};<-2.7mm,-5.2mm>*{^1}**@{},
   <0mm,-0.8mm>*{};<2.7mm,-5.2mm>*{^2}**@{},
\end{xy}
\  - \
\begin{xy}
 <0mm,-0.8mm>*{};<-3.2mm,1.2mm>*{}**@{-},
 <-3.2mm,1.2mm>*{};<0mm,3mm>*{}**@{-},
 <0mm,0mm>*{};<0mm,-3.2mm>*{}**@{-},
 <0mm,3mm>*{};<0mm,6.0mm>*{}**@{-},
<0mm,3mm>*{\bullet};
 <0mm,-0.7mm>*{\bullet};
 <0mm,-0.8mm>*{};<4.4mm,3mm>*{}**@{-},
 <0mm,3mm>*{};<4.4mm,-0.8mm>*{}**@{-},
 <0mm,0mm>*{};<4.7mm,3.4mm>*{^2}**@{},
     <0mm,0mm>*{};<0.1mm,6.3mm>*{^1}**@{},
   <0mm,0mm>*{};<4.7mm,-3.1mm>*{^2}**@{},
   <0mm,-0.8mm>*{};<0.2mm,-5.4mm>*{^1}**@{},
    \end{xy} \ =\ 0,
    \ \ \ , \ \ \
\begin{xy}
 <0mm,2.47mm>*{};<0mm,0.12mm>*{}**@{-},
 <0.5mm,3.5mm>*{};<2.2mm,5.2mm>*{}**@{-},
 <-0.48mm,3.48mm>*{};<-2.2mm,5.2mm>*{}**@{-},
 <0mm,3mm>*{\bullet};<0mm,3mm>*{}**@{},
  <0mm,-0.8mm>*{\bullet};<0mm,-0.8mm>*{}**@{},
<-0.39mm,-1.2mm>*{};<-2.2mm,-3.5mm>*{}**@{-},
 <0.39mm,-1.2mm>*{};<2.2mm,-3.5mm>*{}**@{-},
     <0.5mm,3.5mm>*{};<2.8mm,5.7mm>*{^2}**@{},
     <-0.48mm,3.48mm>*{};<-2.8mm,5.7mm>*{^1}**@{},
   <0mm,-0.8mm>*{};<-2.7mm,-5.2mm>*{^1}**@{},
   <0mm,-0.8mm>*{};<2.7mm,-5.2mm>*{^2}**@{},
\end{xy}
\  - \
\begin{xy}
 <0mm,-0.8mm>*{};<3.2mm,1.2mm>*{}**@{-},
 <3.2mm,1.2mm>*{};<0mm,3mm>*{}**@{-},
 <0mm,0mm>*{};<0mm,-3.2mm>*{}**@{-},
 <0mm,3mm>*{};<0mm,6.0mm>*{}**@{-},
<0mm,3mm>*{\bullet};
 <0mm,-0.7mm>*{\bullet};
 <0mm,-0.8mm>*{};<-4.4mm,3mm>*{}**@{-},
 <0mm,3mm>*{};<-4.4mm,-0.8mm>*{}**@{-},
 <0mm,0mm>*{};<-4.7mm,3.4mm>*{^1}**@{},
     <0mm,0mm>*{};<-0.1mm,6.3mm>*{^2}**@{},
   <0mm,0mm>*{};<-4.7mm,-3.1mm>*{^1}**@{},
   <0mm,-0.8mm>*{};<-0.2mm,-5.4mm>*{^2}**@{},
    \end{xy} \ =\ 0.
$$
It is clear that $\Frob_+^+(m,n)=\bK[\Sy_m]\ot \bK[\Sy_n]$. Analogously one defines two other versions of $\Frob$,
$$
\Frob_+=\left\{\Frob_+(m,n)= \Id_m\ot \bK[\Sy_n]\right\} \ \ \ \mbox{and} \ \ \
\Frob^+=\left\{\Frob_+(m,n)= \bK[\Sy_m]\ot \Id_n\right\},
$$
with comultiplication (resp.\ multiplication) commutative but multiplication (resp.\ comultiplication)
noncommutative. Then applying  again
 Theorem~\ref{comor}(i) or Theorem~\ref{Convolutionproperadinfty}
to $\Qo$ being $(\Frob_+^+)^*$,  $(\Frob_+)^*$  or  $(\Frob^+)^*$
and $\caD$ being $\cF^c(s\bar{\Po})$ we conclude that the vector
spaces,
$$
\Hm\left((\Frob_+^+)^*, \Po \right)= \bigoplus_{m,n}\Po(m,n),
$$
$$
\Hm\left((\Frob_+)^*, \Po \right)=
\bigoplus_{m,n}\Po(m,n)^{\Sy_m},
$$
$$
\Hm\left((\Frob^+)^*, \Po \right)=
\bigoplus_{m,n}\Po(m,n)^{\Sy_n},
$$
admit canonically $L_\infty$-structures proving thereby Claim~(i)-(iii).

\sip

Finally, the natural commutative diagram of morphisms of properads,
$$
\xymatrix{
& \Frob_+  \ar[dr] & \\
  \Frob \ar[ur] \ar[dr]  & &
  \Frob^+_+  \\
& \Frob^+ \ar[ur] &
}
$$
proves claim (v).
\end{proo}

\bip

The prop(erad) $\Frob^+_+$ was designed so that it is generated as an $\Sy$-module by arbitrary $(m,n)$-corollas,
 and the comultiplication $\Delta_{(1,1)}$ in its dual, $(\Frob^+_+)^*$  splits such a corolla into {\em all}\, possible two
 vertex $(m,n)$-graphs. Hence the $L_\infty$-structure claimed in Theorem~\ref{main} is exactly the same as in
 Theorem~\ref{Prop->Homo Lie}.

 \sip

 The $L_\infty$ structures on the direct sum $\oplus P$ and its subspaces of invariants constructed in the proof
 of Theorem~\ref{main} are the most natural ones to consider as they involve {\em all}\, possible compositions in $\Po$.
 However, they are
 by no means  unique in the case of prop(erad)s
(as opposite to the case of operads). For example, the part of
prop(erad) compositions which correspond to so called
$\frac{1}{2}$-propic  graphs or dioperadic graphs, that is graphs
of genus $0$, (see pictures (\ref{12}) below) also combine into a
$L_\infty$-structure on $\oplus P$ (and its subspaces of
invariants) as the following argument shows.

 \bip

 For a prop(erad) $\Po=\{\Po(m,n)\}$ we denote by $\Po^\dag=\{\Po^\dag(m,n)\}$ the associated
``flow reversed" prop(erad) with $\Po^\dag(m,n):=\Po(n,m)$. Let
$\Ass$ be the operad of associative algebras and define the
properad, $\Ass^\dag \hspace{-0.5mm}\bullet \hspace{-0.5mm}\Ass=
\{\Ass^\dag \hspace{-0.5mm}\bullet \hspace{-0.5mm}\Ass(m,n)\}$ by
setting
$$
Ass^\dag  \hspace{-0.5mm}\bullet \hspace{-0.5mm}\Ass(m,n):=
\Ass^\dag(m)\ot \Ass(n)\simeq \bK[\bS_m]\ot \bK[\bS_n]
$$
and defining the compositions $\mu_{(1,1)}$ to be non-zero only on
decorated graphs of the form \Beq
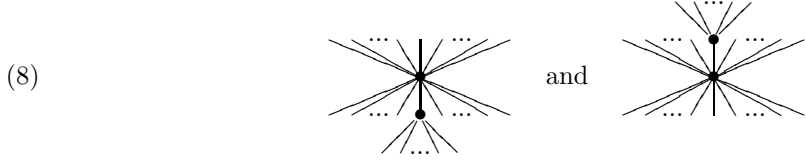
\label{12}
\begin{xy}
<0mm,0mm>*{\bullet}, <0mm,-5mm>*{}**@{-}, <0mm,5mm>*{}**@{-},
<-12mm,-5mm>*{}**@{-}, <-9mm,-5mm>*{}**@{-}, <-3mm,-5mm>*{}**@{-},
<-5.5mm,-5mm>*{...}**@{}, <12mm,-5mm>*{}**@{-},
<9mm,-5mm>*{}**@{-}, <3mm,-5mm>*{}**@{-}, <5.5mm,-5mm>*{...}**@{},
<-12mm,5mm>*{}**@{-}, <-9mm,5mm>*{}**@{-}, <-3mm,5mm>*{}**@{-},
<-5.5mm,5mm>*{...}**@{}, <12mm,5mm>*{}**@{-}, <9mm,5mm>*{}**@{-},
<3mm,5mm>*{}**@{-}, <5.5mm,5mm>*{...}**@{},
<0mm,-5mm>*{\bullet}; <-5mm,-10mm>*{}**@{-},
<-2.5mm,-10mm>*{}**@{-}, <0mm,-10mm>*{...}**@{},
<2.5mm,-10mm>*{}**@{-}, <5mm,-10mm>*{}**@{-},
\end{xy}
\ \ \ \ \mbox{and}\ \ \ \
\begin{xy}
<0mm,0mm>*{\bullet}, <0mm,5mm>*{}**@{-}, <0mm,-5mm>*{}**@{-},
<-12mm,-5mm>*{}**@{-}, <-9mm,-5mm>*{}**@{-}, <-3mm,-5mm>*{}**@{-},
<-5.5mm,-5mm>*{...}**@{}, <12mm,-5mm>*{}**@{-},
<9mm,-5mm>*{}**@{-}, <3mm,-5mm>*{}**@{-}, <5.5mm,-5mm>*{...}**@{},
<-12mm,5mm>*{}**@{-}, <-9mm,5mm>*{}**@{-}, <-3mm,5mm>*{}**@{-},
<-5.5mm,5mm>*{...}**@{}, <12mm,5mm>*{}**@{-}, <9mm,5mm>*{}**@{-},
<3mm,5mm>*{}**@{-}, <5.5mm,5mm>*{...}**@{},
<0mm,5mm>*{\bullet}; <-5mm,10mm>*{}**@{-}, <-2.5mm,10mm>*{}**@{-},
<0mm,10mm>*{...}**@{}, <2.5mm,10mm>*{}**@{-}, <5mm,10mm>*{}**@{-},
\end{xy}
\Eeq on which it is equal to the operadic compositions in $\Ass$.
The properad $\Ass^\dag \hspace{-0.5mm}\bullet
\hspace{-0.5mm}\Ass$ corresponds to the $\frac{1}{2}$-prop
$\mathcal{U}^{properad}_{\frac{1}{2}\textrm{-prop}}(\mathcal{F}rob^+_+)$.

\sip

Let $\Com$ be the properad of commutative algebras, and define the
properads $\Com^\dag \hspace{-0.5mm}\bullet \hspace{-0.5mm} \Ass$,
$\Ass^\dag \hspace{-0.5mm}\bullet \hspace{-0.5mm} \Com$, and
$\Com^\dag \hspace{-0.5mm}\bullet \hspace{-0.5mm} \Com$ by analogy
to  $\Ass^\dag \hspace{-0.5mm}\bullet \hspace{-0.5mm} \Ass$.
Similarly, they correspond to the $\frac{1}{2}$-props
$\mathcal{U}^{properad}_{\frac{1}{2}\textrm{-prop}}(\mathcal{F}rob^+)$,
$\mathcal{U}^{properad}_{\frac{1}{2}\textrm{-prop}}(\mathcal{F}rob_+)$
and
$\mathcal{U}^{properad}_{\frac{1}{2}\textrm{-prop}}(\mathcal{F}rob)$.
Applying Theorem~\ref{comor}(i) to $\Qo$ being coproperads
$(\Ass^\dag \hspace{-0.5mm}\bullet \hspace{-0.5mm} \Ass)^*$,
$(\Ass^\dag \hspace{-0.5mm}\bullet \hspace{-0.5mm} \Com)^*$,
$(\Com^\dag \hspace{-0.5mm}\bullet \hspace{-0.5mm} \Ass)^*$, or
$(\Com^\dag \hspace{-0.5mm}\bullet \hspace{-0.5mm} \Com)^*$, we
conclude that the vector spaces,
$$
\bigoplus_{m,n}P(m,n), \
\bigoplus_{m,n}P(m,n)^{\Sy_m}, \
\bigoplus_{m,n}P(m,n)^{\Sy_n}, \ \bigoplus_{m,n}P(m,n)^{\Sy_m\times \Sy_n},
$$
admit canonically $L_\infty$-structures encoding
$\frac{1}{2}$-prop compositions of the form (\ref{12}). The
natural morphism of operads,
$$
\Ass \lon \Com,
$$
implies that these $L_\infty$-structures are related to each other
via the same commutative diagram of $L_\infty$ morphisms as in
Theorem~\ref{main}(v). In the case of operads the constructed
$L_\infty$-structures are exactly the same as in
Theorem~\ref{main} but for prop(erad)s they are different.


\section{Deformation theory of morphisms of
prop(erad)s}\label{deformation theory}

In this section, we define the deformation theory of morphisms of prop(erad)s. We follow the conceptual method proposed by Quillen in \cite{Quillen67, Quillen70}.

\subsection{Basic Definition} Let $(\Po, d_\Po) \xrightarrow{\varphi} (\Qo, d_\Qo)$
be a morphism of dg prop(erad)s. We would like to define a chain
complex with which we could study the deformation theory of this
map. Following Quillen \cite{Quillen70}, the conceptual method is
to take the total right derived functor of the functor $\Der$ of
derivations from the category of prop(erad)s above $\Qo$ (see also
\cite{Markl96,VanderLaan02}). That is, we consider a cofibrant
replacement $(\mathcal{R}, \p)$ of $\Po$ is the category of dg
prop(erad)s
$$
\xymatrix{\mathcal{R} \ar[r]^{\varepsilon} \ar[dr]_{\gamma}&
\Po \ar[d]^{\varphi} \\ & \Qo.}
$$

Recall that $\Qo$ has an infinitesimal $\Po$-bimodule
(respectively infinitesimal $\mathcal{R}$-bimodule) structure
given by $\varphi$ (respectively $\gamma$).

\begin{lem}
Let $(\R,\, \p)$ be a resolution of $\Po$ and let $f$ be a
homogenous derivation of degree $n$ in $\Der_n \left( \R ,\,
\Qo\right)$, the derivative $D(f)=d_\Qo \circ f - (-1)^{|f|}
f\circ \p$ is a derivation of degree $n-1$ of $\Der_{n-1} \left(
\R,\, \Qo \right)$.
\end{lem}

\begin{proo}
The degree of $D(f)$ is $n-1$. It remains to show that it is a
derivation.  For every pair $r_1$ and $r_2$ of homogenous elements
of $\R$, $D(f)\left( \mu_\R (r_1\bbc r_2)\right)$ is equal to
\begin{eqnarray*}
D(f)\left( \mu_\R (r_1\bbc r_2)\right) &=& (d_\Qo \circ f -
(-1)^{n} f\circ \p)   \left( \mu_\R (r_1\bbc r_2)\right)\\
&=& d_\Qo \left( \mu_\Qo \big(f(r_1) \bbc
\gamma(r_2)+(-1)^{n|r_1|}\gamma(r_1)\bbc f(r_2) \big) \right)\\
&& -(-1)^n f\left(\mu_\R\big(\p(r_1)\bbc r_2 +(-1)^{|r_1|}r_1\bbc
\p(r_2)\big)\right) \\
&=& \mu_\Qo\Big( (d_\Qo\circ f)(r_1) \bbc \gamma(r_2) +(-1)^{n+|r_1|} f(r_1)\bbc (d_\Qo \circ \gamma)(r_2)    \\
&& +(-1)^{n|r_1|}(d_\Qo\circ \gamma)(r_1) \bbc f(r_2)
+(-1)^{|r_1|(n-1)} \gamma(r_1)\bbc
(d_\Qo \circ f)(r_2) \Big) \\
 && -(-1)^n \mu_{\Qo}\Big( (f\circ \p)
(r_1)\bbc \gamma (r_2) + (-1)^{n(|r_1|-1)}(\gamma \circ
\p)(r_1))\bbc
f(r_2)  \\
&&  +(-1)^{|r_1|}f(r_1)\bbc (\gamma \circ \p)(r_2))+
(-1)^{(n-1)|r_1|}\gamma (r_1) \bbc (f\circ \p)(r_2)\Big).
\end{eqnarray*}
Since $\gamma$ is morphism of dg prop(erad)s, it commutes with the
differentials, that is $\gamma \circ \p = d_\Qo \circ \gamma$.
This gives
\begin{eqnarray*}
D(f)\left( \mu_\R (r_1\bbc r_2)\right) &=&
 \mu_\Qo\Big( (d_\Qo\circ f)(r_1) \bbc \gamma(r_2)  -(-1)^n  (f\circ \p)
(r_1)\bbc \gamma (r_2)    \\
&& +(-1)^{|r_1|(n-1)}\big( \gamma(r_1)\bbc (d_\Qo \circ f)(r_2)
-(-1)^n\gamma (r_1) \bbc (f\circ \p)(r_2)\big)\Big)\\
&=& \mu_\Qo\Big(D(f)(r_1)\bbc \gamma(r_2)  +(-1)^{|r_1|(n-1)}
\gamma(r_1)\bbc D(f)(r_2)\Big).
\end{eqnarray*}
\end{proo}

In other words, the space of derivations $\Der(\mathcal{R}, \Qo)$
 is a sub-dg-module of the space of morphisms $\Hom^{\Sy}(\mathcal{R},
 \Qo)$. We define the deformation complex of the morphism $\varphi$ by
$C_\bullet(\varphi):=\big( \Der_\bullet(\mathcal{R}, \Qo),
D\big)$. By Theorem~\ref{Exists QF regular} and Theorem~\ref{Exists QF}, there always exists a
quasi-free cofibrant resolutions. For instance, we can consider the
bar-cobar resolution by Theorem~\ref{BarCobarResolution}. This
will produce an explicit but huge complex which is difficult to
compute. Instead of that, we will work with  the chain complex
obtained from a minimal model of $\Po$ when it exists. Its size is
much smaller but its differential can be not so easy to make
explicit. In this sequel, our main example be the deformation
theory of representations of $\Po$ of the form $\Qo=\End_X$, that is
$\Po$-gebras.

\subsection{Deformation theory of representations of
prop(erad)s.}\label{Deformation theory}

Let $(\Po, d_\Po)$ be a dg prop(erad), let $(X, d_X)$ an arbitrary dg
$\Po$-gebra and let
$(\Po_\infty:= \Omega(\Co), \p)$ be a cofibrant quasi-free resolution of $\Po$ and
$$\xymatrix{\Omega(\Co) \ar[r]^{\varepsilon}
\ar[dr]_{\gamma}& \Po \ar[d]^{\varphi} \\ & \End_X.}$$

\begin{dei}[Deformation complex]
We define the \emph{deformation complex of the $\Po$-gebra
structure of $X$} by $C_\bullet(\Po, X):=\big(
\Der_\bullet(\Omega(\Co), \End_X), D \big)$.
\end{dei}

\begin{thm}\label{boundary=Q}
The deformation complex $\big( \Der_\bullet(\Omega(\Co), \Qo), D
\big)$ is isomorphic to $\Hom_\bullet^\Sy(\bar{\Co}, \Qo)$ with
$D=Q^\gamma$ for $\gamma=\varphi\circ \varepsilon_{|\bar{\Co}}$.
\end{thm}

\begin{proo}
Lemma~\ref{uniqueDer} proves the identification between the two
spaces. Since $\gamma$ is a morphism of dg prop(erad)s from a
quasi-free prop(erad), it is a solution of the Maurer-Cartan
equation $Q(\gamma)=0$ in the convolution
$\textrm{L}_\infty$-algebra $\Hom^\Sy_\bullet(\bar{\Co}, \Qo)$ by
Theorem~\ref{mor}. Let $f$ be an element of $\Hom^\Sy_n(\bar{\Co},
\Qo)$. Following Lemma~\ref{uniqueDer}, we denote by $\p_f$ the
unique derivation of $\Der_n(\Omega(\Co), \Qo)$ induced by $f$. We
have to show that $D(\p_f)_{s^{-1}\bar{\Co}}=Q^\gamma(f)$. For an
element $s^{-1}c \in s^{-1}\Co$, we use the Sweedler type notation
for $\p(s^{-1}c)=\sum_\G \G(s^{-1}c_1, \ldots, s^{-1}c_n)$. By
Lemma~\ref{uniqueDer}, we have
\begin{eqnarray*}
&& \p_f\big(\G(s^{-1}c_1, \ldots, s^{-1}c_n)\big)=\\
&&\sum_{i=1}^n (-1)^{n(|c_1|+\cdots +|c_{i-1}|+i-1)} \mu_\Qo\big(
 \G(\gamma(s^{-1}c_1), \ldots,\gamma(s^{-1}c_{i-1}),f(s^{-1}c_i),\gamma(s^{-1}c_{i+1}),\ldots,
 \gamma(s^{-1}c_n)\big).
 \end{eqnarray*}
Therefore, $D(\p_f)_{s^{-1}\bar{\Co}}$ is
 equal to
\begin{eqnarray*}
&& D(\p_f)(s^{-1}c)=(d_\Qo \circ \p_f-(-1)^n\p_f\circ
\p)(s^{-1}c) =d_\Qo (f(s^{-1}c)) -  \\
&&(-1)^n \sum_\G \sum_{i=1}^n (-1)^{n(|c_1|+\cdots
+|c_{i-1}|+i-1)} \mu_\Qo\big(
 \G(\gamma(s^{-1}c_1), \ldots,\gamma(s^{-1}c_{i-1}),f(s^{-1}c_i),\gamma(s^{-1}c_{i+1}),\ldots,
 \gamma(s^{-1}c_n)\big)\\
 &&= Q^\gamma(f).
\end{eqnarray*}
\end{proo}

In order words, the deformation complex is equal to the
convolution $L_\infty$-algebra $\Hom_\bullet^\Sy(\bar{\Co}, \Qo)$
twisted by the Maurer-Cartan element $\gamma$.


\begin{Rq}  It is natural to consider the augmentation of this chain
complex by $\Hom^{\Sy}(\II, \Qo)$, that is $\Hom_\bullet^\Sy(\Co,
\Qo)$ .
\end{Rq}

In summary, by Theorem~\ref{mor} the vector space
$\Hom^\Sy_\bullet(\Co,\Qo)$ has a canonical filtered
$L_\infty$-structure, $Q$ whose Maurer-Cartan elements are
morphisms of dg prop(erad)s, $ \cP_\infty \to \Qo$, that is
representations of $\cP_\infty$ in $\Qo$. Then let $\ga$ be one of
these morphisms, and let $Q^\ga$ be the associated twisting of the
canonical $L_\infty$-algebra by $\ga$ (see \S \ref{filtered}).
This defines the deformation complex of $\gamma$.

\begin{dei}[Deformation Complex]\label{deftheory}
The deformation complex of a morphism of prop(erad)s $\gamma \,
:\, \Po_\infty \to \Qo$ is the twisted $L_\infty$-algebra
$(\Hom^\Sy_\bullet(\Co,\Qo), Q^\gamma$).
\end{dei}

This definition extends to the case of prop(erad)s the deformation complex of
algebras over operads introduced in \cite{KontsevichSoibelman00,VanderLaan02}.

\sip

With the results on the model category structure on prop(erad)s (see Appendix),
we can now prove the independence of this construction in the
homotopy  category of homotopy prop(erad)s and homotopy Lie algebras.

\begin{thm}\label{independence of homotopy convolution}
Let
$\Omega(\Co_1)$ and $\Omega(\Co_2)$ be two quasi-free cofibrant resolutions
of a dg prop(erad) $\Po$. For any
dg prop(erad) $\Qo$, the homotopy convolution prop(erad)s
$\Hom(\Co_1, \Qo)$ and $\Hom(\Co_2, \Qo)$ are linked by two
quasi-isomorphisms of homotopy prop(erad)s
$$\Hom(\Co_1, \Qo) \leftrightarrows \Hom(\Co_2, \Qo),$$
and the natural maps $$(\Hom^\Sy_\bullet(\Co_1,\Qo), Q^{\gamma_1}) \leftrightarrows   (\Hom^\Sy_\bullet(\Co_2,\Qo), Q^{\gamma_2})$$
is a quasi-isomorphism of homotopy Lie algebras.
\end{thm}

\begin{proo}
 We
apply the left lifting property in the model category of
dg prop(erad)s to the following diagram~:
$$
\xymatrix{0 \ar[r] \ar@{>->}[d]  & \Omega(\Co_1) \ar@{->>}[d]^{\sim}  \ar@{..>}@<0.5ex>[dl] \\
\Omega(\Co_2) \ar@{..>}@<0.5ex>[ur] \ar@{->>}[r]^{\sim}&  \Po, }
$$
to get the two dotted quasi-isomorphisms of prop(erad)s. By Proposition~\ref{preserves EQ}, they induce quasi-isomorphisms of the level of the homotopy coprop(erad)s $\Co_1  \leftrightarrows \Co_2$. We
conclude by Theorem~\ref{functoriality convolution properad} and by Corollary~\ref{functoriality convolution L-infinite}.
\end{proo}

 The homology
groups, $$ H^\ga_\bullet(\Qo):=H_\bullet(\Hom^\Sy(\Co, \Qo), Q_1^\ga),
$$ are independent of the choice of a cofibrant quasi-free resolution of $\cP$
and are called {\em homology groups of the
$\cP_\infty$-representation  $\Qo$}. In the case where
$\Qo=\End_X$, they are called the homology groups of the
$\Po_\infty$-gebra $(X, \gamma)$.

\begin{pro} The Maurer-Cartan elements, $\Gamma$, of $Q^\ga$ are in
one-to-one correspondence with those $\cP_\infty$-structures,
 on $X$,
$$
\rho: (\cP_\infty, \p) \lon (\End_X, d),
$$
whose restrictions to the generating space, $s^{-1}\Co$, of
$\cP_\infty$ are equal precisely  to
  the sum $\ga+\Gamma$.
\end{pro}

This proposition justifies the name `deformation complex' because
the $L_\infty$-algebra $Q^\ga$ {\em controls the deformations of
$\ga$ in the class of homotopy $\cP$-structures.} When applied to
$\Qo=\End_X$ and $\gamma \, : \, \Po_\infty \to \End_X$, this
defines the deformation complex of the $\Po_\infty$-gebra
structure $\gamma$ on $X$. (Some author call this the ``cohomology
of $X$ with coefficients into itself'' but we are reluctant to
make this choice and prefer to view it as a deformation complex).
This definition applies to any homotopy algebra over an operad
(associative algebras, Lie algebras, commutative algebras, PreLie
algebras, Poisson or Gerstenhaber algebras, etc ...) as well as to
any homotopy (bial)gebra over a properad (Lie bialgebras,
associative bialgebras, etc ...) in order to give, for the first
time, a cohomology theory for homotopy
$\Po$-(bial)gebras.

\subsection{Koszul case and cohomology operations}\label{Koszul case}

In Theorem~\ref{Koszul=quadratic model}, we have seen that a
properad $\Po$ is Koszul if and only if it admits a quadratic
model $\Omega(\Po^{\ac})\xrightarrow{\sim} \Po$, where $\Po^{\ac}$
the Koszul dual (strict) coproperad. In this case, by
Theorem~\ref{boundary=Q}, the deformation complex of a
$\Po_\infty$-gebra $\Hom^\Sy(\Po^{\ac}, \End_X)$ is dg Lie algebra
where the boundary map is equal to the twisted differential
$D(f)=d(f)+[\gamma,
f]$.\\

The first definition of this kind of preLie operation appeared in
the seminal paper of M. Gerstenhaber \cite{Gerstenhaber63} in the
case of the cohomology of associative algebras. In the case
treated by M. Gerstenhaber, the cooperad $\Co$ is the Koszul dual
cooperad $\A^{\ac}$ of the operad $\A$ coding associative algebras
and the operad $\Po$ is the endomorphism operad $\End_A$. The
induced Lie bracket is the \emph{intrinsic Lie bracket} of
Stasheff \cite{Stasheff93}. It is equal to the Lie bracket of
Gerstenhaber \cite{Gerstenhaber63} on Hochschild cochain complex
of associative algebras, the Lie bracket of Nijenhuis-Richardson
\cite{NijenhuisRichardson67} on Chevalley-Eilenberg cochain
complex of Lie algebras and the Lie bracket of Stasheff on
Harrison cochain complex of commutative algebras. It is proven by
Balavoine in \cite{Balavoine} that  the deformation complex of
algebras over any Koszul operad admits a Lie structure. This
statement was made more precise by Markl, Shnider and Stasheff in
Section~$3.9$ Part~II of \cite{MarklShniderStasheff02} where they
proves that this Lie bracket comes from a Prelie product. This
result on the level of operads was proved using the space of
coderivations of the cofree $\Po^{\ac}$-coalgebra, which is shown
to be a PreLie algebra. Such a method is impossible to generalize
to prop(erad)s simply because there exists no notion of (co)free
gebra.\\

 As explained here, one has to work with convolution
prop(erad) to prove a similar result. Actually, this method gives
a stronger statement.

\begin{thm}\label{Cohomology operations}
Let $\Po$ be a Koszul properad and let $\varphi\, :\, \Po \to \Qo$ be a morphism of properads, the deformation complex of $\varphi$ is a LR-algebra.

In the non-symmetric case, when $\Po$ and $\Qo$ are non-symmetric properads, the deformation complex is a non-symmetric properad.
\end{thm}

\begin{proo}
It is direct consequence of the definition of the deformation complex and Theorem~\ref{Convolution LR and Lie algebra}. In the non-symmetric case, the deformation complex is directly a non-symmetric convolution properad, since it is not restricted to invariant elements.
\end{proo}

This result provides higher braces or LR-operations (see
Section~\ref{operations on convolution properad}). Recall that
non-symmetric braces play a fundamental role in the proof of
Deligne's conjecture for associative algebras (see
\cite{Tamarkin98, Voronov00, KontsevichSoibelman00,
McClureSmith02, BergerFresse04}) and in the extension of it to
other kind of algebras (see \cite{Vallette06Preprint}
Section~$5.5$). From this rich structure, we derive a
Lie-admissible bracket and then a Lie bracket which can be used to
study the deformations of $X$. We expect the LR-operations to be
used in the future for a better understanding of deformation
theory (in the context of a Deligne conjecture for associative
bialgebras and Gerstenhaber-Schack bicomplex, for
instance).\\

Notice that this Lie bracket was found by hand in one example of
gebras over a properad before this general theory. The properad of
Lie bialgebras is Koszul. Therefore, on the deformation
(bi)complex of Lie bialgebras, there is a Lie bracket. The
construction of this Lie bracket was given by Kosmann-Schwarzbach
in \cite{KosmannSchwarzbach91}. (See also Ciccoli-Guerra
\cite{CiccoliGuerra03} for the interpretation of this bicomplex in
terms of deformations.)

\subsection{Definition \`a la Quillen}\label{Def a la Quillen}
In the previous sections, we defined the deformation complex of
representations of a prop(erad) $\Po$ that admits a quasi-free
model and proved the independence of this definition in the
categories of homotopy prop(erad)s (and homotopy Lie algebras). In
this section, we generalize the definition of the deformation
theory of a morphism of commutative rings due to Quillen
\cite{Quillen70} to the case of prop(erad)s. (See L. Illusie
\cite{Illusie71} for a generalization in the context of topoi and
schemes). Hence, it defines a (relative) deformation complex for
representations of any prop(erad). It also gives rise to the
\emph{cotangent complex} associated to any morphism of
prop(erad)s.

 Since a commutative algebra is
an associative algebra, an associative algebra an operad and an
operad a prop(erad), this generalization of Quillen theory can be
seen as a way to extend results of (commutative) algebraic
geometry to non-commutative non-linear geometry. It is non-linear
because the monoidal product $\bc$ defining prop(erad)s is neither
linear on the left nor on the right, contrary to the tensor product $\otimes$ of vector spaces.\\

Let $\I$ be a `ground' prop(erad) (to recover the previous
section, consider $\I=I$, the unit of the monoidal category of
$\Sy$-bimodules). We look at prop(erad)s $\Po$ under $\I$, $\I \to
\Po$. And for such  a prop(erad) $\Po$, we consider the category
of prop(erad)s over $\Po$, that is

$$\xymatrix@R=25pt@C=25pt{ & \X \ar[d]^{f} \\ \I \ar[r] \ar[ru]^{u} & \Po.} $$
We denote this category by $\textrm{Prop(erad)}/\Po$. Let $M$ be
an infinitesimal bimodule over $\Po$ (see Section~\ref{Infinitesimal bimodule}). The infinitesimal
$\Po$-bimodule $M$ is also an infinitesimal bimodule over any
prop(erad) $\X$ over $\Po$, by pulling back along $\X \to \Po$.
Hence, we can consider the space of $\I$-derivations from $\X$ to
$M$, that is derivations from $\X$ to $M$ which vanish on $\I$. We
denoted this space by $\Der_\I(\X, M)$.

We aim now to represent this bifunctor on the left and on the
right. To represent it on the left, we introduce the
\emph{square-zero (or infinitesimal) extension of $\Po$ by $M$}~:
$\Po \ltimes M := \Po \oplus M$ with the following structure of
prop(erad) over $\Po$. The monoidal product $(\Po \oplus
M)\boxtimes(\Po \oplus M)$ is equal to
$$\Po \boxtimes \Po \, \bigoplus\,  \Po \boxtimes (\Po \oplus
\underbrace{M}_1)\,  \bigoplus \, (\Po \oplus \underbrace{M}_1)
\boxtimes \Po\,  \bigoplus\,   \mathcal{M}, $$ where $\mathcal{M}$
is the sub-$\Sy$-bimodule of $(\Po \oplus M)\boxtimes(\Po \oplus
M)$ composed by at least two elements from $M$. On the first
component $\Po \boxtimes \Po$, the product of $\Po \ltimes M$ is
defined by the product of $\Po$. On the second component, it is
defined by the left action of $\Po$ on $M$. On the third one, it
is defined by the right action of $\Po$ on $M$. Finally, the
product on $\mathcal{M}$ is null.

\begin{lem}\label{infinitesimal extension}
For any prop(erad) $\Po$ and any infinitesimal $\Po$-bimodule $M$,
the infinitesimal extension $\Po\ltimes M$ is a prop(erad).
\end{lem}

\begin{proo}
The definition of infinitesimal $\Po$-bimodule directly implies
the associativity of the prop(erad)ic composition of $\Po \ltimes
M$.
\end{proo}

The purpose of this definition is in the following result, which
states that infinitesimal $\Po$-bimodules are abelian group
objects in the category of prop(erad)s over $\Po$.

\begin{pro}\label{Der-inf. Extension}
There is a natural bijection
$$\Hom_{\textrm{Prop(erad)}/\Po}(\X, \Po \ltimes M)\cong \Der_\I(\X, M),$$
where $\Der_\I(\X, M)$ is the space of $\I$- derivations from $\X$
to $M$.
\end{pro}

\begin{proo}
Let us denote by $f\, : \, \X \to \Po$. Any morphism $\X \to \Po
\ltimes M=\Po\oplus M$  the category of prop(erad)s over $\Po$ is
the sum of $f$ with its component on $M$, which we denote by $D$.
Finally, $f\oplus D\, : \, \X \to \Po\ltimes M$ is a morphism of
prop(erad)s if and only if $D$ is a derivation $\X \to M$.
\end{proo}

To represent the space of derivations on the right, we introduce
the \emph{module of K\"ahler differentials of a prop(erad)}. It is
a quotient of the free infinitesimal $\X$-bimodule over $\I$ on
$\X$ by suitable relations. We recall from Section~$2.5$ of
\cite{Vallette03} that the relative composition product is defined
by the following coequalizer
$$\xymatrix@C=30pt{M \bc \Po \bc N  \ar@<-0.5ex>[r]_(0.55){\rho \bc N}
\ar@<+0.5ex>[r]^(0.55){M \bc \lambda} & M \bc N  \ar@{->>}[r] & M
{\bc}_{\Po} N},$$ where $\lambda$ is the left action of $\Po$ on
$N$, $\Po \bc N \to N$ and $\rho$ the right action of $\Po$ on
$M$, $M \bc \Po  \to M$.\\

Let $f\, :\, \Po \to \Qo$ be a morphism of prop(erad)s. There is a
natural functor from the category of infinitesimal $\Qo$-bimodules
to the category of infinitesimal $\Po$-bimodules by pulling back
along $f$. We denote it by $f^*\, : \textrm{Inf.}\,
\Qo\textrm{-biMod} \to \textrm{Inf.}\, \Po\textrm{-biMod}$.

\begin{pro}\label{free inf P bimod}
The functor $f^*\, : \textrm{Inf.}\, \Qo\textrm{-biMod} \to
\textrm{Inf.}\, \Po\textrm{-biMod}$  admits a left adjoint
$$f_!\ :\  \textrm{Inf.}\, \Po\textrm{-biMod}
\leftrightharpoons \textrm{Inf.}\, \Qo\textrm{-biMod}\ : \ f^*,$$
which is explicitly given by $f_!(M)=\Qo\bc_{\Po}
\underbrace{M}_{1}\bc_\Po \Qo$, for any infinitesimal
$\Po$-bimodule $M$. The $\Sy$-bimodule $\Qo\bc_{\Po}
\underbrace{M}_{1}\bc_\Po \Qo$ is the coequalizer
$$\xymatrix@C=70pt{(\Qo \bc \Po) \bc \underbrace{M}_1 \bc (\Po \bc \Qo)
\ar@<2ex>[r]^(0.55){\rho_\Qo \bc M \bc \lambda_\Qo}
\ar[r]_(0.6){\Qo \bc(\lambda \circ (\Po \bc \rho)) \bc \Qo} & \Qo
\bc \underbrace{M}_1 \bc \Qo \ar@<1ex>@{->>}[r] & \Qo \bc_\Po
\underbrace{M}_1 \bc_\Po \Qo,  }$$ where the notation $\Qo \bc
\underbrace{M}_1 \bc \Qo$ stands for 3-levels graphs with only one
element of $M$ labelling a vertex on the second level and such
that every element of $Q$ on the first and third level have a
common internal edge with this element of $M$. (The action of
$\Po$ on $\Qo$ is given by the morphism $f$.)
\end{pro}

\begin{proo}
We have the natural bijection
$$
\Hom_{\textrm{Inf.} \, \Qo \textrm{-biMod}}(\Qo\bc_{\Po}
\underbrace{M}_{1}\bc_\Po \Qo, N)\cong \Hom_{\textrm{Inf.}\,
\Po\textrm{-biMod}}(M, f^*(N)).$$ Let $\Phi \, : \,\Qo\bc_{\Po}
\underbrace{M}_{1}\bc_\Po \Qo\to N$ be a morphism a infinitesimal
$\Qo$-bimodules. It is characterized by the image of the
projection of the element of $I \bc \underbrace{M}_1 \bc I$ in
$\Qo\bc_{\Po} \underbrace{M}_{1}\bc_\Po \Qo$. Let us call $\varphi
\, :\, M \to N$ this map. It is then easy to see that $\varphi$ is
a morphism of infinitesimal $\Po$-bimodules.
\end{proo}

\begin{Ex}
If we apply the preceding proposition to the unit $u \, : \, I \to
\Po$ of a prop(erad) $\Po$, the functor $u^* \, :\, {\textrm{Inf.}
\, \Po \textrm{-biMod}} \to \Sy\textrm{-biMod}$ is the classical
forgetful functor. Hence $u_!(M)=\Po \bc \underbrace{M}_{1} \bc
\Po$ is the free infinitesimal $\Po$-bimodule associated to any
$\Sy$-bimodule $M$.
\end{Ex}

\begin{dei}[Module of K\"ahler differentials]
Let us denote by $u \, : \, \I \to \X$. The module of K\"ahler
differentials of a prop(erad) $\X$ over $\I$ is the quotient of
the free infinitesimal $\X$-bimodule over the infinitesimal
$\I$-bimodule $u^*(\X)$, that is
$$u_!(u^*(\X))=\X \bc_\I \underbrace{u^*(\X)}_1
\bc_\I \X,$$
 by the relations
$$\pi\big(I\bc \mu_\X(x_1{\bc}_{(1,1)} x_2) \bc I - I \bc x_1{\bc}_{(1,1)}
x_2 - (-1)^{|x_1|} x_1{\bc}_{(1,1)} x_2 {\bc} I\big),$$ where
$\pi$ is the canonical projection of $\X \bc \underbrace{\X}_1 \bc
\X$ on the coequalizer $\X \bc_\I \underbrace{u^*(\X)}_1 \bc_\I
\X$.
 We denote it by $\Omega_{\X/\I}$.
\end{dei}

We define the universal derivation $D \, : \, \X
\to\Omega_{\X/\I}$ by $D(x)$ equal to the class of $I \bc x \bc I$
in $\Omega_{\X/\I}$. Like in the case of commutative algebras (see
J.-L. Loday \cite{LodayBOOK} Section~$1.3$) or associative
algebras (see \cite{Connes85, Karoubi87, LodayBOOK}, the module of
K\"ahler differentials represents the derivations.

\begin{pro}\label{Der-Kahler}
There is a natural bijection
$$\Der_\I(\X, M)\cong \Hom_{\textrm{Inf.} \,
\X\textrm{-biMod}}(\Omega_{\X/\I}, M).$$
\end{pro}

\begin{proo}
Let $d$ be a derivation in $\Der_\I(\X, M)$. There is a unique
morphism of infinitesimal $\X$-bimodules $\theta \, :\,
\Omega_{\X/\I} \to M$ such that the following diagram commutes
$$\xymatrix@R=30pt@C=30pt{\X \ar[r]^{D} \ar[rd]_{d}   &
\Omega_{\X/\I} \ar@{-->}[d]^{\theta} \\ & M.} $$

The image of the class of $I\bc x\bc I$ in $\Omega_{\X/\I}$ under
 $\theta$ is defined by $d(x)$. It extends freely to the
 infinitesimal $\X$-bimodule $\X \bc_\I \underbrace{u^*(\X)}_1 \bc_\I
\X$ and then passes to the quotient thanks to the Leibniz relation
verified by $d$.
\end{proo}

The module of K\"ahler differentials of an associative algebra is
the non-commutative analog of classical differential forms (see A.
Connes \cite{Connes85}). Since operads and prop(erad)s can also be
used to encode geometry (see \cite{Merkulov05, Merkulov06}), the
module of K\"ahler differentials for prop(erad)s seems a promising
tool to study non-linear properties in non-commutative geometry.

\begin{thm}\label{AdjunctionHomDerHom}
For any infinitesimal $\Po$-bimodule $M$,  the following
adjunction holds
$$\Hom_{\textrm{Prop(erad)}/\Po}(\X, \Po \ltimes M)\cong \Der_\I(\X, M) \cong
\Hom_{\textrm{Inf.} \, \Po\textrm{-biMod}}(\Po \bc_\X
\underbrace{\Omega_{\X/\I}}_1 \bc_\X \Po, M)  .$$
\end{thm}

\begin{proo}
It is a direct corollary of Proposition~\ref{Der-inf. Extension}
and Proposition~\ref{Der-Kahler}. The last natural bijection
$$  \Hom_{\textrm{Inf.}\,
\X\textrm{-biMod}}(\Omega_{\X/\I}, f^*(M)) \cong
\Hom_{\textrm{Inf.} \, \Po \textrm{-biMod}}(f_!(\Omega_{\X/\I}),
N)$$ is provided by Proposition~\ref{free inf P bimod} applied to
the morphism $f \, :\, \X \to \Po$.
\end{proo}

In other words, the  following functors form a pair of adjoint
functors
$$\Po \bc_- \underbrace{\Omega_{-/\I}}_1 \bc_- \Po \ : \
\textrm{Prop(erad)}/\Po   \ \rightleftharpoons \ \textrm{Inf.} \,
\Po\textrm{-biMod} \ :   \Po\ltimes - .$$

The model category structure on prop(erad)s induces a model
category structure on $\textrm{Prop(erad)}/\Po$.

\begin{lem}
The category of infinitesimal $\Po$-bimodules is endowed with a
cofibrantly generated model category structure.
\end{lem}

\begin{proo}
We use the  same arguments as in Appendix~\ref{model category},
that is the Transfer Theorem~\ref{TransferTheorem} along the free
infinitesimal $\Po$-bimodule functor $\eta_!\, : \,
\Sy\textrm{-biMod} \to \textrm{Inf.} \Po\textrm{-biMod}$. The
forgetful functor $\eta^*$ creates limits and colimits which
proves $(1)$ and $(2)$. A relative $\eta_!(J)$-cell complex has
the form $A_0 \to A_0 \oplus \Po \boxtimes \underbrace{(
\oplus_{i\ge 0}D^{k_i}_{m_i, n_i}) }_1 \boxtimes \Po$, which is a
quasi-isomorphism of dg $\Sy$-bimodules since the right hand term
$\Po \boxtimes \underbrace{( \oplus_{i\ge 0}D^{k_i}_{m_i, n_i})
}_1 \boxtimes \Po$ is acyclic.
\end{proo}

\begin{pro}
The pair of adjoint functors $$\Po \bc_-
\underbrace{\Omega_{-/\I}}_1 \bc_- \Po \ : \
\textrm{Prop(erad)}/\Po   \ \rightleftharpoons \ \textrm{Inf.} \,
\Po\textrm{-biMod} \ :   \Po\ltimes - $$ form a Quillen
adjunction.
\end{pro}

\begin{proo}
By Lemma~$1.3.4$ of \cite{Hovey99}, it is enough to prove that the
right adjoint $\Po \ltimes -$ preserves fibrations and acyclic
fibrations. Let $f\, : \, M \epi M'$ be a fibration (resp. acyclic
fibration) between two infinitesimal $\Po$-bimodules, that is $f$
is degreewise surjective (resp. and a quasi-isomorphism). Since
$\Po \ltimes (f)$ is the morphism of properads on  $Id_{\Po}
\oplus f \, :\, \Po \oplus M \to \Po \oplus M'$, it is degreewise
surjective (resp. and a quasi-isomorphism), which concludes the
proof.
\end{proo}

Thereofre, we can derive them in the associated homotopy categories.\\

This proves that the homology of $\Der_\I(\mathcal{R}, M)$ is
independent of the choice of the cofibrant resolution of $\Po$
because it is well defined in the homotopy category of prop(erad)s
over $\Po$ and in the homotopy category infinitesimal
$\Po$-bimodules.
$$\Hom_{\mathsf{Ho}(\textrm{Prop(erad)}/\Po)}(\X, \Po \ltimes M)\cong \Der_\I(\X, M) \cong
\Hom_{\mathsf{Ho}(\textrm{Inf.} \, \Po\textrm{-biMod})}(\Po \bc_\X
\underbrace{\Omega_{\X/\I}}_1 \bc_\X \Po, M)  .$$

\begin{dei}[Cotangent complex]
The \emph{cotangent complex} of $\Po$ is the total left derived
functor of the right adjoint, that is
$$\mathbb{L}_{\Po/\I}:=\Po \bc_\mathcal{R}
\underbrace{\Omega_{\mathcal{R}/\I}}_1 \bc_\mathcal{R} \Po,$$ for
$\mathcal{R}$ a cofibrant resolution of $\Po$.
\end{dei}

Since on the homology of the cotangent complex, in the classical
case of commutative rings, one can read the properties of the
morphism $\I \to \Po$ (smooth, locally complete intersection,etc
...), we expect to be able to read such properties on the
generalized version defined here. In the same way, transitivity
and flat base change theorems should be proved for this cotangent
complex but it is not our aim here and will be studied in a future
work.

\begin{Rq}
This section is written in the category of dg-prop(erad)s since we
work in this paper over a field of characteristic $0$. Therefore,
to explicit the cotangent complex and the (co)homology of
prop(erad)s, we have to use cofibrant resolutions in the category
of dg prop(erad)s, for instance quasi-free resolutions (Koszul or
homotopy Koszul). One can extend this section and the Appendix
when the characteristic of the ground ring is not $0$. In this
case, one has to use simplicial resolutions like in M. Andr\'e
\cite{Andre74} and D. Quillen \cite{Quillen70}.
\end{Rq}

\section{Examples of deformation theories.}\label{examples}

In this section, we show that the conceptual deformation theory
defined here coincide to well known theories in the case of
associative algebras, Lie algebras, commutative algebras, Poisson
algebras. As a corollary, we get classical Lie brackets on these
cohomology theory as well as classical Lie brackets in
differential geometry. More surprisingly, we make deformation
theory explicit in the case of associative bialgebras and show
that it corresponds to Gerstenhaber-Schack type bicomplex.

\subsection{Associative algebras.} If $P$ is the properad, $\ass$, of associative algebras, it is generated
by a non-symmetric operad still denoted by $\ass$. This operad is
Koszul that is, its minimal resolution exists and is generated by
the (strict) cooperad $\ass^{\ac}$ with
$$
\ass^{\ac}(m,n)=\left\{\Ba{cc}
s^{n-2}\bK[\bS_1] \ot \bK[\bS_n] & \mbox{for}\ m=1, n\geq 2 \\
0 & \mbox{otherwise.} \Ea \right.
$$
We represent the generating element of $\ass(1,n)$ by a corolla $
\begin{xy}
<0mm,4mm>*{\bullet}; <0mm,4mm>*{};<0mm,8mm>*{}**@{-},
 <0mm,4mm>*{};<-3mm,0mm>*{}**@{-},
 <0mm,4mm>*{};<-5mm,0mm>*{}**@{-},
 <0mm,4mm>*{};<3mm,0mm>*{}**@{-},
 <0mm,4mm>*{};<5mm,0mm>*{}**@{-},
 <0mm,4mm>*{};<0mm,0mm>*{...}**@{},
<0mm,4mm>*{};<0mm,-1mm>*{\underbrace{\ \ \ \ \ \ \ \ \ \ }}**@{},
<0mm,4mm>*{};<0mm,-4mm>*{^n}**@{}, \,
\end{xy}
$. The partial coproduct of this cooperad is given by the formula

$$
\Delta_{(1,1)}\left( \begin{xy}
 <0mm,0mm>*{\bullet};<0mm,0mm>*{}**@{},
 <0mm,0mm>*{};<0mm,5mm>*{}**@{-},
 <0mm,0mm>*{};<-8mm,-5mm>*{}**@{-},
 <0mm,0mm>*{};<-4.5mm,-5mm>*{}**@{-},
 <0mm,0mm>*{};<-1mm,-5mm>*{\ldots}**@{},
 <0mm,0mm>*{};<4.5mm,-5mm>*{}**@{-},
 <0mm,0mm>*{};<8mm,-5mm>*{}**@{-},
   <0mm,0mm>*{};<-8.5mm,-6.9mm>*{^1}**@{},
   <0mm,0mm>*{};<-5mm,-6.9mm>*{^2}**@{},
   <0mm,0mm>*{};<4.5mm,-6.9mm>*{^{n\hspace{-0.5mm}-\hspace{-0.5mm}1}}**@{},
   <0mm,0mm>*{};<9.0mm,-6.9mm>*{^n}**@{},
 \end{xy}\right)=
\sum_{k=0}^{n-1}\sum_{l=1}^{n-k} (-1)^{(l-1)(n-k-l)}
\begin{xy}
<0mm,0mm>*{\bullet}, <4mm,-7mm>*{^{1\ \ \dots\ k\qquad \ \ k+l+1\
\ \dots \ \ n}}, <-14mm,-5mm>*{}**@{-}, <-6mm,-5mm>*{}**@{-},
<20mm,-5mm>*{}**@{-}, <8mm,-5mm>*{}**@{-}, <0mm,-5mm>*{}**@{-},
<0mm,-5mm>*{\bullet}; <-5mm,-10mm>*{}**@{-},
<-2mm,-10mm>*{}**@{-}, <2mm,-10mm>*{}**@{-}, <5mm,-10mm>*{}**@{-},
<0mm,-12mm>*{_{k+1\dots k+l}},
 <0mm,0mm>*{};<0mm,5mm>*{}**@{-},
\end{xy}.$$

\begin{pro}
Let $\ass \xrightarrow{\varphi} \Qo$ be a map of (non-symmetric)
operads. The deformation complex of this map is isomorphism to
$\Qo$ up to the following shift of degree
$$C_\bullet^\varphi(\ass, \Qo)(n)=s^{-1}\Hom_\bullet^\Sy(\ass^{\ac}, \Qo)(n)=s^{-n}\Qo_\bullet(n).$$
The boundary map is given by
$$D(q)=d(q)+ \mu_\Qo(\varphi(\nu); I, q)+ \sum_{i=1}^n (-1)^i \mu_\Qo(q; I, \ldots, I,
\underbrace{\varphi(\nu)}_{i}, I, \ldots,I) +(-1)^{n+1}
\mu_\Qo(\varphi(\nu); q,I) ,$$ for $q\in \Qo(n)$ if we denote by
$\nu$ the generating binary operation of $\ass(2)$.
\end{pro}

\begin{proo}
There is a one-to-one correspondence between $\Sy_n$-equivariant
maps from $\ass^{\ac}(n)$ to $\Qo(n)$ and elements of $\Qo(n)$.
Let us denote by $f_q$ the unique map determined by $q\in\Qo(n)$.
Since $\ass^{\ac}$ is a coooperad and $\Qo$ is an operad, the convolution operad
$\Hom(\ass^{\ac}, \Qo)$ is preLie algebra with product denote $\star$ (see Section~\ref{Convolution prop(erad)}).
By Theorem~\ref{boundary=Q} and Section~\ref{Koszul case}, we have
$D(f_q)=d(f_q)+ \gamma\star f_q - (-1)^{|f_q|}f_q\star \gamma$.
Since $f_q$ vanishes on $\ass^{\ac}(m)$ for $m\neq n$ and since
$\gamma$ vanishes on $\ass^{\ac}(m)$ for $m\neq 2$, the only
non-vanishing component of $\gamma\star f_q=\mu_\Qo \circ (\gamma
\bc_{(1,1)} f_q ) \circ \Delta_{(1,1)}$ is $-\mu_\Qo(\varphi(\nu);
I, q) +(-1)^{n} \mu_\Qo(\varphi(\nu); q,I)$ on
$\Hom_\bullet^\Sy(\ass^{\ac}, \Qo)$. And the only non-vanishing
component of $f_q \star \gamma=\mu_\Qo \circ (f_q \bc_{(1,1)}
\gamma) \circ \Delta_{(1,1)}$ is $\sum_{i=1}^n (-1)^{i+1}
\mu_\Qo(q; I, \ldots, I, \underbrace{\varphi(\nu)}_{i}, I,
\ldots,I) $ on $\Hom_\bullet^\Sy(\ass^{\ac}, \Qo)$, which
concludes the proof.
\end{proo}

This deformation complex appears in many places in the literature
under different names. When $\Qo=\End_X$ with $X$ an associative
algebra, it is the Hochschild (co)chain complex of $X$ (\emph{with
coefficient in $X$}): $s^{-1}\Hom(V, \End_X)=\oplus_{n\geq
2}s^{1-n}\Hom(X^{\ot n}, X)$. The induced $L_\infty$-algebra, $Q$
on it is strict since the operad $\ass$ is Koszul. It is precisely
the Gerstenhaber Lie algebra \cite{Gerstenhaber63} and $Q^\ga$ is
the Hochschild dg Lie algebra controlling deformations of a
particular associative algebra structure, $\ga:\ass \rar \End_X$, on a vector space $X$.\\

In the work of McClure-Smith on Deligne's conjecture
\cite{McClureSmith02}, an operad $\Qo$ with a morphism of operads
$\ass \to \Qo$ is called a \emph{multiplicative operad}. The
simplicial complex that they define on such an operad is exactly
the deformation complex of this map. For the operad $\Qo=Poisson$,
this complex is related to the homology of long knots (see
\cite{Tourtchine04}). More generally, Maxim Kontsevich proposed
the conjecture that the deformation complex of $\ass \to \End_X$
is a $d+1$-algebra when $X$ is a $d$-algebra in
\cite{Kontsevich99}. This conjecture was proved by Tamarkin in
\cite{Tamarkin00}, see also  Hu, Kriz and Voronov
\cite{HuKrizVoronov06}. In this context, this chain complex is
often called the \emph{Hochschild
complex of $\Qo$}.\\

Since this (co)chain complex comes from the general theory of
(co)homology of Quillen, it would be better to call its
(co)homology the \emph{cohomology of $\ass$ with coefficients in
$\Qo$} or the chain complex, the \emph{deformation complex of the
map $\varphi$}.

\sip

Analogously one recovers other classical examples --- Harrison
complex/cohomology and Chevalley-Eilenberg complex/cohomology ---
from the operads of commutative algebras and, respectively, Lie
algebras.

\subsection{Poisson structures.}\label{Poisson structure}
 A {\em  Lie
1-bialgebra}\,  is, by definition, a graded vector space $V$
together with two linear maps,
$$
\Ba{rccc}
\delta: &  V & \lon & \wedge^2 V \\
       & a    & \lon & \sum a_{1}\wedge a_{2}
\Ea \ \ \ \ \  , \ \ \ \ \Ba{rccc}
[\, \bullet\, ]: & \odot^2 V& \lon & V \\
       & a\ot b    & \lon & (-1)^{|a|}[a\bullet b]
\Ea
$$
of degrees $0$ and $-1$ respectively which satisfy the identities,
\begin{itemize}
\item[(i)] $(\delta\ot\Id)\delta a + \tau (\delta\ot\Id)\delta a+
\tau^2 (\delta\ot\Id)\delta a =0$, where $\tau$ is the cyclic
permutation $(123)$ represented naturally on $V\ot V \ot V$
(co-Jacobi identity); \item[(ii)] $[[a\bullet b]\bullet
c]=[a\bullet[b\bullet c]] +
(-1)^{|b||a|+|b|+|a|}[b\bullet[a\bullet c]]$ (Jacobi identity);
\item[(iii)]
 $\delta [a\bullet b]=\sum a_1\wedge [a_2\bullet b] - (-1)^{|a_1||a_2|} a_2\wedge
[a_1\bullet b] +  [a\bullet b_1]\wedge b_2 -
(-1)^{|b_1||b_2|}[a\bullet b_2]\wedge b_1$ (Leibniz type
identity).
\end{itemize}
This notion of Lie 1-bialgebras is similar to the well-known
notion of Lie bialgebras except that in the latter case both
operations, Lie and co-Lie brackets, have degree $0$.

\sip

Let $\Lieb$ be the properad whose representations are Lie
1-bialgebras. It is Koszul contractible, that is its minimal
resolution, $(\Lieb_\infty, \delta)$,
 exists and is generated by the $\bS$-bimodule $V=\{V(m,n)\}_{m,n\geq1, m+n\geq 3}$ with
$$
V(m,n):=s^{m-2}sgn_m\ot \one_n =\mbox{span}\left\langle
\begin{xy}
 <0mm,0mm>*{\bullet};<0mm,0mm>*{}**@{},
 <-0.6mm,0.44mm>*{};<-8mm,5mm>*{}**@{-},
 <-0.4mm,0.7mm>*{};<-4.5mm,5mm>*{}**@{-},
 <0mm,0mm>*{};<-1mm,5mm>*{\ldots}**@{},
 <0.4mm,0.7mm>*{};<4.5mm,5mm>*{}**@{-},
 <0.6mm,0.44mm>*{};<8mm,5mm>*{}**@{-},
   <0mm,0mm>*{};<-8.5mm,5.5mm>*{^1}**@{},
   <0mm,0mm>*{};<-5mm,5.5mm>*{^2}**@{},
   <0mm,0mm>*{};<4.5mm,5.5mm>*{^{m\hspace{-0.5mm}-\hspace{-0.5mm}1}}**@{},
   <0mm,0mm>*{};<9.0mm,5.5mm>*{^m}**@{},
 <-0.6mm,-0.44mm>*{};<-8mm,-5mm>*{}**@{-},
 <-0.4mm,-0.7mm>*{};<-4.5mm,-5mm>*{}**@{-},
 <0mm,0mm>*{};<-1mm,-5mm>*{\ldots}**@{},
 <0.4mm,-0.7mm>*{};<4.5mm,-5mm>*{}**@{-},
 <0.6mm,-0.44mm>*{};<8mm,-5mm>*{}**@{-},
   <0mm,0mm>*{};<-8.5mm,-6.9mm>*{^1}**@{},
   <0mm,0mm>*{};<-5mm,-6.9mm>*{^2}**@{},
   <0mm,0mm>*{};<4.5mm,-6.9mm>*{^{n\hspace{-0.5mm}-\hspace{-0.5mm}1}}**@{},
   <0mm,0mm>*{};<9.0mm,-6.9mm>*{^n}**@{},
 \end{xy}
\right\rangle,
$$
where $sign_m$ stands for the sign representation of $\bS_m$ and
$\one_n$ for the trivial representation of $\bS_n$. The
differential is given on generators by \cite{Merkulov06}
$$
\delta \begin{xy}
 <0mm,0mm>*{\bullet};<0mm,0mm>*{}**@{},
 <0mm,0mm>*{};<-8mm,5mm>*{}**@{-},
 <0mm,0mm>*{};<-4.5mm,5mm>*{}**@{-},
 <0mm,0mm>*{};<-1mm,5mm>*{\ldots}**@{},
 <0mm,0mm>*{};<4.5mm,5mm>*{}**@{-},
 <0mm,0mm>*{};<8mm,5mm>*{}**@{-},
   <0mm,0mm>*{};<-8.5mm,5.5mm>*{^1}**@{},
   <0mm,0mm>*{};<-5mm,5.5mm>*{^2}**@{},
   <0mm,0mm>*{};<4.5mm,5.5mm>*{^{m\hspace{-0.5mm}-\hspace{-0.5mm}1}}**@{},
   <0mm,0mm>*{};<9.0mm,5.5mm>*{^m}**@{},
 <0mm,0mm>*{};<-8mm,-5mm>*{}**@{-},
 <0mm,0mm>*{};<-4.5mm,-5mm>*{}**@{-},
 <0mm,0mm>*{};<-1mm,-5mm>*{\ldots}**@{},
 <0mm,0mm>*{};<4.5mm,-5mm>*{}**@{-},
 <0mm,0mm>*{};<8mm,-5mm>*{}**@{-},
   <0mm,0mm>*{};<-8.5mm,-6.9mm>*{^1}**@{},
   <0mm,0mm>*{};<-5mm,-6.9mm>*{^2}**@{},
   <0mm,0mm>*{};<4.5mm,-6.9mm>*{^{n\hspace{-0.5mm}-\hspace{-0.5mm}1}}**@{},
   <0mm,0mm>*{};<9.0mm,-6.9mm>*{^n}**@{},
 \end{xy}
\ \ = \ \
 \sum_{I_1\sqcup I_2=(1,\ldots,m)\atop {J_1\sqcup J_2=(1,\ldots,n)\atop
 {|I_1|\geq 0, |I_2|\geq 1 \atop
 |J_1|\geq 1, |J_2|\geq 0}}
}\hspace{0mm} (-1)^{\sigma(I_1\sqcup I_2) + |I_1||I_2|}
 \begin{xy}
 <0mm,0mm>*{\bullet};<0mm,0mm>*{}**@{},
 <0mm,0mm>*{};<-8mm,5mm>*{}**@{-},
 <0mm,0mm>*{};<-4.5mm,5mm>*{}**@{-},
 <0mm,0mm>*{};<0mm,5mm>*{\ldots}**@{},
 <0mm,0mm>*{};<4.5mm,5mm>*{}**@{-},
 <0mm,0mm>*{};<13mm,5mm>*{}**@{-},
     <0mm,0mm>*{};<-2mm,7mm>*{\overbrace{\ \ \ \ \ \ \ \ \ \ \ \ }}**@{},
     <0mm,0mm>*{};<-2mm,9mm>*{^{I_1}}**@{},
 <0mm,0mm>*{};<-8mm,-5mm>*{}**@{-},
 <0mm,0mm>*{};<-4.5mm,-5mm>*{}**@{-},
 <0mm,0mm>*{};<-1mm,-5mm>*{\ldots}**@{},
 <0mm,0mm>*{};<4.5mm,-5mm>*{}**@{-},
 <0mm,0mm>*{};<8mm,-5mm>*{}**@{-},
      <0mm,0mm>*{};<0mm,-7mm>*{\underbrace{\ \ \ \ \ \ \ \ \ \ \ \ \ \ \
      }}**@{},
      <0mm,0mm>*{};<0mm,-10.6mm>*{_{J_1}}**@{},
 <13mm,5mm>*{};<13mm,5mm>*{\bullet}**@{},
 <13mm,5mm>*{};<5mm,10mm>*{}**@{-},
 <13mm,5mm>*{};<8.5mm,10mm>*{}**@{-},
 <13mm,5mm>*{};<13mm,10mm>*{\ldots}**@{},
 <13mm,5mm>*{};<16.5mm,10mm>*{}**@{-},
 <13mm,5mm>*{};<20mm,10mm>*{}**@{-},
      <13mm,5mm>*{};<13mm,12mm>*{\overbrace{\ \ \ \ \ \ \ \ \ \ \ \ \ \ }}**@{},
      <13mm,5mm>*{};<13mm,14mm>*{^{I_2}}**@{},
 <13mm,5mm>*{};<8mm,0mm>*{}**@{-},
 <13mm,5mm>*{};<12mm,0mm>*{\ldots}**@{},
 <13mm,5mm>*{};<16.5mm,0mm>*{}**@{-},
 <13mm,5mm>*{};<20mm,0mm>*{}**@{-},
     <13mm,5mm>*{};<14.3mm,-2mm>*{\underbrace{\ \ \ \ \ \ \ \ \ \ \ }}**@{},
     <13mm,5mm>*{};<14.3mm,-4.5mm>*{_{J_2}}**@{},
 \end{xy}
$$
where $\sigma(I_1\sqcup I_2)$ is the sign of the shuffle
$I_1\sqcup I_2=(1,\ldots, m)$.

\sip

Hence, for an arbitrary dg vector space $X$,
$$
s^{-1}\Hom(V, \End_X)=\bigoplus_{m,n\geq 1}s^{1-m}\wedge^m X\ot
\odot^n X\simeq \wedge^\bullet T_X,
$$
where $\wedge^\bullet T_X$ is the vector space of formal germs of
polyvector fields at $0\in X$ when we view $X$ as a formal graded
manifold. It is not hard to show using the above explicit formula
for the differential $\delta$ that the canonically induced, in
accordance with Theorem~\ref{mor}(i), $L_\infty$-structure on
$s^{-1}\Hom(V, \End_X)$ is precisely the classical Schouten Lie
algebra structure on polyvector fields. Thus our theory applied to
Lie 1-bialgebras reproduces deformation theory of Poisson
structures, and $\Lieb$-homology is precisely Poisson homology.

\sip

In a similar way one can check that our construction of
$L_\infty$-algebras applied to the minimal resolution of so called
pre-Lie$^2$-algebras \cite{Merkulov05} gives rise to another
classical geometric object --- the Fr\"olicher-Nijenhuis Lie
brackets on the sheaf, $T_X\ot \Omega^\bullet_X$, of tangent
vector bundle valued differential forms. Thus the associated
deformation theory describes deformations of integrable Nijenhuis
structures.

\subsection{Associative bialgebras} \label{Associative
bialgebras}

In this section, we make explicit the deformation theory of
representation of the properad $\ab$ of associative bialgebras. As
this example has never been rigorously treated in the literature
before, we show full details here.

\sip

As the properad $\ab$ is homotopy Koszul (see Section~\ref{SectionHoKoszul})
it admits a minimal resolution, $\ab_\infty=(\cF(\Co), \p)$, which is
 generated by a relatively small $\bS$-bimodule
 $\Co=\{\Co(m,n)\}_{m,n\geq 1, m+n\geq 3}$,
\[
\Co(m,n):= s^{m+n-3}\bK[\bS_m]\ot
\bK[\bS_n]=\mbox{span}\left\langle
\begin{xy}
 <0mm,0mm>*{\bullet};<0mm,0mm>*{}**@{},
 <0mm,0mm>*{};<-8mm,5mm>*{}**@{-},
 <0mm,0mm>*{};<-4.5mm,5mm>*{}**@{-},
 <0mm,0mm>*{};<-1mm,5mm>*{\ldots}**@{},
 <0mm,0mm>*{};<4.5mm,5mm>*{}**@{-},
 <0mm,0mm>*{};<8mm,5mm>*{}**@{-},
   <0mm,0mm>*{};<-8.5mm,5.5mm>*{^1}**@{},
   <0mm,0mm>*{};<-5mm,5.5mm>*{^2}**@{},
   <0mm,0mm>*{};<4.5mm,5.5mm>*{^{m\hspace{-0.5mm}-\hspace{-0.5mm}1}}**@{},
   <0mm,0mm>*{};<9.0mm,5.5mm>*{^m}**@{},
 <0mm,0mm>*{};<-8mm,-5mm>*{}**@{-},
 <0mm,0mm>*{};<-4.5mm,-5mm>*{}**@{-},
 <0mm,0mm>*{};<-1mm,-5mm>*{\ldots}**@{},
 <0mm,0mm>*{};<4.5mm,-5mm>*{}**@{-},
 <0mm,0mm>*{};<8mm,-5mm>*{}**@{-},
   <0mm,0mm>*{};<-8.5mm,-6.9mm>*{^1}**@{},
   <0mm,0mm>*{};<-5mm,-6.9mm>*{^2}**@{},
   <0mm,0mm>*{};<4.5mm,-6.9mm>*{^{n\hspace{-0.5mm}-\hspace{-0.5mm}1}}**@{},
   <0mm,0mm>*{};<9.0mm,-6.9mm>*{^n}**@{},
 \end{xy}
\right\rangle .
\]
 The differential $\p$ in
$\ab_\infty$ is neither quadratic nor of genus $0$. The derivation
$\p$ on $\F(\Co)$ is equivalent to a structure of
homotopy coproperad on $s^{-1}\Co$. 
 The values of $\p$ on $(1,n)$- and $(m,1)$-corollas
are given, of course, by the well-known $A_\infty$-formulae, while
\Beqrn \p
\begin{xy}
 <0mm,0mm>*{\bullet};<0mm,0mm>*{}**@{},
 <0mm,0mm>*{};<-4.5mm,5mm>*{}**@{-},
 <0mm,0mm>*{};<4.5mm,5mm>*{}**@{-},
   <0mm,0mm>*{};<-5mm,5.5mm>*{^1}**@{},
   <0mm,0mm>*{};<4.5mm,5.5mm>*{^2}**@{},
 <0mm,0mm>*{};<-8mm,-5mm>*{}**@{-},
 <0mm,0mm>*{};<-4.5mm,-5mm>*{}**@{-},
 <0mm,0mm>*{};<-1mm,-5mm>*{\ldots}**@{},
 <0mm,0mm>*{};<4.5mm,-5mm>*{}**@{-},
 <0mm,0mm>*{};<8mm,-5mm>*{}**@{-},
   <0mm,0mm>*{};<-8.5mm,-6.9mm>*{^1}**@{},
   <0mm,0mm>*{};<-5mm,-6.9mm>*{^2}**@{},
   <0mm,0mm>*{};<4.5mm,-6.9mm>*{^{n\hspace{-0.5mm}-\hspace{-0.5mm}1}}**@{},
   <0mm,0mm>*{};<9.0mm,-6.9mm>*{^n}**@{},
 \end{xy}&=&
\sum_{k=0}^{n-2}\sum_{l=2}^{n-k} (-1)^{k+l(n-k-l)+1}
\begin{xy}
<0mm,0mm>*{\bullet}, <4mm,-7mm>*{^{1\ \ \dots\ k\qquad \ \ k+l+1\
\ \dots \ \ n}}, <-14mm,-5mm>*{}**@{-}, <-6mm,-5mm>*{}**@{-},
<20mm,-5mm>*{}**@{-}, <8mm,-5mm>*{}**@{-}, <0mm,-5mm>*{}**@{-},
<0mm,-5mm>*{\bullet}; <-5mm,-10mm>*{}**@{-},
<-2mm,-10mm>*{}**@{-}, <2mm,-10mm>*{}**@{-}, <5mm,-10mm>*{}**@{-},
<0mm,-12mm>*{_{k+1\dots k+l}},
 <0mm,0mm>*{};<-4.5mm,5mm>*{}**@{-},
 <0mm,0mm>*{};<4.5mm,5mm>*{}**@{-},
   <0mm,0mm>*{};<-5mm,5.5mm>*{^1}**@{},
   <0mm,0mm>*{};<4.5mm,5.5mm>*{^2}**@{},
\end{xy}
\\
&& -\ \sum_{k=2}^\infty \sum_{r_1+\cdots +r_k=n}(-1)^s
\begin{xy}
<18mm,0mm>*{};<-18mm,0mm>*{}**@{-},
<-8mm,9mm>*{\Delta^{SU}}; <0mm,9mm>*{\bullet};
<0mm,9mm>*{};<0mm,13mm>*{}**@{-},
 <0mm,9mm>*{};<-3mm,5mm>*{}**@{-},
 <0mm,9mm>*{};<-5mm,5mm>*{}**@{-},
 <0mm,9mm>*{};<3mm,5mm>*{}**@{-},
 <0mm,9mm>*{};<5mm,5mm>*{}**@{-},
 <0mm,9mm>*{};<0mm,5mm>*{...}**@{},
<0mm,9mm>*{};<0mm,4mm>*{\underbrace{\ \ \ \ \ \ \ \ \ \ }}**@{},
<0mm,9mm>*{};<0mm,1mm>*{^k}**@{},
<-12mm,-5mm>*{\bullet}; <-12mm,-5mm>*{};<-9.5mm,-2mm>*{}**@{-},
<-12mm,-5mm>*{};<-14.5mm,-2mm>*{}**@{-},
<-12mm,-5mm>*{};<-13mm,-8mm>*{}**@{-},
<-12mm,-5mm>*{};<-15mm,-8mm>*{}**@{-},
<-12mm,-5mm>*{};<-8mm,-8mm>*{}**@{-}, <-11mm,-8mm>*{_{\ldots}};
<-12mm,-5mm>*{};<-11.8mm,-9.5mm>*{_{\underbrace{ }}}**@{},
<-11.8mm,-12.5mm>*{_{r_1}};
<-3mm,-5mm>*{\bullet}; <-3mm,-5mm>*{};<-0.5mm,-2mm>*{}**@{-},
<-3mm,-5mm>*{};<-5.5mm,-2mm>*{}**@{-},
<-3mm,-5mm>*{};<-4mm,-8mm>*{}**@{-},
<-3mm,-5mm>*{};<-6mm,-8mm>*{}**@{-},
<-3mm,-5mm>*{};<1mm,-8mm>*{}**@{-}, <-2mm,-8mm>*{_{\ldots}};
<-12mm,-5mm>*{};<-2.8mm,-9.5mm>*{_{\underbrace{ }}}**@{},
<-2.8mm,-12.5mm>*{_{r_2}};
<12mm,-5mm>*{\bullet}; <12mm,-5mm>*{};<9.5mm,-2mm>*{}**@{-},
<12mm,-5mm>*{};<14.5mm,-2mm>*{}**@{-},
<12mm,-5mm>*{};<13mm,-8mm>*{}**@{-},
<12mm,-5mm>*{};<15mm,-8mm>*{}**@{-},
<12mm,-5mm>*{};<8mm,-8mm>*{}**@{-}, <11mm,-8mm>*{_{\ldots}};
<12mm,-5mm>*{};<11.8mm,-9.5mm>*{_{\underbrace{ }}}**@{},
<11.8mm,-12.5mm>*{_{r_k}};
<5mm,-5mm>*{...}**@{},
\end{xy}
\Eeqrn where
$$
s=(k-1)(r_1-1) + (k-2)(r_2-1) +\ldots + 1(r_k-1),
$$
$\Delta^{SU}$ is the  Saneblidze-Umble diagonal, and the
horizontal line means
 {\em fraction}\, composition from \cite{Markl06}.
The meaning of this part of the differential is clear: it
describes $A_\infty$-morphisms between an $A_\infty$-structure on
$X$ and the associated Saneblidze-Umble diagonal
$A_\infty$-structure on $X\ot X$. Explicitly, this formula is
obtained by first considering the quasi-free resolution of the
2-colored operad coding two associative algebras and a morphism
between them. While the resolution of the associative operad is
given by the associahedra, this resolution is given by the
multiplihedra. This resolution gives the relaxed notion of
$A_\infty$-algebra and morphism of $A_\infty$-algebras at the same
time. Then, to get the formula above, we applied this resolution
to the $A_\infty$-algebra $X$ and to $X\ot X$ with the
$A_\infty$-algebra structure induced by the Saneblidze-Umble
diagonal.\\

The values of $\p$ on corollas of the form
$$
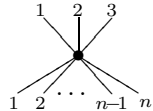
\begin{xy}
 <0mm,0mm>*{\bullet};<0mm,0mm>*{}**@{},
 <0mm,0mm>*{};<0mm,5mm>*{}**@{-},
 <0mm,0mm>*{};<-4.5mm,5mm>*{}**@{-},
 <0mm,0mm>*{};<4.5mm,5mm>*{}**@{-},
   <0mm,0mm>*{};<-5mm,5.5mm>*{^1}**@{},
   <0mm,0mm>*{};<4.5mm,5.5mm>*{^3}**@{},
   <0mm,0mm>*{};<0mm,5.5mm>*{^2}**@{},
 <0mm,0mm>*{};<-8mm,-5mm>*{}**@{-},
 <0mm,0mm>*{};<-4.5mm,-5mm>*{}**@{-},
 <0mm,0mm>*{};<-1mm,-5mm>*{\ldots}**@{},
 <0mm,0mm>*{};<4.5mm,-5mm>*{}**@{-},
 <0mm,0mm>*{};<8mm,-5mm>*{}**@{-},
   <0mm,0mm>*{};<-8.5mm,-6.9mm>*{^1}**@{},
   <0mm,0mm>*{};<-5mm,-6.9mm>*{^2}**@{},
   <0mm,0mm>*{};<4.5mm,-6.9mm>*{^{n\hspace{-0.5mm}-\hspace{-0.5mm}1}}**@{},
   <0mm,0mm>*{};<9.0mm,-6.9mm>*{^n}**@{},
 \end{xy}
$$
describe a {\em homotopy}\, between two natural
$A_\infty$-morphisms from $X$ to $X\ot X \ot X$, values on
corollas with 4 output legs ---  homotopies between homotopies
etc. We conjecture that $(\ab_\infty, \p)$ is a one coloured
version of a certain $\bN$-coloured properad describing
$A_\infty$-algebras, morphisms of $A_\infty$-algebras, homotopies
between morphisms of $A_\infty$-algebras, homotopies of homotopies
etc.,
 and we hope to describe it in a future publication.

\sip

It was proven in \cite{Markl06, MarklVoronov03}   that there
exists a minimal model $(\ab_\infty, \p)$ such that the
differential preserves Kontsevich's path grading of $\ab_\infty$
and has  the form $\p=\p_0 + \p_{pert}$, where $\p_0$ describes
the minimal resolution, $\frac{1}{2}\asb_\infty$ of the prop of
$\frac{1}{2}$-bialgebras (these facts follow also immediately from
our Corollary~\ref{generators}). The perturbation part,
$\p_{pert}$, is a linear combination of so called {\em
fractions}\, and their compositions. We shall assume from now on
that $\p$ has all these properties. By checking genus of these
fractions (or by referring to our proof of homotopy Koszulness of
$\ab$ in Section~\ref{SectionHoKoszul}) one can easily obtain the
following useful (for our purposes)

\begin{fact}\label{delta_0}
The differential $\p_0$ is precisely the quadratic part of $\p$,
i.e.\ it is equal to the composition,
$$
\p_0: \Co \stackrel{\p}{\lon} \cF(\Co) \stackrel{proj}{\lon}
\cF(\Co)^{(2)}.
$$
\end{fact}

\bip

Let $\Qo$ be a dg properad. By Theorem~\ref{mor}, the vector
space,
$$
s^{-1}\Hom^\Sy(\Co, \Qo)=\bigoplus_{n,m\geq 1\atop m+n\geq
3}s^{2-m-n}\Qo(m,n)=: \fg_{GS}(\Qo),
$$
has a canonical homotopy non-symmetric properad and
$L_\infty$-structure $Q$, whose Maurer-Cartan elements are
morphisms of properads $\F(\Co)
\to \Qo$. \\

If $\ga: \ab \rar \Qo$ is a representation of $\ab$, or more
generally of $\ab_\infty$~: $\ab_\infty \rar \Qo$, then, by
Definition~\ref{deftheory}, there exists an associated twisted
$L_\infty$-structures, $Q^\ga=\{Q_n^\ga\}_{n\geq 1}$, on
$\fg_{GS}(\Qo)$ which controls deformations of $\ga$ in the class
of representations of $\ab_\infty$. An explicit formula for the
differential $\p$ would induce an explicit $L_\infty$-structure.
Once again, our main example of this deformation theory is given
by $\Qo=\End_X$. In this case, the complex above is the
deformation complex of associative bialgebra, or more generally of
$\ab_\infty$-gebra, structure on $X$.\\

When $X$ is an associative bialgebra, Gerstenhaber and Schack
defined in \cite{GerstenhaberSchack90} a bicomplex whose homology has nice properties with
respect to the deformations of the associative bialgebra
structure (see also \cite{LazarevMovshev91}). Let us first extend this definition to any properad
$\Qo$ and not only $\End_X$.

\begin{dei}
Let $\gamma \, :\, \ab \to \Qo$ be a representation of $\ab$. We
define the \emph{Gerstenhaber-Schack bicomplex of $\gamma$} by
$C^{m,n}:=\Qo(m,n)$ and the differentials by \Beqr \nonumber d_h \
\begin{xy}
 <0mm,0mm>*{\bullet};<0mm,0mm>*{}**@{},
 <0mm,0mm>*{};<-8mm,5mm>*{}**@{-},
 <0mm,0mm>*{};<-4.5mm,5mm>*{}**@{-},
 <0mm,0mm>*{};<-1mm,5mm>*{\ldots}**@{},
 <0mm,0mm>*{};<4.5mm,5mm>*{}**@{-},
 <0mm,0mm>*{};<8mm,5mm>*{}**@{-},
   <0mm,0mm>*{};<-8.5mm,5.5mm>*{^1}**@{},
   <0mm,0mm>*{};<-5mm,5.5mm>*{^2}**@{},
   <0mm,0mm>*{};<4.5mm,5.5mm>*{^{m\hspace{-0.5mm}-\hspace{-0.5mm}1}}**@{},
   <0mm,0mm>*{};<9.0mm,5.5mm>*{^m}**@{},
 <0mm,0mm>*{};<-8mm,-5mm>*{}**@{-},
 <0mm,0mm>*{};<-4.5mm,-5mm>*{}**@{-},
 <0mm,0mm>*{};<-1mm,-5mm>*{\ldots}**@{},
 <0mm,0mm>*{};<4.5mm,-5mm>*{}**@{-},
 <0mm,0mm>*{};<8mm,-5mm>*{}**@{-},
   <0mm,0mm>*{};<-8.5mm,-6.9mm>*{^1}**@{},
   <0mm,0mm>*{};<-5mm,-6.9mm>*{^2}**@{},
   <0mm,0mm>*{};<4.5mm,-6.9mm>*{^{n\hspace{-0.5mm}-\hspace{-0.5mm}1}}**@{},
   <0mm,0mm>*{};<9.0mm,-6.9mm>*{^n}**@{},
 \end{xy}
 &:=&
\sum_{i=0}^{n-2}(-1)^{i+1}
 \begin{xy}
 <0mm,0mm>*{\bullet};<0mm,0mm>*{}**@{},
 <0mm,0mm>*{};<-8mm,-5mm>*{}**@{-},
 <0mm,0mm>*{};<-3.5mm,-5mm>*{}**@{-},
 <0mm,0mm>*{};<-6mm,-5mm>*{..}**@{},
 <0mm,0mm>*{};<0mm,-5mm>*{}**@{-},
  <0mm,-5mm>*{\bullet};
  <0mm,-5mm>*{};<-2.5mm,-8mm>*{}**@{-},
  <0mm,-5mm>*{};<2.5mm,-8mm>*{}**@{-},
  <0mm,-5mm>*{};<-3mm,-10mm>*{^{i\hspace{-0.2mm}+\hspace{-0.5mm}1}}**@{},
  <0mm,-5mm>*{};<3mm,-10mm>*{^{i\hspace{-0.2mm}+\hspace{-0.5mm}2}}**@{},
<0mm,0mm>*{};<8mm,-5mm>*{}**@{-},
<0mm,0mm>*{};<3.5mm,-5mm>*{}**@{-},
 <0mm,0mm>*{};<6mm,-5mm>*{..}**@{},
   <0mm,0mm>*{};<-8.5mm,-6.9mm>*{^1}**@{},
   <0mm,0mm>*{};<-4mm,-6.9mm>*{^i}**@{},
   <0mm,0mm>*{};<9.0mm,-6.9mm>*{^n}**@{},
 <0mm,0mm>*{};<-8mm,5mm>*{}**@{-},
 <0mm,0mm>*{};<-4.5mm,5mm>*{}**@{-},
 <0mm,0mm>*{};<-1mm,5mm>*{\ldots}**@{},
 <0mm,0mm>*{};<4.5mm,5mm>*{}**@{-},
 <0mm,0mm>*{};<8mm,5mm>*{}**@{-},
   <0mm,0mm>*{};<-8.5mm,5.5mm>*{^1}**@{},
   <0mm,0mm>*{};<-5mm,5.5mm>*{^2}**@{},
   <0mm,0mm>*{};<4.5mm,5.5mm>*{^{m\hspace{-0.5mm}-\hspace{-0.5mm}1}}**@{},
   <0mm,0mm>*{};<9.0mm,5.5mm>*{^m}**@{},
 \end{xy}
  \ \   +\ \ \
\begin{xy}
<18mm,0mm>*{};<-18mm,0mm>*{}**@{-},
<8mm,-9mm>*{\bullet}; <8mm,-9mm>*{};<5mm,-13mm>*{}**@{-},
 <8mm,-9mm>*{};<3mm,-13mm>*{}**@{-},
 <8mm,-9mm>*{};<11mm,-13mm>*{}**@{-},
 <8mm,-9mm>*{};<13mm,-13mm>*{}**@{-},
 <8mm,-9mm>*{};<8mm,-13mm>*{...}**@{},
 <8mm,-9mm>*{};<8mm,-15mm>*{\, ^2\ ^3 \  \ ^{n\hspace{-0.5mm}-\hspace{-0.5mm}1}\, ^n }**@{},
 <8mm,-9mm>*{};<5mm,-5mm>*{}**@{-},
 <8mm,-9mm>*{};<3mm,-5mm>*{}**@{-},
 <8mm,-9mm>*{};<11mm,-5mm>*{}**@{-},
 <8mm,-9mm>*{};<13mm,-5mm>*{}**@{-},
 <8mm,-9mm>*{};<8mm,-5mm>*{...}**@{},
<8mm,-9mm>*{};<8mm,-3.8mm>*{\overbrace{\ \ \ \ \ \ \ \ \ \
}}**@{}, <8mm,-9mm>*{};<8mm,-1.5mm>*{_m}**@{},
<-8mm,-10mm>*{\bullet}; <-8mm,-10mm>*{};<-8mm,-13mm>*{}**@{-},
<-8mm,-10mm>*{};<-8mm,-14.5mm>*{_1}**@{},
<-8mm,-10mm>*{};<-11.5mm,-7.4mm>*{}**@{-},
<-12mm,-7mm>*{};<-15.5mm,-4.5mm>*{}**@{-},
<-8mm,-10mm>*{};<-1mm,-5mm>*{}**@{-}, <-10mm,-8.5mm>*{\bullet};
<-10mm,-8.5mm>*{};<-5mm,-5mm>*{}**@{-},
<-13mm,-6mm>*{\bullet}; <-13mm,-6mm>*{};<-11mm,-4.6mm>*{}**@{-},
<-13mm,-6mm>*{};<-8.4mm,-5.2mm>*{...}**@{},
<8mm,-9mm>*{};<-8mm,-3.8mm>*{\overbrace{\ \ \ \ \ \ \ \ \ \ \ \ \
\ }}**@{}, <8mm,-9mm>*{};<-8mm,-1.5mm>*{_m}**@{},
<-12mm,5mm>*{\bullet}; <-12mm,5mm>*{};<-9.5mm,2mm>*{}**@{-},
<-12mm,5mm>*{};<-14.5mm,2mm>*{}**@{-},
<-12mm,5mm>*{};<-12mm,8mm>*{}**@{-},
<-12mm,5mm>*{};<-11.8mm,9.5mm>*{_1}**@{},
<-6mm,5mm>*{\bullet}; <-6mm,5mm>*{};<-3.5mm,2mm>*{}**@{-},
<-6mm,5mm>*{};<-8.5mm,2mm>*{}**@{-},
<-6mm,5mm>*{};<-6mm,8mm>*{}**@{-},
<-6mm,5mm>*{};<-5.8mm,9.5mm>*{_2}**@{},
<12mm,5mm>*{\bullet}; <12mm,5mm>*{};<9.5mm,2mm>*{}**@{-},
<12mm,5mm>*{};<14.5mm,2mm>*{}**@{-},
<12mm,5mm>*{};<12mm,8mm>*{}**@{-},
<12mm,5mm>*{};<6.4mm,9.5mm>*{_{m\hspace{-0.5mm}-\hspace{-0.5mm}1}}**@{},
<6mm,5mm>*{\bullet}; <6mm,5mm>*{};<3.5mm,2mm>*{}**@{-},
<6mm,5mm>*{};<8.5mm,2mm>*{}**@{-},
<6mm,5mm>*{};<6mm,8mm>*{}**@{-},
<6mm,5mm>*{};<12.1mm,9.5mm>*{_m}**@{},
<6mm,5mm>*{};<0mm,6mm>*{...}**@{},
\end{xy}
 \\
 &&
\nonumber+\ \ (-1)^{n+1}\begin{xy}
<18mm,0mm>*{};<-18mm,0mm>*{}**@{-},
<-8mm,-9mm>*{\bullet}; <-8mm,-9mm>*{};<-5mm,-13mm>*{}**@{-},
 <-8mm,-9mm>*{};<-3mm,-13mm>*{}**@{-},
 <-8mm,-9mm>*{};<-11mm,-13mm>*{}**@{-},
 <-8mm,-9mm>*{};<-13mm,-13mm>*{}**@{-},
 <-8mm,-9mm>*{};<-8mm,-13mm>*{...}**@{},
 <-8mm,-9mm>*{};<-8mm,-15mm>*{\, ^1\ ^2 \ \ \ \ \ \ ^{n\hspace{-0.5mm}-\hspace{-0.5mm}1} }**@{},
 <-8mm,-9mm>*{};<-5mm,-5mm>*{}**@{-},
 <-8mm,-9mm>*{};<-3mm,-5mm>*{}**@{-},
 <-8mm,-9mm>*{};<-11mm,-5mm>*{}**@{-},
 <-8mm,-9mm>*{};<-13mm,-5mm>*{}**@{-},
 <-8mm,-9mm>*{};<-8mm,-5mm>*{...}**@{},
<-8mm,-9mm>*{};<-8mm,-3.8mm>*{\overbrace{\ \ \ \ \ \ \ \ \ \
}}**@{}, <-8mm,-9mm>*{};<-8mm,-1.5mm>*{_m}**@{},
<8mm,-10mm>*{\bullet}; <8mm,-10mm>*{};<8mm,-13mm>*{}**@{-},
<8mm,-10mm>*{};<8mm,-14.5mm>*{_n}**@{},
<8mm,-10mm>*{};<11.5mm,-7.4mm>*{}**@{-},
<12mm,-7mm>*{};<15.5mm,-4.5mm>*{}**@{-},
<8mm,-10mm>*{};<1mm,-5mm>*{}**@{-}, <10mm,-8.5mm>*{\bullet};
<10mm,-8.5mm>*{};<5mm,-5mm>*{}**@{-},
<13mm,-6mm>*{\bullet}; <13mm,-6mm>*{};<11mm,-4.6mm>*{}**@{-},
<13mm,-6mm>*{};<8.4mm,-5.2mm>*{...}**@{},
<8mm,-9mm>*{};<8mm,-3.8mm>*{\overbrace{\ \ \ \ \ \ \ \ \ \ \ \ \ \
}}**@{}, <8mm,-9mm>*{};<8mm,-1.5mm>*{_m}**@{},
<-12mm,5mm>*{\bullet}; <-12mm,5mm>*{};<-9.5mm,2mm>*{}**@{-},
<-12mm,5mm>*{};<-14.5mm,2mm>*{}**@{-},
<-12mm,5mm>*{};<-12mm,8mm>*{}**@{-},
<-12mm,5mm>*{};<-11.8mm,9.5mm>*{_1}**@{},
<-6mm,5mm>*{\bullet}; <-6mm,5mm>*{};<-3.5mm,2mm>*{}**@{-},
<-6mm,5mm>*{};<-8.5mm,2mm>*{}**@{-},
<-6mm,5mm>*{};<-6mm,8mm>*{}**@{-},
<-6mm,5mm>*{};<-5.8mm,9.5mm>*{_2}**@{},
<12mm,5mm>*{\bullet}; <12mm,5mm>*{};<9.5mm,2mm>*{}**@{-},
<12mm,5mm>*{};<14.5mm,2mm>*{}**@{-},
<12mm,5mm>*{};<12mm,8mm>*{}**@{-},
<12mm,5mm>*{};<6.4mm,9.5mm>*{_{m\hspace{-0.5mm}-\hspace{-0.5mm}1}}**@{},
<6mm,5mm>*{\bullet}; <6mm,5mm>*{};<3.5mm,2mm>*{}**@{-},
<6mm,5mm>*{};<8.5mm,2mm>*{}**@{-},
<6mm,5mm>*{};<6mm,8mm>*{}**@{-},
<6mm,5mm>*{};<12.1mm,9.5mm>*{_m}**@{},
<6mm,5mm>*{};<0mm,6mm>*{...}**@{},
\end{xy}
\Eeqr \Beqr d_v \
\begin{xy}
 <0mm,0mm>*{\bullet};<0mm,0mm>*{}**@{},
 <0mm,0mm>*{};<-8mm,5mm>*{}**@{-},
 <0mm,0mm>*{};<-4.5mm,5mm>*{}**@{-},
 <0mm,0mm>*{};<-1mm,5mm>*{\ldots}**@{},
 <0mm,0mm>*{};<4.5mm,5mm>*{}**@{-},
 <0mm,0mm>*{};<8mm,5mm>*{}**@{-},
   <0mm,0mm>*{};<-8.5mm,5.5mm>*{^1}**@{},
   <0mm,0mm>*{};<-5mm,5.5mm>*{^2}**@{},
   <0mm,0mm>*{};<4.5mm,5.5mm>*{^{m\hspace{-0.5mm}-\hspace{-0.5mm}1}}**@{},
   <0mm,0mm>*{};<9.0mm,5.5mm>*{^m}**@{},
 <0mm,0mm>*{};<-8mm,-5mm>*{}**@{-},
 <0mm,0mm>*{};<-4.5mm,-5mm>*{}**@{-},
 <0mm,0mm>*{};<-1mm,-5mm>*{\ldots}**@{},
 <0mm,0mm>*{};<4.5mm,-5mm>*{}**@{-},
 <0mm,0mm>*{};<8mm,-5mm>*{}**@{-},
   <0mm,0mm>*{};<-8.5mm,-6.9mm>*{^1}**@{},
   <0mm,0mm>*{};<-5mm,-6.9mm>*{^2}**@{},
   <0mm,0mm>*{};<4.5mm,-6.9mm>*{^{n\hspace{-0.5mm}-\hspace{-0.5mm}1}}**@{},
   <0mm,0mm>*{};<9.0mm,-6.9mm>*{^n}**@{},
 \end{xy}
 &:=&
\ \sum_{i=0}^{n-2}(-1)^{i+1}
 \begin{xy}
 <0mm,0mm>*{\bullet};<0mm,0mm>*{}**@{},
 <0mm,0mm>*{};<-8mm,5mm>*{}**@{-},
 <0mm,0mm>*{};<-3.5mm,5mm>*{}**@{-},
 <0mm,0mm>*{};<-6mm,5mm>*{..}**@{},
 <0mm,0mm>*{};<0mm,5mm>*{}**@{-},
  <0mm,5mm>*{\bullet};
  <0mm,5mm>*{};<-2.5mm,8mm>*{}**@{-},
  <0mm,5mm>*{};<2.5mm,8mm>*{}**@{-},
  <0mm,5mm>*{};<-3mm,9mm>*{^{i\hspace{-0.2mm}+\hspace{-0.5mm}1}}**@{},
  <0mm,5mm>*{};<3mm,9mm>*{^{i\hspace{-0.2mm}+\hspace{-0.5mm}2}}**@{},
<0mm,0mm>*{};<8mm,5mm>*{}**@{-},
<0mm,0mm>*{};<3.5mm,5mm>*{}**@{-},
 <0mm,0mm>*{};<6mm,5mm>*{..}**@{},
   <0mm,0mm>*{};<-8.5mm,5.5mm>*{^1}**@{},
   <0mm,0mm>*{};<-4mm,5.5mm>*{^i}**@{},
   <0mm,0mm>*{};<9.0mm,5.5mm>*{^m}**@{},
 <0mm,0mm>*{};<-8mm,-5mm>*{}**@{-},
 <0mm,0mm>*{};<-4.5mm,-5mm>*{}**@{-},
 <0mm,0mm>*{};<-1mm,-5mm>*{\ldots}**@{},
 <0mm,0mm>*{};<4.5mm,-5mm>*{}**@{-},
 <0mm,0mm>*{};<8mm,-5mm>*{}**@{-},
   <0mm,0mm>*{};<-8.5mm,-6.9mm>*{^1}**@{},
   <0mm,0mm>*{};<-5mm,-6.9mm>*{^2}**@{},
   <0mm,0mm>*{};<4.5mm,-6.9mm>*{^{n\hspace{-0.5mm}-\hspace{-0.5mm}1}}**@{},
   <0mm,0mm>*{};<9.0mm,-6.9mm>*{^n}**@{},
 \end{xy}
 \
\ \ +\ \ \
\begin{xy}
<18mm,0mm>*{};<-18mm,0mm>*{}**@{-},
<8mm,9mm>*{\bullet}; <8mm,9mm>*{};<5mm,13mm>*{}**@{-},
 <8mm,9mm>*{};<3mm,13mm>*{}**@{-},
 <8mm,9mm>*{};<11mm,13mm>*{}**@{-},
 <8mm,9mm>*{};<13mm,13mm>*{}**@{-},
 <8mm,9mm>*{};<8mm,13mm>*{...}**@{},
 <8mm,9mm>*{};<8mm,14mm>*{\, ^2\ ^3 \  \ ^{m\hspace{-0.5mm}-\hspace{-0.5mm}1}\, ^m }**@{},
 <8mm,9mm>*{};<5mm,5mm>*{}**@{-},
 <8mm,9mm>*{};<3mm,5mm>*{}**@{-},
 <8mm,9mm>*{};<11mm,5mm>*{}**@{-},
 <8mm,9mm>*{};<13mm,5mm>*{}**@{-},
 <8mm,9mm>*{};<8mm,5mm>*{...}**@{},
<8mm,9mm>*{};<8mm,4mm>*{\underbrace{\ \ \ \ \ \ \ \ \ \ }}**@{},
<8mm,9mm>*{};<8mm,1mm>*{^n}**@{},
<-8mm,10mm>*{\bullet};<0mm,0mm>*{}**@{},
<-8mm,10mm>*{};<-8mm,13mm>*{}**@{-},
<-8mm,10mm>*{};<-8mm,14mm>*{^1}**@{},
<-8mm,10mm>*{};<-11.5mm,7.4mm>*{}**@{-},
<-12mm,7mm>*{};<-15.5mm,4.5mm>*{}**@{-},
<-8mm,10mm>*{};<-1mm,5mm>*{}**@{-}, <-10mm,8.5mm>*{\bullet};
<-10mm,8.5mm>*{};<-5mm,5mm>*{}**@{-},
<-13mm,6mm>*{\bullet}; <-13mm,6mm>*{};<-11mm,4.6mm>*{}**@{-},
<-13mm,6mm>*{};<-8.4mm,5.2mm>*{...}**@{},
<8mm,9mm>*{};<-8mm,4mm>*{\underbrace{\ \ \ \ \ \ \ \ \ \ \ \ \ \
}}**@{}, <8mm,9mm>*{};<-8mm,1mm>*{^n}**@{},
<-12mm,-5mm>*{\bullet}; <-12mm,-5mm>*{};<-9.5mm,-2mm>*{}**@{-},
<-12mm,-5mm>*{};<-14.5mm,-2mm>*{}**@{-},
<-12mm,-5mm>*{};<-12mm,-8mm>*{}**@{-},
<-12mm,-5mm>*{};<-11.8mm,-9.5mm>*{_1}**@{},
<-6mm,-5mm>*{\bullet}; <-6mm,-5mm>*{};<-3.5mm,-2mm>*{}**@{-},
<-6mm,-5mm>*{};<-8.5mm,-2mm>*{}**@{-},
<-6mm,-5mm>*{};<-6mm,-8mm>*{}**@{-},
<-6mm,-5mm>*{};<-5.8mm,-9.5mm>*{_2}**@{},
<12mm,-5mm>*{\bullet}; <12mm,-5mm>*{};<9.5mm,-2mm>*{}**@{-},
<12mm,-5mm>*{};<14.5mm,-2mm>*{}**@{-},
<12mm,-5mm>*{};<12mm,-8mm>*{}**@{-},
<12mm,-5mm>*{};<6.4mm,-9.5mm>*{_{n\hspace{-0.5mm}-\hspace{-0.5mm}1}}**@{},
<6mm,-5mm>*{\bullet}; <6mm,-5mm>*{};<3.5mm,-2mm>*{}**@{-},
<6mm,-5mm>*{};<8.5mm,-2mm>*{}**@{-},
<6mm,-5mm>*{};<6mm,-8mm>*{}**@{-},
<6mm,-5mm>*{};<12.1mm,-9.5mm>*{_n}**@{},
<6mm,-5mm>*{};<0mm,-8mm>*{...}**@{},
\end{xy}
\vspace{6mm} \nonumber
\\
&& \ \ +\ \ (-1)^{m+1}
\begin{xy}
<18mm,0mm>*{};<-18mm,0mm>*{}**@{-},
<-8mm,9mm>*{\bullet}; <-8mm,9mm>*{};<-5mm,13mm>*{}**@{-},
 <-8mm,9mm>*{};<-3mm,13mm>*{}**@{-},
 <-8mm,9mm>*{};<-11mm,13mm>*{}**@{-},
 <-8mm,9mm>*{};<-13mm,13mm>*{}**@{-},
 <-8mm,9mm>*{};<-8mm,13mm>*{...}**@{},
 <-8mm,9mm>*{};<-8mm,14mm>*{\, ^1\ ^1 \ \ \ \, \ ^{m\hspace{-0.5mm}-\hspace{-0.5mm}1}}**@{},
 <-8mm,9mm>*{};<-5mm,5mm>*{}**@{-},
 <-8mm,9mm>*{};<-3mm,5mm>*{}**@{-},
 <-8mm,9mm>*{};<-11mm,5mm>*{}**@{-},
 <-8mm,9mm>*{};<-13mm,5mm>*{}**@{-},
 <-8mm,9mm>*{};<-8mm,5mm>*{...}**@{},
<-8mm,9mm>*{};<-8mm,4mm>*{\underbrace{\ \ \ \ \ \ \ \ \ \ }}**@{},
<-8mm,9mm>*{};<-8mm,1mm>*{^n}**@{},
<8mm,10mm>*{\bullet};<0mm,0mm>*{}**@{},
<8mm,10mm>*{};<8mm,13mm>*{}**@{-},
<8mm,10mm>*{};<8mm,14mm>*{^m}**@{},
<8mm,10mm>*{};<11.5mm,7.4mm>*{}**@{-},
<12mm,7mm>*{};<15.5mm,4.5mm>*{}**@{-},
<8mm,10mm>*{};<1mm,5mm>*{}**@{-}, <10mm,8.5mm>*{\bullet};
<10mm,8.5mm>*{};<5mm,5mm>*{}**@{-},
<13mm,6mm>*{\bullet}; <13mm,6mm>*{};<11mm,4.6mm>*{}**@{-},
<13mm,6mm>*{};<8.4mm,5.2mm>*{...}**@{},
<8mm,9mm>*{};<8mm,4mm>*{\underbrace{\ \ \ \ \ \ \ \ \ \ \ \ \ \
}}**@{}, <8mm,9mm>*{};<8mm,1mm>*{^n}**@{},
<-12mm,-5mm>*{\bullet}; <-12mm,-5mm>*{};<-9.5mm,-2mm>*{}**@{-},
<-12mm,-5mm>*{};<-14.5mm,-2mm>*{}**@{-},
<-12mm,-5mm>*{};<-12mm,-8mm>*{}**@{-},
<-12mm,-5mm>*{};<-11.8mm,-9.5mm>*{_1}**@{},
<-6mm,-5mm>*{\bullet}; <-6mm,-5mm>*{};<-3.5mm,-2mm>*{}**@{-},
<-6mm,-5mm>*{};<-8.5mm,-2mm>*{}**@{-},
<-6mm,-5mm>*{};<-6mm,-8mm>*{}**@{-},
<-6mm,-5mm>*{};<-5.8mm,-9.5mm>*{_2}**@{},
<12mm,-5mm>*{\bullet}; <12mm,-5mm>*{};<9.5mm,-2mm>*{}**@{-},
<12mm,-5mm>*{};<14.5mm,-2mm>*{}**@{-},
<12mm,-5mm>*{};<12mm,-8mm>*{}**@{-},
<12mm,-5mm>*{};<6.4mm,-9.5mm>*{_{n\hspace{-0.5mm}-\hspace{-0.5mm}1}}**@{},
<6mm,-5mm>*{\bullet}; <6mm,-5mm>*{};<3.5mm,-2mm>*{}**@{-},
<6mm,-5mm>*{};<8.5mm,-2mm>*{}**@{-},
<6mm,-5mm>*{};<6mm,-8mm>*{}**@{-},
<6mm,-5mm>*{};<12.1mm,-9.5mm>*{_n}**@{},
<6mm,-5mm>*{};<0mm,-8mm>*{...}**@{},
\end{xy}
\vspace{6mm} \nonumber \Eeqr where the general $(m,n)$-corollas
have to be understand as elements of $\Qo(m,n)$. The binary
corollas are the image under $\gamma$ of the generating product
and coproduct of $\ab$. Finally, this pictures represent the
composition of all these elements in $\Qo$.
\end{dei}

Let us compare these with the Gerstenhaber-Schack differential,
$d_{GS}$, in the bicomplex $C_{GS}^{m,n}:=\Hom(X^{\ot n}, X^{\ot
m})$ which is  defined by \cite{GerstenhaberSchack90}
$$
d_{GS}= d_h + d_v,
$$
with $d_h: \Hom(X^{\ot n}, X^{\ot m})\rar \Hom(X^{\ot n+1}, X^{\ot
m})$ given on an arbitrary $f\in  \Hom(X^{\ot n}, X^{\ot m})$ by
\Beqrn (d_hf)(a_0, a_1, \ldots, a_n)&:=&\Delta^m(a_0)\Box f(a_1,
a_2, \ldots, a_n) -
\sum_{i=0}^{n-1} (-1)^if(a_1, \dots, a_ia_{i+1}, \ldots, a_n)\\
&& + (-1)^{n+1}f(a_1, a_2, \ldots, a_n)\Box\Delta^m(a_n)\ \ \
\forall\ \ a_i\in X . \Eeqrn Here the multiplication in $X$ is
denoted by juxtaposition, the induced multiplication in the
algebra $X^{\ot m}$ by $\Box$, the comultiplication in $X$ by
$\Delta$, and
$$
\Delta^n: (\Delta\ot \Id^{\ot m-2 })\circ (\Delta\ot \Id^{\ot
m-3})\circ \ldots \circ \Delta: X \rar X^{\ot m}.
$$
The expression for $d_v$ is an obvious ``dual" analogue of $d_h$.
Now let us represent $d_h$ in graphical terms by associating the
graphs
$$
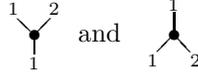
\begin{xy}
 <0mm,-0.55mm>*{};<0mm,-2.5mm>*{}**@{-},
 <0.5mm,0.5mm>*{};<2.2mm,2.2mm>*{}**@{-},
 <-0.48mm,0.48mm>*{};<-2.2mm,2.2mm>*{}**@{-},
 <0mm,0mm>*{\bullet};<0mm,0mm>*{}**@{},
 <0mm,-0.55mm>*{};<0mm,-3.8mm>*{_1}**@{},
 <0.5mm,0.5mm>*{};<2.7mm,2.8mm>*{^2}**@{},
 <-0.48mm,0.48mm>*{};<-2.7mm,2.8mm>*{^1}**@{},
 \end{xy}
\ \ \mbox{and}\ \ \
 \begin{xy}
 <0mm,0.66mm>*{};<0mm,3mm>*{}**@{-},
 <0.39mm,-0.39mm>*{};<2.2mm,-2.2mm>*{}**@{-},
 <-0.35mm,-0.35mm>*{};<-2.2mm,-2.2mm>*{}**@{-},
 <0mm,0mm>*{\bullet};<0mm,0mm>*{}**@{},
   <0mm,0.66mm>*{};<0mm,3.4mm>*{^1}**@{},
   <0.39mm,-0.39mm>*{};<2.9mm,-4mm>*{^2}**@{},
   <-0.35mm,-0.35mm>*{};<-2.8mm,-4mm>*{^1}**@{},
\end{xy}
$$
to comultiplication and, respectively, multiplication while the
corolla
$$
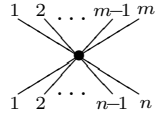
\begin{xy}
 <0mm,0mm>*{\bullet};<0mm,0mm>*{}**@{},
 <0mm,0mm>*{};<-8mm,5mm>*{}**@{-},
 <0mm,0mm>*{};<-4.5mm,5mm>*{}**@{-},
 <0mm,0mm>*{};<-1mm,5mm>*{\ldots}**@{},
 <0mm,0mm>*{};<4.5mm,5mm>*{}**@{-},
 <0mm,0mm>*{};<8mm,5mm>*{}**@{-},
   <0mm,0mm>*{};<-8.5mm,5.5mm>*{^1}**@{},
   <0mm,0mm>*{};<-5mm,5.5mm>*{^2}**@{},
   <0mm,0mm>*{};<4.5mm,5.5mm>*{^{m\hspace{-0.5mm}-\hspace{-0.5mm}1}}**@{},
   <0mm,0mm>*{};<9.0mm,5.5mm>*{^m}**@{},
 <0mm,0mm>*{};<-8mm,-5mm>*{}**@{-},
 <0mm,0mm>*{};<-4.5mm,-5mm>*{}**@{-},
 <0mm,0mm>*{};<-1mm,-5mm>*{\ldots}**@{},
 <0mm,0mm>*{};<4.5mm,-5mm>*{}**@{-},
 <0mm,0mm>*{};<8mm,-5mm>*{}**@{-},
   <0mm,0mm>*{};<-8.5mm,-6.9mm>*{^1}**@{},
   <0mm,0mm>*{};<-5mm,-6.9mm>*{^2}**@{},
   <0mm,0mm>*{};<4.5mm,-6.9mm>*{^{n\hspace{-0.5mm}-\hspace{-0.5mm}1}}**@{},
   <0mm,0mm>*{};<9.0mm,-6.9mm>*{^n}**@{},
 \end{xy}
$$
to $f$. Then the r.h.s\ of the formula for $d_h$ reads,
$$
\begin{xy}
<18mm,0mm>*{};<-18mm,0mm>*{}**@{-},
<8mm,-9mm>*{\bullet}; <8mm,-9mm>*{};<5mm,-13mm>*{}**@{-},
 <8mm,-9mm>*{};<3mm,-13mm>*{}**@{-},
 <8mm,-9mm>*{};<11mm,-13mm>*{}**@{-},
 <8mm,-9mm>*{};<13mm,-13mm>*{}**@{-},
 <8mm,-9mm>*{};<8mm,-13mm>*{...}**@{},
 <8mm,-9mm>*{};<8mm,-15mm>*{\, ^{a_1}\ \ \  \ \ \ \ ^{a_n} }**@{},
 <8mm,-9mm>*{};<5mm,-5mm>*{}**@{-},
 <8mm,-9mm>*{};<3mm,-5mm>*{}**@{-},
 <8mm,-9mm>*{};<11mm,-5mm>*{}**@{-},
 <8mm,-9mm>*{};<13mm,-5mm>*{}**@{-},
 <8mm,-9mm>*{};<8mm,-5mm>*{...}**@{},
<8mm,-9mm>*{};<8mm,-3.8mm>*{\overbrace{\ \ \ \ \ \ \ \ \ \
}}**@{}, <8mm,-9mm>*{};<8mm,-1.5mm>*{_m}**@{},
<-8mm,-10mm>*{\bullet}; <-8mm,-10mm>*{};<-8mm,-13mm>*{}**@{-},
<-8mm,-10mm>*{};<-8mm,-14.5mm>*{_{a_0}}**@{},
<-8mm,-10mm>*{};<-11.5mm,-7.4mm>*{}**@{-},
<-12mm,-7mm>*{};<-15.5mm,-4.5mm>*{}**@{-},
<-8mm,-10mm>*{};<-1mm,-5mm>*{}**@{-}, <-10mm,-8.5mm>*{\bullet};
<-10mm,-8.5mm>*{};<-5mm,-5mm>*{}**@{-},
<-13mm,-6mm>*{\bullet}; <-13mm,-6mm>*{};<-11mm,-4.6mm>*{}**@{-},
<-13mm,-6mm>*{};<-8.4mm,-5.2mm>*{...}**@{},
<8mm,-9mm>*{};<-8mm,-3.8mm>*{\overbrace{\ \ \ \ \ \ \ \ \ \ \ \ \
\ }}**@{}, <8mm,-9mm>*{};<-8mm,-1.5mm>*{_m}**@{},
<-12mm,5mm>*{\bullet}; <-12mm,5mm>*{};<-9.5mm,2mm>*{}**@{-},
<-12mm,5mm>*{};<-14.5mm,2mm>*{}**@{-},
<-12mm,5mm>*{};<-12mm,8mm>*{}**@{-},
<-12mm,5mm>*{};<-11.8mm,9.5mm>*{_1}**@{},
<-6mm,5mm>*{\bullet}; <-6mm,5mm>*{};<-3.5mm,2mm>*{}**@{-},
<-6mm,5mm>*{};<-8.5mm,2mm>*{}**@{-},
<-6mm,5mm>*{};<-6mm,8mm>*{}**@{-},
<-6mm,5mm>*{};<-5.8mm,9.5mm>*{_2}**@{},
<12mm,5mm>*{\bullet}; <12mm,5mm>*{};<9.5mm,2mm>*{}**@{-},
<12mm,5mm>*{};<14.5mm,2mm>*{}**@{-},
<12mm,5mm>*{};<12mm,8mm>*{}**@{-},
<12mm,5mm>*{};<6.4mm,9.5mm>*{_{m\hspace{-0.5mm}-\hspace{-0.5mm}1}}**@{},
<6mm,5mm>*{\bullet}; <6mm,5mm>*{};<3.5mm,2mm>*{}**@{-},
<6mm,5mm>*{};<8.5mm,2mm>*{}**@{-},
<6mm,5mm>*{};<6mm,8mm>*{}**@{-},
<6mm,5mm>*{};<12.1mm,9.5mm>*{_m}**@{},
<6mm,5mm>*{};<0mm,6mm>*{...}**@{},
\end{xy}
 \
-
\ \sum_{i=0}^{n-1}(-1)^i
 \begin{xy}
 <0mm,0mm>*{\bullet};<0mm,0mm>*{}**@{},
 <0mm,0mm>*{};<-8mm,-5mm>*{}**@{-},
 <0mm,0mm>*{};<-3.5mm,-5mm>*{}**@{-},
 <0mm,0mm>*{};<-6mm,-5mm>*{..}**@{},
 <0mm,0mm>*{};<0mm,-5mm>*{}**@{-},
  <0mm,-5mm>*{\bullet};
  <0mm,-5mm>*{};<-2.5mm,-8mm>*{}**@{-},
  <0mm,-5mm>*{};<2.5mm,-8mm>*{}**@{-},
  <0mm,-5mm>*{};<-3mm,-10mm>*{^{a_i}}**@{},
  <0mm,-5mm>*{};<3mm,-10mm>*{^{a_{i\hspace{-0.2mm}+\hspace{-0.5mm}1}}}**@{},
<0mm,0mm>*{};<8mm,-5mm>*{}**@{-},
<0mm,0mm>*{};<3.5mm,-5mm>*{}**@{-},
 <0mm,0mm>*{};<6mm,-5mm>*{..}**@{},
   <0mm,0mm>*{};<-8.5mm,-6.9mm>*{^{a_0}}**@{},
   <0mm,0mm>*{};<-4mm,-6.9mm>*{}**@{},
   <0mm,0mm>*{};<9.0mm,-6.9mm>*{^n}**@{},
 <0mm,0mm>*{};<-8mm,5mm>*{}**@{-},
 <0mm,0mm>*{};<-4.5mm,5mm>*{}**@{-},
 <0mm,0mm>*{};<-1mm,5mm>*{\ldots}**@{},
 <0mm,0mm>*{};<4.5mm,5mm>*{}**@{-},
 <0mm,0mm>*{};<8mm,5mm>*{}**@{-},
   <0mm,0mm>*{};<-8.5mm,5.5mm>*{^1}**@{},
   <0mm,0mm>*{};<-5mm,5.5mm>*{^2}**@{},
   <0mm,0mm>*{};<4.5mm,5.5mm>*{^{m\hspace{-0.5mm}-\hspace{-0.5mm}1}}**@{},
   <0mm,0mm>*{};<9.0mm,5.5mm>*{^m}**@{},
 \end{xy}
+\  (-1)^{n+1}
\begin{xy}
<18mm,0mm>*{};<-18mm,0mm>*{}**@{-},
<-8mm,-9mm>*{\bullet}; <-8mm,-9mm>*{};<-5mm,-13mm>*{}**@{-},
 <-8mm,-9mm>*{};<-3mm,-13mm>*{}**@{-},
 <-8mm,-9mm>*{};<-11mm,-13mm>*{}**@{-},
 <-8mm,-9mm>*{};<-13mm,-13mm>*{}**@{-},
 <-8mm,-9mm>*{};<-8mm,-13mm>*{...}**@{},
 <-8mm,-9mm>*{};<-8mm,-15mm>*{\, ^{a_0}\ \ \ \ \ \ \ \ ^{a_{n\hspace{-0.5mm}-\hspace{-0.5mm}1}} }**@{},
 <-8mm,-9mm>*{};<-5mm,-5mm>*{}**@{-},
 <-8mm,-9mm>*{};<-3mm,-5mm>*{}**@{-},
 <-8mm,-9mm>*{};<-11mm,-5mm>*{}**@{-},
 <-8mm,-9mm>*{};<-13mm,-5mm>*{}**@{-},
 <-8mm,-9mm>*{};<-8mm,-5mm>*{...}**@{},
<-8mm,-9mm>*{};<-8mm,-3.8mm>*{\overbrace{\ \ \ \ \ \ \ \ \ \
}}**@{}, <-8mm,-9mm>*{};<-8mm,-1.5mm>*{_{m}}**@{},
<8mm,-10mm>*{\bullet}; <8mm,-10mm>*{};<8mm,-13mm>*{}**@{-},
<8mm,-10mm>*{};<8mm,-14.5mm>*{_{a_n}}**@{},
<8mm,-10mm>*{};<11.5mm,-7.4mm>*{}**@{-},
<12mm,-7mm>*{};<15.5mm,-4.5mm>*{}**@{-},
<8mm,-10mm>*{};<1mm,-5mm>*{}**@{-}, <10mm,-8.5mm>*{\bullet};
<10mm,-8.5mm>*{};<5mm,-5mm>*{}**@{-},
<13mm,-6mm>*{\bullet}; <13mm,-6mm>*{};<11mm,-4.6mm>*{}**@{-},
<13mm,-6mm>*{};<8.4mm,-5.2mm>*{...}**@{},
<8mm,-9mm>*{};<8mm,-3.8mm>*{\overbrace{\ \ \ \ \ \ \ \ \ \ \ \ \ \
}}**@{}, <8mm,-9mm>*{};<8mm,-1.5mm>*{_m}**@{},
<-12mm,5mm>*{\bullet}; <-12mm,5mm>*{};<-9.5mm,2mm>*{}**@{-},
<-12mm,5mm>*{};<-14.5mm,2mm>*{}**@{-},
<-12mm,5mm>*{};<-12mm,8mm>*{}**@{-},
<-12mm,5mm>*{};<-11.8mm,9.5mm>*{_1}**@{},
<-6mm,5mm>*{\bullet}; <-6mm,5mm>*{};<-3.5mm,2mm>*{}**@{-},
<-6mm,5mm>*{};<-8.5mm,2mm>*{}**@{-},
<-6mm,5mm>*{};<-6mm,8mm>*{}**@{-},
<-6mm,5mm>*{};<-5.8mm,9.5mm>*{_2}**@{},
<12mm,5mm>*{\bullet}; <12mm,5mm>*{};<9.5mm,2mm>*{}**@{-},
<12mm,5mm>*{};<14.5mm,2mm>*{}**@{-},
<12mm,5mm>*{};<12mm,8mm>*{}**@{-},
<12mm,5mm>*{};<6.4mm,9.5mm>*{_{m\hspace{-0.5mm}-\hspace{-0.5mm}1}}**@{},
<6mm,5mm>*{\bullet}; <6mm,5mm>*{};<3.5mm,2mm>*{}**@{-},
<6mm,5mm>*{};<8.5mm,2mm>*{}**@{-},
<6mm,5mm>*{};<6mm,8mm>*{}**@{-},
<6mm,5mm>*{};<12.1mm,9.5mm>*{_m}**@{},
<6mm,5mm>*{};<0mm,6mm>*{...}**@{},
\end{xy}
$$
which is precisely the first three summands in the previous
definition. The other three terms correspond to $d_v$. Therefore,
when $\Qo=\End_X$ and $\ga \, :\, \ab \to \End_X$ is an
associative bialgebra structure on $X$, the preceding bicomplex is
exactly Gerstenhaber-Schack bicomplex
\cite{GerstenhaberSchack90}.\\

However, even without an explicit minimal model of $\ab$, we can
show the following general result.

\begin{thm}\label{GS}
Let $(\ab_\infty, \p)\stackrel{\pi}{\lon} \ab$, be a minimal model
of the properad of bialgebras and $\ga:\ab \rar \Qo$ an arbitrary
representation of $\ab$.  Then  the differential,
$$
Q_1^\ga=Q\circ e^{\ga\odot}
$$
associated to this minimal model in the twisted
$L_\infty$-structure, $Q^\ga$, on $\fg_{GS}$,
 is isomorphic to the Gerstenhaber-Schack differential.
Hence the deformation complex of representation of $\ab$ is
isomorphic to the Gerstenhaber-Schack bicomplex.
\end{thm}

\begin{proo}
Let $(\ab_\infty=\cF(\Co), \p)$ be a minimal model of the properad
of bialgebras, and let $I$ be the ideal in $\cF(\Co)$ generated by
graphs in $\cF(\Co)^{(\geq 2)}$ with at least two non-binary
(i.e.\ neither $\begin{xy}
 <0mm,-0.55mm>*{};<0mm,-2.5mm>*{}**@{-},
 <0.5mm,0.5mm>*{};<2.2mm,2.2mm>*{}**@{-},
 <-0.48mm,0.48mm>*{};<-2.2mm,2.2mm>*{}**@{-},
 <0mm,0mm>*{\bullet};<0mm,0mm>*{}**@{},
 \end{xy}$ nor
\begin{xy}
 <0mm,0.66mm>*{};<0mm,3mm>*{}**@{-},
 <0.39mm,-0.39mm>*{};<2.2mm,-2.2mm>*{}**@{-},
 <-0.35mm,-0.35mm>*{};<-2.2mm,-2.2mm>*{}**@{-},
 <0mm,0mm>*{\bullet};<0mm,0mm>*{}**@{},
\end{xy})
vertices, and let
$$
{\hsB}:= \frac{\ab_\infty}{(I, \p I)}
$$
be the associated quotient {\em dg }\, properad. The induced
differential in ${\hsB}$ we denote by ${\hp}$. It is precisely
this quotient part, $\hp$, of the total  differential $\p$ which
completely determines the $L_\infty$-differential differential
$Q_1^\ga$. Thus our plan is the following: in the next Lemma we
present an explicit, up to an automorphism, form of the
differential $\hp$ (despite the fact that $\p$ is not explicit !)
 and thereafter compare the resulting $Q_1^\ga$ with the
Gerstenhaber-Schack definition.

\sip

The major step in the proof is the following Lemma (in its
formulation we use fraction notations again).

\begin{lem}\label{keylemma}
(i)  The derivation, $d$, of\, ${\hsB}$ given on generators by,
\Beq\label{d-1} d\begin{xy}
 <0mm,-0.55mm>*{};<0mm,-2.5mm>*{}**@{-},
 <0.5mm,0.5mm>*{};<2.2mm,2.2mm>*{}**@{-},
 <-0.48mm,0.48mm>*{};<-2.2mm,2.2mm>*{}**@{-},
 <0mm,0mm>*{\bullet};<0mm,0mm>*{}**@{},
 <0mm,-0.55mm>*{};<0mm,-3.8mm>*{_1}**@{},
 <0.5mm,0.5mm>*{};<2.7mm,2.8mm>*{^2}**@{},
 <-0.48mm,0.48mm>*{};<-2.7mm,2.8mm>*{^1}**@{},
 \end{xy}
 =0
 \ \ , \  \ \ \
 d
 \begin{xy}
 <0mm,0.66mm>*{};<0mm,3mm>*{}**@{-},
 <0.39mm,-0.39mm>*{};<2.2mm,-2.2mm>*{}**@{-},
 <-0.35mm,-0.35mm>*{};<-2.2mm,-2.2mm>*{}**@{-},
 <0mm,0mm>*{\bullet};<0mm,0mm>*{}**@{},
   <0mm,0.66mm>*{};<0mm,3.4mm>*{^1}**@{},
   <0.39mm,-0.39mm>*{};<2.9mm,-4mm>*{^2}**@{},
   <-0.35mm,-0.35mm>*{};<-2.8mm,-4mm>*{^1}**@{},
\end{xy}
=0, \Eeq

\Beq\label{d-2}
 d
 \begin{xy}
 <0.39mm,-0.39mm>*{};<2.2mm,-2.2mm>*{}**@{-},
 <-0.35mm,-0.35mm>*{};<-2.2mm,-2.2mm>*{}**@{-},
 <0.39mm,0.39mm>*{};<2.2mm,2.2mm>*{}**@{-},
 <-0.35mm,0.35mm>*{};<-2.2mm,2.2mm>*{}**@{-},
 <0mm,0mm>*{\bullet};<0mm,0mm>*{}**@{},
   <0.39mm,-0.39mm>*{};<2.9mm,-4mm>*{^2}**@{},
   <-0.35mm,-0.35mm>*{};<-2.8mm,-4mm>*{^1}**@{},
    <0mm,0.66mm>*{};<2.8mm,3.4mm>*{^2}**@{},
<0mm,0.66mm>*{};<-2.8mm,3.4mm>*{^1}**@{},
\end{xy}
=
 \begin{xy}
 <0mm,2.47mm>*{};<0mm,-0.5mm>*{}**@{-},
 <0.5mm,3.5mm>*{};<2.2mm,5.2mm>*{}**@{-},
 <-0.48mm,3.48mm>*{};<-2.2mm,5.2mm>*{}**@{-},
 <0mm,3mm>*{\bullet};<0mm,3mm>*{}**@{},
  <0mm,-0.8mm>*{\bullet};<0mm,-0.8mm>*{}**@{},
<0mm,-0.8mm>*{};<-2.2mm,-3.5mm>*{}**@{-},
 <0mm,-0.8mm>*{};<2.2mm,-3.5mm>*{}**@{-},
     <0.5mm,3.5mm>*{};<2.8mm,5.7mm>*{^2}**@{},
     <-0.48mm,3.48mm>*{};<-2.8mm,5.7mm>*{^1}**@{},
   <0mm,-0.8mm>*{};<-2.7mm,-5.2mm>*{^1}**@{},
   <0mm,-0.8mm>*{};<2.7mm,-5.2mm>*{^2}**@{},
\end{xy}
\ - \
\begin{xy}
 <0mm,0mm>*{\bullet};<0mm,0mm>*{}**@{},
 <0mm,-0.49mm>*{};<0mm,-3.0mm>*{}**@{-},
 <-0.5mm,0.5mm>*{};<-3mm,2mm>*{}**@{-},
 <-3mm,2mm>*{};<0mm,4mm>*{}**@{-},
 <0mm,4mm>*{\bullet};<-2.3mm,2.3mm>*{}**@{},
 <0mm,4mm>*{};<0mm,7.4mm>*{}**@{-},
<0mm,0mm>*{};<2.2mm,1.5mm>*{}**@{-},
 <6mm,0mm>*{\bullet};<0mm,0mm>*{}**@{},
 <6mm,4mm>*{};<3.8mm,2.5mm>*{}**@{-},
 <6mm,4mm>*{};<6mm,7.4mm>*{}**@{-},
 <6mm,4mm>*{\bullet};<-2.3mm,2.3mm>*{}**@{},
 <0mm,4mm>*{};<6mm,0mm>*{}**@{-},
<6mm,4mm>*{};<9mm,2mm>*{}**@{-}, <6mm,0mm>*{};<9mm,2mm>*{}**@{-},
<6mm,0mm>*{};<6mm,-3mm>*{}**@{-},
   <-1.8mm,2.8mm>*{};<0mm,7.8mm>*{^1}**@{},
   <-2.8mm,2.9mm>*{};<0mm,-4.3mm>*{_1}**@{},
<-1.8mm,2.8mm>*{};<6mm,7.8mm>*{^2}**@{},
   <-2.8mm,2.9mm>*{};<6mm,-4.3mm>*{_2}**@{},
 \end{xy}
\Eeq and, for all other generators with $m+n\geq 4$, by
\Beqr\label{d-3} d
\begin{xy}
 <0mm,0mm>*{\bullet};<0mm,0mm>*{}**@{},
 <0mm,0mm>*{};<-8mm,5mm>*{}**@{-},
 <0mm,0mm>*{};<-4.5mm,5mm>*{}**@{-},
 <0mm,0mm>*{};<-1mm,5mm>*{\ldots}**@{},
 <0mm,0mm>*{};<4.5mm,5mm>*{}**@{-},
 <0mm,0mm>*{};<8mm,5mm>*{}**@{-},
   <0mm,0mm>*{};<-8.5mm,5.5mm>*{^1}**@{},
   <0mm,0mm>*{};<-5mm,5.5mm>*{^2}**@{},
   <0mm,0mm>*{};<4.5mm,5.5mm>*{^{m\hspace{-0.5mm}-\hspace{-0.5mm}1}}**@{},
   <0mm,0mm>*{};<9.0mm,5.5mm>*{^m}**@{},
 <0mm,0mm>*{};<-8mm,-5mm>*{}**@{-},
 <0mm,0mm>*{};<-4.5mm,-5mm>*{}**@{-},
 <0mm,0mm>*{};<-1mm,-5mm>*{\ldots}**@{},
 <0mm,0mm>*{};<4.5mm,-5mm>*{}**@{-},
 <0mm,0mm>*{};<8mm,-5mm>*{}**@{-},
   <0mm,0mm>*{};<-8.5mm,-6.9mm>*{^1}**@{},
   <0mm,0mm>*{};<-5mm,-6.9mm>*{^2}**@{},
   <0mm,0mm>*{};<4.5mm,-6.9mm>*{^{n\hspace{-0.5mm}-\hspace{-0.5mm}1}}**@{},
   <0mm,0mm>*{};<9.0mm,-6.9mm>*{^n}**@{},
 \end{xy}
 &=&
\sum_{i=0}^{n-2}(-1)^{i+1}
 \begin{xy}
 <0mm,0mm>*{\bullet};<0mm,0mm>*{}**@{},
 <0mm,0mm>*{};<-8mm,-5mm>*{}**@{-},
 <0mm,0mm>*{};<-3.5mm,-5mm>*{}**@{-},
 <0mm,0mm>*{};<-6mm,-5mm>*{..}**@{},
 <0mm,0mm>*{};<0mm,-5mm>*{}**@{-},
  <0mm,-5mm>*{\bullet};
  <0mm,-5mm>*{};<-2.5mm,-8mm>*{}**@{-},
  <0mm,-5mm>*{};<2.5mm,-8mm>*{}**@{-},
  <0mm,-5mm>*{};<-3mm,-10mm>*{^{i\hspace{-0.2mm}+\hspace{-0.5mm}1}}**@{},
  <0mm,-5mm>*{};<3mm,-10mm>*{^{i\hspace{-0.2mm}+\hspace{-0.5mm}2}}**@{},
<0mm,0mm>*{};<8mm,-5mm>*{}**@{-},
<0mm,0mm>*{};<3.5mm,-5mm>*{}**@{-},
 <0mm,0mm>*{};<6mm,-5mm>*{..}**@{},
   <0mm,0mm>*{};<-8.5mm,-6.9mm>*{^1}**@{},
   <0mm,0mm>*{};<-4mm,-6.9mm>*{^i}**@{},
   <0mm,0mm>*{};<9.0mm,-6.9mm>*{^n}**@{},
 <0mm,0mm>*{};<-8mm,5mm>*{}**@{-},
 <0mm,0mm>*{};<-4.5mm,5mm>*{}**@{-},
 <0mm,0mm>*{};<-1mm,5mm>*{\ldots}**@{},
 <0mm,0mm>*{};<4.5mm,5mm>*{}**@{-},
 <0mm,0mm>*{};<8mm,5mm>*{}**@{-},
   <0mm,0mm>*{};<-8.5mm,5.5mm>*{^1}**@{},
   <0mm,0mm>*{};<-5mm,5.5mm>*{^2}**@{},
   <0mm,0mm>*{};<4.5mm,5.5mm>*{^{m\hspace{-0.5mm}-\hspace{-0.5mm}1}}**@{},
   <0mm,0mm>*{};<9.0mm,5.5mm>*{^m}**@{},
 \end{xy}
 \ \   +\ \ \
\begin{xy}
<18mm,0mm>*{};<-18mm,0mm>*{}**@{-},
<8mm,-9mm>*{\bullet}; <8mm,-9mm>*{};<5mm,-13mm>*{}**@{-},
 <8mm,-9mm>*{};<3mm,-13mm>*{}**@{-},
 <8mm,-9mm>*{};<11mm,-13mm>*{}**@{-},
 <8mm,-9mm>*{};<13mm,-13mm>*{}**@{-},
 <8mm,-9mm>*{};<8mm,-13mm>*{...}**@{},
 <8mm,-9mm>*{};<8mm,-15mm>*{\, ^2\ ^3 \  \ ^{n\hspace{-0.5mm}-\hspace{-0.5mm}1}\, ^n }**@{},
 <8mm,-9mm>*{};<5mm,-5mm>*{}**@{-},
 <8mm,-9mm>*{};<3mm,-5mm>*{}**@{-},
 <8mm,-9mm>*{};<11mm,-5mm>*{}**@{-},
 <8mm,-9mm>*{};<13mm,-5mm>*{}**@{-},
 <8mm,-9mm>*{};<8mm,-5mm>*{...}**@{},
<8mm,-9mm>*{};<8mm,-3.8mm>*{\overbrace{\ \ \ \ \ \ \ \ \ \
}}**@{}, <8mm,-9mm>*{};<8mm,-1.5mm>*{_m}**@{},
<-8mm,-10mm>*{\bullet}; <-8mm,-10mm>*{};<-8mm,-13mm>*{}**@{-},
<-8mm,-10mm>*{};<-8mm,-14.5mm>*{_1}**@{},
<-8mm,-10mm>*{};<-11.5mm,-7.4mm>*{}**@{-},
<-12mm,-7mm>*{};<-15.5mm,-4.5mm>*{}**@{-},
<-8mm,-10mm>*{};<-1mm,-5mm>*{}**@{-}, <-10mm,-8.5mm>*{\bullet};
<-10mm,-8.5mm>*{};<-5mm,-5mm>*{}**@{-},
<-13mm,-6mm>*{\bullet}; <-13mm,-6mm>*{};<-11mm,-4.6mm>*{}**@{-},
<-13mm,-6mm>*{};<-8.4mm,-5.2mm>*{...}**@{},
<8mm,-9mm>*{};<-8mm,-3.8mm>*{\overbrace{\ \ \ \ \ \ \ \ \ \ \ \ \
\ }}**@{}, <8mm,-9mm>*{};<-8mm,-1.5mm>*{_m}**@{},
<-12mm,5mm>*{\bullet}; <-12mm,5mm>*{};<-9.5mm,2mm>*{}**@{-},
<-12mm,5mm>*{};<-14.5mm,2mm>*{}**@{-},
<-12mm,5mm>*{};<-12mm,8mm>*{}**@{-},
<-12mm,5mm>*{};<-11.8mm,9.5mm>*{_1}**@{},
<-6mm,5mm>*{\bullet}; <-6mm,5mm>*{};<-3.5mm,2mm>*{}**@{-},
<-6mm,5mm>*{};<-8.5mm,2mm>*{}**@{-},
<-6mm,5mm>*{};<-6mm,8mm>*{}**@{-},
<-6mm,5mm>*{};<-5.8mm,9.5mm>*{_2}**@{},
<12mm,5mm>*{\bullet}; <12mm,5mm>*{};<9.5mm,2mm>*{}**@{-},
<12mm,5mm>*{};<14.5mm,2mm>*{}**@{-},
<12mm,5mm>*{};<12mm,8mm>*{}**@{-},
<12mm,5mm>*{};<6.4mm,9.5mm>*{_{m\hspace{-0.5mm}-\hspace{-0.5mm}1}}**@{},
<6mm,5mm>*{\bullet}; <6mm,5mm>*{};<3.5mm,2mm>*{}**@{-},
<6mm,5mm>*{};<8.5mm,2mm>*{}**@{-},
<6mm,5mm>*{};<6mm,8mm>*{}**@{-},
<6mm,5mm>*{};<12.1mm,9.5mm>*{_m}**@{},
<6mm,5mm>*{};<0mm,6mm>*{...}**@{},
\end{xy}
 \\
 &&
+\ \ (-1)^{n+1}\begin{xy} <18mm,0mm>*{};<-18mm,0mm>*{}**@{-},
<-8mm,-9mm>*{\bullet}; <-8mm,-9mm>*{};<-5mm,-13mm>*{}**@{-},
 <-8mm,-9mm>*{};<-3mm,-13mm>*{}**@{-},
 <-8mm,-9mm>*{};<-11mm,-13mm>*{}**@{-},
 <-8mm,-9mm>*{};<-13mm,-13mm>*{}**@{-},
 <-8mm,-9mm>*{};<-8mm,-13mm>*{...}**@{},
 <-8mm,-9mm>*{};<-8mm,-15mm>*{\, ^1\ ^2 \ \ \ \ \ \ ^{n\hspace{-0.5mm}-\hspace{-0.5mm}1} }**@{},
 <-8mm,-9mm>*{};<-5mm,-5mm>*{}**@{-},
 <-8mm,-9mm>*{};<-3mm,-5mm>*{}**@{-},
 <-8mm,-9mm>*{};<-11mm,-5mm>*{}**@{-},
 <-8mm,-9mm>*{};<-13mm,-5mm>*{}**@{-},
 <-8mm,-9mm>*{};<-8mm,-5mm>*{...}**@{},
<-8mm,-9mm>*{};<-8mm,-3.8mm>*{\overbrace{\ \ \ \ \ \ \ \ \ \
}}**@{}, <-8mm,-9mm>*{};<-8mm,-1.5mm>*{_m}**@{},
<8mm,-10mm>*{\bullet}; <8mm,-10mm>*{};<8mm,-13mm>*{}**@{-},
<8mm,-10mm>*{};<8mm,-14.5mm>*{_n}**@{},
<8mm,-10mm>*{};<11.5mm,-7.4mm>*{}**@{-},
<12mm,-7mm>*{};<15.5mm,-4.5mm>*{}**@{-},
<8mm,-10mm>*{};<1mm,-5mm>*{}**@{-}, <10mm,-8.5mm>*{\bullet};
<10mm,-8.5mm>*{};<5mm,-5mm>*{}**@{-},
<13mm,-6mm>*{\bullet}; <13mm,-6mm>*{};<11mm,-4.6mm>*{}**@{-},
<13mm,-6mm>*{};<8.4mm,-5.2mm>*{...}**@{},
<8mm,-9mm>*{};<8mm,-3.8mm>*{\overbrace{\ \ \ \ \ \ \ \ \ \ \ \ \ \
}}**@{}, <8mm,-9mm>*{};<8mm,-1.5mm>*{_m}**@{},
<-12mm,5mm>*{\bullet}; <-12mm,5mm>*{};<-9.5mm,2mm>*{}**@{-},
<-12mm,5mm>*{};<-14.5mm,2mm>*{}**@{-},
<-12mm,5mm>*{};<-12mm,8mm>*{}**@{-},
<-12mm,5mm>*{};<-11.8mm,9.5mm>*{_1}**@{},
<-6mm,5mm>*{\bullet}; <-6mm,5mm>*{};<-3.5mm,2mm>*{}**@{-},
<-6mm,5mm>*{};<-8.5mm,2mm>*{}**@{-},
<-6mm,5mm>*{};<-6mm,8mm>*{}**@{-},
<-6mm,5mm>*{};<-5.8mm,9.5mm>*{_2}**@{},
<12mm,5mm>*{\bullet}; <12mm,5mm>*{};<9.5mm,2mm>*{}**@{-},
<12mm,5mm>*{};<14.5mm,2mm>*{}**@{-},
<12mm,5mm>*{};<12mm,8mm>*{}**@{-},
<12mm,5mm>*{};<6.4mm,9.5mm>*{_{m\hspace{-0.5mm}-\hspace{-0.5mm}1}}**@{},
<6mm,5mm>*{\bullet}; <6mm,5mm>*{};<3.5mm,2mm>*{}**@{-},
<6mm,5mm>*{};<8.5mm,2mm>*{}**@{-},
<6mm,5mm>*{};<6mm,8mm>*{}**@{-},
<6mm,5mm>*{};<12.1mm,9.5mm>*{_m}**@{},
<6mm,5mm>*{};<0mm,6mm>*{...}**@{},
\end{xy}
 \nonumber\\
  &&
+\ \sum_{i=0}^{n-2}(-1)^{i+1}
 \begin{xy}
 <0mm,0mm>*{\bullet};<0mm,0mm>*{}**@{},
 <0mm,0mm>*{};<-8mm,5mm>*{}**@{-},
 <0mm,0mm>*{};<-3.5mm,5mm>*{}**@{-},
 <0mm,0mm>*{};<-6mm,5mm>*{..}**@{},
 <0mm,0mm>*{};<0mm,5mm>*{}**@{-},
  <0mm,5mm>*{\bullet};
  <0mm,5mm>*{};<-2.5mm,8mm>*{}**@{-},
  <0mm,5mm>*{};<2.5mm,8mm>*{}**@{-},
  <0mm,5mm>*{};<-3mm,9mm>*{^{i\hspace{-0.2mm}+\hspace{-0.5mm}1}}**@{},
  <0mm,5mm>*{};<3mm,9mm>*{^{i\hspace{-0.2mm}+\hspace{-0.5mm}2}}**@{},
<0mm,0mm>*{};<8mm,5mm>*{}**@{-},
<0mm,0mm>*{};<3.5mm,5mm>*{}**@{-},
 <0mm,0mm>*{};<6mm,5mm>*{..}**@{},
   <0mm,0mm>*{};<-8.5mm,5.5mm>*{^1}**@{},
   <0mm,0mm>*{};<-4mm,5.5mm>*{^i}**@{},
   <0mm,0mm>*{};<9.0mm,5.5mm>*{^m}**@{},
 <0mm,0mm>*{};<-8mm,-5mm>*{}**@{-},
 <0mm,0mm>*{};<-4.5mm,-5mm>*{}**@{-},
 <0mm,0mm>*{};<-1mm,-5mm>*{\ldots}**@{},
 <0mm,0mm>*{};<4.5mm,-5mm>*{}**@{-},
 <0mm,0mm>*{};<8mm,-5mm>*{}**@{-},
   <0mm,0mm>*{};<-8.5mm,-6.9mm>*{^1}**@{},
   <0mm,0mm>*{};<-5mm,-6.9mm>*{^2}**@{},
   <0mm,0mm>*{};<4.5mm,-6.9mm>*{^{n\hspace{-0.5mm}-\hspace{-0.5mm}1}}**@{},
   <0mm,0mm>*{};<9.0mm,-6.9mm>*{^n}**@{},
 \end{xy}
 \
\ \ +\ \ \
\begin{xy}
<18mm,0mm>*{};<-18mm,0mm>*{}**@{-},
<8mm,9mm>*{\bullet}; <8mm,9mm>*{};<5mm,13mm>*{}**@{-},
 <8mm,9mm>*{};<3mm,13mm>*{}**@{-},
 <8mm,9mm>*{};<11mm,13mm>*{}**@{-},
 <8mm,9mm>*{};<13mm,13mm>*{}**@{-},
 <8mm,9mm>*{};<8mm,13mm>*{...}**@{},
 <8mm,9mm>*{};<8mm,14mm>*{\, ^2\ ^3 \  \ ^{m\hspace{-0.5mm}-\hspace{-0.5mm}1}\, ^m }**@{},
 <8mm,9mm>*{};<5mm,5mm>*{}**@{-},
 <8mm,9mm>*{};<3mm,5mm>*{}**@{-},
 <8mm,9mm>*{};<11mm,5mm>*{}**@{-},
 <8mm,9mm>*{};<13mm,5mm>*{}**@{-},
 <8mm,9mm>*{};<8mm,5mm>*{...}**@{},
<8mm,9mm>*{};<8mm,4mm>*{\underbrace{\ \ \ \ \ \ \ \ \ \ }}**@{},
<8mm,9mm>*{};<8mm,1mm>*{^n}**@{},
<-8mm,10mm>*{\bullet};<0mm,0mm>*{}**@{},
<-8mm,10mm>*{};<-8mm,13mm>*{}**@{-},
<-8mm,10mm>*{};<-8mm,14mm>*{^1}**@{},
<-8mm,10mm>*{};<-11.5mm,7.4mm>*{}**@{-},
<-12mm,7mm>*{};<-15.5mm,4.5mm>*{}**@{-},
<-8mm,10mm>*{};<-1mm,5mm>*{}**@{-}, <-10mm,8.5mm>*{\bullet};
<-10mm,8.5mm>*{};<-5mm,5mm>*{}**@{-},
<-13mm,6mm>*{\bullet}; <-13mm,6mm>*{};<-11mm,4.6mm>*{}**@{-},
<-13mm,6mm>*{};<-8.4mm,5.2mm>*{...}**@{},
<8mm,9mm>*{};<-8mm,4mm>*{\underbrace{\ \ \ \ \ \ \ \ \ \ \ \ \ \
}}**@{}, <8mm,9mm>*{};<-8mm,1mm>*{^n}**@{},
<-12mm,-5mm>*{\bullet}; <-12mm,-5mm>*{};<-9.5mm,-2mm>*{}**@{-},
<-12mm,-5mm>*{};<-14.5mm,-2mm>*{}**@{-},
<-12mm,-5mm>*{};<-12mm,-8mm>*{}**@{-},
<-12mm,-5mm>*{};<-11.8mm,-9.5mm>*{_1}**@{},
<-6mm,-5mm>*{\bullet}; <-6mm,-5mm>*{};<-3.5mm,-2mm>*{}**@{-},
<-6mm,-5mm>*{};<-8.5mm,-2mm>*{}**@{-},
<-6mm,-5mm>*{};<-6mm,-8mm>*{}**@{-},
<-6mm,-5mm>*{};<-5.8mm,-9.5mm>*{_2}**@{},
<12mm,-5mm>*{\bullet}; <12mm,-5mm>*{};<9.5mm,-2mm>*{}**@{-},
<12mm,-5mm>*{};<14.5mm,-2mm>*{}**@{-},
<12mm,-5mm>*{};<12mm,-8mm>*{}**@{-},
<12mm,-5mm>*{};<6.4mm,-9.5mm>*{_{n\hspace{-0.5mm}-\hspace{-0.5mm}1}}**@{},
<6mm,-5mm>*{\bullet}; <6mm,-5mm>*{};<3.5mm,-2mm>*{}**@{-},
<6mm,-5mm>*{};<8.5mm,-2mm>*{}**@{-},
<6mm,-5mm>*{};<6mm,-8mm>*{}**@{-},
<6mm,-5mm>*{};<12.1mm,-9.5mm>*{_n}**@{},
<6mm,-5mm>*{};<0mm,-8mm>*{...}**@{},
\end{xy}
\vspace{6mm} \nonumber
\\
&& \ \ +\ \ (-1)^{m+1}
\begin{xy}
<18mm,0mm>*{};<-18mm,0mm>*{}**@{-},
<-8mm,9mm>*{\bullet}; <-8mm,9mm>*{};<-5mm,13mm>*{}**@{-},
 <-8mm,9mm>*{};<-3mm,13mm>*{}**@{-},
 <-8mm,9mm>*{};<-11mm,13mm>*{}**@{-},
 <-8mm,9mm>*{};<-13mm,13mm>*{}**@{-},
 <-8mm,9mm>*{};<-8mm,13mm>*{...}**@{},
 <-8mm,9mm>*{};<-8mm,14mm>*{\, ^1\ ^1 \ \ \ \, \ ^{m\hspace{-0.5mm}-\hspace{-0.5mm}1}}**@{},
 <-8mm,9mm>*{};<-5mm,5mm>*{}**@{-},
 <-8mm,9mm>*{};<-3mm,5mm>*{}**@{-},
 <-8mm,9mm>*{};<-11mm,5mm>*{}**@{-},
 <-8mm,9mm>*{};<-13mm,5mm>*{}**@{-},
 <-8mm,9mm>*{};<-8mm,5mm>*{...}**@{},
<-8mm,9mm>*{};<-8mm,4mm>*{\underbrace{\ \ \ \ \ \ \ \ \ \ }}**@{},
<-8mm,9mm>*{};<-8mm,1mm>*{^n}**@{},
<8mm,10mm>*{\bullet};<0mm,0mm>*{}**@{},
<8mm,10mm>*{};<8mm,13mm>*{}**@{-},
<8mm,10mm>*{};<8mm,14mm>*{^m}**@{},
<8mm,10mm>*{};<11.5mm,7.4mm>*{}**@{-},
<12mm,7mm>*{};<15.5mm,4.5mm>*{}**@{-},
<8mm,10mm>*{};<1mm,5mm>*{}**@{-}, <10mm,8.5mm>*{\bullet};
<10mm,8.5mm>*{};<5mm,5mm>*{}**@{-},
<13mm,6mm>*{\bullet}; <13mm,6mm>*{};<11mm,4.6mm>*{}**@{-},
<13mm,6mm>*{};<8.4mm,5.2mm>*{...}**@{},
<8mm,9mm>*{};<8mm,4mm>*{\underbrace{\ \ \ \ \ \ \ \ \ \ \ \ \ \
}}**@{}, <8mm,9mm>*{};<8mm,1mm>*{^n}**@{},
<-12mm,-5mm>*{\bullet}; <-12mm,-5mm>*{};<-9.5mm,-2mm>*{}**@{-},
<-12mm,-5mm>*{};<-14.5mm,-2mm>*{}**@{-},
<-12mm,-5mm>*{};<-12mm,-8mm>*{}**@{-},
<-12mm,-5mm>*{};<-11.8mm,-9.5mm>*{_1}**@{},
<-6mm,-5mm>*{\bullet}; <-6mm,-5mm>*{};<-3.5mm,-2mm>*{}**@{-},
<-6mm,-5mm>*{};<-8.5mm,-2mm>*{}**@{-},
<-6mm,-5mm>*{};<-6mm,-8mm>*{}**@{-},
<-6mm,-5mm>*{};<-5.8mm,-9.5mm>*{_2}**@{},
<12mm,-5mm>*{\bullet}; <12mm,-5mm>*{};<9.5mm,-2mm>*{}**@{-},
<12mm,-5mm>*{};<14.5mm,-2mm>*{}**@{-},
<12mm,-5mm>*{};<12mm,-8mm>*{}**@{-},
<12mm,-5mm>*{};<6.4mm,-9.5mm>*{_{n\hspace{-0.5mm}-\hspace{-0.5mm}1}}**@{},
<6mm,-5mm>*{\bullet}; <6mm,-5mm>*{};<3.5mm,-2mm>*{}**@{-},
<6mm,-5mm>*{};<8.5mm,-2mm>*{}**@{-},
<6mm,-5mm>*{};<6mm,-8mm>*{}**@{-},
<6mm,-5mm>*{};<12.1mm,-9.5mm>*{_n}**@{},
<6mm,-5mm>*{};<0mm,-8mm>*{...}**@{},
\end{xy}
\vspace{6mm} \nonumber \Eeqr is a differential.

\bip

(ii) The dg properads $(\hsB,\hp)$ and $(\hsB, d)$ are isomorphic.
\end{lem}
\begin{proo}
(i) It is easy to see that among for 2-vertex connected binary
graphs\footnote{Equivalence classes of graphs in  $\hsB$ we call
simply graphs for shortness.} attached to any other graph in $\sB$
the bialgebra relations,
$$
\begin{xy}
 <0mm,0mm>*{\bullet};<0mm,0mm>*{}**@{},
 <0mm,-0.49mm>*{};<0mm,-3.0mm>*{}**@{-},
 <0.49mm,0.49mm>*{};<1.9mm,1.9mm>*{}**@{-},
 <-0.5mm,0.5mm>*{};<-1.9mm,1.9mm>*{}**@{-},
 <-2.3mm,2.3mm>*{\bullet};<-2.3mm,2.3mm>*{}**@{},
 <-1.8mm,2.8mm>*{};<0mm,4.9mm>*{}**@{-},
 <-2.8mm,2.9mm>*{};<-4.6mm,4.9mm>*{}**@{-},
   <0.49mm,0.49mm>*{};<2.7mm,2.3mm>*{^3}**@{},
   <-1.8mm,2.8mm>*{};<0.4mm,5.3mm>*{^2}**@{},
   <-2.8mm,2.9mm>*{};<-5.1mm,5.3mm>*{^1}**@{},
 \end{xy}
\ - \
\begin{xy}
 <0mm,0mm>*{\bullet};<0mm,0mm>*{}**@{},
 <0mm,-0.49mm>*{};<0mm,-3.0mm>*{}**@{-},
 <0.49mm,0.49mm>*{};<1.9mm,1.9mm>*{}**@{-},
 <-0.5mm,0.5mm>*{};<-1.9mm,1.9mm>*{}**@{-},
 <2.3mm,2.3mm>*{\bullet};<-2.3mm,2.3mm>*{}**@{},
 <1.8mm,2.8mm>*{};<0mm,4.9mm>*{}**@{-},
 <2.8mm,2.9mm>*{};<4.6mm,4.9mm>*{}**@{-},
   <0.49mm,0.49mm>*{};<-2.7mm,2.3mm>*{^1}**@{},
   <-1.8mm,2.8mm>*{};<0mm,5.3mm>*{^2}**@{},
   <-2.8mm,2.9mm>*{};<5.1mm,5.3mm>*{^3}**@{},
 \end{xy}
=0, \ \ \
 \begin{xy}
 <0mm,0mm>*{\bullet};<0mm,0mm>*{}**@{},
 <0mm,0.69mm>*{};<0mm,3.0mm>*{}**@{-},
 <0.39mm,-0.39mm>*{};<2.4mm,-2.4mm>*{}**@{-},
 <-0.35mm,-0.35mm>*{};<-1.9mm,-1.9mm>*{}**@{-},
 <-2.4mm,-2.4mm>*{\bullet};<-2.4mm,-2.4mm>*{}**@{},
 <-2.0mm,-2.8mm>*{};<0mm,-4.9mm>*{}**@{-},
 <-2.8mm,-2.9mm>*{};<-4.7mm,-4.9mm>*{}**@{-},
    <0.39mm,-0.39mm>*{};<3.3mm,-4.0mm>*{^3}**@{},
    <-2.0mm,-2.8mm>*{};<0.5mm,-6.7mm>*{^2}**@{},
    <-2.8mm,-2.9mm>*{};<-5.2mm,-6.7mm>*{^1}**@{},
 \end{xy}
\ - \
 \begin{xy}
 <0mm,0mm>*{\bullet};<0mm,0mm>*{}**@{},
 <0mm,0.69mm>*{};<0mm,3.0mm>*{}**@{-},
 <0.39mm,-0.39mm>*{};<2.4mm,-2.4mm>*{}**@{-},
 <-0.35mm,-0.35mm>*{};<-1.9mm,-1.9mm>*{}**@{-},
 <2.4mm,-2.4mm>*{\bullet};<-2.4mm,-2.4mm>*{}**@{},
 <2.0mm,-2.8mm>*{};<0mm,-4.9mm>*{}**@{-},
 <2.8mm,-2.9mm>*{};<4.7mm,-4.9mm>*{}**@{-},
    <0.39mm,-0.39mm>*{};<-3mm,-4.0mm>*{^1}**@{},
    <-2.0mm,-2.8mm>*{};<0mm,-6.7mm>*{^2}**@{},
    <-2.8mm,-2.9mm>*{};<5.2mm,-6.7mm>*{^3}**@{},
 \end{xy}
=0, \ \ \
 \begin{xy}
 <0mm,2.47mm>*{};<0mm,-0.5mm>*{}**@{-},
 <0.5mm,3.5mm>*{};<2.2mm,5.2mm>*{}**@{-},
 <-0.48mm,3.48mm>*{};<-2.2mm,5.2mm>*{}**@{-},
 <0mm,3mm>*{\bullet};<0mm,3mm>*{}**@{},
  <0mm,-0.8mm>*{\bullet};<0mm,-0.8mm>*{}**@{},
<0mm,-0.8mm>*{};<-2.2mm,-3.5mm>*{}**@{-},
 <0mm,-0.8mm>*{};<2.2mm,-3.5mm>*{}**@{-},
     <0.5mm,3.5mm>*{};<2.8mm,5.7mm>*{^2}**@{},
     <-0.48mm,3.48mm>*{};<-2.8mm,5.7mm>*{^1}**@{},
   <0mm,-0.8mm>*{};<-2.7mm,-5.2mm>*{^1}**@{},
   <0mm,-0.8mm>*{};<2.7mm,-5.2mm>*{^2}**@{},
\end{xy}
\ - \
\begin{xy}
 <0mm,0mm>*{\bullet};<0mm,0mm>*{}**@{},
 <0mm,-0.49mm>*{};<0mm,-3.0mm>*{}**@{-},
 <-0.5mm,0.5mm>*{};<-3mm,2mm>*{}**@{-},
 <-3mm,2mm>*{};<0mm,4mm>*{}**@{-},
 <0mm,4mm>*{\bullet};<-2.3mm,2.3mm>*{}**@{},
 <0mm,4mm>*{};<0mm,7.4mm>*{}**@{-},
<0mm,0mm>*{};<2.2mm,1.5mm>*{}**@{-},
 <6mm,0mm>*{\bullet};<0mm,0mm>*{}**@{},
 <6mm,4mm>*{};<3.8mm,2.5mm>*{}**@{-},
 <6mm,4mm>*{};<6mm,7.4mm>*{}**@{-},
 <6mm,4mm>*{\bullet};<-2.3mm,2.3mm>*{}**@{},
 <0mm,4mm>*{};<6mm,0mm>*{}**@{-},
<6mm,4mm>*{};<9mm,2mm>*{}**@{-}, <6mm,0mm>*{};<9mm,2mm>*{}**@{-},
<6mm,0mm>*{};<6mm,-3mm>*{}**@{-},
   <-1.8mm,2.8mm>*{};<0mm,7.8mm>*{^1}**@{},
   <-2.8mm,2.9mm>*{};<0mm,-4.3mm>*{_1}**@{},
<-1.8mm,2.8mm>*{};<6mm,7.8mm>*{^2}**@{},
   <-2.8mm,2.9mm>*{};<6mm,-4.3mm>*{_2}**@{},
 \end{xy}
 =0,
 $$
hold. Using this fact it is an easy and straightforward
calculation  to check that $d^2=0$. We omit the details. (In fact
we shall show below that $d$ is essentially a graph  encoding of
the Gerstenhaber-Schack differential $d_{GS}$ so  this calculation
is essentially identical to the one which establishes
$d_{GS}^2=0$.)

\bip

(ii)  We begin our proof of Lemma~\ref{keylemma}(ii) with the
following

\sip

{\bf Claim 1}. {\em The natural projection $p: (\hsB, d)\rar \ab$
is a quasi-isomorphism.}

\sip

Indeed, the dg properad $(\hsB, d)$ has a natural increasing and
bounded above filtration\footnote{One might prove
Claim~1 using another filtration of $\hsB$ by the number of vertices and Fact~\ref{delta_0} provided one assumes (without
any losses) that $\hsB$ is completed with respect to this filtration.}, $\{\mathsf F_{-p}\hsB\}_{p \geq 0}$, with
$\mathsf F_{-p}\hsB$ being the span of equivalence classes of
graphs which admit a representative in $\ab_\infty^{(\geq p)}$.
As the differential $d$ is
connected and preserves
the induced path gradation, the associated spectral sequence $(\sE_r, \sd_r)$ converges to $H^\bullet(\hsB, d)$.  The $0$th
term $(\sE_0,\sd_0)$  has the differential given on generators by
\Beq\label{d_1} \sd_0
\begin{xy}
 <0mm,-0.55mm>*{};<0mm,-2.5mm>*{}**@{-},
 <0.5mm,0.5mm>*{};<2.2mm,2.2mm>*{}**@{-},
 <-0.48mm,0.48mm>*{};<-2.2mm,2.2mm>*{}**@{-},
 <0mm,0mm>*{\bullet};<0mm,0mm>*{}**@{},
 <0mm,-0.55mm>*{};<0mm,-3.8mm>*{_1}**@{},
 <0.5mm,0.5mm>*{};<2.7mm,2.8mm>*{^2}**@{},
 <-0.48mm,0.48mm>*{};<-2.7mm,2.8mm>*{^1}**@{},
 \end{xy}
 =0
 \ \ , \  \ \ \
 \sd_0
 \begin{xy}
 <0mm,0.66mm>*{};<0mm,3mm>*{}**@{-},
 <0.39mm,-0.39mm>*{};<2.2mm,-2.2mm>*{}**@{-},
 <-0.35mm,-0.35mm>*{};<-2.2mm,-2.2mm>*{}**@{-},
 <0mm,0mm>*{\bullet};<0mm,0mm>*{}**@{},
   <0mm,0.66mm>*{};<0mm,3.4mm>*{^1}**@{},
   <0.39mm,-0.39mm>*{};<2.9mm,-4mm>*{^2}**@{},
   <-0.35mm,-0.35mm>*{};<-2.8mm,-4mm>*{^1}**@{},
\end{xy}
=0\ \ , \ \ \
 \sd_0
 \begin{xy}
 <0.39mm,-0.39mm>*{};<2.2mm,-2.2mm>*{}**@{-},
 <-0.35mm,-0.35mm>*{};<-2.2mm,-2.2mm>*{}**@{-},
 <0.39mm,0.39mm>*{};<2.2mm,2.2mm>*{}**@{-},
 <-0.35mm,0.35mm>*{};<-2.2mm,2.2mm>*{}**@{-},
 <0mm,0mm>*{\bullet};<0mm,0mm>*{}**@{},
   <0.39mm,-0.39mm>*{};<2.9mm,-4mm>*{^2}**@{},
   <-0.35mm,-0.35mm>*{};<-2.8mm,-4mm>*{^1}**@{},
    <0mm,0.66mm>*{};<2.8mm,3.4mm>*{^2}**@{},
<0mm,0.66mm>*{};<-2.8mm,3.4mm>*{^1}**@{},
\end{xy}
=
 \begin{xy}
 <0mm,2.47mm>*{};<0mm,-0.5mm>*{}**@{-},
 <0.5mm,3.5mm>*{};<2.2mm,5.2mm>*{}**@{-},
 <-0.48mm,3.48mm>*{};<-2.2mm,5.2mm>*{}**@{-},
 <0mm,3mm>*{\bullet};<0mm,3mm>*{}**@{},
  <0mm,-0.8mm>*{\bullet};<0mm,-0.8mm>*{}**@{},
<0mm,-0.8mm>*{};<-2.2mm,-3.5mm>*{}**@{-},
 <0mm,-0.8mm>*{};<2.2mm,-3.5mm>*{}**@{-},
     <0.5mm,3.5mm>*{};<2.8mm,5.7mm>*{^2}**@{},
     <-0.48mm,3.48mm>*{};<-2.8mm,5.7mm>*{^1}**@{},
   <0mm,-0.8mm>*{};<-2.7mm,-5.2mm>*{^1}**@{},
   <0mm,-0.8mm>*{};<2.7mm,-5.2mm>*{^2}**@{},
\end{xy}\ ,
\Eeq and, for all other generators with $m+n\geq 4$,
\Beq\label{d_11} \sd_0
\begin{xy}
 <0mm,0mm>*{\bullet};<0mm,0mm>*{}**@{},
 <0mm,0mm>*{};<-8mm,5mm>*{}**@{-},
 <0mm,0mm>*{};<-4.5mm,5mm>*{}**@{-},
 <0mm,0mm>*{};<-1mm,5mm>*{\ldots}**@{},
 <0mm,0mm>*{};<4.5mm,5mm>*{}**@{-},
 <0mm,0mm>*{};<8mm,5mm>*{}**@{-},
   <0mm,0mm>*{};<-8.5mm,5.5mm>*{^1}**@{},
   <0mm,0mm>*{};<-5mm,5.5mm>*{^2}**@{},
   <0mm,0mm>*{};<4.5mm,5.5mm>*{^{m\hspace{-0.5mm}-\hspace{-0.5mm}1}}**@{},
   <0mm,0mm>*{};<9.0mm,5.5mm>*{^m}**@{},
 <0mm,0mm>*{};<-8mm,-5mm>*{}**@{-},
 <0mm,0mm>*{};<-4.5mm,-5mm>*{}**@{-},
 <0mm,0mm>*{};<-1mm,-5mm>*{\ldots}**@{},
 <0mm,0mm>*{};<4.5mm,-5mm>*{}**@{-},
 <0mm,0mm>*{};<8mm,-5mm>*{}**@{-},
   <0mm,0mm>*{};<-8.5mm,-6.9mm>*{^1}**@{},
   <0mm,0mm>*{};<-5mm,-6.9mm>*{^2}**@{},
   <0mm,0mm>*{};<4.5mm,-6.9mm>*{^{n\hspace{-0.5mm}-\hspace{-0.5mm}1}}**@{},
   <0mm,0mm>*{};<9.0mm,-6.9mm>*{^n}**@{},
 \end{xy}
 =
 \sum_{i=0}^{n-2}(-1)^i
 \begin{xy}
 <0mm,0mm>*{\bullet};<0mm,0mm>*{}**@{},
 <0mm,0mm>*{};<-8mm,5mm>*{}**@{-},
 <0mm,0mm>*{};<-3.5mm,5mm>*{}**@{-},
 <0mm,0mm>*{};<-6mm,5mm>*{..}**@{},
 <0mm,0mm>*{};<0mm,5mm>*{}**@{-},
  <0mm,5mm>*{\bullet};
  <0mm,5mm>*{};<-2.5mm,8mm>*{}**@{-},
  <0mm,5mm>*{};<2.5mm,8mm>*{}**@{-},
  <0mm,5mm>*{};<-3mm,9mm>*{^{i\hspace{-0.2mm}+\hspace{-0.5mm}1}}**@{},
  <0mm,5mm>*{};<3mm,9mm>*{^{i\hspace{-0.2mm}+\hspace{-0.5mm}2}}**@{},
<0mm,0mm>*{};<8mm,5mm>*{}**@{-},
<0mm,0mm>*{};<3.5mm,5mm>*{}**@{-},
 <0mm,0mm>*{};<6mm,5mm>*{..}**@{},
   <0mm,0mm>*{};<-8.5mm,5.5mm>*{^1}**@{},
   <0mm,0mm>*{};<-4mm,5.5mm>*{^i}**@{},
   <0mm,0mm>*{};<9.0mm,5.5mm>*{^m}**@{},
 <0mm,0mm>*{};<-8mm,-5mm>*{}**@{-},
 <0mm,0mm>*{};<-4.5mm,-5mm>*{}**@{-},
 <0mm,0mm>*{};<-1mm,-5mm>*{\ldots}**@{},
 <0mm,0mm>*{};<4.5mm,-5mm>*{}**@{-},
 <0mm,0mm>*{};<8mm,-5mm>*{}**@{-},
   <0mm,0mm>*{};<-8.5mm,-6.9mm>*{^1}**@{},
   <0mm,0mm>*{};<-5mm,-6.9mm>*{^2}**@{},
   <0mm,0mm>*{};<4.5mm,-6.9mm>*{^{n\hspace{-0.5mm}-\hspace{-0.5mm}1}}**@{},
   <0mm,0mm>*{};<9.0mm,-6.9mm>*{^n}**@{},
 \end{xy}
 +
\sum_{i=0}^{n-2}(-1)^i
 \begin{xy}
 <0mm,0mm>*{\bullet};<0mm,0mm>*{}**@{},
 <0mm,0mm>*{};<-8mm,-5mm>*{}**@{-},
 <0mm,0mm>*{};<-3.5mm,-5mm>*{}**@{-},
 <0mm,0mm>*{};<-6mm,-5mm>*{..}**@{},
 <0mm,0mm>*{};<0mm,-5mm>*{}**@{-},
  <0mm,-5mm>*{\bullet};
  <0mm,-5mm>*{};<-2.5mm,-8mm>*{}**@{-},
  <0mm,-5mm>*{};<2.5mm,-8mm>*{}**@{-},
  <0mm,-5mm>*{};<-3mm,-10mm>*{^{i\hspace{-0.2mm}+\hspace{-0.5mm}1}}**@{},
  <0mm,-5mm>*{};<3mm,-10mm>*{^{i\hspace{-0.2mm}+\hspace{-0.5mm}2}}**@{},
<0mm,0mm>*{};<8mm,-5mm>*{}**@{-},
<0mm,0mm>*{};<3.5mm,-5mm>*{}**@{-},
 <0mm,0mm>*{};<6mm,-5mm>*{..}**@{},
   <0mm,0mm>*{};<-8.5mm,-6.9mm>*{^1}**@{},
   <0mm,0mm>*{};<-4mm,-6.9mm>*{^i}**@{},
   <0mm,0mm>*{};<9.0mm,-6.9mm>*{^n}**@{},
 <0mm,0mm>*{};<-8mm,5mm>*{}**@{-},
 <0mm,0mm>*{};<-4.5mm,5mm>*{}**@{-},
 <0mm,0mm>*{};<-1mm,5mm>*{\ldots}**@{},
 <0mm,0mm>*{};<4.5mm,5mm>*{}**@{-},
 <0mm,0mm>*{};<8mm,5mm>*{}**@{-},
   <0mm,0mm>*{};<-8.5mm,5.5mm>*{^1}**@{},
   <0mm,0mm>*{};<-5mm,5.5mm>*{^2}**@{},
   <0mm,0mm>*{};<4.5mm,5.5mm>*{^{m\hspace{-0.5mm}-\hspace{-0.5mm}1}}**@{},
   <0mm,0mm>*{};<9.0mm,5.5mm>*{^m}**@{},
 \end{xy}\\
\Eeq We want to compute homology, $\sE_1=H^\bullet(\sE_0, \sd_0)$,
of this complex and show that $\sE_1\simeq \ab$. For this purpose
consider a 3-step filtration, $0\subset F_{-2}\subset F_{-1}\subset F_{0}=\sE_0$, of
the complex $(\sE_0, \sd_0)$ with
$$
F_{-1}:=\mbox{span}\left\langle
\begin{xy}
 <0mm,-0.55mm>*{};<0mm,-2.5mm>*{}**@{-},
 <0.5mm,0.5mm>*{};<2.2mm,2.2mm>*{}**@{-},
 <-0.48mm,0.48mm>*{};<-2.2mm,2.2mm>*{}**@{-},
 <0mm,0mm>*{\bullet};<0mm,0mm>*{}**@{},
 \end{xy}
 \ , \
 \begin{xy}
 <0mm,0.66mm>*{};<0mm,3mm>*{}**@{-},
 <0.39mm,-0.39mm>*{};<2.2mm,-2.2mm>*{}**@{-},
 <-0.35mm,-0.35mm>*{};<-2.2mm,-2.2mm>*{}**@{-},
 <0mm,0mm>*{\bullet};<0mm,0mm>*{}**@{},
\end{xy}
\  ,  \
 \begin{xy}
 <0.39mm,-0.39mm>*{};<2.2mm,-2.2mm>*{}**@{-},
 <-0.35mm,-0.35mm>*{};<-2.2mm,-2.2mm>*{}**@{-},
 <0.39mm,0.39mm>*{};<2.2mm,2.2mm>*{}**@{-},
 <-0.35mm,0.35mm>*{};<-2.2mm,2.2mm>*{}**@{-},
 <0mm,0mm>*{\bullet};<0mm,0mm>*{}**@{},
\end{xy}
\right\rangle, \ \ \
F_{-2}:=\mbox{span}\left\langle
\begin{xy}
 <0mm,-0.55mm>*{};<0mm,-2.5mm>*{}**@{-},
 <0.5mm,0.5mm>*{};<2.2mm,2.2mm>*{}**@{-},
 <-0.48mm,0.48mm>*{};<-2.2mm,2.2mm>*{}**@{-},
 <0mm,0mm>*{\bullet};<0mm,0mm>*{}**@{},
 \end{xy}
 \ , \
 \begin{xy}
 <0mm,0.66mm>*{};<0mm,3mm>*{}**@{-},
 <0.39mm,-0.39mm>*{};<2.2mm,-2.2mm>*{}**@{-},
 <-0.35mm,-0.35mm>*{};<-2.2mm,-2.2mm>*{}**@{-},
 <0mm,0mm>*{\bullet};<0mm,0mm>*{}**@{},
\end{xy}
\right\rangle,
$$
and let $(\cE_r, \p_r)$ be the associated spectral sequence. The
differential $\p_0$ is zero on the generators of $F_{-1}$ and is
equal to $\sd_0$ on all the other generators. Thus, modulo shifts
of gradings, actions of finite groups and tensor products by
trivial (i.e. with zero differential) complexes,
 the complex $(\cE_0, \p_0)$ is isomorphic to the tensor product of
two isomorphic operadic complexes (one with ``time" flow reversed
upside down relative to another) which were studied on page 40 of
\cite{MarklMerkulovShadrin06} and which have the differential (in
notations of that paper) given by
$$
d_1\begin{xy}
 <0mm,0mm>*{\bullet};<0mm,0mm>*{}**@{},
 <0mm,0mm>*{};<-8mm,-5mm>*{}**@{-},
 <0mm,0mm>*{};<-4.5mm,-5mm>*{}**@{-},
 <0mm,0mm>*{};<0mm,-4mm>*{\ldots}**@{},
 <0mm,0mm>*{};<4.5mm,-5mm>*{}**@{-},
 <0mm,0mm>*{};<8mm,-5mm>*{}**@{-},
   <0mm,0mm>*{};<-11mm,-7.9mm>*{^{1}}**@{},
   <0mm,0mm>*{};<-4mm,-7.9mm>*{^{2}}**@{},
   <0mm,0mm>*{};<10.0mm,-7.9mm>*{^n}**@{},
 <0mm,0mm>*{};<0mm,5mm>*{}**@{-},
 \end{xy}=
 \sum_{i=0}^{n-2}(-1)^{i+1}
\begin{xy}
 <0mm,0mm>*{\bullet};<0mm,0mm>*{}**@{},
 <0mm,0mm>*{};<-8mm,-5mm>*{}**@{-},
 <0mm,0mm>*{};<-4mm,-5mm>*{}**@{-},
 <0mm,0mm>*{};<-5.6mm,-5mm>*{..}**@{},
 <0mm,0mm>*{};<2.6mm,-5mm>*{..}**@{},
 <0mm,0mm>*{};<4.5mm,-5mm>*{}**@{-},
 <0mm,0mm>*{};<-0.5mm,-5mm>*{}**@{-},
 <0mm,0mm>*{};<8mm,-5mm>*{}**@{-},
   <0mm,0mm>*{};<-11mm,-7.9mm>*{^1}**@{},
   <0mm,0mm>*{};<-4.5mm,-7.9mm>*{^i}**@{},
   <0mm,0mm>*{};<10.0mm,-7.9mm>*{^n}**@{},
 <0mm,0mm>*{};<0mm,5mm>*{}**@{-},
  <-0.5mm,-5.4mm>*{\bullet};
<-0.7mm,-5.4mm>*{};<-3.6mm,-9mm>*{}**@{-},
<-0.7mm,-5.4mm>*{};<3mm,-9mm>*{}**@{-},
<-0.7mm,-5.4mm>*{};<-4mm,-11mm>*{_{i+1}}**@{},
<-0.7mm,-5.4mm>*{};<4mm,-11mm>*{_{i+2}}**@{},
 \end{xy}
$$
It is shown in \cite{MarklMerkulovShadrin06} that the cohomology
of this complex is concentrated in degree $0$ and is isomorphic to
the operad of associative algebras. In our context this result
immediately implies that $(\cE_1, \p_1)$ is isomorphic to $F_{-1}$
with the differential $\p_1$ given on generators by (\ref{d_1}).
Its cohomology is obviously concentrated in degree 0 (and is
equal, in fact, to  the properad, $\frac{1}{2}\sB$, of
infinitesimal bialgebras). Hence
 $H^\bullet(\hsB, d)$ is
concentrated in degree $0$ proving  Claim~1.

\bip {\bf Claim 2}. {\em The natural projection ${\pi}: (\hsB,
\hp)\rar \ab$ is a quasi-isomorphism.}

\sip Indeed, the defined above filtration, $\{\mathsf
F_{-p}\hsB\}_{p \geq 0}$, by the number of vertices is also
compatible with the differential $\hp$. Let $({E}_r, {d}_r)$ be
the associated spectral sequence. Its first nontrivial term,
$({E}_1, {d}_1)$ is, by Fact~\ref{delta_0}, isomorphic to the
complex $(\sE_1, \sd_1)$ above. Hence we can apply the same
reasoning as in the proof of Claim~1.

\bip

{\bf Claim 3}. {\em The exits a morphism of dg properads $\Phi$
making the diagram
\[
 \xymatrix{
  & (\hsB, d)
  \ar[d]^{p} \\
 (\ab_\infty, \p)\ar[ur]^\Phi\ar[r]_{\pi} &
 (\ab, 0)
 }
\]
commutative}. \sip

Since $\ab$ is a properad concentrated in degree $0$, the map $p$
is surjective. Since $p$ is a quasi-isomorphism by Claim~2, it is
an acyclic fibration in the model category of properads (see
Appendix~\ref{model category}). By Corollary~\ref{quasi-free ->
cofibrant}, $\ab_\infty$ is a cofibrant properad. Finally, the
morphism $\Phi$ is given by the left lifting property in the model
category of properads. Hence, the existence of $\Phi$ is clear but
we need to make it more precise. We construct it as follows and
refine it in Claim~$4$.\\

As $\ab_\infty=\cF(\Co)$ is a free properad, a morphism $\Phi$ is
completely determined by its values on the generating
$(m,n)$-corollas which span the vector space $\Co$, and one can
construct $\Phi$ by a simple induction \footnote{this induction is
a straightforward analogue of the Whitehead lifting trick in the
theory of $CW$-complexes in algebraic topology.}
 on the degree, $r:=m+n-3\geq 0$, of such corollas. For $r=0$ we set $\Phi$ to be identity, i.e.
 $$
\Phi\left(\begin{xy}
 <0mm,-0.55mm>*{};<0mm,-2.5mm>*{}**@{-},
 <0.5mm,0.5mm>*{};<2.2mm,2.2mm>*{}**@{-},
 <-0.48mm,0.48mm>*{};<-2.2mm,2.2mm>*{}**@{-},
 <0mm,0mm>*{\bullet};<0mm,0mm>*{}**@{},
 <0mm,-0.55mm>*{};<0mm,-3.8mm>*{_1}**@{},
 <0.5mm,0.5mm>*{};<2.7mm,2.8mm>*{^2}**@{},
 <-0.48mm,0.48mm>*{};<-2.7mm,2.8mm>*{^1}**@{},
 \end{xy}\right)
=
\begin{xy}
 <0mm,-0.55mm>*{};<0mm,-2.5mm>*{}**@{-},
 <0.5mm,0.5mm>*{};<2.2mm,2.2mm>*{}**@{-},
 <-0.48mm,0.48mm>*{};<-2.2mm,2.2mm>*{}**@{-},
 <0mm,0mm>*{\bullet};<0mm,0mm>*{}**@{},
 <0mm,-0.55mm>*{};<0mm,-3.8mm>*{_1}**@{},
 <0.5mm,0.5mm>*{};<2.7mm,2.8mm>*{^2}**@{},
 <-0.48mm,0.48mm>*{};<-2.7mm,2.8mm>*{^1}**@{},
 \end{xy}
 \ \ , \ \
\Phi\left(
 \begin{xy}
 <0mm,0.66mm>*{};<0mm,3mm>*{}**@{-},
 <0.39mm,-0.39mm>*{};<2.2mm,-2.2mm>*{}**@{-},
 <-0.35mm,-0.35mm>*{};<-2.2mm,-2.2mm>*{}**@{-},
 <0mm,0mm>*{\bullet};<0mm,0mm>*{}**@{},
   <0mm,0.66mm>*{};<0mm,3.4mm>*{^1}**@{},
   <0.39mm,-0.39mm>*{};<2.9mm,-4mm>*{^2}**@{},
   <-0.35mm,-0.35mm>*{};<-2.8mm,-4mm>*{^1}**@{},
\end{xy}
\right)=
 \begin{xy}
 <0mm,0.66mm>*{};<0mm,3mm>*{}**@{-},
 <0.39mm,-0.39mm>*{};<2.2mm,-2.2mm>*{}**@{-},
 <-0.35mm,-0.35mm>*{};<-2.2mm,-2.2mm>*{}**@{-},
 <0mm,0mm>*{\bullet};<0mm,0mm>*{}**@{},
   <0mm,0.66mm>*{};<0mm,3.4mm>*{^1}**@{},
   <0.39mm,-0.39mm>*{};<2.9mm,-4mm>*{^2}**@{},
   <-0.35mm,-0.35mm>*{};<-2.8mm,-4mm>*{^1}**@{},
\end{xy}.
$$
Assume we constructed values of $\Phi$ on all corollas of degree
$r\leq N$. Let $e$ be a generating corolla of non-zero weight
$r=N+1$. Note that $\delta e$ is a linear combination of graphs
whose vertices are decorated by corollas of weight $\leq N$ (as
differential $\delta$ has degee $-1$). Then, by induction,
$\Phi(\delta e)$ is a well-defined element in $\hsB$. As
$\pi(e)=0$, the element,
$$
\Phi(\delta e)
$$
is a closed element in $\hsB$ which projects under $p$ to zero. By
Claim 1, the surjection $p$ is a quasi-isomorphism. Hence this
element is exact and there exists $\mathfrak e\in \hsB$ such that
$$
d {\mathfrak e}= \Phi(\delta e).
$$
We set $\Phi(e):= {\mathfrak e}$ completing thereby inductive
construction of $\Phi$.

\bip {\bf Claim 4}. {\em A morphism $\Phi$ can be chosen so that}
$$
\Phi\left(
\begin{xy}
 <0mm,0mm>*{\bullet};<0mm,0mm>*{}**@{},
 <0mm,0mm>*{};<-8mm,5mm>*{}**@{-},
 <0mm,0mm>*{};<-4.5mm,5mm>*{}**@{-},
 <0mm,0mm>*{};<-1mm,5mm>*{\ldots}**@{},
 <0mm,0mm>*{};<4.5mm,5mm>*{}**@{-},
 <0mm,0mm>*{};<8mm,5mm>*{}**@{-},
   <0mm,0mm>*{};<-8.5mm,5.5mm>*{^1}**@{},
   <0mm,0mm>*{};<-5mm,5.5mm>*{^2}**@{},
   <0mm,0mm>*{};<4.5mm,5.5mm>*{^{m\hspace{-0.5mm}-\hspace{-0.5mm}1}}**@{},
   <0mm,0mm>*{};<9.0mm,5.5mm>*{^m}**@{},
 <0mm,0mm>*{};<-8mm,-5mm>*{}**@{-},
 <0mm,0mm>*{};<-4.5mm,-5mm>*{}**@{-},
 <0mm,0mm>*{};<-1mm,-5mm>*{\ldots}**@{},
 <0mm,0mm>*{};<4.5mm,-5mm>*{}**@{-},
 <0mm,0mm>*{};<8mm,-5mm>*{}**@{-},
   <0mm,0mm>*{};<-8.5mm,-6.9mm>*{^1}**@{},
   <0mm,0mm>*{};<-5mm,-6.9mm>*{^2}**@{},
   <0mm,0mm>*{};<4.5mm,-6.9mm>*{^{n\hspace{-0.5mm}-\hspace{-0.5mm}1}}**@{},
   <0mm,0mm>*{};<9.0mm,-6.9mm>*{^n}**@{},
 \end{xy}
\right) =
\begin{xy}
 <0mm,0mm>*{\bullet};<0mm,0mm>*{}**@{},
 <0mm,0mm>*{};<-8mm,5mm>*{}**@{-},
 <0mm,0mm>*{};<-4.5mm,5mm>*{}**@{-},
 <0mm,0mm>*{};<-1mm,5mm>*{\ldots}**@{},
 <0mm,0mm>*{};<4.5mm,5mm>*{}**@{-},
 <0mm,0mm>*{};<8mm,5mm>*{}**@{-},
   <0mm,0mm>*{};<-8.5mm,5.5mm>*{^1}**@{},
   <0mm,0mm>*{};<-5mm,5.5mm>*{^2}**@{},
   <0mm,0mm>*{};<4.5mm,5.5mm>*{^{m\hspace{-0.5mm}-\hspace{-0.5mm}1}}**@{},
   <0mm,0mm>*{};<9.0mm,5.5mm>*{^m}**@{},
 <0mm,0mm>*{};<-8mm,-5mm>*{}**@{-},
 <0mm,0mm>*{};<-4.5mm,-5mm>*{}**@{-},
 <0mm,0mm>*{};<-1mm,-5mm>*{\ldots}**@{},
 <0mm,0mm>*{};<4.5mm,-5mm>*{}**@{-},
 <0mm,0mm>*{};<8mm,-5mm>*{}**@{-},
   <0mm,0mm>*{};<-8.5mm,-6.9mm>*{^1}**@{},
   <0mm,0mm>*{};<-5mm,-6.9mm>*{^2}**@{},
   <0mm,0mm>*{};<4.5mm,-6.9mm>*{^{n\hspace{-0.5mm}-\hspace{-0.5mm}1}}**@{},
   <0mm,0mm>*{};<9.0mm,-6.9mm>*{^n}**@{},
 \end{xy}
+ \mbox{terms with $\geq 2$ number of vertices}.
$$

\sip

Indeed, the differential $d$ in $\hsB$ has the form,
$$
d=\sd_1 + d_{rest},
$$
where $\sd_1$ is the quadratic {\em differential}\, in $\hsB$
defined by (\ref{d_11}) and   the part $d_{part}$ corresponds to
graphs lying in $F_{-3}\sB$. We shall prove Claim~4 by induction
on the degree $r=m+n-3$ of the generating $(m,n)$-corollas in
$\ab_\infty$ (cf.\ proof of Theorem~43 in \cite{Markl06}). For
$r=0$ the Claim is true. Assume we have already constructed $\Phi$
such that the claim is true for values of $\Phi$ on corollas with
non-zero degree $\leq N$ and consider a generating corolla, $e$,
of degree $N+1$. The value, ${\mathfrak e}:= \Phi(e)$, is a
solution of the equation, \Beq\label{rest} \sd_1 {\mathfrak e} +
d_{rest} {\mathfrak e} = \Phi(\p_0 e) + \Phi(\p_{pert}e). \Eeq Let
$\pi_1$ and $\pi_2$ denote projections in $\hsB$ to the subspaces
spanned by equivalence classes of graphs with $1$ and,
respectively, 2 vertices. Then
 equation (\ref{rest}) implies,
$$
\pi_2\circ \sd_1 ({\mathfrak e})= \sd_1\circ \pi_1(\mathfrak e)  =
\pi_2\circ\Phi(\p_0 e),
$$ as both $d_{rest} {\mathfrak e}$ and $ \Phi(\p_{pert}e)$ are
spanned by graphs lying in $F_{-3}\sB$. Using now the explicit
form for the differential $\p_0$ (given, e.g., by formula (14) in
\cite{Markl06}) and the induction assumption we immediately
conclude that
$$
\pi_1 ({\mathfrak e})=\begin{xy}
 <0mm,0mm>*{\bullet};<0mm,0mm>*{}**@{},
 <0mm,0mm>*{};<-8mm,5mm>*{}**@{-},
 <0mm,0mm>*{};<-4.5mm,5mm>*{}**@{-},
 <0mm,0mm>*{};<-1mm,5mm>*{\ldots}**@{},
 <0mm,0mm>*{};<4.5mm,5mm>*{}**@{-},
 <0mm,0mm>*{};<8mm,5mm>*{}**@{-},
   <0mm,0mm>*{};<-8.5mm,5.5mm>*{^1}**@{},
   <0mm,0mm>*{};<-5mm,5.5mm>*{^2}**@{},
   <0mm,0mm>*{};<4.5mm,5.5mm>*{^{m\hspace{-0.5mm}-\hspace{-0.5mm}1}}**@{},
   <0mm,0mm>*{};<9.0mm,5.5mm>*{^m}**@{},
 <0mm,0mm>*{};<-8mm,-5mm>*{}**@{-},
 <0mm,0mm>*{};<-4.5mm,-5mm>*{}**@{-},
 <0mm,0mm>*{};<-1mm,-5mm>*{\ldots}**@{},
 <0mm,0mm>*{};<4.5mm,-5mm>*{}**@{-},
 <0mm,0mm>*{};<8mm,-5mm>*{}**@{-},
   <0mm,0mm>*{};<-8.5mm,-6.9mm>*{^1}**@{},
   <0mm,0mm>*{};<-5mm,-6.9mm>*{^2}**@{},
   <0mm,0mm>*{};<4.5mm,-6.9mm>*{^{n\hspace{-0.5mm}-\hspace{-0.5mm}1}}**@{},
   <0mm,0mm>*{};<9.0mm,-6.9mm>*{^n}**@{},
 \end{xy}
 $$
 completing the proof of Claim~4.

\bip {\bf Claim 5}. {\em The morphism $\Phi$ induces a dg
isomorphism  $(\hsB, \hp) \rar (\hsB, d)$.}

\sip

Indeed, $\Phi$ sends the ideal $I$ to zero. Since $\Phi$ respects
differentials, it sends the ideal $(I, \p I)$ to zero as well and
hence induces, by Claims 3 and 4, a required isomorphism. This
completes proof of Lemma~\ref{keylemma}.
\end{proo}

\bip Now we continue with the proof of Theorem~\ref{GS}. The
differential $Q_1^\ga$ in the graded vector space
$\fg_{GS}(\Qo)=\oplus_{m,n}s^{2-m-n}\Qo$ is completely determined
by the quotient differential, $\hp$, of the full differential $\p$
in $\ab_\infty$. By Lemma~\ref{keylemma}, this quotient
differential is given, up to automorphims, by formulae
(\ref{d-1})-(\ref{d-3}). The proof of the Theorem~\ref{GS} is
completed.
\end{proo}

As a direct corollary, we have

\begin{cor}
The Gerstenhaber-Schack bicomplex of an associative bialgebra $X$
is a homotopy non-symmetric properad and a twisted
$L_\infty$-algebra whose Maurer-Cartan elements are deformations
of the first structure.
\end{cor}

The homotopy non-symmetric properad structure induces, on this
chain complex, (homotopy) LR-operations which play the same role
than the non-symmetric braces for Hochschild cochain complex. They
are are expected to be used in the proof of a Deligne conjecture
for associative bialgebras.

\subsection{Twisted $L_\infty$-algebras and dg prop(erad)s}
 For any quasi-free prop(erad)
$(\cP_\infty=\cF(s^{-1}\cC), \p_\cP)$ and any prop(erad) $\Qo$
there exists, in accordance with Theorem~\ref{mor},
 a canonical $L_\infty$ structure, $Q$,
on the graded vector space $s^{-1}\Hm(\Co,\Qo)$ whose
Maurer-Cartan elements are in one-to-one correspondence with
representation of $P_\infty$  in $\Qo$. If $\ga$ is any particular
representation of $P_\infty$, then the associated twisted
$L_\infty$-algebra, $Q^\ga$, describes deformation theory of $\ga$
within the class of representation of $\cP_\infty$ (see \S 7.2).
Remarkably, {\em  there always exists a quasi-free prop(erad)
$(\cP^{(2)}_\infty, \p)$ whose representations in $\Qo$ are in
one-to-one correspondence with pairs, $(\ga, \Ga)$, where $\ga$ is
a representation of $\cP_\infty$ on $\Qo$ and $\Ga$ is a MC
element in $Q^\ga$}. Thus the dg prop(erad) $\cP^{(2)}_\infty$
gives a complete description of the deformation theory of a
generic representation of $\cP_\infty$. In fact, this
constructions can be obviously iterated giving rise to quasi-free
prop(erad)s $\cP^{(3)}_\infty$, $\cP^{(4)}_\infty$  etc.

\sip

By definition, $\cP^{(2)}_\infty$ is a free prop(erad) on the
$\bS$-bimodule $s^{-1}\cC\oplus s^{-1}\cC$ but the differential,
$\p$, in $\cP^{(2)}_\infty$ is {\em not}\, a direct sum,
$\p_\cP\oplus \p_\cP$, of differentials in $\cP_\infty$. We
illustrate the above claim in the case of $\cP=\Ass$, the operad
of associative algebras, before giving the general definition.

\sip

Let $\Ass^{(2)}_\infty$ be a quasi-free operad generated by an
$\bS$-module,
$$
s^{n-2}\bK[\Sigma_n]\, \oplus  s^{n-2}\bK[\Sigma_n] \simeq  \mbox{span}
\left(
\begin{xy}
 <0mm,0mm>*{\bullet};<0mm,0mm>*{}**@{},
 <0mm,0mm>*{};<-8mm,-5mm>*{}**@{-},
 <0mm,0mm>*{};<-4.5mm,-5mm>*{}**@{-},
 <0mm,0mm>*{};<0mm,-4mm>*{\ldots}**@{},
 <0mm,0mm>*{};<4.5mm,-5mm>*{}**@{-},
 <0mm,0mm>*{};<8mm,-5mm>*{}**@{-},
   <0mm,0mm>*{};<-11mm,-7.9mm>*{^{\sigma(1)}}**@{},
   <0mm,0mm>*{};<-4mm,-7.9mm>*{^{\sigma(2)}}**@{},
   <0mm,0mm>*{};<10.0mm,-7.9mm>*{^{\sigma(n)}}**@{},
 <0mm,0mm>*{};<0mm,5mm>*{}**@{-},
 \end{xy}\ ,
 \
\begin{xy}
 <0mm,0mm>*{\blacktriangledown};<0mm,0mm>*{}**@{},
 <0mm,0mm>*{};<-8mm,-5mm>*{}**@{-},
 <0mm,0mm>*{};<-4.5mm,-5mm>*{}**@{-},
 <0mm,0mm>*{};<0mm,-4mm>*{\ldots}**@{},
 <0mm,0mm>*{};<4.5mm,-5mm>*{}**@{-},
 <0mm,0mm>*{};<8mm,-5mm>*{}**@{-},
   <0mm,0mm>*{};<-11mm,-7.9mm>*{^{\sigma(1)}}**@{},
   <0mm,0mm>*{};<-4mm,-7.9mm>*{^{\sigma(2)}}**@{},
   <0mm,0mm>*{};<10.0mm,-7.9mm>*{^{\sigma(n)}}**@{},
 <0mm,0mm>*{};<0mm,5mm>*{}**@{-},
 \end{xy}
\right)_{\sigma\in \bS_n}
$$
and equipped with a differential given on generators by,
$$
\p\left( \begin{xy}
 <0mm,0mm>*{\bullet};<0mm,0mm>*{}**@{},
 <0mm,0mm>*{};<0mm,5mm>*{}**@{-},
 <0mm,0mm>*{};<-8mm,-5mm>*{}**@{-},
 <0mm,0mm>*{};<-4.5mm,-5mm>*{}**@{-},
 <0mm,0mm>*{};<-1mm,-5mm>*{\ldots}**@{},
 <0mm,0mm>*{};<4.5mm,-5mm>*{}**@{-},
 <0mm,0mm>*{};<8mm,-5mm>*{}**@{-},
   <0mm,0mm>*{};<-8.5mm,-6.9mm>*{^1}**@{},
   <0mm,0mm>*{};<-5mm,-6.9mm>*{^2}**@{},
   <0mm,0mm>*{};<4.5mm,-6.9mm>*{^{n\hspace{-0.5mm}-\hspace{-0.5mm}1}}**@{},
   <0mm,0mm>*{};<9.0mm,-6.9mm>*{^n}**@{},
 \end{xy}\right)=
\sum_{k=0}^{n-1}\sum_{l=1}^{n-k} (-1)^{(l-1)(n-k-l)}
\begin{xy}
<0mm,0mm>*{\bullet}, <4mm,-7mm>*{^{1\ \ \dots\ k\qquad \ \ k+l+1\
\ \dots \ \ n}}, <-14mm,-5mm>*{}**@{-}, <-6mm,-5mm>*{}**@{-},
<20mm,-5mm>*{}**@{-}, <8mm,-5mm>*{}**@{-}, <0mm,-5mm>*{}**@{-},
<0mm,-5mm>*{\bullet}; <-5mm,-10mm>*{}**@{-},
<-2mm,-10mm>*{}**@{-}, <2mm,-10mm>*{}**@{-}, <5mm,-10mm>*{}**@{-},
<0mm,-12mm>*{_{k+1\dots k+l}},
 <0mm,0mm>*{};<0mm,5mm>*{}**@{-},
\end{xy},
$$
\begin{eqnarray*}
\p\left( \begin{xy}
 <0mm,0mm>*{\blacktriangledown};<0mm,0mm>*{}**@{},
 <0mm,0mm>*{};<0mm,5mm>*{}**@{-},
 <0mm,0mm>*{};<-8mm,-5mm>*{}**@{-},
 <0mm,0mm>*{};<-4.5mm,-5mm>*{}**@{-},
 <0mm,0mm>*{};<-1mm,-5mm>*{\ldots}**@{},
 <0mm,0mm>*{};<4.5mm,-5mm>*{}**@{-},
 <0mm,0mm>*{};<8mm,-5mm>*{}**@{-},
   <0mm,0mm>*{};<-8.5mm,-6.9mm>*{^1}**@{},
   <0mm,0mm>*{};<-5mm,-6.9mm>*{^2}**@{},
   <0mm,0mm>*{};<4.5mm,-6.9mm>*{^{n\hspace{-0.5mm}-\hspace{-0.5mm}1}}**@{},
   <0mm,0mm>*{};<9.0mm,-6.9mm>*{^n}**@{},
 \end{xy}\right)&=&
\sum_{k=0}^{n-1}\sum_{l=1}^{n-k} (-1)^{(l-1)(n-k-l)}\left(
\begin{xy}
<0mm,0mm>*{\blacktriangledown}, <4mm,-7mm>*{^{1\ \ \dots\ k\qquad \ \ k+l+1\
\ \dots \ \ n}}, <-14mm,-5mm>*{}**@{-}, <-6mm,-5mm>*{}**@{-},
<20mm,-5mm>*{}**@{-}, <8mm,-5mm>*{}**@{-}, <0mm,-5mm>*{}**@{-},
<0mm,-5mm>*{\bullet}; <-5mm,-10mm>*{}**@{-},
<-2mm,-10mm>*{}**@{-}, <2mm,-10mm>*{}**@{-}, <5mm,-10mm>*{}**@{-},
<0mm,-12mm>*{_{k+1\dots k+l}},
 <0mm,0mm>*{};<0mm,5mm>*{}**@{-},
\end{xy} + \right.\\
&&\left. +
\begin{xy}
<0mm,0mm>*{\bullet}, <4mm,-7mm>*{^{1\ \ \dots\ k\qquad \ \ k+l+1\
\ \dots \ \ n}}, <-14mm,-5mm>*{}**@{-}, <-6mm,-5mm>*{}**@{-},
<20mm,-5mm>*{}**@{-}, <8mm,-5mm>*{}**@{-}, <0mm,-5mm>*{}**@{-},
<0mm,-5mm>*{\blacktriangledown}; <-5mm,-10mm>*{}**@{-},
<-2mm,-10mm>*{}**@{-}, <2mm,-10mm>*{}**@{-}, <5mm,-10mm>*{}**@{-},
<0mm,-12mm>*{_{k+1\dots k+l}},
 <0mm,0mm>*{};<0mm,5mm>*{}**@{-},
\end{xy}
+
\begin{xy}
<0mm,0mm>*{\blacktriangledown}, <4mm,-7mm>*{^{1\ \ \dots\ k\qquad \ \ k+l+1\
\ \dots \ \ n}}, <-14mm,-5mm>*{}**@{-}, <-6mm,-5mm>*{}**@{-},
<20mm,-5mm>*{}**@{-}, <8mm,-5mm>*{}**@{-}, <0mm,-5mm>*{}**@{-},
<0mm,-5mm>*{\blacktriangledown}; <-5mm,-10mm>*{}**@{-},
<-2mm,-10mm>*{}**@{-}, <2mm,-10mm>*{}**@{-}, <5mm,-10mm>*{}**@{-},
<0mm,-12mm>*{_{k+1\dots k+l}},
 <0mm,0mm>*{};<0mm,5mm>*{}**@{-},
\end{xy}
\right)
\end{eqnarray*}

\sip

\begin{pro}
There is a one-to-one correspondence between representations of
the dg operad $\Ass_\infty^{(2)}$ in an operad $\Qo$ and degree
$-1$ elements, $\ga$ and $\Ga$, in the deformation complex
(Hochschild complex) $\displaystyle \oplus_{n\geq 2}
s^{1-n}\Qo(n)$ such that,
\begin{eqnarray*}
[d+ \ga, d+\ga]_G &=& 0   \\
\left[d+ \ga, \Ga\right]_G + \frac{1}{2}\left[\Ga,\Ga\right]_G &=& 0,
\end{eqnarray*}
where $[\ ,\ ]_G$ stands for the Gerstenhaber brackets.
\end{pro}

Proof is obvious and hence is omitted. The data $(d,\ga)$
describes a representation of $\Ass_\infty$ on $\Qo$, and the data
$\Ga$ describes a deformation of this representation.\\

In general, $\Po_\infty^{(2)}$ is the prop(erad) given by
$\F(s^{-1}\Co^\bullet \oplus s^{-1}\Co^\blacktriangledown)$. We
denote by $\p_\Po(c)=\sum \G(c_1, \ldots, c_n)$ the image under
the differential $\p_\Po$ of an element $c$ of $s^{-1}\Co$, with
$c_1, \ldots, c_n\in \Co$. The differential $\p$ of
$\Po_\infty^{(2)}$ is defined by
\begin{eqnarray*}
\p(c^\bullet):=\sum \G(c^\bullet_1,
\ldots, c^\bullet_n) &\textrm{for}& c^\bullet \in s^{-1}\Co^\bullet, \\
\p(c^\blacktriangledown):=\sum \G(c^{i_1}, \ldots, c^{i_n}_n)
&\textrm{for}& c^{\blacktriangledown} \in
s^{-1}\Co^{\blacktriangledown},
\end{eqnarray*}
where the $i_1, \ldots, i_n$ are in $\{ \bullet,\,
\blacktriangledown\}$ with at least one equal to
$\blacktriangledown$. It is easy to see that $\p^2=0$.

\begin{proo}
The formula for $\p(c^\bullet)$ gives the first relation. With the
formula for $\p(c^\blacktriangledown)$, it gives the second one.
\end{proo}

\begin{pro}
There is a one-to-one correspondence between representations of
the dg prop(erad) $\Po_\infty^{(2)}$ in an prop(erad) $\Qo$ and
degree $-1$ elements, $\ga$ and $\Ga$, in the deformation complex
$s^{-1}\Hm(\Co,\Qo)$ such that
\begin{eqnarray*}
Q(\ga) &=& 0   \\
Q^\ga(\Ga)=Q(\ga+\Ga) &=& 0,
\end{eqnarray*}
where $Q$ stands for the $L_\infty$-algebra structure.
\end{pro}

In the following proposition, we interpret $\Ass^{(2)}_\infty$ as
the Koszul resolution of a new operad, denoted by $\Ass^{(2)}$.

\begin{pro}\label{Ass infty 2}
The dg operad $(\Ass_\infty^{(2)}, \p)$ is a quadratic resolution
of a quadratic operad $\Ass^{(2)}$ defined as the quotient of the
free operad on the $\bS$-module
$$
A(n):= \left\{ \Ba{ll} \bK[\bS_2]\oplus \bK[\bS_2]=\mbox{span}\ \left(
\ \
\begin{xy}
<0mm,0mm>*{\bullet},
<0mm,5mm>*{}**@{-},
<-3mm,-5mm>*{}**@{-},
<3mm,-5mm>*{}**@{-},
<-4mm,-7mm>*{_{\sigma(1)}},
<4mm,-7mm>*{_{\sigma(2)}},
\end{xy}
\
,
\
\begin{xy}
<0mm,0mm>*{\blacktriangledown},
<0mm,5mm>*{}**@{-},
<-3mm,-5mm>*{}**@{-},
<3mm,-5mm>*{}**@{-},
<-4mm,-7mm>*{_{\sigma(1)}},
<4mm,-7mm>*{_{\sigma(2)}},
\end{xy}
\right)_{\sigma\in \Sigma_2} & \mbox{for}\ n=2 \vspace{3mm}\\
0 & \mbox{otherwise},
\Ea
\right.
$$
modulo the ideal  generated by
relations,
$$
\begin{xy}
<0mm,0mm>*{\bullet},
<0mm,5mm>*{}**@{-},
<-3mm,-5mm>*{}**@{-},
<3mm,-5mm>*{}**@{-},
<-3mm,-5mm>*{\bullet};
<-6mm,-10mm>*{}**@{-},
<0mm,-10mm>*{}**@{-},
<-6mm,-12mm>*{_{\sigma(1)}},
<0mm,-12mm>*{_{\sigma(2)}},
<3mm,-7mm>*{_{\sigma(3)}},
\end{xy}
-
\begin{xy}
<0mm,0mm>*{\bullet},
<0mm,5mm>*{}**@{-},
<-3mm,-5mm>*{}**@{-},
<3mm,-5mm>*{}**@{-},
<3mm,-5mm>*{\bullet};
<6mm,-10mm>*{}**@{-},
<0mm,-10mm>*{}**@{-},
<6mm,-12mm>*{_{\sigma(3)}},
<0mm,-12mm>*{_{\sigma(2)}},
<-3mm,-7mm>*{_{\sigma(1)}},
\end{xy}
=0,  \ \ \ \ \forall \sigma\in \bS_3.
$$
and
$$
\begin{xy}
<0mm,0mm>*{\blacktriangledown},
<0mm,5mm>*{}**@{-},
<-3mm,-5mm>*{}**@{-},
<3mm,-5mm>*{}**@{-},
<-3mm,-5mm>*{\bullet};
<-6mm,-10mm>*{}**@{-},
<0mm,-10mm>*{}**@{-},
<-6mm,-12mm>*{_{\sigma(1)}},
<0mm,-12mm>*{_{\sigma(2)}},
<3mm,-7mm>*{_{\sigma(3)}},
\end{xy}
-
\begin{xy}
<0mm,0mm>*{\blacktriangledown},
<0mm,5mm>*{}**@{-},
<-3mm,-5mm>*{}**@{-},
<3mm,-5mm>*{}**@{-},
<3mm,-5mm>*{\bullet};
<6mm,-10mm>*{}**@{-},
<0mm,-10mm>*{}**@{-},
<6mm,-12mm>*{_{\sigma(3)}},
<0mm,-12mm>*{_{\sigma(2)}},
<-3mm,-7mm>*{_{\sigma(1)}},
\end{xy}
+
\begin{xy}
<0mm,0mm>*{\bullet},
<0mm,5mm>*{}**@{-},
<-3mm,-5mm>*{}**@{-},
<3mm,-5mm>*{}**@{-},
<-3mm,-5mm>*{\blacktriangledown};
<-6mm,-10mm>*{}**@{-},
<0mm,-10mm>*{}**@{-},
<-6mm,-12mm>*{_{\sigma(1)}},
<0mm,-12mm>*{_{\sigma(2)}},
<3mm,-7mm>*{_{\sigma(3)}},
\end{xy}
-
\begin{xy}
<0mm,0mm>*{\bullet},
<0mm,5mm>*{}**@{-},
<-3mm,-5mm>*{}**@{-},
<3mm,-5mm>*{}**@{-},
<3mm,-5mm>*{\blacktriangledown};
<6mm,-10mm>*{}**@{-},
<0mm,-10mm>*{}**@{-},
<6mm,-12mm>*{_{\sigma(3)}},
<0mm,-12mm>*{_{\sigma(2)}},
<-3mm,-7mm>*{_{\sigma(1)}},
\end{xy}
+
\begin{xy}
<0mm,0mm>*{\blacktriangledown},
<0mm,5mm>*{}**@{-},
<-3mm,-5mm>*{}**@{-},
<3mm,-5mm>*{}**@{-},
<-3mm,-5mm>*{\blacktriangledown};
<-6mm,-10mm>*{}**@{-},
<0mm,-10mm>*{}**@{-},
<-6mm,-12mm>*{_{\sigma(1)}},
<0mm,-12mm>*{_{\sigma(2)}},
<3mm,-7mm>*{_{\sigma(3)}},
\end{xy}
-
\begin{xy}
<0mm,0mm>*{\blacktriangledown},
<0mm,5mm>*{}**@{-},
<-3mm,-5mm>*{}**@{-},
<3mm,-5mm>*{}**@{-},
<3mm,-5mm>*{\blacktriangledown};
<6mm,-10mm>*{}**@{-},
<0mm,-10mm>*{}**@{-},
<6mm,-12mm>*{_{\sigma(3)}},
<0mm,-12mm>*{_{\sigma(2)}},
<-3mm,-7mm>*{_{\sigma(1)}},
\end{xy}
=0,  \ \ \ \ \forall \sigma\in \bS_3.
$$

\end{pro}

\begin{proo}
Let $F_{-p}(\Ass^{(2)}_\infty)$ be the subspace of
$\Ass^{(2)}_\infty$ spanned by trees with a least $p$ internal
edges between one vertex labelled by $\bullet$ and the other one
labelled by $\blacktriangledown$. This defines an increasing
filtration with is bounded on $\Ass^{(2)}_\infty(n)$ for each $n$.
Therefore it converges to the homology of $\Ass^{(2)}_\infty(n)$
by the Classical Convergence Theorem~$5.5.1$ of \cite{Weibel}. The
first term $E^0_{pq}$ is equal to the subspace of
$\Ass^{(2)}_\infty$ spanned by trees with exactly $p$ internal
edges between one vertex labelled by $\bullet$ and the other one
labelled by $\blacktriangledown$. And the differential $d^0$ is
equal to the sum of the differentials
$\p_{\Ass^\bullet_\infty}+\p_{\Ass^\blacktriangledown_\infty}$,
that it is splits $\bullet$ and $\blacktriangledown$ vertices into
pure $\bullet$ and $\blacktriangledown$ trees. Hence $(E^0, d^0)$
is the coproduct
$\Ass^\bullet_\infty\vee\Ass^\blacktriangledown_\infty$ (see
Section~\ref{LimColimProp(erad)s}) of two resolutions of $\Ass$,
which is acyclic. Finally, we have $E^1_{pq}=0$ for $p+q\neq 0$
and $\bigoplus_{p\ge 0}E^1_{-p\, p}=\Ass^\bullet \vee
\Ass^\blacktriangledown$. The spectral sequence collapses and the
homology of $\Ass^{(2)}_\infty$ is concentrated in degree $0$.
Another presentation of this homology group is given by the
quotient of the free operad on degree $0$ elements, namely the two
binary products $\bullet$ and $\blacktriangledown$, by ideal
generated by the image under $\p$ of the degree $1$ elements of
$\Ass^{(2)}_\infty$.
\end{proo}

In other words, the operad $\Ass^{(2)}$ is Koszul. A
representation of $\Ass^{(2)}$ in a vector space $X$ is equivalent
to a pair of linear maps $\mu: X\ot X\rar X$ and $\nu: X\ot X\rar
X$ such that both $(X, \mu)$ and $(X, \mu+\nu)$ are associative
algebras. As a corollary, we get the following isomorphism of
$\Sy$-modules $\Ass^{(2)}\cong \Ass \vee \Ass$.

\begin{Rq}
The example of $\Ass^{(2)}$ is also interesting from the viewpoint
of Koszul operad. It comes from a set theoretic operad. It is
Koszul whereas the method of \cite{Vallette06} cannot be applied
because $\Ass^{(2)}$ is not basic set, that is the composition of
operations is not injective. The product $\blacktriangledown$ has
an ``absorbing'' effect.
\end{Rq}

In the same way, we define the operad $\Li^{(2)}$ by
$\F([\,,\,]_\bullet\oplus
[\,,\,]_\blacktriangledown)/(\textrm{Jac}^\bullet_\bullet\oplus
(\textrm{Jac}^\bullet_\blacktriangledown+\textrm{Jac}^\blacktriangledown_\bullet+\textrm{Jac}^\blacktriangledown_\blacktriangledown))$,
where $[\,,\,]_\bullet$ and $[\,,\,]_\blacktriangledown$ stand for
two skew-symmetric brackets and where $\textrm{Jac}_a^b$ stands
for the ``Jacobi'' relation
$[[X,Y]_a,Z]_b+[[Y,Z]_a,X]_b+[[Z,X]_a,Y]_b=0$. This operad enjoys
the same properties with $\Li$ than the non-symmetric operad
 $\Ass^{(2)}$ with $\Ass$ explained above. More generally, to any
 binary quadratic operad $\Po$ (eventually non-symmetric) with its minimal
 model $\Po_\infty$, we can associate an operad $\Po^{(2)}$ such
 that $\Po^{(2)}_\infty$ is its minimal model.

 We summary this
 result with the explicit form of $\Po^{(2)}$ in terms of Manin
 products in the following theorem.
 Let $\Po=\F(V)/(R)$ and $\Qo=\F(W)/(S)$ be two binary quadratic non-symmetric operads. There exists a morphism of
 $\psi \, :\, \F(V)(3)\ot  \F(W)(3) \to \F(V\ot W)(3)$. \emph{Manin's black square product of $\Po$ and $\Qo$} is equal to  $\Po\blacksquare \Qo:=\F(V\ot W)/(\psi(R\ot S))$. In the symmetric case, the definition is similar but the morphism $\psi$ is more involved and requires the use of signature representations  (see \cite{Vallette06Preprint} for more details).

\begin{thm}
For any binary quadratic non-symmetric operad $\Po$ which admits a
minimal model $\Po_\infty$, the non-symmetric operad
$\Po^{(2)}_\infty$ is a minimal model (resolution) of $\Po\,
\blacksquare\,  \Ass^{(2)}$, which is isomorphic, as a graded
module, to $\Po \vee \Po$, where the coproduct has to be taken in
the category of non-symmetric operads .

For any binary quadratic operad $\Po$ which admits a minimal model
$\Po_\infty$, the operad $\Po^{(2)}_\infty$ is a minimal model
(resolution) of $\Po\, \bullet\,  \Li^{(2)}$which is isomorphic,
as an $\Sy$-module to $\Po \vee \Po$.
\end{thm}

\begin{proo}
By the same argument as in Proposition~\ref{Ass infty 2} above,
the homology of $\Po_\infty^{(2)}$ is concentrated in degree $0$.
If we denote the quadratic operad $\Po$ by $F(V)/(\mathcal{R})$,
this non-trivial homology group is equal to $F(V_\bullet \oplus
V_\blacktriangledown)/(\mathcal{R}^\bullet_\bullet \oplus
\mathcal{R}^\bullet_\blacktriangledown+\mathcal{R}^\blacktriangledown_\bullet+
\mathcal{R}^\blacktriangledown_\blacktriangledown)$, which is
equal to the black product of $\Po$ with $\Li^{(2)}$.
\end{proo}



\appendix
\section{Model category structure for prop(erad)s}
\label{model category}

In this appendix, we prove that the categories of props and
properads have a cofibrantly generated model category structure.
We make precise coproducts, pushouts, cofibrations and cofibrant
objects. We refer the reader to the book \cite{Hovey99} of M.
Hovey for a comprehensive treatment of model
categories. (In order to be self-contained in this appendix, with we will not avoid the prefix dg here.) \\

Let us denote by \emph{dg Properads} the category of dg properads
and by \emph{dg Props} the category of dg props. It means either
the category of reduced dg prop(erad)s ($\Po(m,n)=0$ for $n=0$ or
$m=0$) or the category of dg prop(erad)s over a field of
characteristic $0$. By default, we work over unbounded chain
complexes but the following proofs hold over bounded chain
complexes as well. We transfer the cofibrantly generated model
category structure of $\Sy$-bimodules to the category of
prop(erad)s via the free prop(erad) functor.

\subsection{Model category structure on
$\Sy$-bimodules}\label{CofGenModCatSbiMod} The category of dg
$\Sy$-bimodules is endowed with a cofibrantly generated model
category structure coming from the cofibrantly
generated model category structure on dg $\KK$-modules.\\

Recall from \cite{Hovey99} Theorem~$2.3.11$ that the category of
dg $R$-modules has a cofibrantly generated model category
structure for any ring $R$. Quasi-isomorphisms form the class of
weak equivalences and degreewise surjective maps form the class of
fibrations. Let us make explicit the (generating) acyclic
cofibrations. The model category of dg $\KK$-modules is
cofibrantly generated by the acyclic cofibrations $J^k\, :\, 0 \to
D^k$, where $D^k$ is the
 chain complex
$$ \cdots \to 0 \to \underbrace{\KK}_{k} \xrightarrow{\Id} \underbrace{\KK}_{k-1} \to 0 \to
\cdots$$ and by the cofibrations $I^k\, :\, S^{k-1} \to D^k$,
where $S^{k-1}$ is the following
 chain complex
$$ \cdots \to 0 \to \underbrace{\KK}_{k-1}  \to 0 \to
\cdots.$$

For any $m,n \in \NN$, the category of left $\Sy_m$ and right
$\Sy_n$-bimodules is the category of dg modules over the group
ring $\KK[\Sy_m^{\textrm{op}}\times \Sy_n]$. By
 Theorem~$2.3.11$ of \cite{Hovey99} the preceding
theorem, it has a cofibrantly generated model category structure,
where the generating acyclic cofibrations are the maps
$J^k_{m,n}\, :\, 0 \to D^k_{m,n}$, where $D^k_{m,n}$ is the
acyclic dg $\KK[\Sy_m^{\textrm{op}}\times \Sy_n]$-module
$$ \cdots \to 0 \to \underbrace{\KK[\Sy_m^{\textrm{op}}\times \Sy_n]}_{k}
\xrightarrow{\Id} \underbrace{\KK[\Sy_m^{\textrm{op}}\times
\Sy_n]}_{k-1} \to 0 \to \cdots$$ and where the generating
cofibrations are the maps $I^k_{m,n}\, :\, S^{k-1}_{m,n} \to
D^k_{m,n}$, with $S^{k-1}_{m,n}$ the following  dg
$\KK[\Sy_m^{\textrm{op}}\times \Sy_n]$-module
$$ \cdots \to 0 \to \underbrace{\KK[\Sy_m^{\textrm{op}}\times \Sy_n]}_{k-1}
 \to 0 \to \cdots.$$

Since the category of dg $\Sy$-bimodules is the product over
$(m,n) \in \NN^2$ of the model categories of left $\Sy_m$ and
right dg $\Sy_n$-bimodules, it is naturally endowed with a
term-by-term cofibrantly generated model category structure. The
set of generating acyclic cofibrations can be chosen to be
$J=\lbrace \widetilde{J}^k_{m,n}; \, k\in \ZZ, m,n \in \NN
\rbrace$, where $\widetilde{J}^k_{m,n}$ is equal to ${J}^k_{m,n}
\, :\, 0 \to D^k_{m,n}$ in arity $(m,n)$ and $0$ elsewhere.
Similarly, the set of generating cofibrations can be chosen to be
$I=\lbrace \widetilde{I}^k_{m,n}; \, k\in \ZZ, m,n \in \NN
\rbrace$, where $\widetilde{I}^k_{m,n}$ is equal to ${I}^k_{m,n}
\, :\, S^{k-1}_{m,n} \to  D^k_{m,n}$ in arity $(m,n)$ and $0$
elsewhere. Notice that the domains of elements of $I$ or $J$ are
sequentially small with respect to any map in the category of dg
$\Sy$-bimodules.

\subsection{Transfer Theorem}
In the section, we recall the theorem of transfer, mainly due to
Quillen \cite{Quillen67} Section $\textrm{II}.4$ (see also S. E.
Crans \cite{Crans95} Theorem~$3.3$ and M. Hovey \cite{Hovey99}
Proposition~$2.1.19$). We will use it to endow the category of dg
prop(erad)s with a model category structure.

\begin{dei}[Relative I-cell complexes]
For every class I of maps of a category, a \emph{relative I-cell
complexes} is a sequential colimit of pushouts of maps of I.
\end{dei}

Let us make explicit this type of morphisms. A relative I-cell
complex is a map $A_0 \xrightarrow{\varphi} A_\infty$ which comes
from a sequential colimit
$$\xymatrix{A_0 \ar[r]^{\i_0} \ar[d]^{\varphi} & A_1 \ar[r]^{\i_1}
\ar[dl] & \cdots\ar[r]  & A_n\ar[r]^{\i_n} \ar[dlll]  & A_{n+1}
\ar[dllll]
\ar[r] & \cdots \\
A_\infty:=\mathop{Colim}_\NN \, A_n,  &&&&&  } $$ where each map
$A_n \xrightarrow{i_n} A_{n+1}$ is defined by a pushout
$$\xymatrix{ \vee_\alpha S_\alpha \ar[d]^{\vee_\alpha j_\alpha}  \ar[r] & A_n \ar[d]^{i_n} \\
\vee_\alpha T_\alpha \ar[r]&  A_{n+1},}$$ with $j_\alpha \in I$. As usual, we
denote the collection of relative I-cell complexes by I-cell.

\begin{thm}[\cite{Quillen67} Section $\textrm{II}.4$, \cite{Crans95} Theorem~$3.3$, \cite{Hovey99}
Proposition~$2.1.19$]\label{TransferTheorem} Let $\Co$ be a
cofibrantly generated model category with $I$ as the set of
generating cofibrations and $J$ as the set of generating acyclic
cofibrations. Let $F \, : \Co \rightleftharpoons \Do \, : \, U$ be
an adjunction, where $F$ is the left adjoint and $U$ the right
adjoint. Suppose that \draftnote{ici tout est est
"petit-sequentiel" verifier que l'on peut en effet faire comme
cela.}
\begin{enumerate}
\item $\Do$ has finite limits and colimits,

\item the functor $U$ preserves filtered colimits,

\item the image under $U$ of any relative $F(J)$-cell complex is a
weak equivalence in $\Co$.
\end{enumerate}

A map $f$ in $\Do$ is defined to be a weak equivalence (resp.
fibration) if the associated map $U(f)$ is a weak equivalence
(resp. fibration) in $\Co$. The class of cofibrations in $\Do$ is
the class of map that verify the left lifting property (LLP) with respect to acyclic fibrations. \\

These three classes of maps provide the category $\Do$ with a
model category structure cofibrantly generated by $F(I)$ as the
the set of generating cofibrations and $F(J)$ as the set of
generating acyclic cofibrations.
\end{thm}

We also refer the reader to Section~$2.5$ of
\cite{BergerMoerdijk03} for the application of this Theorem with
stronger and sometimes more convenient hypotheses.  Remark that Transfer
Theorem~\ref{TransferTheorem} was used (and rephrased) by V.
Hinich in \cite{Hinich97} to provide a model category structure to
the category of operads over unbounded chain complexes (see
Theorem~$2.2.1$ of \cite{Hinich97} and the corrected version of
Theorem~$6.6.1$ in \cite{Hinich03}). M. Spitzweck also applied
this theorem to prove a general result about model category
structures on categories of algebras over a triple (Theorem~$1$ of
\cite{Spitzweck01}).

\subsection{Limits and Colimits of prop(erad)s}\label{LimColimProp(erad)s} In this section, we prove that
the category of prop(erad)s has all limits and finite colimits. We
also make explicit the coproducts and pushouts of prop(erad)s.

\begin{pro}
The category of prop(erad)s has all limits.
\end{pro}

\begin{proo}
We recall from D. Borisov and Y.I. Manin \cite{BorisovManin06}
 that the free prop(erad) functor induces
a triple $\F \, : \Sy\textrm{-biMod} \to \Sy\textrm{-biMod}$ such
that an algebra over it is a prop(erad). Since the underlying
category of $\Sy$-bimodules has limits, the category of
prop(erad)s has all limits (Section~$1.5$ of
\cite{GetzlerJones94}).
\end{proo}

To prove that the category of prop(erad)s has finite colimits, we
first make explicit coproducts and pushouts. This section is the
generalization of Section~$1.5$ of \cite{GetzlerJones94} from
operads to prop(erad)s. Once again, the situation is more subtle
for prop(erad)s than for operads since it requires the notion of
adjacent vertices of a graph
(see Section~\ref{admissible subgraph}).\\

Let $\Po$ and $\Qo$ be two prop(erad)s. The coproduct of $\Po$ and
$\Qo$ is given by a quotient of the free prop(erad) on their sum
$\F(\Po \oplus \Qo)$. On this space, we define an equivalence
relation by the following generating relation : if a graph $g$,
with vertices indexed by elements of $\Po$ and $\Qo$ has two
adjacent vertices indexed elements of $\Po$ (or $\Qo$), it is
equivalent to the same graph, where the two adjacent vertices are
contracted and the new vertex is labelled by the composition in
$\Po$ (or $\Qo$) of the two associated elements of $\Po$ (or
$\Qo$). The quotient of $\F(\Po \oplus \Qo)$ by this relation is
the coproduct of $\Po$ and $\Qo$. We denote it by $\Po \vee \Qo$.
This $\Sy$-bimodule has the following basis. It can be represented
by the sum over (connected) graphs with vertices indexed by
elements of $\Po$ and $\Qo$ such that no adjacent pair of vertices
are labelled by the same kind of elements (see
Figure~\ref{coproduct prop(erad)}).

\begin{figure}[h]
$$\xymatrix@R=12pt@C=8pt{ \ar@{-}[rd] & & \ar@{-}[ld] &  \\
&  *+[F-,]{\ \Qo \ } \ar@{=}[d]   & &  \\
& *+[F-,]{\ \Po \ }\ar@{-}[dd] \ar@{-}[rd] & \ar@{-}[d]& \ar@{-}[ld] \\
& & *+[F-,]{\ \Qo \ } \ar@{-}[rd] \ar@{=}[ld] &  \\
& *+[F-,]{\ \Po \ } \ar@{-}[d]\ar@{-}[rd]\ar@{-}[ld] & &  \\
& & &  }  $$ \caption{Element of the coproduct $\Po\vee \Qo$}
\label{coproduct prop(erad)}
\end{figure}

Let $\Po$, $\Qo$ and $\Ro$ be three prop(erad)s. Let $f\, : \, \Po
\to \Qo$ and $g\, : \, \Po \to \Ro$ be two morphisms of
prop(erad)s. Their pushout is isomorphic to the quotient of
$\Qo\vee \Ro$ by the ideal generated by $\lbrace f(p)-g(p),\, p\in
\Po \rbrace$. (We refer to Appendix~B of \cite{Vallette06Preprint}
for the notion of ideal of a prop(erad). The notion of ideal
generated by a sub-$\Sy$-bimodule is also made explicit there.)
The pushout $\Qo\vee_\Po \Ro$ is represented by labelled
(connected) graphs as above but further quotient by the following
relation~: if a vertex is labelled by an element of the form
$f(p)$ for $p \in \Po$, it can be replaced by the same vertex
labelled by the corresponding element $g(p)$ and vice-versa. When
this operation generates two adjacent vertices indexed by elements
of the same prop(erad), there are to be composed.

\begin{pro}
The category of prop(erad)s has finite colimits.
\end{pro}

\begin{proo}
This result can be proved with two methods.

First, recall that the free properad on an $\Sy$-bimodule $V$ is
given by the sum on (connected) graphs without level whose
vertices are coherently labelled by elements of $V$ (see
Section~$2.7$ of \cite{Vallette03}). We denote it by
$$\F(V)=\left( \bigoplus_{g\in \mathcal{G}_c} \bigotimes_{\nu \in
\mathcal{N}(g)} V(|Out(\nu)|,\, |In(\nu)|)\right) \Bigg/ \approx \
,$$ on Theorem~$2.3$ of \cite{Vallette03}, where $\mathcal{N}(g)$
is the set of vertices of a graph $g$. Since the tensor product of
dg $\Sy$-bimodules preserves colimits, the functor
\begin{eqnarray*}
\Sy\textrm{-biMod} &\to& \Sy\textrm{-biMod} \\
V &\mapsto& \bigotimes_{\nu \in \mathcal{N}(g)} V(|Out(\nu)|,\,
|In(\nu)|),
\end{eqnarray*}
associated to any graph $g$, preserves filtered colimits (see
Lemma~$1.14$ of \cite{GetzlerJones94}). Then the triple $\F \,: \,
\Sy\textrm{-biMod} \to \Sy\textrm{-biMod} $ associated to the free
prop(erad) functor preserves filtered colimits. The argument of
Page~$16$ of \cite{GetzlerJones94} proves that the category of
prop(erad)s has filtered colimits. Since it has pushouts and
filtered colimits, it has finite colimits by Chapter~IX of
\cite{MacLane}.

We can also construct coequalizers in this category. Since it is
an additive category, it is enough to construct cokernels. Let
$f\, :\, \Po \to \Qo$ be a morphism of prop(erad)s. Its cokernel
is given by the quotient of $\Qo$ with the ideal generated by the
image of $f$. Since it has coproducts and coequalizers, this
category has finite colimits by Theorem~$2.1$ Chapter~$V$ of
\cite{MacLane}.
\end{proo}

\subsection{Model category structure}
In this section, we apply the Transfer
Theorem~\ref{TransferTheorem} to provide a cofibrantly generated
model category structure on the category of prop(erad)s.\\

We consider the free prop(erad) adjunction $\F \, : \textrm{dg} \
\Sy\textrm{-biMod} \rightleftharpoons  \textrm{dg Prop(erad)s} \,
: \, U$. We proved in \ref{CofGenModCatSbiMod} that the category
on the left hand side is a cofibrantly generated model category.
We apply the Transfer Theorem~\ref{TransferTheorem} to this
adjunction as follows. The generating acyclic cofibrations are
$\F(J)=\lbrace I \to \F(D^k_{m,n})\rbrace$ and the generating
cofibrations are $\F(I)=\lbrace \F(S^{k-1}_{m,n}) \to
\F(D^k_{m,n})\rbrace$.

\begin{lem}\label{F(I),F(J)-cell}
A morphism of dg properads is a relative $\F(J)$-cell complex if and only is it is a   map $\Po \to \Po\vee \F(D)$,
where $D=\bigoplus_{d\ge 1} D_{i}$ is an acyclic dg $\Sy$-bimodule whose components are
free $\Sy$-bimodules with each $D_{i}$ equal to a direct sum of dg $\Sy$-bimodules $D^k_{m,n}$.

A morphism of dg properads is a relative $\F(I)$-cell complex if and only if it is a map
$\Po \to \Po\vee \F(S)$, where $S$ is a dg $\Sy$-bimodule, whose
components are free $\Sy$-bimodules, endowed with an exhaustive filtration $$S_0=\lbrace 0\rbrace \subset S_1\subset
S_2\subset \cdots \subset \textrm{Colim}_i S_i =S$$ such that $d \, : \, S_i \to \F(S_{i-1})$ and such that $S_{i-1}\mono S_i$ are split monomorphisms of dg $\Sy$-bimodules with cokernels isomorphic to a free $\Sy$-bimodule.
\end{lem}

\begin{proo}
Pushouts of elements of $\F(J)$ are as follows:
$$\xymatrix@H=16pt{I \ar@{>->}[d]^{\vee_\alpha \F(J^\alpha)}  \ar[r] & \Po \ar[d] \\
\bigvee_\alpha \F(D^\alpha) \ar[r]&  \Po \vee \big(\bigvee_\alpha \F(D^\alpha)\big),}$$
with each $D^\alpha$ equal to a $D^k_{m,n}$. Since the
coproduct of free prop(erad)s is the free prop(erad) on the sum of
their generating spaces,  $\F(V)\vee\F(V')\cong \F(V\oplus V')$,
the composite of two such maps is equal to $\Po \to \Po \vee
\F(\bigoplus_\alpha D^\alpha \oplus \bigoplus_\beta D^\beta)$. Hence a sequential
colimit of such pushouts has the form $\Po \to \Po\vee\F(D)$, with
$D=\bigoplus_{d\ge 1} D_{i}$ an acyclic dg $\Sy$-bimodule whose components are free
$\Sy$-bimodules.

A pushout of an element of $\F(I)$ is
$$\xymatrix@H=16pt{\vee_\alpha \F(S^\alpha) \ar@{>->}[d]^{\vee_\alpha \F(I^\alpha)}  \ar[r]^(0.6){f} & \Po \ar[d] \\
\vee_\alpha \F(D^\alpha) \ar[r]&  \Qo,}$$ with each $S^\alpha$ equal to an $S^k_{m,n}$ and  $D^\alpha$ equal to a $D^k_{m,n}$. We denote by $z$ the image under
$f$ of the generating element of $S^{k-1}_{m,n}$. Notice that $z$
is a cycle in $\Po$. If we denote by $\xi$ and $d\xi$ the
generating elements of $D^{k}_{m,n}$, the pushout $\Qo$ is equal
to $\Po\vee \F(\bigoplus_\alpha \xi^\alpha.\KK[\Sy_m^{\textrm{op}}\times \Sy_n] )$ with $d\xi=z$. Therefore a relative
$\F(I)$-cell complex is a map $\Po \to \Po\vee\F(S)$, with $S$ a
dg $\Sy$-bimodule whose components are free
$\Sy$-bimodules. Since a relative $\F(I)$-cell complex is a sequential colimit of such pushouts, the filtration of $S$ is given by this sequential guing of cells.
\end{proo}

\begin{thm}\label{CofGenModelCategoryOnProperads}
The category of prop(erad)s has a cofibrantly model category
structure provided by the following three classes of morphisms. A
map $\Po \xrightarrow{f} \Qo$ is a

\begin{itemize}
\item weak equivalence if and only if it is a quasi-isomorphism of
dg $\Sy$-bimodules, that is a quasi-isomorphism in any arity,

\item fibration if and only if it is a degreewise surjection in
any arity,

\item cofibration if and only if it has the left lifting property
with respect to acyclic fibrations.
\end{itemize}
The generating cofibrations are the maps $\F(I)=\lbrace
\F(S^{k-1}_{m,n}) \to \F(D^k_{m,n})\rbrace$ and the generating
acyclic cofibrations are the maps $\F(J)=\lbrace I \to
\F(D^k_{m,n})\rbrace$.
\end{thm}

\begin{proo}
The category of prop(erad)s has finite limits and colimits $(1)$
by the preceding section. To any dg $\Sy$-bimodule $M$, we can
consider the trivial (abelian) prop(erad) structure on $I\oplus
M$, that is the composite product is zero on $M$. So, it is easy
to check that the forgetful functor preserves filtered colimits
$(2)$. Recall from \ref{LimColimProp(erad)s} that the coproduct
$\Po\vee \F(D)$ admits a basis composed by (connected) graphs with
vertices indexed elements of $\Po$ and $D$ such that there is no
pair of adjacent vertices indexed by two elements of $\Po$.
Therefore, $\Po \vee \F(D)$ is equal to the directed sum $\Po
\oplus X$, where $X$ has a basis given by graphs indexed by
elements coming from $\Po$ and at least one element from $D$. The
map $\Po \to \Po\vee \F(D)$ is the inclusion of $\Po$ into the
first summand so that it is enough to prove that $X$ is an acyclic
chain complex. For every graph $g$ indexed by elements of $\Po$
and at least one element of $D$, the resulting chain complex is
isomorphic to a quotient by the action of some symmetric groups of
tensor products of $\Po$ and at least one $D$. Since $D$ is an
acyclic chain complex made of free $\KK[\Sy_m^{\textrm{op}}\times
\Sy_n]$-modules, it is an acyclic projective chain complex over
any ring of symmetric subgroup . Hence the chain complex
associated to any graph $g$ indexed by elements of $\Po$ and at
least one element of $D$ is acyclic, which proves hypothesis $(3)$
of Transfer Theorem~\ref{TransferTheorem}.
\end{proo}

\subsection{Cofibrations and Cofibrant objects} In this section, we make explicit the
cofibrations and the cofibrant objects in the model category of dg
prop(erad)s. We refer to the Appendix of \cite{Fresse04} for the case of operads.

\begin{pro}\label{Cofibrations}
A map $\xymatrix{f\, :\,   \Po\  \ar@{>->}[r] & \Qo}$ is a cofibration in
the model category of dg prop(erad)s if and only if it is a retract of a map $\Po \to \Po
\vee \F(S)$, with isomorphisms on domains, where $S$ is a dg
$\Sy$-bimodule whose components are  free $\Sy$-bimodules, endowed with an exhaustive filtration $$S_0=\lbrace 0\rbrace \subset S_1\subset
S_2\subset \cdots \subset \textrm{Colim}_i S_i =S$$ such that $d \, : \, S_i \to \F(S_{i-1})$ and such that $S_{i-1}\mono S_i$ are split monomorphisms of dg $\Sy$-bimodules with cokernels isomorphic to a free $\Sy$-bimodule.

A map  $\xymatrix{f\, :\,   \Po\
\ar@{>->}[r]^(0.6){\sim} & \Qo}$  is an acyclic cofibration in the model category of
dg prop(erad)s  if and only if it is  a retract of a map $\Po \to \Po \vee \F(D)$, with
isomorphisms on domains, where $D=\bigoplus_{d\ge 1} D_{i}$ is an acyclic dg $\Sy$-bimodule whose components are
free $\Sy$-bimodules with each $D_{i}$ equal to a direct sum of dg $\Sy$-bimodules $D^k_{m,n}$.
\end{pro}

\begin{proo}
The proposition follows from general results on the (acyclic)
cofibrations of cofibrantly generated model categories.
Explicitely, we apply Proposition~$2.1.18$ of \cite{Hovey99} to
the cofibrantly generated model category of prop(erad)s. This
proposition gives explicitly that (acyclic) cofibrations of
prop(erad)s are retracts of relative $\F(I)$-cell complexes
(relative $\F(J)$-cell complexes). We conclude by
Lemma~\ref{F(I),F(J)-cell}.
\end{proo}

Applied to $\Po=I$, this proposition gives to following  corollary.

\begin{pro}\label{criterion cofibrant}
A dg prop(erad) is  cofibrant for this model category structure if and only if it is  a retract of a quasi-free
prop(erad) $\F(S)$, where the components of $S$ are free
$\Sy$-bimodules, endowed with an exhaustive filtration $$S_0=\lbrace 0\rbrace \subset S_1\subset
S_2\subset \cdots \subset \textrm{Colim}_i S_i =S$$ such that $d \, : \, S_i \to \F(S_{i-1})$ and such that $S_{i-1}\mono S_i$ are split monomorphisms of dg $\Sy$-bimodules with cokernels isomorphic to a free $\Sy$-bimodule.
\end{pro}

\begin{Rq}
In the model category of dg prop(erads) on non-negatively graded
dg $\Sy$-bimodules, a dg prop(erad) is cofibrant if and only if it
is retract of a quasi-free prop(erad) $\F(S)$ whose components are
free $\Sy$-bimodules. The extra assumption on the filtration is automatically given by the homological degree.
\end{Rq}

Recall that we are working over a field of characteristic $0$.

\begin{lem}\label{retract of free}
Any quasi-free prop(erad) $\F(X)$ is a retract of a quasi-free
 prop(erad) $\F(S)$, where the components of $S$ are free
 $\Sy$-bimodules. Moreover, if $X$ is endowed with an exhaustive filtration $$X_0=\lbrace 0\rbrace \subset X_1\subset
X_2\subset \cdots \subset \textrm{Colim}_i X_i =X$$ such that $d \, : \, X_i \to \F(X_{i-1})$ and such that $X_{i-1}\mono X_i$ are split monomorphisms of dg $\Sy$-bimodules, then $S$ can be chosen with the same property and such that the cokernels of the $S_{i-1}\mono S_i$ are free $\Sy$-bimodules.
\end{lem}

\begin{proo}
Let $\bar{X}(m, n)$ denote the set of equivalence classes under
the action of $\Sy_m^{\textrm{op}}\times \Sy_n$. For simplicity,
 we use the generic notation  $\bar{X}$. We choose a set of
 representatives $\lbrace \bar{x}_i \rbrace_{i\in \mathcal{I}}$ of $\bar{X}$.
Let $S$ be the free $\Sy$-bimodule generated by the $\lbrace
\bar{x}_i \rbrace_{i\in \mathcal{I}}$. The generator associated to
$\bar{x}_i$ will be denoted by $s_i$. For any $x$ in $X$, we
consider the sub-group $\Sy_{x}:=\lbrace \sigma \in
\Sy_m^{\textrm{op}}\times \Sy_n  \, |\,  x.\sigma=\chi(\sigma)x,
\, \chi(\sigma)\in \KK\rbrace$. In this case, $\chi$ is a
character of $\Sy_x$. We define the following element of $S$~:
$$N(\bar{x}_i):=\frac{1}{|\Sy_{\bar{x}_i}|} \sum
\chi(\sigma^{-1}).s_i\sigma,$$ where the sum runs over $\sigma \in
\Sy_{\bar{x}_i}$.
 The image under the boundary map
$\partial$ of an $\bar{x}_i$ is a sum of graphs of the form $\sum
\mathcal{G}(\bar{x}_{i_1}, \ldots, \bar{x}_{i_k})$. We define the
boundary map $\partial'$ on $\F(S)$ by
$$\partial'(s_i):=\sum \frac{1}{|\Sy_{\bar{x}_i}|} \sum
\chi(\sigma^{-1}). \mathcal{G}(N(\bar{x}_{i_1}), \ldots,
N(\bar{x}_{i_k}))\sigma,$$ where the second sum runs over $\sigma
\in \Sy_{\bar{x}_i}$. Finally, we define the maps of dg
prop(erad)s $\F(S) \to \F(X)$ by $s_i \mapsto \bar{x}_i$ and
$\F(X) \to \F(S)$ by $\bar{x}_i \mapsto N(\bar{x}_i)$. They form a
deformation retract, which preserves the filtration of $X$ when it exists.
\end{proo}

\begin{cor}\label{quasi-free -> cofibrant}
In the model category of dg prop(erads), any quasi-free properad $\F(X)$, where
$X$ is  endowed with an exhaustive filtration $$X_0=\lbrace 0\rbrace \subset X_1\subset
X_2\subset \cdots \subset \textrm{Colim}_i X_i =X$$ such that $d \, : \, X_i \to \F(X_{i-1})$ and such that the $X_{i-1}\mono X_i$ are split monomorphisms of dg $\Sy$-bimodules is cofibrant.
\end{cor}

\begin{Rq}
In the non-negatively graded case, any quasi-free prop(erad) is cofibrant.
\end{Rq}

\begin{proo}
It is a direct corollary of Proposition~\ref{criterion cofibrant}
and Lemma~\ref{retract of free}.
\end{proo}

\begin{thm}\label{Exists QF regular}
Any dg properad $\Qo$ admits a cofibrant replacement of the form $\xymatrix@C=20pt{\F(S) \ar@{->>}[r]^{\sim}& Q}$,
where the components of $S$ are free
$\Sy$-bimodules, endowed with an exhaustive filtration $$S_0=\lbrace 0\rbrace \subset S_1\subset
S_2\subset \cdots \subset \textrm{Colim}_i S_i =S$$ such that $d \, : \, S_i \to \F(S_{i-1})$ and such that $S_{i-1}\mono S_i$ are split monomorphisms of dg $\Sy$-bimodules with cokernels isomorphic to a free $\Sy$-bimodule.
\end{thm}

\begin{proo}
Any dg properad $\Qo$ admits a cofibrant replacement $\xymatrix@C=20pt{I\  \ar@{>->}[r]& \Po \ar@{->>}[r]^{\sim}& Q}$. Since $\Po$ is cofibrant, it is retract $\xymatrix@C=20pt{\Po \ar[r]^{\sim} & \F(S) \ar@{->>}[r]^{\sim}& \Po}$ of such an $\F(S)$ by Proposition~\ref{criterion cofibrant}.
\end{proo}

We can simply such a cofibrant replacement as follows.

\begin{thm}\label{Exists QF}
A quasi-free cofibrant replacement $\xymatrix@C=20pt{\F(S) \ar@{->>}[r]^{\sim}& Q}$ induces a quasi-free cofibrant replacement $\xymatrix@C=20pt{\F(X) \ar@{->>}[r]^{\sim}& Q}$,
where the action of the symmetric groups on the components of $X$ is the same then the action on their image in $\Qo$. Moreover, $X$ is endowed with an exhaustive filtration $$X_0=\lbrace 0\rbrace \subset X_1\subset
X_2\subset \cdots \subset \textrm{Colim}_i X_i =X$$ such that $d \, : \, X_i \to \F(X_{i-1})$ and such that $X_{i-1}\mono X_i$ are split monomorphisms of dg $\Sy$-bimodules.
\end{thm}

\begin{proo}
 The dg $\Sy$-bimodule which generates the quasi-free cofibrant replacement $\xymatrix@C=20pt{\F(S) \ar@{->>}[r]^{\sim}& Q}$ is a free $\Sy$-bimodule. Let us denote by $s_\alpha$ the generators and $q_\alpha$ their image in $\Qo$. We define $X$ to be  the $\Sy$-bimodule generated by the $q_\alpha$ and we consider the free properad $\F(X)$ on $X$. In $\F(S)$, the image of $s_\alpha$ under the differential map $d$ is equal to $d(s_\alpha)=\sum
\mathcal{G}(s_{\alpha_1}, \ldots, s_{\alpha_k})$. We define the differential map of $\F(X)$ by
$d(q_\alpha)=\sum
\mathcal{G}(q_{\alpha_1}, \ldots, q_{\alpha_k})$. The map $\xymatrix{\F(S) \ar@{->>}[r]^{\sim}& Q}$ factors through
$\xymatrix@C=20pt{\F(S) \ar@{->>}[r]^{\sim}& \F(X) \ar@{->>}[r]^{\sim}&Q}$. Finally, $\F(X)$ is cofibrant by Corollary~\ref{quasi-free -> cofibrant}.
\end{proo}

The difference between resolution $\F(S)$ and $\F(X)$ is that in $\F(S)$, the symmetry of the operations of $\Qo$ is deformed up to homotopy whereas in $\F(X)$ only the relations are deformed up to homotopy. (The same phenomenon appears for resolutions of the operad $Com$ of commutative algebras where the former corresponds to $E_\infty$ operads and the later to $C_\infty$).

We can now choose to work with such cofibrant models. The extra filtration on the space of generators, which appears conceptually here, is similar to the one used by Sullivan  \cite{Sullivan77} in rational homotopy theory  and by Markl in \cite{Markl96M} for operads.\\

Let $\Po$ be a dg properad. Its space of \emph{indecomposable elements} is the cokernel of the composite map with non-trivial elements, $\mu \, :\, \oPo \bc \oPo \to \oPo$. The space of indecomposable elements inherits a differential map from the one of $\Po$ which makes it into a dg $\Sy$-bimodule. The associated functor $\textrm{Indec} \, :\, \textrm{dg} \ \textrm{properads} \to \textrm{dg} \Sy\textrm{-bimodules}$ is left adjoint to the augmentation functor $M \mapsto M\oplus I$, where the properad structure on $M\oplus I$ is the trivial one.

The following last result will allow us to proof that the deformation complex defined in ... does not depend on the quasi-free model chosen to make it explicit.

\begin{pro}\label{preserves EQ}
Any weak equivalence (quasi-isomorphism) between two quasi-free cofibrant dg properads $\F(X) \xrightarrow{\sim} \F(Y)$ induces a weak equivalence (quasi-isomorphism) between the spaces of indecomposable elements $X \xrightarrow{\sim} Y$.
\end{pro}

\begin{proo}
The two categories of dg properads and dg $\Sy$-bimodules have  model categories structures. Since the augmentation functor preserves fibrations and acyclic fibrations, by Lemma~$1.3.4$ of \cite{Hovey99} the indecomposable functors, being its left adjoint, preserves cofibrations and acyclic cofibrations. And by Brown's Lemma (Lemma~$1.1.12$ of \cite{Hovey99}), it preserves weak equivalences between cofibrant objects.
\end{proo}

\bip


\begin{center}
\textsc{Acknowledgements}
\end{center}
This work was partially supported by the G\"oran Gustafsson
foundation. The first author expresses his thanks to the Max
Planck Institute in Mathematics in Bonn and \'Ecole Normale
Sup\'erieure in Paris for stimulating conditions during the work
on this project. The second author would like to thank Jean-Louis
Loday, Dennis Sullivan, Benoit Fresse, Clemens Berger, Bernhard Keller and
Domenico Fiorenza   for useful discussions. The second author
expresses his thanks to Northwestern University for the invitation
in Spring quarter 2007. B.V. is supported by the ANR OBTH.

\bibliographystyle{amsalpha}
\bibliography{bib}

\medskip

{\small \textsc{Sergei Merkulov, Department of
Mathematics, Stockholm University,}\\
\textsc{10691 Stockholm, Sweden}\\
E-mail address : \texttt{sm@math.su.se}\\
URL : \texttt{http://www.math.su.se/$\sim$sm/}}\\

{\small \textsc{Bruno Vallette, Laboratoire J.A. Dieudonn\'e,
Universit\'e de Nice Sophia-Antipolis,} \\
\textsc{Parc Valrose, 06108 Nice
Cedex 02, France}\\
E-mail address : \texttt{brunov@math.unice.fr}\\
URL : \texttt{http://math.unice.fr/$\sim$brunov}}

\end{document}